\pdfoutput=1

\RequirePackage{fix-cm}

\documentclass[leqno, fleqn, centertags, 12pt]{article}

\usepackage{latexsym}
\usepackage{amsmath}
\usepackage{amssymb}
\usepackage{verbatim}

\usepackage[T1]{fontenc}

\usepackage[upright, widespace]{fourier}

\usepackage{fullpage}

\usepackage{tikz-cd}

\tikzcdset{row sep/special/.initial=12ex}
\tikzcdset{row sep/spe/.initial=10.5ex}
\tikzcdset{column sep/spec/.initial=5em}
\tikzcdset{column sep/specc/.initial=8em}
\tikzcdset{column sep/speccc/.initial=10em}

\usepackage{comment}

\newskip\nineskipamount \nineskipamount=9pt plus 0pt minus 0pt
\newskip\zeroskipamount \zeroskipamount=0pt plus 0pt minus 0pt
\usepackage[noorphans,vskip=\nineskipamount]{quoting}

\makeatletter
\renewcommand{\@makefntext}[1]{\vspace*{0.5ex}\parindent=0em
\hspace*{-0.4em}
\hbox to 0.4em{\hss\@makefnmark}\hspace*{0.4em}{#1}
}
\makeatother

\newcounter{mysectionnumber}
\newcommand{\mysection}[2]{\setcounter{footnote}{0}
\setcounter{myparnum}{0}
\refstepcounter{mysectionnumber}
\vspace{21pt}{\Large {\themysectionnumber.} {#1}}\label{#2}\vspace*{15pt}}

\newcommand{\mynonumbersection}[1]{
\vspace{21pt}{\Large {#1}}\vspace*{15pt}}

\newcommand{\myit}[1]{\textbf{\textit{#1}}\hspace{0.0em}}

\newcounter{myparnum}[mysectionnumber]
\renewcommand{\themyparnum}{\arabic{mysectionnumber}.\arabic{myparnum}}
\newcommand{\mypar}[2]{\refstepcounter{myparnum}{\vspace{\medskipamount}\textbf{{\themyparnum. #1}\label{#2}}\hspace{0.5em}}}

\newcounter{mylemmanum}[myparnum]
\newcommand{\myuppar}[1]{\vspace{\medskipamount}\textbf{#1}\hspace*{0.5em}}

\newcounter{myappendnumber}
\newcounter{myaparnum}[myappendnumber]
\newcommand{\myappend}[2]{\setcounter{footnote}{0}
\setcounter{myaparnum}{0}
\setcounter{myparnum}{0}
\refstepcounter{myappendnumber}
\vspace{21pt}{\Large A\dff.{\themyappendnumber.}\oss {#1}}\label{#2}\vspace*{15pt}}

\newcommand{\proof}{\vspace{\medskipamount}{\textbf{{\emph{Proof}.}}\hspace*{1em}}}
\newcommand{\prooftitle}[1]{\vspace{\medskipamount}{\textbf{{\emph{#1}.}}\hspace*{1em}}}

\newcommand{\eproof}{ $\blacksquare$}
\newcommand{\esubproof}{ $\square$}

\newcommand{\dis}{\displaystyle}

\def\sss{\hspace{0.05em}\ }
\def\dss{\hspace{0.1em}\ }
\def\trs{\hspace{0.15em}\ }
\def\qss{\hspace{0.2em}\ }
\def\pss{\hspace{0.3em}\ }
\def\oss{\hspace{0.4em}\ }

\def\halfff{\hspace*{0.025em}}
\def\fff{\hspace*{0.05em}}
\def\dff{\hspace*{0.1em}}
\def\trf{\hspace*{0.15em}}
\def\qff{\hspace*{0.2em}}
\def\pff{\hspace*{0.3em}}
\def\off{\hspace*{0.4em}}

\newcommand{\hnsp}{\hspace*{-0.05em}}
\newcommand{\nsp}{\hspace*{-0.1em}}
\newcommand{\nnsp}{\hspace*{-0.15em}}
\newcommand{\snsp}{\hspace*{-0.175em}}
\newcommand{\dnsp}{\hspace*{-0.2em}}

\renewcommand{\leq}{\leqslant}
\renewcommand{\geq}{\geqslant}
\newcommand{\id}{\mathop{\mbox{id}}\nolimits}

\newcommand{\zzz}{\mathbb{Z}}
\newcommand{\rrr}{\mathbb{R}}

\newcommand{\ftwo}{\mathbb{F}_{\dff 2}}

\newcommand{\num}[1]{|\qff #1 \qff|}
\newcommand{\hclass}[1]{[\dff #1 \dff]}
\newcommand{\fclass}[1]{[\snsp [\dff #1\dff]\snsp ]}

\newcommand{\conv}[1]{\langle\dff #1 \dff\rangle}

\newcommand{\ttoo}{\hspace*{0.2em}\longrightarrow\hspace*{0.2em}}

\newcommand{\bd}{\operatorname{bd}}

\begin{document}

\setlength{\baselineskip}{12pt plus 0pt minus 0pt}
\setlength{\parskip}{12pt plus 0pt minus 0pt}
\setlength{\abovedisplayskip}{12pt plus 0pt minus 0pt}
\setlength{\belowdisplayskip}{12pt plus 0pt minus 0pt}

\newskip\smallskipamount \smallskipamount=3pt plus 0pt minus 0pt
\newskip\medskipamount   \medskipamount  =6pt plus 0pt minus 0pt
\newskip\bigskipamount   \bigskipamount =12pt plus 0pt minus 0pt

\setlength{\skip\footins}{6pt}

\author{Nikolai\qss V.\qss Ivanov}
\title{Scarf's\pss theorems,\oss simplices,\oss and\qss oriented\qss matroids}
\date{}

\footnotetext{\hspace*{-0.65em}\copyright\oss 
Nikolai\qss V.\qss Ivanov,\oss 2019,\oss 2022.\oss
With\sss the exception of\dss the preface,\oss 
this\sss paper was completed\sss in\qss 2019.\oss}

\footnotetext{\hspace*{-0.65em}The author\dss is\dss grateful\dss to\qss V.{\dff}I.\dss Danilov\trs
for attracting\dss his attention\dss to\sss the works of\qss H.\dss Scarf.}

\maketitle

\vspace*{12ex}

{\renewcommand{\baselinestretch}{1}
\selectfont

\myit{\hspace*{0em}\large Contents}\vspace*{1.5ex} \\ 
\hbox to 0.8\textwidth{\myit{Preface} \hfil 2}\hspace*{0.5em} \vspace*{1.5ex}\\
\hbox to 0.8\textwidth{\myit{\phantom{1}1.}\hspace*{0.5em} Pseudo-simplices \hfil 5}\hspace*{0.5em} 
\vspace*{0.25ex}\\
\hbox to 0.8\textwidth{\myit{\phantom{1}2.}\hspace*{0.5em} Combinatorics of\trs families 
of\trs linear orders \hfil 12}\hspace*{0.5em} \vspace*{0.25ex}\\
\hbox to 0.8\textwidth{\myit{\phantom{1}3.}\hspace*{0.5em} Scarf's\qss proof\dss
of\qss Brouwer's\qss fixed\dss point\dss theorem \hfil 16}\hspace*{0.5em} \vspace*{0.25ex}\\
\hbox to 0.8\textwidth{\myit{\phantom{1}4.}\hspace*{0.5em} An example\fff:\pss integer\dss points in\sss
a simplex  \hfil 18}\hspace*{0.5em}\vspace*{0.25ex}\\
\hbox to 0.8\textwidth{\myit{\phantom{1}5.}\hspace*{0.5em} Oriented\dss matroids \hfil 27}\hspace*{0.5em} \vspace*{0.25ex}\\
\hbox to 0.8\textwidth{\myit{\phantom{1}6.}\hspace*{0.5em} Oriented\dss matroid colorings\fff:\oss the non-degenerate case \hfil 31}\hspace*{0.5em} \vspace*{0.25ex}\\
\hbox to 0.8\textwidth{\myit{\phantom{1}7.}\hspace*{0.5em} Vector colorings \hfil 39}\hspace*{0.5em} \vspace*{0.25ex}\\
\hbox to 0.8\textwidth{\myit{\phantom{1}8.}\hspace*{0.5em} Oriented\dss matroid colorings\fff:\oss the general case \hfil 42}\hspace*{0.5em} \vspace*{0.25ex}\\
\hbox to 0.8\textwidth{\myit{\phantom{1}9.}\hspace*{0.5em} Scarf's\qss proof\dss 
of\qss Kakutani's\qss fixed\dss point\dss theorem \hfil 46}\hspace*{0.5em} \vspace*{0.25ex}\\
\hbox to 0.8\textwidth{\myit{10.}\hspace*{0.5em} Chains in $\rrr^{\dff n}${\nnsp} and\dss 
their\sss applications \hfil 50}\hspace*{0.5em} \vspace*{1.5ex}\\
\myit{Appendices}\hspace*{0.5em}  \hspace*{0.5em} \vspace*{1.5ex}\\
\hbox to 0.8\textwidth{\myit{A.1.}\hspace*{0.5em} Todd's\dss theorem\hfil 59}\hspace*{0.5em} \vspace*{0.25ex}\\
\hbox to 0.8\textwidth{\myit{A.2.}\hspace*{0.5em} Extensions of\dss oriented\dss matroids\hfil 61}\hspace*{0.5em}  \vspace*{1.5ex}\\
\hbox to 0.8\textwidth{\myit{References}\hspace*{0.5em}\hfil 62}\hspace*{0.5em}  \vspace*{0.25ex}  

}

\renewcommand{\baselinestretch}{1}
\selectfont

\newpage
\mynonumbersection{Preface}

\myuppar{Scarf\dss proof\dss of\qss Brouwer's\qss fixed\dss point\dss theorem
and\dss Scarf\dss combinatorial\dss theorem.}
In\qss 1967\qss H.\dss Scarf\pss \cite{sc2}\qss suggested a new\dss
proof\dss of\qss Brouwer's\qss fixed\dss point\dss theorem.\oss
His paper\dss opens with an outline of\trs the well\dss known\dss proof\dss of\qss 
Brouwer's\qss theorem\dss 
based on\qss Sperner's\trs lemma and\dss
Knaster--Kuratowski--Mazurkiewicz\dss argument\halfff.\oss
As\dss is\dss well\dss known,\oss
this proof\dss begins\sss with a\sss triangulation of\dss the standard simplex\sss
$\Delta^n\qff \subset\qff \rrr^{\dff n\dff +\dff 1}$\sss
and\sss eventually\dss requires\sss triangulations of\sss $\Delta^n$\sss 
into\sss arbitrarily small\sss simplices.\oss

Scarf's\dss proof\qss follows a similar outline,\oss
but\dss triangulations are replaced\sss by sufficiently\dss dense 
finite subsets\sss $X\qff \subset\qff \Delta^n$\dnsp.\oss
In\sss particular,\oss one can\sss take set\sss of\dss vertices of\dss 
a\sss triangulation\sss into small\sss simplices as $X$\nnsp.\oss
As in\dss Sperner's\dss lemma,\oss the points of\dss $X$\sss
are colored\sss into\sss $n\qff +\qff 1$\sss colors.\oss 
Actually,\oss \emph{colorings}\pss is\dss a relatively modern\sss terminology,\oss
neither\dss Sperner,\oss nor\dss Scarf\trs used\dss it\halfff.\oss
Scarf's\sss key\dss result\dss is\dss an analogue of\qss Sperner's\trs lemma.\oss
It\dss claims\sss that\sss under assumptions similar\sss to\sss the assumptions of\qss
Sperner's\trs lemma\sss there exists\sss
$\sigma\qff \subset\qff X$\sss such\sss that\sss the colors of\dss different\sss
elements of\sss $\sigma$ are different\sss and\sss a\sss technical\sss condition holds.\oss
The\sss latter ensures\sss that\sss $\sigma$\sss is\dss small\dss 
if\dss $X$\sss is\dss sufficiently dense.\oss
In contrast\sss with\dss Sperner's\dss lemma $\sigma$ may contain\sss
less\sss than\sss $n\qff +\qff 1$\sss points.\oss
Replacing\sss sets of\dss vertices of\dss triangulations
by sufficiently dense,\oss but\sss otherwise arbitrary,\oss finite subsets of\dss $\Delta^n$\sss
requires a subtle modification of\pss 
Knaster--Kuratowski--Mazurkiewicz\qss coloring\sss rule\halfff:\oss
the inequalities must\dss be reversed.

While\dss Scarf\trs uses a very\sss geometric\sss language,\oss
his analogue of\trs Sperner's\dss lemma\dss is\dss a\sss truly\sss combinatorial\sss
result\sss and can\sss be reformulated\sss in\sss terms of\dss finite collections of\dss
linear orders on a finite set.\oss 
See\qss \emph{Scarf's\dss combinatorial\dss theorem}\pss in\dss Section\qss \ref{orders}.\oss
The origins of\qss this\sss theorem and\dss its\sss proof\trs belong\dss
to\sss the\sss linear\dss programming\dss and\dss the game\sss theory\halfff.\oss
Actually\halfff,\oss Scarf's\qss proof\dss of\qss Brouwer's\qss theorem\dss is\dss
a\sss byproduct\sss of\qss Scarf's\qss fundamental\dss paper\qss \cite{sc1}\qss
in\dss game\sss theory\halfff.\oss

\myuppar{Vector\sss colorings.}
For\sss the game\sss theory\dss Scarf\trs needed\dss
to consider a much more general\sss notion of\dss colorings.\oss
Namely,\oss he needed\qss \emph{vector\sss colorings},\oss
assigning\sss to each point\sss $x$\sss of\sss $X\qff \subset\qff \Delta^n$\sss
a vector\sss $\varphi\dff(\trf x\trf)\qff \in\qff \rrr^{\dff n\dff +\dff 1}$\dnsp.\oss
The usual\sss colorings are naturally\sss interpreted as vector colorings\sss
taking as values\sss the vertices of\sss $\Delta^n$\dnsp.\oss
Scarf's\dss analogue of\trs Sperner's\dss lemma\sss claims\sss that\sss 
under assumptions about\sss a vector coloring 
$\varphi$ similar\sss to\sss the assumptions of\qss Sperner's\trs lemma\dss
there exists a subset\sss $\sigma\qff \subset\qff X$\sss such\sss that
$\varphi\dff(\trf \sigma\trf)$ together\sss with some\sss vertices of\sss $\Delta^n$\sss
(which ones depends only on $\sigma$\nnsp)\qss
forms a basis of\sss $\rrr^{\dff n\dff +\dff 1}$\sss
and,\oss moreover,\oss the coordinates of\dss the center of\dss the standard simplex\sss $\Delta^n$\sss
with respect\sss to\sss this basis are non-negative.\oss
The\sss latter\sss property\sss replaces\sss the classical\dss property of\dss
having\sss all\sss available colors as colors of\dss vertices of\dss a simplex.\oss
See\qss \emph{Scarf\qss theorem}\pss in\dss Section\qss \ref{vector-colorings}\qss 
for\sss the details.\oss

\myuppar{Matroid colorings.}
This\sss theorem of\qss Scarf\trs appears\sss to be\sss a\sss truly\sss geometric one,\oss
in contrast\sss with\dss Scarf's\dss combinatorial\dss theorem,\oss
and\dss Scarf's\dss proof\trs is\dss based on\sss geometric 
ideas from\sss the\sss linear\sss programming.\oss 
Nevertheless,\oss as we will\sss see,\oss there\dss is\dss a purely combinatorial\sss
generalization of\dss this\sss theorem of\qss Scarf.\oss
Namely,\oss one can\sss replace vector colorings by colorings\sss
taking values in an\qss \emph{oriented\dss matroid},\oss
which may\sss be\sss thought\sss as finite combinatorial\dss models
of\dss vector spaces over\qss \emph{ordered}\pss fields.\oss
We will\sss call\sss such colorings\qss \emph{matroid\sss colorings}.\oss

One of\dss the main\sss goals of\dss the present\sss paper\dss
is\dss to prove a generalization of\qss Scarf\trs theorem\sss to matroid colorings.\oss
See\dss Sections\qss \ref{scarf-oriented-matroids}\qss and\qss \ref{scarf-oriented-matroids-general}.\oss
We attempted\sss to made\sss this generalization accessible\sss to\sss the readers not\sss
familiar even with\sss the definition of\dss oriented\sss matroids.\oss
Section\qss \ref{oriented-matroids}\qss is\dss devoted\sss 
to an\sss introduction\sss to oriented\sss matroids,\oss
with\sss more\sss advanced\sss results of\qss Todd\qss \cite{t}\qss
and\trs Las Vergnas\qss \cite{lv}\qss relegated\sss to\sss two appendices.\oss
In\dss Section\qss \ref{scarf-oriented-matroids}\qss we prove our\sss generalization\sss
in\sss a non-degenerate case,\oss where\sss the main\sss ideas are\sss the most\sss transparent.\oss
In\dss Section\qss \ref{vector-colorings}\qss we prove\sss the original\trs Scarf's\dss
theorem about\sss vector colorings,\oss including\sss the degenerate case.\oss
Finally,\oss in\dss Section\qss \ref{scarf-oriented-matroids-general}\qss
we consider\sss the general\dss matroid colorings.\oss

\myuppar{Scarf'\dss methods and\dss topology.}
At\sss the first\sss sight\trs Scarf's\dss methods appear\sss to be radically\sss
different\dss from\sss the methods of\dss algebraic\qss (or combinatorial\fff)\qss topology.\oss
In any case,\oss this was\sss the first\dss impression of\dss the present\sss author.\oss
For a\sss topologist\sss such as\sss the present\sss author,\oss 
the most\sss striking\sss feature of\qss Scarf's\dss methods\dss
is\dss the absence of\trs triangulations.\oss
Nevertheless,\oss more general\sss simplicial\sss complexes are only\sss hidden,\oss
and already\trs H.\dss Kuhn\qss \cite{ku}\qss related\dss Scarf's\dss combinatorial\dss theorem\sss
with\dss Sperner's\dss lem\-ma and classical\dss topological\sss concepts of\dss
pseudo-manifolds and simplicial\sss subdivision.\oss

In\sss fact,\oss the classical\dss ideas of\dss combinatorial\dss topology\sss
provide\sss the right\sss framework\sss for\dss Scarf's\dss results.\oss
In\sss the spirit\sss of\dss author's papers\qss \cite{i1}\qss --\qss \cite{i3}\qss
we approach\sss the results of\qss Scarf\trs and\sss their\sss generalization\sss
to matroid colorings\sss from\sss the point\sss of\dss view of\dss
the classical\sss combinatorial\sss topology,\oss with simplices,\oss
chains,\oss and cochains being\sss the main\sss tools.\oss

\myuppar{Simplices,\oss chains,\oss and cochains.}
The present\sss paper begins with a generalization and\sss refinement\sss 
of\dss the ideas of\qss Kuhn\qss \cite{ku}.\oss
The first\dss basic notion\dss is\dss that\sss of\dss a\qss \emph{pseudo-simplex},\oss 
a simplicial\sss complex equipped\sss
with a family of\dss subcomplexes imitating\dss the family of\dss faces\qss
(of\dss arbitrary codimension)\qss of\dss a simplex.\oss
For example,\oss a closed\sss pseudo-manifolds\sss 
with a\sss top-dimensional\sss simplex\sss removed,\oss
considered\dss together with\sss the collection of\dss all\dss proper
faces of\dss the removed simplex\dss is\dss a pseudo-simplex.\oss
Pseudo-manifolds with a\sss top-dimensional\sss simplex\sss removed
were used\dss by\trs Kuhn\qss \cite{ku}\qss in\sss his approach\sss to\dss
Scarf's\dss results,\oss but\dss in\sss this situation\sss the finer\sss structure of\dss
faces\dss is\dss trivial,\oss and\sss was ignored\dss by\trs Kuhn.\oss

The next\sss basic notion\dss is\dss that\sss of\dss a\qss \emph{chain-simplex}.\oss 
A\qss \emph{chain-simplex}\qss is\dss defined as a simplicial\sss complex equipped\sss
with a family of\dss subcomplexes imitating\dss the family of\dss faces of\dss a simplex
only at\sss the\sss level\sss of\dss chains with coefficients in\dss the field of\dss two elements $\ftwo$\nsp.\oss
Every\sss pseudo-simplex\dss is\dss a chain-simplex\sss in a natural\sss way.\oss
The notion of\dss a chain-simplex\dss is\dss an axiomatization of\dss the main\sss features
behind\sss the proofs of\trs lemmas of\pss Alexander\dss and\qss Sperner\halfff,\oss
as\sss they are presented\dss in\qss \cite{i2}.\oss
As we will\sss see,\oss the same structure\dss is\dss behind\sss the proof\dss
of\qss Scarf\dss combinatorial\sss theorem.\oss
Pseudo-simplices and\sss chain-simplices are\sss 
the\sss topic of\trs Section\qss \ref{pseudo-simplices}.\oss

\myuppar{Scarf's\dss construction of\dss pseudo-simplices.}
From\sss this point\sss of\dss view a key\sss idea of\qss Scarf\dss
can\sss be interpreted as a construction of\dss a pseudo-simplex\sss 
from a family of\trs linear orders on a finite set.\oss
This\dss is\dss purely combinatorial\sss result.\oss
When\sss combined\sss with\sss the analogue of\trs Alexander's\qss 
(or,\qss what\dss is\dss the same,\oss Sperner's)\qss 
lemma\sss from\dss Section\qss \ref{pseudo-simplices},\oss
this\sss leads\sss to\trs Scarf's\dss combinatorial\dss theorem.\oss
This\dss is\dss the\sss topic of\trs Section\qss \ref{orders}.\oss

In\dss Section\qss \ref{scarf-brouwer}\qss we follow\dss Scarf\trs
and apply\dss Scarf's\dss combinatorial\dss theorem\dss to finite subsets
of\dss $\Delta^n$\sss with orders defined\dss by\sss the barycentric coordinates,\oss
or,\oss equivalently,\oss by\dss the usual\sss coordinates in\sss $\rrr^{\dff n\dff +\dff 1}$\dss
(since\sss $\Delta^n$\sss is\dss the standard simplex\sss in\sss $\rrr^{\dff n\dff +\dff 1}$\dnsp,\oss
they are\sss the same).\oss
Together with\sss the modified\trs
Knaster--Kuratowski--Mazurkiewicz\dss argument\dss
this\sss leads\sss to\dss Scarf's\dss proof\dss
of\qss Brouwer's\qss fixed\dss point\dss theorem.\oss

\myuppar{Scarf's\dss pseudo-simplices and\dss triangulations.}
One may suspect\sss that\sss the simplices of\trs Scarf's\dss pseudo-simplex constructed\sss
from a finite subset\sss of\dss $\Delta^n$ are,\oss
in\sss fact,\oss simplices of\dss some\sss triangulation of\sss $\Delta^n$\sss
with vertices belonging\sss to $X$\nnsp,\oss
and\sss the whole proof\trs is\dss almost\sss the standard\sss proof\dss in disguise.\oss
But,\oss while\sss there\dss is\dss a canonical\sss surjective map 
from\qss ({\fff}the\sss geometric realization of\dff)\qss this pseudo-simplex\sss to $\Delta^n$\dnsp,\oss
in\sss general\sss the images of\dss the interiors of\dss different\sss simplices overlap 
and we don't\sss get\sss a\sss triangulation of\sss $\Delta^n$\dnsp.\oss
Such examples were constructed\dss by\dss Scarf\qss \cite{sc3}.\oss
We do not\sss discuss\sss them\sss in\sss the present\sss paper.\oss

But\sss for some natural\sss choices of\sss $X\qff \subset\qff \Delta^n$\sss 
one indeed\sss gets\sss triangulations of\sss $\Delta^n$\dnsp,\oss
and very\sss beautiful\sss ones.\oss
In\dss Section\qss \ref{example}\qss we discuss\sss them\sss in details  
following\trs Scarf\qss \cite{sc3}.\oss
These\sss triangulations\dss implicitly\dss appeared\sss in\sss the works of\trs Lebesgue.\oss 
See\qss \cite{i3},\oss Section\qss 2\qss for\sss the details.\oss 
Later\sss these\sss triangularions were explicitly\sss
constructed\dss by\trs Freudenthal\qss \cite{f}\qss 
as a solution of\dss a problem of\trs Brouwer.\oss
Simultaneously\sss with\dss Scarf\trs these\sss 
triangulations were rediscovered\sss by\trs Quillen\dss in\sss his work\qss \cite{q}\qss on\sss
the higher algebraic\dss $K$\dnsp-theory.\oss
See\dss Segal\qss \cite{se},\oss Appendix\qss 1.\oss

\myuppar{Scarf's\dss proof\dss of\qss Kakutani\dss fixed\sss point\dss theorem.}
Kakutani\dss fixed\sss point\dss theorem\dss is\dss a natural\sss generalization
of\qss Brouwer's\dss fixed\sss point\dss theorem\dss to multivalued\sss maps\sss
having closed convex sets as\qss (sets of\fff)\qss values.\oss
Kakutani\qss \cite{k}\qss deduced\sss his\sss theorem\sss from\dss Brouwer's.\oss
Scarf\qss \cite{sc3}\qss suggested a more natural\dss proof\dss
and deduced\trs Kakutani\trs theorem directly\sss from\sss
his\sss theorem about\sss vector colorings.\oss
In author's opinion,\oss this proof\dss clarifies both\dss the\dss Kakutani\dss
fixed\sss point\dss theorem and\dss Scarf's\dss theorem about\sss vector colorings.\oss
In\dss Section\qss \ref{kakutani}\qss we present\sss a version of\qss Scarf's\dss proof\dss
making clear\sss that\sss for a vector coloring $\varphi$\sss the vectors\sss
$\varphi\dff(\trf x\trf)\qff -\qff x$\sss should\sss be\sss thought\sss as\sss
tangent\sss vectors\sss to $\Delta^n$\dnsp.\oss

\myuppar{Combinatorial\dss topology of\sss $\rrr^{\dff n}$\nsp\dnsp.}
In\trs Section\qss \ref{scarf-theorem-geometric}\qss 
we approach\sss to\sss this circle of\dss questions\sss from\sss 
the viewpoint\sss of\dss the classical\sss combinatorial\dss topology.\oss
We consider\sss rectilinear chains in\sss $\rrr^{\dff n}$\sss
and\sss prove\sss the basic properties of\dss the intersection\sss numbers\sss
$c\dff \cdot\dff d$\dss of\dss such chains\sss in\sss the case 
when\sss the dimension of\sss $c$ or $d$\sss is\dss equal\sss to $0$ or\sss $1$\nnsp.\oss 
This\sss leads\sss to new proofs of\qss Scarf\trs and\qss Kakutani\trs theorems,\oss
and\sss to a\sss generalization of\qss Kannai's\qss \cite{k}\qss
\emph{Generalized\qss Sperner's\qss lemma}.\oss
But,\oss naturally,\oss these methods are not\sss sufficient\sss
to deal\sss with\dss matroid colorings.\oss

\newpage
\mysection{Pseudo-simplices}{pseudo-simplices}

\myuppar{Simplex-families.}
Let\sss $I$ be a finite set\halfff.\oss
A\qss \emph{simplex-family}\qss based on $I$
is\dss defined as a map\qss
$\mathcal{D}\dff \colon\dff
A\off \longmapsto\off \mathcal{D}\dff(\dff A\dff)$\dss
assigning\dss to every\sss subset\dss $A\qff \subset\pff I$\dss
a\sss simplicial\sss complex $\mathcal{D}\dff(\dff A\dff)$
of\dss dimension\dss $\num{A}\qff -\qff 1$\nnsp.\oss
Note\dss that\dss
$\mathcal{D}\dff(\trf B\trf)$ is\dss \emph{not\sss assumed}\pss to be 
a subcomplex of\qss $\mathcal{D}\dff(\dff A\dff)$\dss
when\qss $B\pff \subset\pff A$\nnsp.\oss
It\dss may\dss be comfortable\sss to\sss think\dss that\sss all\sss complexes
$\mathcal{D}\dff(\dff A\dff)$ are subcomplexes of\dss some universal\sss 
simplicial\sss complex.\oss

For every\dss finite set\sss $A$\sss let\sss $\Delta\dff(\dff A\dff)$
be\sss the abstract\sss simplicial\sss complex\sss having\dss $A$\dss
as its set\sss of\dss vertices and\sss all\dss subsets of\dss $A$ as its simplices.\oss
In other\sss words,\pss $\Delta\dff(\dff A\dff)$ is\dss the simplex\sss having\sss $A$
as\sss the set\sss of\dss vertices and considered as a simplicial\sss complex.\oss
A basic example of\dss a\sss simplex-family\dss based on\dss $I$\dss is\dss the map\qss
$\Delta\dff \colon\dff
A\off \longmapsto\off \Delta\dff(\dff A\dff)$\nnsp,\oss
where\dss $A\qff \subset\pff I$\nnsp.\oss

Let\dss $d\dff(\dff A\dff)$\dss be\sss the dimension of\dss
$\Delta\dff(\dff A\dff)$\nnsp,\oss
and\dss let\dss $e\dff(\dff A\dff)\off =\off d\dff(\dff A\dff)\qff -\qff 1$\nnsp.\oss
So,\pss $d\dff(\dff A\dff)\off =\off \num{A}\qff -\qff 1$\dss
and\dss $e\dff(\dff A\dff)\off =\off \num{A}\qff -\qff 2$\nnsp.\oss
So,\oss
for a simplex-family\sss $\mathcal{D}$\sss based on\dss $I$\dss and a subset\dss
$A\qff \subset\pff I$\dss the dimension of\trs the simplicial\sss complex\sss $\mathcal{D}\dff(\dff A\dff)$
is\dss equal\dss to $d\dff(\dff A\dff)$\nnsp.\oss
We will\dss be mostly\dss interested\dss in
simplices of\trs $\mathcal{D}\dff(\dff A\dff)$ of\dss dimensions
$d\dff(\dff A\dff)$ and $e\dff(\dff A\dff)$\nnsp.\oss

\myuppar{Pseudo-simplices.}
Suppose\sss that\sss $\mathcal{D}$ is\dss a\sss simplex-family\dss based on $I$\nnsp.\oss
Let\dss $A\qff \subset\pff I$\dss and\dss let\sss $\sigma$\sss be an
$e\dff(\dff A\dff)$\dnsp-simplex of\dss either\dss $\mathcal{D}\dff(\dff A\dff)$\dss
or\dss $\mathcal{D}\dff(\trf B\dff)$\dss for some
$e\dff(\dff A\dff)$\dnsp-subset\dss $B\qff \subset\qff A$\nnsp.\oss
Let\dss us denote by\dss $r_{\dff A}\dff(\dff \sigma\dff)$\dss the number of\dss
$d\dff(\dff A\dff)$\dnsp-simplices of\trs $\mathcal{D}\dff(\dff A\dff)$\dss
having $\sigma$ as a face,\oss
and\dss by\dss $s_{\dff A}\dff(\dff \sigma\dff)$\dss the number of\dss 
$e\dff(\dff A\dff)$\dnsp-subsets\dss $B\qff \subset\qff A$\dss such\dss that\sss $\sigma$\sss
is\dss a\sss simplex\dss of\qss $\mathcal{D}\dff(\trf B\dff)$\nnsp.\oss

The simplex-family\dss $\mathcal{D}$\sss is\dss said\dss to
be a\dss \emph{pseudo-simplex}\pss if\qss\vspace{1.5pt}
\begin{equation}
\label{sum-2}
\quad
r_{\dff A}\dff(\dff \sigma\dff)\qff +\qff s_{\dff A}\dff(\dff \sigma\dff)
\off =\off
2
\end{equation}

\vspace{-10.5pt}
for every\sss $A\qff \subset\pff I$ and $e\dff(\dff A\dff)$\dnsp-simplex $\sigma$
of\dss either\dss  $\mathcal{D}\dff(\dff A\dff)$
or\dss $\mathcal{D}\dff(\trf B\dff)$ for an
$e\dff(\dff A\dff)$\dnsp-subset\dss $B\qff \subset\qff A$\nnsp.

Every\dss triangulation\sss $T$\sss of\dss a geometric simplex $\delta$\sss leads\sss
to a\sss pseudo-simplex.\oss
Let\sss $I$\sss be\dss the set\sss of\dss vertices of\dss
$\delta$\nnsp.\oss
For\qss $A\qff \subset\pff I$\qss let\sss $\delta_{\dff A}$\dss be\sss
the face of\dss $\delta$\sss having\sss $A$\sss as\sss its set\sss of\dss vertices.\oss
In\dss particular\halfff,\pss $\delta\off =\off \delta_{\trf I}$\nsp.\oss
The\sss triangulation\sss $T$\sss induces a\sss triangulation\sss $T_{\fff A}$\sss of\dss
$\delta_{\dff A}$\nnsp.\oss
Let\dss $\mathcal{D}_{\qff T}\dff(\dff A\trf)$\dss
be\sss the abstract\sss simplicial\sss complex associated\dss with\sss $T_{\fff A}$\nnsp.\oss
The non-branching\sss property\sss of\trs $T$\sss implies\sss that\dss
the map\qss
$\mathcal{D}_{\qff T}\dff \colon\dff
A\off \longmapsto\off \mathcal{D}_{\qff T}\dff(\dff A\dff)$\dss
is\dss a\sss pseudo-simplex.\oss 
In\dss fact\halfff,\oss even more\dss is\dss true.\oss
Namely\halfff,\oss every\sss complex\dss $\mathcal{D}_{\qff T}\dff(\dff A\trf)$\dss
is\dss a\sss pseudo-manifold of\dss dimension $d\dff(\dff A\dff)$\nnsp.\oss
In\dss particular\halfff,\pss $\mathcal{D}_{\qff T}\dff(\dff A\trf)$\dss is\dss
dimensionally\dss homogeneous,\oss
i.e.\qss every\sss simplex\dss is\dss a\sss face of\dss a simplex
of\dss dimension $d\dff(\dff A\dff)$\nnsp.\oss
There\-fore\dss 
$r_{\dff A}\dff(\dff \sigma\dff)\qff \geq\qff 1$\dss
for every\sss simplex $\sigma$ of\trs  $\mathcal{D}_{\qff T}\dff(\dff A\trf)$\nnsp.\oss
In\dss the next\sss section we will\sss deal\dss with\dss pseudo-simplices\sss
for\dss which\dss $r_{\dff A}\dff(\dff \sigma\dff)\qff \geq\qff 1$\dss
does\dss not\dss necessarily\dss hold.\oss

\myuppar{Dimension $-\qff 1$\nnsp.}
By\sss a standard convention,\oss 
for every\sss simplicial\sss complex\dss
the empty\sss set\sss $\varnothing$\sss 
is\dss the only\sss simplex of\dss dimension
$-\qff 1$\nnsp.\oss
In an agreement\dss with\dss this,\oss the boundary\dss of\dss 
every\sss $0$\dnsp-simplex\dss is\dss interpreted
as\sss the unique\sss $-\qff 1$\dnsp-dimensional\sss simplex $\varnothing$\nnsp.\oss

The following\dss lemma exploits\sss these conventions.\oss
Perhaps,\oss some readers\dss will\dss prefer\dss to\sss take its conclusion
as a\sss part\sss of\trs the definition of\dss a\sss pseudo-simplex.\oss

\mypar{Lemma.}{empty}
\emph{If\pss $\mathcal{D}$\dss is\dss a\sss pseudo-simplex,\oss
then\dss the complex\dss $\mathcal{D}\dff(\dff A\dff)$\dss
has only\sss one vertex\dss for every\dss $1$\dnsp-element\dss subset\qss
$A\qff \subset\pff I$\nnsp.\oss}

\proof
{\dff}If\qss $\num{A}\off =\off 1$\nnsp,\oss then\dss $e\dff(\dff A\dff)\off =\off -\qff 1$\dss
and\dss hence $\varnothing$\sss there\dss is\dss the only\sss 
$e\dff(\dff A\dff)$\dnsp-subset\sss of\trs $A$\nnsp.\oss
The dimension of\qss $\mathcal{D}\dff(\dff \varnothing\dff)$\sss is\dss $-\qff 1$\dss
and\dss hence\dss $\varnothing$\dss is\dss the only\sss simplex of\trs
$\mathcal{D}\dff(\dff \varnothing\dff)$\nnsp.\oss
It\dss follows\dss that\trs
$s_{\dff A}\dff(\dff \varnothing\dff)\off =\off 1$\nnsp.\oss
In view of\qss (\ref{sum-2})\qss this implies\sss that\trs 
$r_{\dff A}\dff(\dff \varnothing\dff)\off =\off 1$\nnsp.\oss
But\sss every\sss simplex\dss has $\varnothing$ as a face and\dss hence\dss
$r_{\dff A}\dff(\dff \varnothing\dff)$\sss is\dss the number of\dss
$0$\dnsp-simplices of\trs $\mathcal{D}\dff(\dff A\dff)$\nnsp.\oss
The lemma follows.\oss  \eproof

\myuppar{Chain-simplices.}
Let\sss $\mathcal{D}$ be a simplex-family\dss based on $I$\nnsp.\oss
For every\sss  
$A\qff \subset\pff I$\qss
let\dss $\mathcal{D}\dff\fclass{A\fff}$\dss be\sss the formal\sss sum of\trs all\sss
$d\dff(\dff A\trf)$\dnsp-simplices of\dss $\mathcal{D}\dff(\dff A\trf)$\dss
considered as a chain\dss with coefficients in\dss $\ftwo$\nsp.\oss
In\dss particular\halfff,\pss 
$\mathcal{D}\dff\fclass{\varnothing}\off =\off \varnothing$\nnsp,\oss
where $\varnothing$ in\dss the right\dss hand side\dss is\dss 
considered as\sss $-\qff 1$\dnsp-dimensional\sss simplex.\oss
When\dss there\dss is\dss no danger of\dss confusion,\oss
the notation\sss $\mathcal{D}\dff\fclass{A\fff}$\sss will\dss be abbreviated\dss
to\sss $\fclass{A}$\nnsp.\oss
The simplex-family\dss $\mathcal{D}$\sss is\dss said\dss to be
a\qss \emph{chain-simplex}\qss based\dss on\dss $I$\dss 
if\trs for every\dss $A\qff \subset\pff I$\qss\vspace{2.25pt}
\begin{equation}
\label{chain-simplex-def}
\quad
\partial\qff
\mathcal{D}\dff\fclass{A\fff}
\off =\off
\sum\nolimits\qff \mathcal{D}\dff\fclass{\fff B\fff}
\qff,
\end{equation}

\vspace{-9.75pt}
where the sum\sss is\dss taken over\dss all\dss 
simplices $B$ of\dss $\Delta\dff(\dff A\trf)$ of\dss dimension\dss $e\dff(\dff A\dff)$\nnsp,\oss
i.e.\qss over all\sss subsets\sss $B$ of\dss $A$ with\dss
$\num{B}\off =\off \num{A\fff}\qff -\qff 1$\nnsp.\oss
The notion of\dss a chain-simplex\dss is\dss an axiomatization of\dss the main homological\dss relations
underlying\dss proofs of\trs lemmas of\pss Alexander\dss and\qss Sperner\halfff.\oss

\mypar{Lemma.}{pseudo-is-chain}
\emph{Every\dss pseudo-simplex\dss is\dss a\sss chain-simplex.\oss}

\proof
Let\dss us\dss compute\sss the both sides of\qss (\ref{chain-simplex-def}).\oss
By\dss the definition of\trs the boundary\sss operator $\partial$\vspace{2.25pt}\vspace{0.125pt}
\begin{equation}
\label{left-sum}
\quad
\partial\qff
\mathcal{D}\dff\fclass{A\fff}
\off =\qff\off
\sum\nolimits_{\qff \sigma}\qff r_{\dff A}\dff(\dff \sigma\dff)\qff \sigma
\qff,
\end{equation}

\vspace{-9.75pt}
where\sss the sum\dss is\dss taken over all\sss 
$e\dff(\dff A\dff)$\dnsp-simplices $\sigma$ of\qss $\mathcal{D}\dff(\dff A\dff)$\nnsp.\oss
On\dss the other\dss hand\dss\vspace{2.25pt}\vspace{0.125pt}
\begin{equation}
\label{right-sum}
\quad
\sum\nolimits\qff \mathcal{D}\dff\fclass{\fff B\fff}
\off =\off
\sum\nolimits_{\qff \sigma}\qff s_{\dff A}\dff(\dff \sigma\dff)\qff \sigma
\qff,
\end{equation}

\vspace{-9.75pt}
where\sss the sum\dss is\dss taken over all\sss 
$e\dff(\dff A\dff)$\dnsp-simplices $\sigma$ of\dss complexes\dss $\mathcal{D}\dff(\trf B\dff)$\dss
with\dss $B\qff \subset\qff A$\dss and\dss $\num{B}\off =\off \num{A}\qff -\qff 1$\nnsp.\oss

If\qss an $e\dff(\dff A\dff)$\dnsp-simplex $\sigma$ is\dss not\sss a simplex of\qss
$\mathcal{D}\dff(\dff A\dff)$\nnsp,\oss
then\dss $r_{\dff A}\dff(\dff \sigma\dff)\off =\off 0$\nnsp.\oss
Sim\-i\-lar\-ly\halfff,\oss
if\qss an $e\dff(\dff A\dff)$\dnsp-simplex $\sigma$ is\dss not\sss a simplex of\qss
$\mathcal{D}\dff(\trf B\dff)$\dss for any\sss $e\dff(\dff A\dff)$\dnsp-subset\dss 
$B\qff \subset\qff A$\nnsp,\oss
then\dss $s_{\dff A}\dff(\dff \sigma\dff)\off =\off 0$\nnsp.\oss
Therefore,\oss the sums in\dss the right\dss hand sides of\qss (\ref{left-sum})\qss
and\qss (\ref{right-sum})\qss can\dss be both\dss taken over all\sss
$e\dff(\dff A\dff)$\dnsp-simplices $\sigma$ 
of\dss either\dss $\mathcal{D}\dff(\dff A\dff)$\dss
or\dss $\mathcal{D}\dff(\trf B\dff)$\dss for some
$e\dff(\dff A\dff)$\dnsp-subset\dss $B\qff \subset\qff A$\nnsp.\oss
Since\dss $2\off =\off 0$\dss in\dss $\ftwo$\nsp,\oss
this observation\dss together\dss with\qss (\ref{sum-2})\qss implies\sss that\dss these
sums are equal.\oss
It\dss follows\dss that\dss the left\dss hand sides of\qss (\ref{left-sum})\qss
and\qss (\ref{right-sum})\qss are equal\sss and\dss hence\qss (\ref{chain-simplex-def})\qss
holds.\oss  \eproof

\myuppar{Envelopes of\dss simplex-families.}
Let\sss $\mathcal{D}$ be a simplex-family\dss based on $I$\nnsp.\oss
The\dss \emph{envelope}\dss of\trs $\mathcal{D}$ 
is\dss the simplex-family\dss $\mathcal{E}$\dss based on\dss $I$\dss and defined as follows.\oss

Let\sss
$V_{\fff \mathcal{D}}$\sss be\sss the union of\trs the sets of\dss vertices of\dss complexes\dss
$\mathcal{D}\dff(\dff A\dff)$\dss with\dss $A\qff \subset\pff I$\nnsp.\oss
We may assume\sss that\dss $V_{\fff \mathcal{D}}$\dss and\dss $I$\dss are disjoint\halfff.\oss
In any case,\oss nobody\sss expects\sss the\sss sets\dss 
$V_{\fff \mathcal{D}}\fff,\pff I$\dss to intersect\halfff.\oss
If\trs $A$\sss is\dss a\sss proper subset\sss of\trs $I$\nnsp,\oss
then\dss 
$\mathcal{E}\dff (\dff A\trf)\off =\off \Delta\dff(\dff A\trf)$\nnsp.\oss
The simplicial\sss complex\dss $\mathcal{E}\dff (\trf I\trf)$\dss
has\sss the union\dss $V_{\fff \mathcal{D}}\qff \cup\pff I$\dss as\sss the set\sss of\dss vertices.\oss
For every\qss $A\qff \subset\pff I$\dss
the unions of\trs the form
$\sigma\qff \cup\pff K$\nnsp,\oss
where $\sigma$ is\dss a\sss simplex of\dss $\mathcal{D}\dff (\dff A\trf)$\dss
and\trs $K$\dss is\dss proper\dss subset\sss
of\trs $I$\sss disjoint\dss from $A$\nnsp,\oss
are simplices of\dss $\mathcal{E}\dff (\trf I\trf)$\nnsp.\oss
Of\dss course,\oss if\qss $A\off \neq\off \varnothing$\nnsp,\oss
then\dss $K$\dss is\dss automatically\sss a\sss proper subset\halfff.\oss
There are no other simplices.\oss

\myuppar{\dnsp$*$\dnsp-product\sss of\dss simplices.}
Suppose\sss that\dss $A\qff \subset\pff I$\dss and\dss $A\off \neq\off \varnothing$\nnsp.\oss
If\trs $\sigma$\sss is\dss a simplex of\dss $\mathcal{D}\dff (\dff A\dff)$\dss
and\dss
$K\qff \subset\pff I\pff \smallsetminus\qff A$\nnsp,\oss
then\dss $\sigma\qff \cup\qff K$\dss is\dss a\sss simplex of\dss
$\mathcal{E}\dff (\trf I\trf)$\nnsp.\oss
It\dss will\dss be convenient\dss to denote\sss the simplex\dss
$\sigma\qff \cup\qff K$\dss also\sss by\dss $\sigma\dff *\trf K$\nnsp.\oss
The operation\dss $(\dff \sigma\fff,\pff K\trf)\off \longmapsto\off \sigma\dff *\trf K$\dss
extends\sss to a\dss bilinear\dss product\sss $\alpha\dff *\trf \Gamma$\dss
between\dss the chains $\alpha$ of\dss $\mathcal{D}\dff(\dff A\dff)$\dss
and\dss the chains\dss $\Gamma$\sss of\trs 
$\Delta\dff(\trf I\pff \smallsetminus\qff A\trf)$\nnsp.\oss

\mypar{Lemma.}{star-leibniz}
$\partial\dff (\dff \sigma\dff *\dff K\trf)
\off =\off\qff
(\trf \partial\dff \tau\dff)\dff *\trf K
\off\qff +\off\qff
\sigma\dff *\qff \partial\trf K$\nsp.\oss

\proof
The boundary\sss $\partial\dff (\dff \sigma\dff *\dff K\trf)$
is\dss the sum of\dss simplices obtained\dss by\dss removing\dss from\dss the simplex
$\sigma\dff *\dff K$ either one vertex\sss of\dss $\sigma$ or\dss 
one vertex\sss of\dss $K$\nnsp.\oss
Clearly\halfff,\pss $(\trf \partial\dff \tau\dff)\dff *\trf K$\sss is\dss 
the
sum of\trs the simplices of\trs the first\dss type and\dss
$\sigma\dff *\dff \partial\qff L$\dss is\dss 
the sum of\trs the simplices of\trs
the second\dss type.\oss  \eproof

\mypar{Lemma.}{top-simplices}
\emph{The dimension of\qss a\sss simplex\sss
of\pss $\mathcal{E}\dff (\trf I\trf)$\qss
is\dss equal\dss to\sss $d\dff(\trf I\trf)$
if\qss and\trs only\qss if\pss it\dss
has\sss the\dss form\dss
$\sigma\dff *\dff (\qff I\qff \smallsetminus\qff A\trf)$\qss
for some non-empty\dss subset\qss $A\qff \subset\pff I$\qss
and\sss some\sss $d\dff(\dff A\dff)$\dnsp-simplex\sss $\sigma$\sss 
of\oss $\mathcal{D}\dff (\dff A\dff)$\nnsp.\oss}

\proof
Let\dss $A\qff \subset\pff I$\nnsp.\oss
Let\sss
$\sigma$ be\dss a\sss simplex of\dss $\mathcal{D}\dff (\dff A\dff)$
and\sss $K$\sss be a\sss proper\sss subset\sss
of\trs $I$\sss disjoint\dss from $C$\nnsp.\oss
Clearly\halfff,\qss  
$\sigma\dff *\dff K$\dss is\dss a $d\dff(\trf I\trf)$\dnsp-simplex\dss
if\trs and\dss only\trs if\qss
$\num{\sigma\qff \cup\pff K}\off =\off \num{\trf I\trf}$\nnsp.\oss
Since\sss the dimension\dss of\dss $\mathcal{D}\dff (\dff A\dff)$\sss is\dss equal\dss to\dss
$d\dff(\dff A\dff)\off =\off \num{A}\qff -\qff 1$\dss
and\dss hence\dss $\num{\sigma}\qff \leq\qff \num{A}$\nnsp,\oss
the\sss latter condition\dss holds\sss
if\trs and\dss only\trs if\qss
$\num{\sigma}\off =\off \num{A}$\dss
and\dss
$K\off =\off I\pff \smallsetminus\qff A$\nnsp.\oss
Since $K$\dss is\dss a\sss proper subset\sss of\trs $I$\nnsp,\oss the set\sss $A$\dss
has\sss to be
non-empty\halfff.\oss
This\sss proves\sss the\qss ``only\trs if\dff''\qss part.\oss
The\qss ``if\dff''\qss part\sss is\dss obvious.\oss  \eproof

\mypar{Theorem.}{envelope-pseudo}
\emph{If\pss $\mathcal{D}$ is\dss a\dss pseudo-simplex\halfff,\oss
then\dss its\dss envelope $\mathcal{E}$ is\dss also a\dss pseudo-simplex\halfff.\oss
Also,\pss $\mathcal{E}\dff(\trf I\dff)$ is\dss a\sss non-branching complex\sss
and\dss its\dss boundary\dss is\dss
equal\dss to\sss the boundary\sss of\qss
$\Delta\dff(\trf I\trf)$\nnsp.}

\proof
Let\sss $A$\sss be\dss a\sss proper subset\sss of\trs $I$\nnsp.\oss
Then $\mathcal{E}\dff (\dff A\dff)\off =\off\Delta\dff(\dff A\dff)$\nnsp.\oss
Every $e\dff(\dff A\dff)$\dnsp-simplex of\trs $\Delta\dff(\dff A\dff)$\sss
is\dss a\sss subset\dss $C\qff \subset\qff A$\dss
such\dss that\trs $\num{C\fff}\off =\off \num{A}\qff -\qff 1$\nnsp.\oss
Clearly\halfff,\pss $A$\sss is\dss the only\sss
$d\dff(\dff A\dff)$\dnsp-simplex of\trs  $\Delta\dff(\dff A\dff)$\dss
having $C$ as a face,\oss
and $C$ is\dss a simplex of\qss
$\mathcal{E}\dff (\trf B\dff)\off =\off \Delta\dff (\trf B\dff)$\dss
for a\sss proper subset\trs $B\qff \subset\qff A$\dss
only\dss if\trs $B\off =\off C$\nnsp.\oss
This proves\qss (\ref{sum-2})\qss for\dss $\mathcal{E}$\dss and\dss proper subsets $A$ of\trs $I$\nnsp.\oss
It\dss remains\sss to deal\dss with\sss $\mathcal{E}\dff(\trf I\trf)$\nnsp.\oss

Let\dss $\tau\dff *\trf K$\trs be\dss an $e\dff(\trf I\trf)$\dnsp-simplex
of\trs $\mathcal{E}\dff(\trf I\trf)$\nnsp,\oss
and\dss let\trs $D\off =\off I\pff \smallsetminus\qff K$\nnsp.\oss
Then\dss $\num{\tau}\off =\off \num{D}\qff -\qff 1$\dss
and\dss
$D\off \neq\off \varnothing$\nnsp.\oss
By\qss Lemma\qss \ref{top-simplices}\qss every
$d\dff(\trf I\trf)$\dnsp-simplex\sss of\dss $\mathcal{E}\dff(\trf I\dff)$\dss
has\sss the\dss form\dss
$\sigma\dff *\dff (\qff I\qff \smallsetminus\qff A\trf)$\qss
for some non-empty\dss subset\qss $A\qff \subset\pff I$\qss
and\sss some\sss $d\dff(\dff A\dff)$\dnsp-simplex $\sigma$\ 
of\qss $\mathcal{D}\dff (\dff A\dff)$\nnsp.\oss
Obviously\halfff,\pss 
$\num{\sigma}
\off =\off
\num{A}$\nnsp.\oss

If\qss $\tau\dff *\trf K$\dss is\dss a\sss face of\qss
$\sigma\dff *\dff (\qff I\qff \smallsetminus\qff A\trf)$\nnsp,\oss
then\dss either\dss $\tau$\sss is\dss a\sss face of\trs $\sigma$\sss
such\dss that\dss
$\num{\sigma}\off =\off \num{\tau}\qff +\qff 1$\dss
and\dss
$K\off =\off I\pff \smallsetminus\qff A$\nnsp,\oss
or\trs
$\tau\off =\off \sigma$\dss and\dss
$K$\sss is\dss a\sss  subset\sss of\trs $I\pff \smallsetminus\qff A$\dss
such\dss that\trs $\num{K}\off =\off \num{\dff I\pff \smallsetminus\qff A}\qff -\qff 1$\nnsp.\oss
In\dss the first\sss case\dss $A\off =\off I\pff \smallsetminus\qff K\off =\off D$\nnsp,\oss
and\dss in\dss the second case\sss $A$ is\dss a\sss subset\sss
of\trs $I\qff \smallsetminus\qff K\off =\off D$\dss
such\dss that\trs
$\num{A}\off =\off \num{D}\qff -\qff 1$\nnsp.\oss

Conversely\halfff,\oss 
if\trs$\tau$\sss is\dss a\sss face of\dss a simplex $\sigma$
of\qss $\mathcal{D}\dff(\trf D\dff)$\sss
such\dss that\dss
$\num{\sigma}\off =\off \num{\tau}\qff +\qff 1$\nnsp,\oss
then\dss the $*$\dnsp-product\trs 
$\sigma\dff *\dff (\qff I\qff \smallsetminus\qff D\trf)
\off =\off
\sigma\dff *\trf K$\dss
is\dss a\sss simplex\dss having\dss $\tau\dff *\trf K$\dss as\dss its\dss face.\oss
The number of\trs such simplices\dss 
$\sigma\dff *\dff (\qff I\qff \smallsetminus\qff D\trf)$\dss
is\dss equal\dss to $r_{\trf D}\dff(\dff \tau\dff)$\nnsp.\oss
If\trs $A$ is\dss a\sss subset\sss
of\trs $D$\dss
such\dss that\trs
$\num{A}\off =\off \num{D}\qff -\qff 1$\dss
and\dss
$A\off \neq\off \varnothing$\nnsp,\oss
then\dss the $*$\dnsp-product\trs 
$\tau\dff *\dff (\qff I\qff \smallsetminus\qff A\trf)$\dss
is\dss a\sss simplex\dss having\dss $\tau\dff *\trf K$\dss as\dss its\dss face.\oss
Since\dss
$\num{A}\off =\off \num{D}\qff -\qff 1\off =\off \num{\tau}$\nnsp,\oss
the condition\dss $A\off \neq\off \varnothing$\dss is\dss equivalent\dss to\dss
$\tau\off \neq\off \varnothing$\nnsp.\oss
Hence\dss if\trs $\tau\off \neq\off \varnothing$\nnsp,\oss
then\dss the number of\trs such simplices\dss
$\tau\dff *\dff (\qff I\qff \smallsetminus\qff A\trf)$\dss
is\dss equal\dss to $s_{\trf D}\dff(\dff \tau\dff)$\nnsp.\oss

It\dss follows\dss that\trs if\trs $\tau\off \neq\off \varnothing$\nnsp,\oss
then\dss the number of\dss
$d\dff(\trf I\trf)$\dnsp-simplices\sss of\dss $\mathcal{E}\dff(\trf I\dff)$\dss
having\dss
$\tau\dff *\trf K$\dss as\dss a\sss face\dss is\dss equal\dss to\dss 
$r_{\trf D}\dff(\dff \tau\dff)\qff +\qff s_{\trf D}\dff(\dff \tau\dff)$\nnsp.\oss
By\qss (\ref{sum-2})\qss this number\dss is\dss equal\dss to\sss $2$\nnsp.\oss
If\trs $\tau\off =\off \varnothing$\nnsp,\oss
then\dss $\num{D}\off =\off 1$\dss and\dss
the only\sss $d\dff(\trf I\trf)$\dnsp-simplices of\trs $\mathcal{E}\dff(\trf I\trf)$\dss
having\dss $\tau\dff *\trf K\off =\off \varnothing\dff *\trf K\off =\off K$\dss 
as\dss a\sss face are\sss simplices of\trs 
the form\dss $\sigma\dff *\trf K$\nnsp,\oss
where $\sigma$ is\dss a $0$\dnsp-simplex of\trs $\mathcal{D}\dff(\trf D\dff)$\nnsp.\oss
By\qss Lemma\qss \ref{empty}\qss the number of\dss such simplices\dss is\dss
equal\dss to\sss $1$\nnsp.\oss

We see\sss that\sss every\sss $e\dff(\trf I\trf)$\dnsp-simplex of\trs
$\mathcal{E}\dff(\trf I\trf)$\dss
is\dss either a face of\dss exactly\dss two $d\dff(\trf I\trf)$\dnsp-simplices,\oss
or\dss is\dss a\sss subset\trs $K\qff \subset\pff I$\dss such\dss that\dss
$\num{K}\off =\off \num{\dff I\dff}\qff -\qff 1$\dss 
and\dss is\dss a face of\dss exactly\dss one $d\dff(\trf I\trf)$\dnsp-simplex.\oss
Therefore\trs $\mathcal{E}\dff(\trf I\trf)$\dss
is\dss non-branching\sss and\qss
(\ref{sum-2})\qss holds for\dss $A\off =\off I$\dss
and\dss $\mathcal{E}$\dss in\dss the role of\trs $\mathcal{D}$\nnsp.\oss  \eproof

\mypar{Theorem.}{envelope-chain}
\emph{If\pss $\mathcal{D}$ is\dss a\dss chain-simplex\halfff,\oss
then $\mathcal{E}$ is\dss also a\dss chain-simplex and\qss
$\partial\qff \mathcal{E}\dff\fclass{\trf I\trf}\off =\off \partial\trf I$\nnsp.}

\proof
With $\mathcal{E}$ in\dss the role of\dss $\mathcal{D}$\nnsp,\oss
the equality\qss (\ref{chain-simplex-def})\qss is\dss trivial\dss if\trs
$A\off \neq\off I$\nnsp.\oss
Let\dss us\dss prove\sss it\dss 
for\dss $A\off =\off I$\nnsp.\oss
Lemma\qss \ref{top-simplices}\qss implies\sss that\vspace{4.5pt} 
\begin{equation*}
\quad
\mathcal{E}\dff\fclass{\trf I\trf}
\off =\qff\off
\sum\nolimits_{\qff C}\off
\sum\nolimits_{\qff \sigma}\off \sigma\dff *\dff (\qff I\qff \smallsetminus\qff C\trf)
\qff,
\end{equation*}

\vspace{-7.5pt}
where $C$\sss runs over all\sss non-empty\sss subsets of\trs $I$\dss and\dss 
$\sigma$ over $d\dff(\trf C\trf)$\dnsp-simplices of\dss
$\mathcal{D}\dff (\trf C\trf)$\dss
for each\sss $C$\nnsp.\oss
Together\dss with\qss Lemma\qss \ref{star-leibniz}\qss 
this\sss implies\sss that\vspace{4.5pt}
\begin{equation}
\label{enlarged-fclass-bd}
\quad
\partial\qff
\mathcal{E}\dff\fclass{\trf I\trf}
\off\qff =\qff\off
\sum\nolimits_{\qff C}\off
\sum\nolimits_{\qff \sigma}\qff \partial\dff \sigma\dff *\dff (\qff I\qff \smallsetminus\qff C\trf)
\off\qff +\off\qff
\sum\nolimits_{\qff C}\off
\sum\nolimits_{\qff \sigma}\qff \sigma\dff *\dff \partial\dff (\qff I\qff \smallsetminus\qff C\trf)
\qff.
\end{equation}

\vspace{-7.5pt}
Let\dss us consider\dss these\sss two double sums separately\halfff.\oss
The first\sss double sum\dss is\dss equal\dss to\vspace{4.5pt}
\[
\quad
\sum\nolimits_{\qff C}\pff
\partial\pff 
\biggl(\pff
\sum\nolimits_{\qff \sigma}\qff \sigma
\qff\biggr)\dff *\dff (\qff I\qff \smallsetminus\qff C\trf)
\off\qff =\off\qff
\sum\nolimits_{\qff C}\off
\bigl(\qff
\partial\pff 
\fclass{C\fff}
\trf\bigr)\dff *\dff (\qff I\qff \smallsetminus\qff C\trf)
\qff,
\]

\vspace{-7.5pt}
where\dss $\fclass{\bullet\halfff}$\dss stands for\dss $\mathcal{D}\dff\fclass{\bullet\halfff}$\nnsp.\oss
Now\qss (\ref{chain-simplex-def})\qss implies\sss
that\dss the first\sss double sum\dss is\dss equal\dss to\vspace{4.5pt}
\begin{equation}
\label{first-sum}
\quad
\sum\nolimits_{\qff C}\off
\sum\nolimits_{\qff D}\off
\fclass{D}
\dff *\dff (\qff I\qff \smallsetminus\qff C\trf)
\qff,
\end{equation}

\vspace{-7.5pt}
where\sss $C$\sss runs over\dss non-empty\sss subsets of\trs $I$\dss and\dss
$D$ runs over subsets of\trs $I$ such\dss that\dss $D\qff \subset\qff C$\dss and\dss
$\num{D}\off =\off \num{C\fff}\qff -\qff 1$\nnsp.\oss
The second double sum\dss is\dss equal\dss to\vspace{4.5pt}
\[
\quad
\sum\nolimits_{\qff C}\off
\sum\nolimits_{\qff \sigma}\qff \sigma\dff *\dff 
\biggl(\pff
\sum\nolimits_{\qff E}\qff (\qff I\qff \smallsetminus\qff E\trf)
\qff\biggr)
\qff,
\]

\vspace{-7.5pt}
where\sss $C$\sss runs over\dss non-empty\sss subsets of\trs $I$\dss and\dss
$E$ runs over subsets of\dss $I$ such\dss that\dss 
$C\qff \subset\pff E$\dss and\dss
$\num{E}\off =\off \num{C\fff}\qff +\qff 1$\nnsp.\oss
By\dss interchanging\dss the order of\dss summation\dss
we see\sss that\dss it\dss is\dss equal\dss to\vspace{6pt}
\[
\quad
\sum\nolimits_{\qff C}\off
\sum\nolimits_{\qff E}\off
\biggl(\pff
\sum\nolimits_{\qff \sigma}\qff \sigma\dff *\dff 
(\qff I\qff \smallsetminus\qff E\trf)
\qff\biggr)
\off\qff =\off\qff
\sum\nolimits_{\qff C}\off
\sum\nolimits_{\qff E}\off
\fclass{C\fff}\dff *\dff 
(\qff I\qff \smallsetminus\qff E\trf)
\qff,
\]

\vspace{-6pt}
where\dss $C\fff,\pff E$\dss are subject\dss to\sss the same conditions as above.\oss
By\dss renaming\dss $C\fff,\pff E$\dss into\dss $D\fff,\pff C$\dss respectively\halfff,\oss
we see\sss that\dss the second double sum\dss is\dss equal\dss to\vspace{6pt}
\begin{equation}
\label{second-sum}
\quad
\sum\nolimits_{\qff D}\off
\sum\nolimits_{\qff C}\off
\fclass{D}\dff *\dff 
(\qff I\qff \smallsetminus\qff C\trf)
\off\qff =\off\qff
\sum\nolimits_{\qff C}\off
\sum\nolimits_{\qff D}\off
\fclass{D}\dff *\dff 
(\qff I\qff \smallsetminus\qff C\trf)
\qff,
\end{equation}

\vspace{-6pt}
where\sss $D$\sss runs over\dss non-empty\sss subsets of\trs $I$\dss and\dss
$C$ runs over subsets of\dss $I$ such\dss that\dss 
$D\qff \subset\pff C$\dss and\dss
$\num{C\fff}\off =\off \num{D}\qff +\qff 1$\nnsp,\oss
or\halfff,\oss what\dss is\dss the same,\pss
$\num{D}\off =\off \num{C\fff}\qff -\qff 1$\nnsp.\oss
The sums\qss (\ref{first-sum})\qss and\qss (\ref{second-sum})\qss
differ only\dss in\dss that\dss in\qss (\ref{first-sum})\qss
the set\dss $D$\sss is\dss allowed\dss to be empty\halfff.\oss
Therefore in\dss the sum of\pss  (\ref{first-sum})\qss and\qss (\ref{second-sum})\qss
all\dss terms with non-empty\dss $D$\sss cancel.\oss
Hence\qss (\ref{enlarged-fclass-bd})\qss implies\dss that\vspace{4.5pt}
\[
\quad
\partial\qff
\mathcal{E}\dff\fclass{\trf I\trf}
\off\qff =\qff\off
\sum\nolimits_{\qff C}\off
\fclass{\varnothing}\dff *\dff 
(\qff I\qff \smallsetminus\qff C\trf)
\]

\vspace{-33pt}
\[
\quad
\phantom{\partial\qff
\mathcal{E}\dff\fclass{\trf I\trf}
\off\qff }=\off\qff
\sum\nolimits_{\qff C}\off
\varnothing\dff *\dff 
(\qff I\qff \smallsetminus\qff C\trf)
\off\qff =\off\qff
\sum\nolimits_{\qff C}\off
I\qff \smallsetminus\qff C
\off\qff =\off\qff
\partial\trf I
\qff,
\]

\vspace{-7.5pt}
where\sss $C$\sss runs over one-element\sss subsets of\trs $I$\nnsp.\oss
It\dss follows\dss that\vspace{4.5pt}
\[
\quad
\partial\qff
\mathcal{E}\dff\fclass{\trf I\trf}
\off\qff =\qff\off
\sum\nolimits_{\qff D}\off D
\off\qff =\off\qff
\sum\nolimits_{\qff D}\off \mathcal{E}\dff\fclass{D}
\qff,
\]

\vspace{-7.5pt}
where\sss $D$\sss runs over\sss subsets of\trs $I$\sss
such\dss that\dss
$\num{D}\off =\off \num{\trf I\trf}\qff -\qff 1$\nnsp.\oss
This proves\qss (\ref{chain-simplex-def})\qss for\qss $A\off =\off I$\qss
and $\mathcal{E}$ in\dss the role of\dss $\mathcal{D}$\nnsp.\oss
The\sss theorem\dss follows.\oss  \eproof

\myuppar{Classical\dss colorings.}
Recall\dss that\sss $V_{\fff \mathcal{D}}$\qss is\dss
the union of\trs the sets of\dss vertices of\dss complexes\dss
$\mathcal{D}\dff(\dff A\dff)$\dss with\dss $A\qff \subset\pff I$\nnsp.\oss
A\qss \emph{classical\dss coloring}\qss of\trs $\mathcal{D}$\dss 
is\dss a simply\sss a map\qss
$V_{\fff \mathcal{D}}\qff \ttoo\qff I$\nnsp.\oss
A classical\sss coloring of\trs $\mathcal{D}$\dss induces
for every\trs $A\qff \subset\pff I$\dss a map from\dss the set\sss
of\dss vertices of\trs $\mathcal{D}\dff(\dff A\dff)$\nnsp,\pss
which we will\sss consider as a simplicial\dss map\dss
$\mathcal{D}\dff(\dff A\dff)\qff \ttoo\qff \Delta\dff(\trf I\trf)$\dnsp.\oss

A classical\sss coloring $c$ is\dss
called 
an\qss \emph{Alexander--Sperner\dss coloring}\pss  
if\dss $c$\sss maps\dss $\mathcal{D}\dff(\dff A\dff)$\dss to\dss $\Delta\dff(\dff A\dff)$\dss
for\sss every\dss $A\qff \subset\pff I$\nnsp,\oss
i.e.\pss if\dss
$c\dff(\dff v\trf)$\dss belongs\sss to\sss $A$\sss for every\dss
subset\dss 
$A\qff \subset\pff I$\qss and
every\dss vertex\sss $v$\dss of\qss $\mathcal{D} (\trf A\trf)$\nnsp.\oss 
A classical\sss coloring\sss $c$ is\dss called\dss a\qss 
\emph{Scarf\qss coloring}\oss 
if\qss
$c\dff(\dff v\trf)\off \not\in\off A$\pss for every\dss
proper subset\dss 
$A$\dss of\trs $I$\qss and
every vertex\sss $v$\dss of\qss $\mathcal{D}\dff (\dff A\trf)$\nnsp.\oss
Cf\halfff.\pss \cite{sc4},\oss Lemma\qss 3.5.\oss

It\dss may\dss happen\dss that\dss there are\sss no\dss Alexander--Sperner\dss colorings.\oss
For example,\oss if\trs a vertex\dss $v$\dss is\dss
a vertex of\dss complexes
$\mathcal{D}\dff (\fff A\trf)$ and\sss $\mathcal{D}\dff (\dff A\fff'\trf)$\dss
for\sss disjoint\sss subsets\qss $A\fff,\pff A'\qff \subset\pff I$\nnsp,\oss
then\dss there are\sss no\dss Alexander--Sperner\dss colorings.\oss
Similarly\halfff,\oss it\dss may\dss happen\dss that\trs there are no\qss Scarf\qss colorings.\oss

\mypar{Theorem.}{as-colorings}
\emph{Let\dss $c$\sss be\dss an\dss Alexander--Sperner\dss coloring\dss
of\dss a\dss chain-simplex\dss $\mathcal{D}$\dss based\sss on\dss $I$\nnsp.\oss
Then\qss $c\trf(\dff \sigma\dff)\off =\off I$\pss for\dss
some\dss $d\dff(\trf I\trf)$\dnsp-simplex\sss $\sigma$ of\pss $\mathcal{D}\dff(\trf I\trf)$\dss
and\dss the\sss number\dss of\qss such\sss $\sigma$ 
is\dss odd.}

\proof
We will\dss use\sss the abbreviated\dss notation\dss
$\fclass{A\fff}\off =\off \mathcal{D}\dff\fclass{A\fff}$\nnsp.\oss
It\dss is sufficient\dss to prove\sss that\qss
$c_{\fff *}\dff 
\fclass{\trf I\trf}
\off =\off
I$\nnsp.\oss
Using an\dss induction\dss by\qss $n\off =\off \num{\fff I\fff}\qff -\qff 1$\nnsp,\oss we can assume\sss
that\dss
$c_{\fff *}\dff\fclass{\trf B\trf}
\off =\off
B$\dss
for every $n$\dnsp-subset\trs $B\qff \subset\pff I$\nnsp.\oss
Then\qss (\ref{chain-simplex-def})\qss implies\sss that\vspace{3pt}
\[
\quad
\partial\trf c_{\fff *}\dff\fclass{\trf I\trf}
\off\dff =\off\qff
c_{\fff *}\fff\bigl(\trf \partial\trf \fclass{\trf I\trf}\trf\bigr)
\off\dff =\off\dff
c_{\fff *}\fff\biggl(\qff \sum\nolimits\pff\fclass{\trf B\trf}\qff\biggr)
\]

\vspace{-30pt}
\[
\quad
\phantom{\partial\trf c_{\fff *}\dff\fclass{\trf I\trf}
\off\dff =\off\qff
c_{\fff *}\fff\bigl(\trf \partial\trf \fclass{\trf I\trf}\trf\bigr)
\off\dff }
=\off\qff
\sum\nolimits\qff c_{\fff *}\dff\fclass{\trf B\trf}
\off\dff =\off\qff
\sum\nolimits\pff B
\off\qff =\off\qff
\partial\qff I
\qff,
\]

\vspace{-6pt}
where\sss the sums\sss are\dss taken over\dss all\sss 
$n$\dnsp-subsets\dss $B$\dss of\dss $I$\nnsp.\oss
Since\sss there is only\sss one $n$\dnsp-chain of\qss
$\Delta\dff(\trf I\trf)$\dss with\dss the boundary\sss equal\dss to\dss
$\partial\dff I$\nnsp,\oss
namely\sss $I$\nnsp,\oss
it\dss follows\sss that\qss
$c_{\fff *}\dff\fclass{\trf I\trf}
\off =\off\dff
I$\nnsp.\oss  \eproof

\myuppar{Extensions of\dss classical\sss colorings.}
Let\trs
$c\dff \colon\dss
V_{\fff \mathcal{D}}\qff \ttoo\qff I$\dss
be\dss a\sss map
and\dss let\dss
$\varphi\qff \colon\dff
V_{\fff \mathcal{D}}\qff \cup\pff I
\qff \ttoo\qff I$\dss
be a map
equal\dss to $c$ on\sss $V_{\fff \mathcal{D}}$\sss
and\dss inducing a\sss bijection\dss $I\qff \ttoo\qff I$\nnsp.\oss
If\dss $\mathcal{E}$\dss is\dss the envelope\sss of\trs $\mathcal{D}$\dss as above,\oss
then\dss $\varphi$ can\dss be considered as
a simplicial\dss map\qss
$\mathcal{E}\dff (\trf I\trf)\qff \ttoo\qff \Delta\dff (\trf I\trf)$\qss
equal\dss to\sss $c$ on each\sss complex\dss 
$\mathcal{D}\dff (\dff A\dff)$\nnsp,\oss where\dss $A\qff \subset\pff I$\nnsp,\oss 
and\dss inducing an automorphism\dss
of\trs the boundary\dss
$\bd\dff \Delta\dff(\trf I\trf)$\nnsp.\oss

\mypar{Lemma.}{d-of-image}
\emph{If\qss $\mathcal{D}$\sss is\dss a\sss chain-simplex,\oss
then\dss
$\varphi_{\dff *}\dff\bigl(\qff \mathcal{E}\dff\fclass{\trf I\trf}\qff\bigr)
\off\dff =\off\qff
I$\nnsp.\oss}

\proof
By\qss Theorem\qss \ref{envelope-chain}\qss the envelope\dss $\mathcal{E}$\dss is\dss also\sss a\sss 
chain-simplex.\oss
By\dss the definition of\trs the envelope,\pss 
$\mathcal{E}\dff\fclass{\trf B\trf}\off =\off B$\dss
for every\dss proper subset\dss $B\qff \subset\qff I$\nnsp.\oss
It\dss follows\dss that\dss 
$\partial\trf \mathcal{E}\dff\fclass{\trf I\trf}
\off =\off
\partial\trf I$\dss and\vspace{4.5pt} 
\[
\quad
\partial\qff \varphi_{\dff *}\dff\bigl(\qff \mathcal{E}\dff\fclass{\trf I\trf}\qff\bigr)
\off\dff =\off\qff
\varphi_{\dff *}\dff\bigl(\qff \partial\trf \mathcal{E}\dff\fclass{\trf I\trf}\qff\bigr)
\off\dff =\off\qff
\varphi_{\dff *}\dff(\trf \partial\trf I \qff)
\qff.
\]

\vspace{-7.5pt}
But\sss $\varphi$ is\dss equal\dss to\sss
the identity\sss on
$\partial\dff \Delta\dff(\trf I\trf)$
and\dss hence\dss
$\varphi_{\dff *}\dff(\trf \partial\trf I \qff)
\off =\off
\partial\trf I$\nnsp.\oss
It\dss follows\dss that\qss\vspace{4.5pt} 
\[
\quad
\partial\qff \varphi_{\dff *}\dff\bigl(\qff \mathcal{E}\dff\fclass{\trf I\trf}\qff\bigr)
\off\dff =\off\qff
\partial\trf I
\qff.
\]

\vspace{-7.5pt}
Since\sss the simplex $I$ is\dss the only $d\dff(\trf I\trf)$\dnsp-chain\dss of\dss
$\Delta\dff(\trf I\trf)$ with\sss the\sss boundary\sss
$\partial\dff I$\nnsp,\pss
this implies\sss that\trs
$\varphi_{\dff *}\dff\bigl(\qff \mathcal{E}\dff\fclass{\trf I\trf}\qff\bigr)
\pff =\off
I$\nnsp.\oss 
This completes\sss the proof\halfff.\oss  \eproof

\mypar{Theorem.}{scarf-arbitrary-colorings}
\emph{If\qss $\mathcal{D}$\sss is\dss a\sss chain-simplex,\oss
then\dss for every\sss classical\sss coloring\qss $c\dff \colon\dff V_{\fff \mathcal{D}}\qff \ttoo\qff I$\qss 
there\dss exist\dss a\sss non-empty\dss subset\pss $A\pff \subset\off I$\qss
and\dss a\sss simplex\dss $\sigma$\sss of\oss $\mathcal{D}\dff (\dff A\dff)$\dss
such\dss that\qss
$c\trf(\dff \sigma\dff)\off =\off A$\dnsp.\oss
Moreover\halfff,\oss the number of\dss such\dss pairs $A\fff,\pff \sigma$ is\dss odd.\oss}

\proof
Lemma\qss \ref{d-of-image}\qss implies\sss that\trs
$\varphi\dff(\trf \rho\dff)
\off =\off
I$\dss
for some $d\dff(\trf I\trf)$\dnsp-simplex\sss $\rho$\sss of\qss $\mathcal{E}\dff (\trf I\trf)$\nnsp.\oss
By\qss Lemma\qss \ref{top-simplices}\qss the simplex\sss $\rho$\sss
has\sss the form\qss
$\rho
\off =\off
\sigma\dff *\dff (\qff I\pff \smallsetminus\qff A\trf)$\qss
for a non-empty\dss subset\qss $A\qff \subset\pff I$\qss
and\sss a\sss $d\dff(\dff A\dff)$\dnsp-simplex\sss $\sigma$\sss of\oss 
$\mathcal{D}\dff (\dff A\dff)$\nnsp.\oss
Since\dss 
$\num{\rho} 
\off =\off \num{I}$\dss
and\qss $\varphi\dff(\trf \rho\dff)
\off =\off
I$\nnsp,\oss
the map\qss $\rho\qff \ttoo\qff I$\qss induced\dss by\sss $\varphi$\sss
is\dss a\sss bijection.\oss
Since $\varphi$ is\dss equal\dss to\sss the identity\sss on\dss $I$\nnsp,\oss
this implies\sss that\qss 
$\varphi\dff(\dff \sigma\dff)
\off =\off
A$\nnsp.\oss
In\dss turn,\oss since $\varphi$ is equal\dss to\sss $c$\sss on\dss 
$\mathcal{D}\dff (\dff A\dff)$\nnsp,\oss
this implies\sss that\qss 
$c\trf(\dff \sigma\dff)
\off =\off
A$\nnsp.\oss
In\dss particular\halfff,\oss there exist\sss simplices $\sigma$\sss with\dss required\dss
properties.\oss

Conversely\halfff,\oss if\trs $A$ is\dss a\sss non-empty\sss
subset\sss of\trs $I$\dss and $\sigma$ is\dss a\sss simplex of\trs the complex\sss
$\mathcal{D}\dff(\dff A\dff)$\sss such\dss that\qss
$c\trf(\dff \sigma\dff)\off =\off A$\nnsp,\oss
then\dss
$\rho
\off =\off
\sigma\dff *\dff (\qff I\pff \smallsetminus\qff A\trf)$\qss
is\dss a\sss $d\dff(\trf I\trf)$\dnsp-simplex of\oss 
$\mathcal{E}\dff (\trf I\trf)$\dss
and\dss
$\varphi\dff(\trf \rho\dff)
\off =\off
I$\nnsp.\oss
Therefore\sss there\dss is\dss a\sss one-to-one correspondence\sss
between\dss pairs\dss $A\fff,\pff \sigma$\dss as in\dss the\sss theorem
and $d\dff(\trf I\trf)$\dnsp-simplices $\rho$ of\oss 
$\mathcal{E}\dff (\trf I\trf)$\dss
such\dss that\trs
$\varphi\dff(\trf \rho\dff)
\off =\off
I$\nnsp.\oss
But\qss $\varphi_{\dff *}\dff\bigl(\qff \mathcal{E}\dff\fclass{\trf I\trf}\qff\bigr)
\pff =\off
I$\qss implies\sss that\dss the number of\dss such simplices $\rho$ is\dss odd.\oss
It\dss follows\dss that\dss the number of\dss 
pairs\dss $A\fff,\pff \sigma$\dss as in\dss the\sss theorem\dss
is\dss also\sss odd.\oss
This completes\sss the proof\halfff.\oss  \eproof

\mypar{Theorem.}{scarf-scarf-colorings}
\emph{If\qss $\mathcal{D}$\sss is\dss a\sss chain-simplex\sss
and\qss $c\dff \colon\dff V_{\fff \mathcal{D}}\qff \ttoo\qff I$\qss is\dss a\qss
Scarf\dss coloring\fff,\oss
then\qss $c\dff(\dff \sigma\dff)\off =\off I$\qss for\dss
some\qss simplex\dss $\sigma$\sss of\oss $\mathcal{D}\dff (\trf I\trf)$
and\qss the\sss number\dss of\qss such\sss simplices $\sigma$ is\dss odd.\oss}

\proof
By\qss Lemma\qss \ref{top-simplices}\qss every $d\dff(\trf I\trf)$\dnsp-simplex\sss $\rho$\sss
of\trs $\mathcal{E}\dff (\trf I\trf)$\dss
has\sss the form\qss
$\rho
\off =\off
\sigma\dff *\dff (\qff I\pff \smallsetminus\qff A\trf)$\qss
for a non-empty\dss subset\qss $A\qff \subset\pff I$\qss
and\sss a\sss $d\dff(\dff A\dff)$\dnsp-simplex\sss $\sigma$\sss 
of\oss $\mathcal{D}\dff (\dff A\dff)$\nnsp.\oss
If\dss $\rho$\dss is\dss not\sss a simplex of\dss
$\mathcal{D}\dff (\trf I\trf)$\nnsp,\oss
then $A$ is a proper subset\sss of\trs $I$\nnsp.\oss
Since $c$ is\dss a\dss Scarf\dss coloring\halfff,\oss
this implies\sss that\dss
$c\dff(\dff \sigma\dff)$\dss is disjoint\dss from\dss
$A$\nnsp,\oss
i.e.\qss
$c\dff(\dff \sigma\dff)\qff \subset\pff I\qff \smallsetminus\qff A$\nnsp.\oss
Since\sss $\varphi$\sss is\dss equal\dss to\sss
the identity\sss on\dss $I$\nnsp,\oss
it\dss follows\sss that\vspace{3pt}
\[
\quad
\varphi\dff(\dff \rho \dff)
\off =\off
\varphi\dff(\qff \sigma\dff *\dff (\qff I\pff \smallsetminus\qff A\trf) \qff)
\off \subset\dff\off 
I\pff \smallsetminus\qff A
\qff.
\]

\vspace{-9pt}
But\sss $A$ is\dss non-empty\sss
and\dss hence\dss $\varphi\dff(\trf \rho \dff)$\dss
is\dss an $m$\dnsp-simplex 
for\sss some\qss $m\qff \leq\pff d\dff(\trf I\trf)\qff -\qff 1$\nnsp.\oss
Therefore\qss
$\varphi_{\dff *}\fff(\dff \rho \dff)
\off =\off
0$\qss
for every\sss $d\dff(\trf I\trf)$\dnsp-simplex $\rho$
which\dss is\dss not\sss a\sss simplex of\dss 
$\mathcal{D}\dff (\trf I\trf)$\nnsp.\oss
It\dss follows\dss that\vspace{4.5pt}
\[
\quad
\varphi_{\dff *}\dff\bigl(\qff \mathcal{D}\dff\fclass{\trf I\trf}\qff\bigr)
\off =\off
\varphi_{\dff *}\dff\bigl(\qff \mathcal{E}\dff\fclass{\trf I\trf}\qff\bigr)
\qff.
\]

\vspace{-7.5pt}
Together\dss with\qss Lemma\qss \ref{d-of-image}\qss
this implies\dss that\qss
$\dis
c_{\dff *}\dff\bigl(\qff \mathcal{D}\dff\fclass{\trf I\trf}\qff\bigr)
\off\dff =\off\qff 
\varphi_{\dff *}\dff\bigl(\qff \mathcal{D}\dff\fclass{\trf I\trf}\qff\bigr)
\off\dff =\off\qff 
I$\nnsp.\oss
Now\dss the\sss theorem\dss follows\dss by\sss a\sss standard\sss argument\halfff.\oss  \eproof

\prooftitle{Another proof\qss of\pss Theorem\qss \ref{as-colorings}}
In contrast\dss with\dss the first\dss proof\halfff,\oss
this proof\dss does not\sss rely\sss on\dss induction\dss by\dss $\num{I}$\nnsp.\oss
Let\dss $f\dff \colon\dff I\qff \ttoo\qff I$\dss be a bijection such\dss that\trs
$f\dff(\dff A\dff)\off \neq\off A$\dss for every\trs $A\qff \subset\pff I$\dss
different\dss from\dss $\varnothing\fff,\pff I$\nnsp.\oss
Equivalently\halfff,\oss let\sss $f$\dss be a cyclic permutation of\trs $I$\dss
(i.e.\qss a permutation consisting of\dss one cycle).\oss
Let\sss $\varphi$\sss be\sss the extension of\dss $c$\sss equal\dss to $f$ on\dss $I$\nnsp.\oss
By\qss Lemma\qss \ref{top-simplices}\qss every\sss
$d\dff(\trf I\trf)$\dnsp-simplex $\rho$ of\dss $\mathcal{E}\dff (\trf I\trf)$
has\sss the form
$\rho
\off =\off
\sigma\dff *\dff (\qff I\pff \smallsetminus\qff A\trf)$
for a non-empty\sss subset\sss $A$ of\trs $I$\dss and a simplex $\sigma$
of\trs $\mathcal{D}\dff (\dff A\dff)$\nnsp.\oss
Suppose\sss that\dss
$\varphi\dff(\trf \rho\dff)
\off =\off
I$\nnsp.\oss
Since $c$ is\dss an\dss Alexander--Sperner\dss coloring,\pss 
$\varphi\dff(\dff \sigma\dff)\qff \subset\qff A$\dss
and\dss hence\dss
$\varphi\dff(\trf \rho\dff)
\off =\off
I$\dss
implies\sss that\dss 
$\varphi\dff(\qff I\pff \smallsetminus\qff A\trf)
\off =\off
I\pff \smallsetminus\qff A$\nnsp.\oss
Since $A$ is\dss non-empty\halfff,\oss by\dss the choice of\dss $\varphi$\sss
this\dss is\dss possible only\dss if\trs $A\off =\off I$\nnsp.\oss
In\dss this case\dss $\rho\off =\off \sigma$\dss and\sss $\sigma$ is\dss a\sss simplex
of\trs $\mathcal{D}\dff (\trf I\trf)$\nnsp.\oss
Therefore $d\dff(\trf I\trf)$\dnsp-simplex\sss $\rho$\sss of\qss $\mathcal{E}\dff (\trf I\trf)$\sss
is\dss actually\sss a simplex of\trs $\mathcal{D}\dff (\trf I\trf)$\nnsp.\oss
By\sss a standard\sss argument\qss Lemma\qss \ref{d-of-image}\qss
implies\sss that\dss the number of\dss such simplices $\rho$ is\dss odd.\oss  \eproof

\newpage
\mysection{Combinatorics\qss of\pss families\qss of\pss linear\qss orders}{orders}

\myuppar{Linear orders\dss and\dss dominant sets.}
Let $X$ be a non-empty\dss finite set\halfff.\oss
Suppose that a family of\dss linear orders $<_{\dff i}$ on\dss $X$\nnsp,\oss
labeled\dss by elements $i$ of\dss a\sss finite set $I$\nnsp,\pss
is given.\oss
For a non-empty\dss subset\qss $\sigma\qff \subset\qff X$\qss let\qss
$\min_{\dff i}\qff \sigma$\qss be the minimal element of $\sigma$ 
with respect to the order\dss $<_{\dff i}$\nnsp.\oss
A subset\qss $\sigma\qff \subset\qff X$\qss is\dss said\dss to\dss be\qss 
\emph{dominant\dss with\dss respect\dss to}\qss
a non-empty subset\dss 
$C$\dss of\qss $I$\qss if\trs\vspace{3pt} 
\begin{equation*}
\quad
\mbox{there\dss is\dss no\dss element}\qff\oss 
y\qff \in\qff X\qff\oss 
\mbox{such\dss that}\qff\oss 
\min\nolimits_{\dff i}\qff \sigma\off <_{\dff i}\off y\qff\oss 
\mbox{for\dss all}\qff\oss
i\qff \in\qff C
\dff.
\end{equation*}

\vspace{-9pt}
By\dss a convention,\qss $\varnothing$\dss
is dominant\dss with respect to every\qss $C\off \neq\off \varnothing$\dnsp,\oss
but\dss not\dss with respect to\dss $\varnothing$\dnsp.

\mypar{Lemma.}{minima}
\emph{If\oss $\sigma\qff \subset\pff X$\qss is dominant with respect to\pss $C\qff \subset\off I$\nnsp,\oss
then}\oss
$\dis
\sigma
\off =\off
\bigl\{\off \min\nolimits_{\dff i}\dff \sigma \off 
\bigl|\bigr.\off 
i\qff \in\qff C \off\bigr\}$.

\proof
Clearly\halfff,\pss $\sigma$ contains all minima\dss $\min_{\fff i}\dff \sigma$\nnsp.\oss
Suppose that\qss $x\qff \in\qff \sigma$\qss is different from all\qss
$\min_{\dff i}\dff \sigma$\qss with\qss $i\qff \in\qff C$\nnsp.\oss
Then\qss $\min_{\dff i}\dff \sigma\qff <_{\dff i}\qff x$\qss for all\qss $i\qff \in\qff C$\nnsp,\oss
contrary to the assumption.\oss  \eproof

\mypar{Corollary\halfff.}{surjectivity}
\emph{
If\dff\oss $\sigma\qff \subset\pff X$\pss 
is dominant with respect to\pss 
$C\qff \subset\pff I$\nnsp,\qff\oss
then\oss $\num{\sigma}\qff \leq\qff \num{C}$\nnsp.\oss}  \eproof

\myuppar{Deleting and adding elements.}
For a set $A$ and an element\qss $a\qff \in\qff A$\qss
we will denote by\qss
$A\qff -\qff a$\qss the set\qss 
$A\qff \smallsetminus\qff \{\qff a \qff\}$\nnsp.\oss
Sim\-i\-lar\-ly\halfff,\oss for\qss $b\qff \not\in A$\qss
we will denote by\qss
$A\qff +\qff b$\qss
the set\qss
$A\qff \cup\qff \{\qff b \qff\}$\nnsp.\oss
The set\qss $A\qff -\qff a$\qss is defined only if\qss $a\qff \in\qff A$\nnsp,\oss
and the set\qss $A\qff +\qff b$\qss is defined only if\qss $b\qff \not\in\qff A$\nnsp.\oss
We\dss will\dss interpret\qss $A\qff -\qff a\qff +\qff b$\qss
as\qss $(\dff A\qff -\qff a \dff)\qff +\qff b$\nnsp,\oss
and\dss interpret\sss similar\sss expressions in\dss the same way\halfff.\oss

\myuppar{Cells and\dss faces.}
Let\dss us\sss fix a subset\qss $C\qff \subset\pff I$\qss
and\sss suppose\sss that\dss $\sigma$\dss is a\sss subset\sss of\qss $X$\nnsp.\oss
We\sss will\dss say\dss that\dss $\sigma$\dss is\dss a\qss \emph{$C$\dnsp-cell}\oss if\qss
$\sigma$\dss is dominant\dss with respect\dss to\dss $C$\dss and\qss
$\num{C}\off =\off \num{\sigma}$\nnsp,\oss
and\dss that\dss $\sigma$\dss is a\qss \emph{$C$\dnsp-face}\oss if\qss
$\sigma$\dss is dominant\dss with respect\dss to\dss $C$\dss and\qss
$\num{C}\off =\off \num{\sigma}\qff +\qff 1$\nnsp.\oss

Let\qss $m_{\trf \sigma}\dff \colon\dff C\qff \ttoo\qff \sigma$\qss 
be\sss the map defined\dss by\dss the rule\vspace{3pt}
\[
\quad
i\off \longmapsto\off \min\nolimits_{\dff i}\qff \sigma
\qff.
\]

\vspace{-9pt}
By\qss Lemma\qss \ref{minima},\oss if\qss $\sigma$\dss is dominant\dss with\sss respect\dss
to\dss $C$\nnsp,\oss then\dss $m_{\trf \sigma}$\dss is a surjection.\oss
Therefore,\oss if\qss $\sigma$\dss is a\dss $C$\dnsp-cell,\oss
then\dss $m_{\trf \sigma}$\dss is a bijection and\dss
and\dss for every\qss $x\qff \in\qff \sigma$\qss
there\dss is\dss a unique element\qss $i\qff \in\qff C$\qss
such\dss that\qss
$x\off =\off \min_{\dff i}\qff \sigma$\nnsp.\oss
We will denote\sss this element\dss $i$\dss by\dss $i_{\trf \sigma}\dff(\dff x\trf)$\nnsp.\pss
In other\dss words,\oss if\qss $\sigma$\dss is a\dss $C$\dnsp-cell,\oss
then\dss $i_{\trf \sigma}$\dss is\dss the inverse of\qss $m_{\trf \sigma}$\nnsp.\oss

\myuppar{From cells\sss to faces.}
Clearly,\oss if\qss $\sigma$\dss is a\dss $C$\dnsp-cell\sss and\qss
$i\qff \not\in\qff C$\nnsp,\oss
then\dss $\sigma$\dss is also a\dss $(\dff C\dff +\dff i\trf)$\dnsp-face.\qff\oss
Removing an element\dss from a\dss $C$\dnsp-cell\dss results in a\sss $C$\dnsp-face.\oss
Indeed,\oss suppose\sss that\dss
$\sigma$\dss is\dss a\dss $C$\dnsp-cell\sss and\qss $x\qff \in\qff \sigma$\nnsp.\oss
Since\qss
$\min\nolimits_{\dff i}\qff \sigma
\off \leq_i\off
\min\nolimits_{\dff i}\off 
(\dff
\sigma\qff -\qff x
\trf)$\qss
for every\sss $i$\nnsp,\oss
the set\qss $\sigma\qff -\qff x$\qss is dominant with respect to\dss $C$\dss together with\dss $\sigma$\nnsp.\oss
Also,\oss $\num{C}\off =\off \num{\sigma\qff -\qff x}\qff +\qff 1$\nnsp.\oss
Hence\qss $\sigma\qff -\qff x$\qss is a\dss $C$\dnsp-face.\oss

\myuppar{From\dss faces\sss to cells.}
Let\qss $\sigma$\dss be\dss a\dss non-empty\dss $C$\dnsp-face.\oss
Then\dss there are\sss two elements\qss
$k\fff,\pff l\qff \in\qff C$\qss such\dss that\qss
$m_{\trf \sigma}\dff(\dff k\trf)\off =\off m_{\trf \sigma}\dff(\dff l\trf)$\qss
and\dss $m_{\qff \sigma}$\dss is\dss injective on\qss $C\qff -\qff k\qff -\qff l$\nnsp.\oss
For\sss every\pss $j\qff \in\qff C$\pss let\oss\vspace{3pt}
\[
\quad
\mathbb{M}_{\fff j}
\off\off =\off\off
\bigl\{\off\dff y
\qff \in\qff X 
\off\qff \bigl|\bigr.\off\dff 
\min\nolimits_{\dff i}\dff \sigma
\off <_{\dff i}\off y
\hspace{1.0em}\mbox{for\dss all}\hspace{1.0em}
i\qff \in\qff C\qff -\qff j
\qff\off\bigr\}\dff.
\]

\vspace{-9pt}
If\pss $\mathbb{M}_{\fff j}\off \neq\off \varnothing$\nnsp,\oss
then\dss we will denote by $m_{\fff j}$ the maximal element of\qss $\mathbb{M}_{\fff j}$\qss
with respect\dss to\qss $<_{\fff j}$\nsp.\oss\vspace{1.375pt}

\mypar{Lemma.}{face-up-cell}
\emph{The set\qss
$\sigma\qff +\qff a$\qss is\dss a\qss $C$\dnsp-cell\pss if\pss and\qss only\qss if\oss
$a\off =\off m_{\fff j}$\pss for\dss some\qss
$j\qff \in\qff \{\qff k\fff,\pff l\pff\}$\qss such\dss that\oss
$\mathbb{M}_{\fff j}\off \neq\off \varnothing$\nnsp.\oss}\vspace{1.375pt}

\proof
To begin with,\oss let us observe that\vspace{1.5pt}
\begin{equation}
\label{minimum-1}
\quad
\min\nolimits_{\dff i}\pff (\dff \sigma\qff +\qff a \trf)
\off =\off 
a
\hspace*{5em}\mbox{if}\hspace*{2.3em}
a\off <_{\dff i}\off \min\nolimits_{\dff i}\qff \sigma\dff,
\hspace*{2.5em}\mbox{and}\hspace*{1.5em}
\end{equation}

\vspace*{-37.5pt}
\begin{equation}
\label{minimum-2}
\quad
\min\nolimits_{\dff i}\pff (\dff \sigma\qff +\qff a \trf)
\off =\off 
\min\nolimits_{\dff i}\qff \sigma
\hspace*{2.3em}\mbox{if}\hspace*{2.3em}
\min\nolimits_{\dff i}\qff \sigma
\off <_{\dff i}\off 
a\dff.
\end{equation}

\vspace{-9pt}
In particular\halfff,\oss
$\min\nolimits_{\dff i}\qff (\dff \sigma\qff +\qff a \trf)
\off =\off
\min\nolimits_{\dff i}\qff \sigma$\oss
or\oss
$a$\oss 
for every\oss $i\qff \in\qff C$\dnsp.\qff\oss\vspace{12pt}

Since\dss $\sigma$\dss is\dss a\dss $C$\dnsp-face,\pss
$\sigma$\dss is\dss dominant\dss with\dss respect\dss to\dss $C$\dss and\dss
Lemma\qss \ref{minima}\qss implies that\oss\vspace{1.5pt}
\begin{equation}
\label{minimum-3}
\quad
\bigl\{\off 
\min\nolimits_{\dff i}\qff \sigma
\off \bigl|\bigr.\off 
i\qff \in\qff C 
\off\bigr\}
\off =\off
\sigma
\qff.
\end{equation}

\vspace{-10.5pt}
If\qss $\sigma\qff +\qff a$\qss is\dss a\dss $C$\dnsp-cell,\oss
then\qss $\sigma\qff +\qff a$\qss is\dss also dominant\dss with\dss respect\sss to\dss $C$\dss
and\qss\vspace{1.5pt}
\begin{equation}
\label{minimum-4}
\quad
\bigl\{\off 
\min\nolimits_{\dff i}\pff  (\dff \sigma\qff +\qff a \dff) 
\off \bigl|\bigr.\off 
i\qff \in\qff C 
\off\bigr\}
\off\off
=\off\off
\sigma
\off +\off
a 
\end{equation}

\vspace{-10.5pt}
by\qss Lemma\qss \ref{minima}.\oss\vspace{12pt}

If\pss (\ref{minimum-3})\qss and\qss  (\ref{minimum-4})\qss hold,\dss 
then\oss\vspace{1.5pt}
\[
\quad
\min\nolimits_{\dff i}\qff (\dff \sigma\qff +\qff a \trf)
\off =\off
\min\nolimits_{\dff i}\qff \sigma
\]

\vspace{-10.5pt}
for all\qss
$i\qff \in\qff C\qff \smallsetminus\qff \{\qff k\fff,\pff l \pff\}$\qss
and\dss for $i$ equal\dss to one of\trs the elements
of\dss the pair\qss $\{\qff k\fff,\pff l \pff\}$\nnsp,\oss 
and\vspace{1.5pt}
\[
\quad
\min\nolimits_{\dff i}\qff (\dff \sigma\qff +\qff a \trf)
\off =\off
a
\]

\vspace{-10.5pt}
for $i$ equal\dss to the other element\dss 
of\qss $\{\qff k\fff,\pff l \pff\}$\nnsp.\oss\vspace{12pt}

Therefore,\oss
if\qss $\sigma\qff +\qff a$\qss is\dss a\dss $C$\dnsp-cell,\oss
then\dss we may assume that\vspace{1.5pt}
\begin{equation}
\label{change-of-minima}
\quad
\min\nolimits_{\dff i}\qff (\dff \sigma\qff +\qff a \trf)
\off =\off
\min\nolimits_{\dff i}\qff \sigma
\hspace*{1.2em}\mbox{for\dss all}\hspace*{1.2em}
i\qff \in\qff C\qff -\qff k
\hspace*{1em}\mbox{and}\hspace*{1em}
\min\nolimits_{\dff k}\qff (\dff \sigma\qff +\qff a \trf)
\off =\off
a\dff.
\end{equation}

\vspace{-9pt}
By\qss (\ref{minimum-1})\qss and\qss (\ref{minimum-2})\qss in this case\oss\vspace{4.5pt}
\[
\quad
\min\nolimits_{\dff i}\qff \sigma
\off <_{\dff i}\off 
a
\hspace*{1.2em}\mbox{for\dss all}\hspace*{1.2em}
i\qff \in\qff C\qff -\qff k
\hspace*{1em}\mbox{and}\hspace*{1em}
a
\off <_{\dff k}\off
\min\nolimits_{\dff k}\qff \sigma
\qff.
\]

\vspace{-7.5pt}
It\trs follows that\qss
$a\qff \in\qff \mathbb{M}_{\dff k}$\nsp.\oss
Since\qss $\sigma\qff +\qff a$\qss is\dss dominant\dss with\dss respect\sss to\dss $C$\nnsp,\oss
the element\dss
$a$\dss is\dss the maximal\dss in\dss
$\mathbb{M}_{\dff k}$\dss
with respect\dss to\qss $<_{\dff k}$\nsp.\oss
In other\dss terms,\pss $a\off =\off m_{\dff k}$\nnsp.\oss

Conversely\halfff,\oss
if\halfff,\oss say\halfff,\oss
$\mathbb{M}_{\dff k}\qff \neq\qff \varnothing$\qss
and\pss
$a\qff \in\qff \mathbb{M}_{\dff k}$\nnsp,\oss
then\vspace{4.5pt}
\[
\quad
\min\nolimits_{\dff i}\qff \sigma
\off <_{\dff i}\off 
a
\hspace*{1.2em}\mbox{for\dss all}\hspace*{1.2em}
i\qff \in\qff C\qff -\qff k\dff.
\]

\vspace{-7.5pt}
If\trs also\qss
$\min\nolimits_{\dff k}\qff \sigma
\off <_{\dff k}\off 
a$\nnsp,\oss
then\dss $\sigma$\dss is not dominant with respect\dss to\dss $C$\dnsp,\oss
contrary to the assumption.\oss
Therefore\qss
$a
\off <_{\dff k}\off
\min\nolimits_{\dff k}\qff \sigma$\dnsp.\oss
By applying\qss 
(\ref{minimum-1})\qss and\qss (\ref{minimum-2})\qss
we see that\qss (\ref{change-of-minima})\qss holds.\oss
It\dss follows that\qss
if\oss
$a\off =\off m_{\dff k}$\nsp,\oss
then\qss
$\sigma\qff +\qff a$\qss is dominant with respect\dss to\trs $C$\qss
and\qss $\sigma\qff +\qff a$\qss
is\dss a\dss $C$\dnsp-cell.\oss  \eproof

\mypar{Lemma.}{two-different}
\emph{If\oss $\mathbb{M}_{\dff k}$\dss
and\pss $\mathbb{M}_{\dff l}$\dss are\dss both\dss non-empty\halfff,\oss
then\qss $m_{\dff k}\off \neq\off m_{\dff l}$\qss
and\dss hence\dss the\dss $C$\dnsp-cells\qss
$\sigma\qff +\qff m_{\dff k}$\qss and\oss $\sigma\qff +\qff m_{\dff l}$\qss
are not\dss equal\halfff.\oss}

\proof
Since\dss $\sigma$\dss is dominant\dss with\dss respect\dss to\dss $C$\nnsp,\oss
the sets\dss $\mathbb{M}_{\dff k}$\dss and\dss $\mathbb{M}_{\dff l}$\dss
are disjoint\halfff.\oss
The lemma follows.\oss  \eproof

\mypar{Lemma.}{face-down-cell}
\emph{The\dss set\dss $\sigma$\sss is\qss a\qss $(\qff C\qff -\dff j\trf)$\dnsp-cell\pss 
if\pss and\qss only\qss if\oss
$j\qff \in\qff \{\qff k\fff,\pff l\pff\}$\qss and\oss
$\mathbb{M}_{\dff j}\off =\off \varnothing$\nsp.}

\proof
Since\dss $\sigma$\dss is a\dss $C$\dnsp-face,\oss
$\num{\sigma}\off =\off \num{C}\qff -\qff 1\off =\off \num{C\qff -\qff j}$\nnsp.\oss
The condition\qss
$\mathbb{M}_{\dff j}\off =\off \varnothing$\qss
is\dss obviously equivalent\dss to\dss $\sigma$\dss being dominant\sss with\dss
respect\dss to\qss $C\qff -\qff j$\nnsp.\oss
Therefore\dss it\dss remains\sss to show\dss that\qss
$\mathbb{M}_{\dff j}\off \neq\off \varnothing$\qss
if\pss
$j\qff \not\in\qff \{\qff k\fff,\pff l\pff\}$\nnsp.\oss
But\qss if\qss
$j\qff \not\in\pff \{\qff k\fff,\pff l\pff\}$\nnsp,\oss
then\vspace{3pt}
\[
\quad
\min\nolimits_{\dff i}\dff \sigma
\off \neq\off
\min\nolimits_{\dff j}\dff \sigma
\]

\vspace{-9pt}
for every\qss $i\off \neq\off j$\nnsp.\qff\oss
It\dss follows\sss that\dss in\dss this case\qss\vspace{3pt}\vspace{-0.125pt}
\[
\quad
\min\nolimits_{\dff i}\dff \sigma
\off <_{\fff i}\off
\min\nolimits_{\dff j}\dff \sigma
\]

\vspace{-9pt}
for every\qss $i\off \neq\off j$\qss
and\dss hence\qss
$\min\nolimits_{\dff j}\dff \sigma\qff \in\qff 
\mathbb{M}_{\dff j}$\nsp.\oss  \eproof

\mypar{Theorem\qss ({\dff}The main combinatorial\dss lemma).}{faces-cells}
\emph{If\qss $\sigma$ is\dss a\dss non-empty\trs $C$\dnsp-face,\oss
then $\sigma$ is contained\dss in\dss no more\sss than\dss two\dss $C$\dnsp-cells.\qff\oss
If\qss $r$\dss is\dss the\dss
number\dss of\oss $C$\dnsp-cells\dss containing\sss $\sigma$\sss
and\sss $s$\sss is\dss the number of\qss subsets\qss
$D\qff \subset\pff C$\qss
such\dss that\dss $\sigma$\sss is\dss a\pss $D$\dnsp-cell,\oss
then\qss $r\qff +\qff s\off =\off 2$\nnsp.\oss}

\proof
This follows from\dss Lemmas\qss \ref{face-up-cell}\qss --\qss \ref{face-down-cell}.\oss
The number\dss $r$\dss is\sss equal\dss to\sss the number of\dss non-empty\sss sets\sss
among\qss
$\mathbb{M}_{\dff k}$\qss
and\qss $\mathbb{M}_{\dff l}$\nsp.\oss  \eproof

\myuppar{The pseudo-simplex associated\dss with a family of\dss orders.}
Theorem\qss \ref{faces-cells}\qss is\dss an analogue of\trs the
non-branching\dss property of\trs triangulations of\dss geometric simplices.\oss
Let\dss us\dss introduce an analogue of\trs the pseudo-simplices\sss defined\dss
by\dss triangulations of\trs geometric simplices\qss 
(see\dss Section\qss \ref{pseudo-simplices}).\oss
For every subset\qss $A\qff \subset\pff I$\qss
let\dss $\mathcal{T} (\dff A\dff)$\dss
be an abstract\sss simplicial complex\sss defined as follows.\oss
The set\sss of\dss vertices of\dss $\mathcal{T} (\dff A\dff)$\dss
is\dss the union of\dss all $A$\dnsp-cells,\pss
and a set\sss of\dss vertices\dss is\dss a\sss simplex\dss if\dss and\dss only\dss if\dss
it\dss is contained\dss in\sss a\dss $A$\dnsp-cell.\oss
Clearly\halfff,\oss
the map\qss
$A\off \longmapsto\off \mathcal{T}(\dff A\dff)$\dss
is\dss a\sss simplex-family\halfff.\oss

\mypar{Theorem.}{scarf-pseudo-simplex}
\emph{The\dss simplex-family\qss $\mathcal{T}$ is\qss a\dss pseudo-simplex.\oss}

\proof
Suppose\sss that\sss $A$\sss be a subset\sss of\trs $I$\nnsp.\oss
Let\dss us\dss prove\qss (\ref{sum-2})\qss
for\qss $\mathcal{D}\off =\off \mathcal{T}$\dnsp\dnsp.\oss
Let\dss us\dss consider some $e\dff(\dff A\dff)$\dnsp-simplex $\sigma$ of\qss
$\mathcal{T} (\dff A\dff)$\nnsp.\oss
By\dss the definition\dss the simplex $\sigma$ is\dss contained\dss in some\dss $A$\dnsp-cell,\oss
i.e.\dss $\sigma\qff +\qff a$\dss is\dss an\dss $A$\dnsp-cell\dss
for some $a$\nnsp.\oss
In\dss particular\halfff,\pss $\sigma$ is\dss an $A$\dnsp-face.\pss
By\trs Theorem\qss \ref{faces-cells}\qss
either $\sigma$ is\dss a face of\dss exactly\dss two $d\dff(\dff A\dff)$\dnsp-simplices,\oss
or $\sigma$ is\dss a face of\dss exactly one $d\dff(\dff A\dff)$\dnsp-simplex and\dss
is\dss a $B$\dnsp-cell\dss for exactly\sss one
subset $B\qff \subset\qff A$\nnsp.\oss
This\sss implies\dss that $\mathcal{T}$ is\dss a\sss pseudo-simplex.\oss  \eproof

\mypar{Corollary.}{scarf-envelope-pseudomanifold}
\emph{Let\qss $\mathcal{S}$\dss be\sss the envelope of\pss $\mathcal{T}$\nnsp.\oss
Then\qss $\mathcal{S}\dff(\trf I\trf)$\dss is\dss a\sss non-branching\dss
dimensionally\dss homogeneous simplicial\dss complex.\oss}

\proof
The envelope $\mathcal{S}\dff(\trf I\trf)$ is\dss dimensionally\dss
homogenous by\dss the construction.\oss
Theorem\qss \ref{envelope-pseudo}\qss implies\sss that\sss
$\mathcal{S}\dff(\trf I\trf)$ is\dss non-branching.\oss  \eproof

\mypar{Corollary.}{scarf-chain-simplex}
\emph{The\dss simplex-family\qss $\mathcal{T}$ is\qss a\dss chain-simplex.\oss}

\proof
It\dss is\dss sufficient\dss to combine\qss Theorem\qss \ref{scarf-pseudo-simplex}\qss
and\qss Lemma\qss \ref{pseudo-is-chain}.\oss  \eproof

\myuppar{Classical\sss colorings.}
In\dss the present\sss context\dss the definitions of\dss colorings from\dss
Section\qss \ref{pseudo-simplices}\qss take\sss the following\dss form.\oss
A\qss \emph{classical\dss coloring}\qss 
is\dss a map\qss
$c\dff \colon\dff
X\qff \ttoo\qff I$\nnsp.\oss 
The map $c$ is\dss
called 
an\qss \emph{Alexander--Sperner\dss coloring}\pss  
if\dss $c$\sss maps\dss $\mathcal{T}(\dff A\dff)$\dss to\dss $\Delta\dff(\dff A\dff)$\dss
for\sss every\dss $A\qff \subset\pff I$\nnsp,\oss
i.e.\pss if\dss
$c\dff(\dff v\trf)\qff \in\qff A$\sss for every\dss
element\sss $v$\sss of\dss every\sss $A$\dnsp-cell.\oss
Similarly\halfff,\pss $c$ is\dss called a\qss 
\emph{Scarf\qss coloring}\oss 
if\qss
$c\dff(\dff v\trf)\off \not\in\off A$\pss 
every element\sss $v$\dss of\dss every\sss $A$\dnsp-cell\dss such\dss that\dss
$A\off \neq\off I$\nnsp.\oss

\myuppar{Scarf's\qss combinatorial\dss theorem.}
\emph{For every\sss classical\dss coloring\qss $c\dff \colon\dff X\ttoo I$\qss 
there exist\sss a\sss non-empty\dss subset\pss $A\qff \subset\pff I$\qss
and\dss an\qss $A$\dnsp-cell\qss $\sigma$\sss
such\dss that\qss
$c\trf(\dff \sigma\dff)\off =\off A$\nnsp.\oss}

\proof
It\dss is\dss sufficient\dss to combine\qss Theorem\qss \ref{scarf-arbitrary-colorings}\qss
with\qss Corollary\qss \ref{scarf-chain-simplex}.\oss  \eproof

\myuppar{Scarf's\qss dual\dss form\sss of\qss Sperner's\qss lemma.}
\emph{If\qss $c\dff \colon\dff X\ttoo I$\qss is\dss an\qss
Scarf\dss coloring\fff,\oss
then\qss $c\dff(\dff \sigma\dff)\off =\off I$\qss for\dss
some\qss $I$\dnsp-cell\qss $\sigma$\dss
and\qss the\sss number\dss of\qss such\qss $I$\dnsp-cells $\sigma$ is\dss odd.\oss}

\proof
It\dss is\dss sufficient\dss to combine\qss Theorem\qss \ref{scarf-scarf-colorings}\qss
with\qss Corollary\qss \ref{scarf-chain-simplex}.\oss  \eproof

\newpage
\mysection{Scarf's\qss proof\pss of\pss Brouwer's\qss fixed\qss point\pss theorem}{scarf-brouwer}

\myuppar{The standard $n$\dnsp-simplex.}
Let\sss $n$ be a natural\dss number 
and\dss 
$I\off =\off
\{\trf 0\fff,\pff 1\fff,\pff \ldots\fff,\pff n\trf\}$\nnsp.\oss
We will\dss number\dss the coordinates\sss in\dss $\rrr^{\fff n\dff +\dff 1}$
by\qss $i\qff \in\qff I$\nnsp.\oss
For a point\sss $x\qff \in\qff \rrr^{\fff n\dff +\dff 1}$\sss
we will denote by $x_{\fff i}$ its $i${\dnsp}th\dss coordinate,\oss
so that\oss
$x\off =\off
(\qff x_{\trf 0}\dff,\pff  x_{\dff 1}\dff,\pff \ldots\dff,\pff  x_{\dff n} \qff)$\dnsp.\oss
Let\dss $\Delta^{n}\off \subset\off \rrr^{\fff n\dff +\dff 1}$\dss
be defined\dss by\dss the equation\vspace{0pt} 
\[
\quad
x_{\trf 0}\qff +\qff x_{\dff 1}\qff +\qff \ldots\qff +\qff x_{\dff n}
\off =\off 
1
\]

\vspace{-12pt}
and\dss the inequalities\qss $x_{\dff i}\qff \geq\qff 0$\qss
for each\qss $i\qff \in\qff I$\nnsp.\oss
Then\dss $\Delta^{n}$ is\trs the standard
$n$\dnsp-simplex\halfff.\oss

\myuppar{Approximating\dss $\Delta^{n}$\dss by\dss finite subsets.}
Let\qss $X\off \subset\off \Delta^{n}$\qss be a finite set\halfff.\oss
If\trs $X$\dss is\dss sufficiently\sss dense in\dss $\Delta^{n}$\dnsp,\oss
then\dss $X$\dss may\sss serve as an approximation\dss to\dss $\Delta^{n}$\dnsp.\oss
The set\dss $X$\dss is\dss not\sss assumed\dss to have any\sss
additional\sss structure.\oss
In\dss particular\halfff,\pss $X$\dss is\dss not\sss assumed\dss to be\sss
the set\sss of\dss vertices of\dss a\sss triangulation of\dss $\Delta^{n}$\dnsp.\oss
Instead,\oss an additional\sss structure\dss is\dss induced on\dss $X$\dss
from\dss $\Delta^{n}$\dnsp.\oss
Let\dss us\sss choose for each\qss
$i\qff \in\qff I$\qss
a linear order\dss $<_{\fff i}$\dss on\dss $X$\dss such\dss that\vspace{1.25pt}
\begin{equation*}
\quad
x_{\dff i}\off <\off y_{\dff i}
\hspace*{1.2em}
\mbox{implies}
\hspace*{1.2em}
x\off <_{\fff i}\off y
\end{equation*}

\vspace{-10.75pt}
for every\qss $x\fff,\qff y\qff \in\qff X$\nnsp.\oss
Obviously\halfff,\pss such orders exist\halfff.\oss
Moreover\halfff,\oss under a mild\dss non-degeneracy\sss assumption\sss
such orders are unique.\oss
In\dss fact\halfff,\oss the order\dss $<_{\fff i}$\dss is\dss 
uniquely\sss determined\dss if\trs and\dss only\trs if\dss
no\sss two different\sss elements of\trs $X$\dss have\sss the same
$i${\dnsp}th\dss coordinate.\oss
This condition can\dss be achieved\dss by\sss a small\dss perturbation of\trs $X$\nnsp,\oss
but\dss it\dss precludes having\dss many\dss points of\trs $X$\dss on\dss faces of\dss
$\Delta^{n}$\dnsp,\oss
which\dss is\dss an\dss important\dss feature of\dss approximations
by\dss the vertices of\dss a\sss triangulation.\oss

\myuppar{Geometric simplices associated\dss with subsets of\trs $X$\nnsp.}
Suppose\sss that\qss $\sigma\off \subset\off X$\qss
and\pss $C\pff \subset\off I$\nnsp.\oss
Then\sss the subsets $\sigma$ and\dss $C$\dss lead\dss to a\sss geometric
simplex\dss $\Delta\dff(\dff \sigma,\off C\trf)$\dss contained\dss
in\dss $\Delta^{n}$\dnsp.\oss
It\dss is\dss defined as follows.\oss
For each\qss $i\qff \in\qff I$\pss
let\qss $m_{\dff i}\off =\off \min_{\fff i}\dff \sigma$\nnsp,\oss
and,\pss
in an agreement\dss with\dss the above notations,\pss
let\dss $m_{\trf i\halfff i}$ be\sss the $i${\nnsp}th\dss coordinate of\dss $m_{\dff i}$\nsp.\oss
Clearly\halfff,\pss $m_{\trf i\halfff i}\qff \leq\qff 1$\qss for every\sss $i$\nnsp.\oss
Let\dss $\Delta\dff(\dff \sigma,\off C\trf)$\dss
be\sss the subset\sss of\trs $\rrr^{\fff n\dff +\dff 1}$\dss
defined\dss by\dss the equation\vspace{0pt}
\[
\quad
x_{\trf 0}\qff +\qff x_{\dff 1}\qff +\qff \ldots\qff +\qff x_{\dff n}
\off =\off 
1
\]

\vspace{-12pt}
and\dss the\sss inequalities\qss
$x_{\dff i}\off \geq\off m_{\trf i\halfff i}$\pss
with\qss $i\qff \in\pff C$\qss
and\qss $x_{\dff i}\off \geq\off 0$\pss
with\qss $i\qff \in\pff I\qff \smallsetminus\pff C$\nnsp.\qff\oss
The definition of\trs the numbers\dss $m_{\trf i\halfff i}$\dss
implies\sss that\qss 
$\sigma\off \subset\off \Delta\dff(\dff \sigma,\off C\trf)$\nnsp.\oss
In\dss particular\halfff,\oss $\Delta\dff(\dff \sigma,\off C\trf)$\dss
is\dss non-empty\halfff.

Let\dss us\dss introduce new\sss coordinates\dss
$y_{\trf 0}\dff,\pff  y_{\dff 1}\dff,\pff \ldots\dff,\pff  y_{\dff n}$\qss
by\dss setting\qss
$y_{\dff i}\off =\off x_{\dff i}\qff -\qff m_{\trf i\halfff i}$\qss
for\qss $i\qff \in\pff C$\qss
and\qss
$y_{\dff i}\off =\off x_{\dff i}$\qss
for\qss $i\qff \in\qff I\qff \smallsetminus\pff C$\nnsp.\oss
In\dss these coordinates\dss $\Delta\dff(\dff \sigma,\off C\trf)$\dss
is\dss defined\dss by\dss the equation\vspace{4.125pt}
\begin{equation}
\label{y-coordinates}
\quad
y_{\trf 0}\qff +\qff y_{\dff 1}\qff +\qff \ldots\qff +\qff y_{\dff n}
\off =\off 
1
\qff -\off 
\sum\nolimits_{\qff i\qff \in\pff C}\qff m_{\trf i\halfff i}
\end{equation}

\vspace{-7.875pt}
and\dss the\sss inequalities\qss
$y_{\dff i}\off \geq\off 0$\pss for\dss all\qss $i\qff \in\pff I$\nnsp.\oss
Since\qss
$\Delta\dff(\dff \sigma,\off C\trf)
\off \neq\off
\varnothing$\nnsp,\oss
the right\dss hand side of\trs the last\sss equation\dss is\qss $\geq\qff 0$\nnsp.\oss
If\trs the right\dss hand side\dss is\dss equal\dss to $0$\nnsp,\oss
then\dss $\Delta\dff(\dff \sigma,\off C\trf)$\dss
is\dss a $1$\dnsp-point\sss set\halfff,\oss
i.e.\qss a $0$\dnsp-simplex.\oss
Since\qss 
$\sigma\off \subset\off \Delta\dff(\dff \sigma,\off C\trf)$\nnsp,\oss
in\dss this case $\sigma$ is\dss also a $1$\dnsp-point\sss set\halfff.\oss
If\trs the right\dss hand side\dss is\qss $>\qff 0$\nnsp,\oss
then\dss $\Delta\dff(\dff \sigma,\off C\trf)$\dss
is\dss a\dss geometric $n$\dnsp-simplex.\oss
Moreover\halfff,\pss in\dss this case $\Delta\dff(\dff \sigma,\off C\trf)$\sss
is\dss similar\sss to\sss $\Delta^{n}$\dss in\dss the sense of\dss
elementary\dss geometry\halfff.\oss
In\dss fact\halfff,\pss up\sss to a\sss translation\dss
$\Delta\dff(\dff \sigma,\off C\trf)$\dss
is\dss homothetic\sss to\sss $\Delta^{n}$\dss
with\dss the ratio equal\dss to\sss the right\dss hand side of\qss (\ref{y-coordinates}).\oss

Clearly\halfff,\oss the $(\fff n\dff -\dff 1\fff)$\dnsp-faces of\qss 
$\Delta\dff(\dff \sigma,\off C\trf)$\sss
are defined\dss by\dss the equations\qss 
$y_{\dff i}\off =\off 0$\nnsp.\oss
Equivalently\halfff,\oss
they\sss are defined\dss by\dss the equations\qss
$x_{\dff i}\off =\off m_{\trf i\halfff i}$\pss
with\qss $i\qff \in\pff C$\qss
and\qss $x_{\dff i}\off =\off 0$\pss
with\qss $i\qff \in\pff I\qff \smallsetminus\pff C$\nnsp.\oss
In\dss particular\halfff,\oss
the intersection of\qss
$\Delta\dff(\dff \sigma,\off C\trf)$\sss
with\dss the face of\dss $\Delta^{n}$\sss
defined\dss by\dss the equations\qss
$x_{\dff i}\off =\off 0$\qss
with\qss $i\qff \in\pff I\qff \smallsetminus\pff C$\qss
is\dss non-empty\halfff.\oss

\mypar{Lemma.}{dominant-dense}
\emph{Let\qss $\varepsilon\qff >\qff 0$\nnsp.\oss
If\oss $\varepsilon'\qff >\qff 0$\qss is\trs sufficiently\trs small\dss
and\qss $X$\dss is\qss $\varepsilon'$\dnsp-dense\dss in\dss $\Delta^{n}$\nnsp,\oss
then\dss the\dss diameter\sss of\pss 
$\Delta\dff(\dff \sigma,\off C\trf)$\dss is\pss $<\pff \varepsilon$\pss
for\dss every\dss subset\qss $\sigma\off \subset\off X$\qss
dominant\trs with\dss respect\dss to\pss
$C\pff \subset\off I$\nnsp.\oss}

\proof
Let\dss $d$\trs be\sss the diameter of\trs the simplex\dss $\Delta^{n}$\sss
and\dss let\sss $r$\sss be\sss the distance from\dss its\sss bary\-center\dss
to its\sss boundary\halfff.\oss
Suppose\sss that\qss $\varepsilon'\qff <\pff \varepsilon\dff r/d$\qss
and\trs that\trs $X$\dss is\dss $\varepsilon'$\dnsp-dense\dss in\dss $\Delta^{n}$\dnsp.\oss

Let\qss $\sigma\off \subset\off X$\qss be a subset\dss
dominant\dss with\dss respect\dss to\qss
$C\pff \subset\off I$\nnsp.\oss
If\dss $\sigma$\dss consists of\dss $1$\sss point\halfff,\oss
then\dss $\Delta\dff(\dff \sigma,\off C\trf)$\sss is\dss also consists of\dss
$1$\dnsp-point\sss and\dss hence its diameter\dss is $0$\nnsp.\oss
Therefore we can\sss assume\sss that\qss
$\num{\sigma}\qff >\qff 1$\qss and\dss hence\dss 
$\Delta\dff(\dff \sigma,\off C\trf)$\sss
is\dss an $n$\dnsp-simplex.\oss
By\dss the definition of\dss  dominant\sss sets,\oss
the simplex\dss $\Delta\dff(\dff \sigma,\off C\trf)$\sss
does not\sss contain\sss elements of\trs $X$\dss in\dss its\dss interior\halfff.\oss

Suppose\sss that\dss the diameter of\dss $\Delta\dff(\dff \sigma,\off C\trf)$\sss
is\dss $\geq\qff \varepsilon$\nnsp.\oss
Since\dss $\Delta\dff(\dff \sigma,\off C\trf)$\sss is\dss similar\dss to\dss  $\Delta^{n}$\dnsp,\oss
in\dss this case\sss the distance from\dss the barycenter of\trs
the simplex\dss $\Delta\dff(\dff \sigma,\off C\trf)$
to\dss its\sss boundary\trs is\qss $\geq\pff \varepsilon\dff r/d$\nnsp.\oss
Since\dss $X$\dss is\dss $\varepsilon'$\dnsp-dense\dss in\dss $\Delta^{n}$\dss
and\qss $\varepsilon'\qff <\pff \varepsilon\dff r/d$\nnsp,\oss
there\dss is\dss a\sss point\qss $x\qff \in\qff X$\qss
with\dss the distance\qss $<\pff \varepsilon\dff r/d$\qss
from\dss the\sss barycenter\halfff.\oss  
Clearly\halfff,\qss $x$\dss is\dss
contained\dss in\dss the interior of\dss $\Delta\dff(\dff \sigma,\off C\trf)$\nnsp,\oss
contrary\dss to\sss the above observation.\oss
The contradiction shows\sss that\dss the diameter of\dss $\Delta\dff(\dff \sigma,\off C\trf)$\sss
is\qss $<\qff \varepsilon$\nnsp.\oss  \eproof

\myuppar{Scarf's\qss modification of\pss KKM\qss colorings.}
Let\dss us\dss consider a continuous map\qss
$f\dff \colon\dff
\Delta^{n}\qff \ttoo\qff \Delta^{n}$\dnsp.\oss
If\qss $x\qff \in\pff X$\qss
and\qss
$y\off\dff =\off f\dff(\dff x \trf)$\nnsp,\oss 
then\vspace{3pt}
\begin{equation}
\label{xy-sums}
\quad
x_{\trf 0}\qff +\qff x_{\dff 1}\qff +\qff \ldots\qff +\qff x_{\dff n}
\off =\off
y_{\dff 0}\qff +\qff y_{\dff 1}\qff +\qff \ldots\qff +\qff y_{\dff n}
\end{equation} 

\vspace{-9pt}
and\dss hence\sss there exists\qss $i\qff \in\qff I$\qss such\dss that\qss
$x_{\dff i}\qff \leq\qff y_{\dff i}$\nnsp.\oss
Let\dss $c\dff(\dff x\trf)$\dss be equal\dss to some\sss such\dss $i$\nnsp.\oss
Then\dss $c$\dss is\dss a\sss map\qss $X\qff \ttoo\qff I$\nnsp.\oss
This\dss is\dss the classical\sss coloring
associated\dss by\qss Scarf\qss 
with\dss $f$\dss and\trs $X$\nnsp.\oss

Note\sss that\qss Scarf's\trs rule defining $c$ is\dss nearly\sss opposite\sss
to\sss the\qss
Knaster--Kuratowski--Mazurkiewicz\qss rule,\oss
which allows color $i$\qss if\pss 
$x_{\dff i}\qff \geq\qff y_{\dff i}$\nnsp.\oss
At\dss the first\sss sight\dss it\sss seems\sss that\dss the difference
between\dss these\sss two rules\dss is\dss insignificant\halfff.\oss
But\dss the\qss KKM\qss argument\sss does not\dss work\dss with\qss Scarf's\trs rule,\oss
and\qss Scarf's\qss proof\dss does not\dss work\dss with\dss the\qss KKM\qss rule.\oss

\prooftitle{Scarf's\qss proof\qss of\pss Brouwer's\trs theorem}
Let\qss
$X_{\dff 1}\fff,\pff X_{\dff 2}\fff,\pff  X_{\dff 3}\dff,\pff \ldots$\qss
be a\sss sequence of\dss finite subsets of\dss $\Delta^{n}$\dnsp.\oss
Every\trs $X_{\dff k}$\dss is\dss equipped\dss with\dss the corresponding
orders\dss $<_{\fff i}$\nnsp,\oss where\qss $i\qff \in\qff I$\nnsp.\oss
If\dss we assume\sss that\dss the above non-degeneracy\sss
assumption holds for every\trs $X_{\dff k}$\nnsp,\oss
then\dss these orders are uniquely\sss determined.\oss
By\qss Scarf's\qss combinatorial\dss theorem\dss there exist\sss subsets\qss
$\sigma_{\dff k}\off \subset\off X_{\dff k}$\qss
and\pss
$C_{\dff k}\off \subset\off I$\qss
such\dss that\dss $\sigma_{\dff k}$\dss is\dss dominant\dss 
with\dss respect\dss to\qss $C_{\dff k}$\dss
and\qss
$c\dff(\dff \sigma_{\dff k}\trf)\pff =\off\dff C_{\dff k}$\nsp.\oss
Let\qss\vspace{2.625pt}
\[
\quad
\Delta_{\dff k}
\off =\off
\Delta\dff(\dff \sigma_{\dff k}\fff,\off C_{\dff k}\trf)
\qff.
\]

\vspace{-9.375pt}
We claim\dss that\dss for each\qss $i\qff \in\qff I$\qss
the simplex\dss $\Delta_{\dff k}$\dss contains a\sss point\dss $x$\dss
such\dss that\qss
$x_{\dff i}\qff \leq\qff y_{\dff i}$\nnsp,\oss
where\qss $y\off\dff =\off f\dff(\dff x \trf)$\nnsp.\oss
Indeed,\oss
if\qss $i\qff \in\pff C_{\dff k}$\nsp,\oss
then\qss
$c\dff(\dff x\trf)\off =\off i$\qss
for some\qss $x\qff \in\qff \sigma_{\dff k}$\nnsp.\oss
Then\qss
$x_{\dff i}\qff \leq\qff y_{\dff i}$\qss
by\dss the choice of\trs 
$c$\nnsp.\oss
If\qss
$i\qff \in\pff I\qff \smallsetminus\pff C_{\dff k}$\nsp,\oss
then\dss the intersection of\dss $\Delta_{\dff k}$
with\dss the face of\dss $\Delta^{n}$\sss
defined\dss by\dss the equation\qss
$x_{\dff i}\off =\off 0$\qss
is\dss non-empty\sss and\dss hence\qss
$x_{\dff i}\off =\off 0\off \leq\off y_{\dff i}$\qss
for any\dss point\sss $x$ in\dss this\sss intersection.\oss
This proves our claim.\oss

Now suppose\sss that\dss the sets\dss $X_{\dff k}$\dss are
chosen\sss in such a\sss way\dss that\dss $X_{\dff k}$\dss
is\dss $\varepsilon_{\dff k}$\hnsp\dnsp-dense\sss in\dss $X_{\dff k}$\dss
for some sequence\qss
$\varepsilon_{\dff k}\qff \ttoo\qff 0$\nnsp.\oss
Lemma\qss \ref{dominant-dense}\qss implies\sss that\dss the diameters
of\trs the simplices\dss $\Delta_{\dff k}$\dss tend\dss to $0$\nnsp.\oss
After\dss passing\dss to a subsequence\sss we can assume\sss
that\sss simplices\sss $\Delta_{\dff k}$\sss converge\sss to a single
point\qss $x\qff \in\qff \Delta^{n}$\dnsp,\oss
i.e.\qss every\sss sequence of\dss points\qss
$x\dff(\dff k\trf)\qff \in\qff \Delta_{\dff k}$\qss
converges\sss to\dss $x$\nnsp.\oss
Since\dss $f$\dss is\dss continuous,\oss
the inequalities of\trs the previous\sss paragraph
survive passing\dss to\sss the limit\halfff.\oss
It\dss follows\dss that\qss
$x_{\dff i}\qff \leq\qff y_{\dff i}$\qss
for every\qss $i\qff \in\qff I$\nnsp,\oss
where\qss $y\off\dff =\off f\dff(\dff x \trf)$\nnsp.\oss
Together\dss with\qss (\ref{xy-sums})\qss this implies\sss
that\qss
$x_{\dff i}\off =\off y_{\dff i}$\qss
for every\qss $i\qff \in\qff I$\nnsp,\qff\oss
i.e.\qss $x\off =\off y\off =\off f\dff(\dff x\trf)$\nnsp.\oss
Hence\dss $x$\dss is\dss a\sss fixed\dss point\sss of\qss $f$\nnsp.\oss  \eproof

\myuppar{Orders and\dss triangulations.}
One may\sss ask\trs if\pss Scarf's\qss proof\dss only\sss avoids mentioning\dss
triangulations of\dss $\Delta^{n}$\dnsp,\oss but\dss nevertheless\sss
they\sss are working\dss in\dss the background.\oss
Indeed,\oss the complex\dss $\mathcal{T}(\trf I\trf)$\dss
from\qss Section\qss \ref{orders}\qss
played a role similar\dss to\sss triangulations.\oss
The set\dss $X$\dss seems\sss to\dss be a natural\sss candidate\sss 
to\sss be\sss the set\sss of\trs the vertices
of\dss a\sss hidden\dss triangulation,\oss
but\dss $X$\dss can\dss be entirely\sss contained\dss in\dss the interior of\dss
$\Delta^{n}$\dss
and\dss hence has not\sss enough\dss points\sss to\sss triangulate\dss $\Delta^{n}$\dnsp.\oss
A natural\dss way\dss to deal\sss with\dss this issue\sss
is\dss to enlarge\dss $\Delta^{n}$\dss to include vertices corresponding\dss
to elements of\trs $I$\nnsp.\oss
Let\dss $E^{\dff n}\off \subset\off \rrr^{\fff n\dff +\dff 1}$\dss
be defined\dss by\dss the same equation as\dss $\Delta^{n}$\dss
and\dss the inequalities\qss
$x_{\dff i}\qff \geq\qff 1\qff -\qff n$\qss
with\qss $i\qff \in\qff I$\nnsp.\oss 
Then\dss $E^{\dff n}$\dss is\dss a\sss simplex containing\dss $\Delta^{n}$\dnsp.\oss
Its vertices are\vspace{1.5pt}
\[
\quad
v_{\fff i}
\off =\off
(\qff 1\dff,\pff \ldots\dff,\pff 1\dff,\pff 
1\qff -\qff n\dff,\pff
1\dff,\pff \ldots\dff,\pff 1 \qff)
\qff,
\]

\vspace{-10.5pt}
where\qss $1\qff -\qff n$\qss stands at\dss the $i${\dnsp}th\dss place
and $1$\nnsp's\dss at\sss all\sss other\sss places.\oss
Let\sss $\mathcal{S}\dff (\trf I\trf)$\dss be\sss the envelope\sss of\dss
$\mathcal{T}(\trf I\trf)$\dss in\dss the sense of\qss Section\qss \ref{pseudo-simplices}.\oss
Let\qss 
$V
\off =\off
\{\qff v_{\fff i}\qff \mid\qff i\qff \in\qff I \pff\}$\qss
and\dss let\dss us\dss identify\sss each vertex\qss $i\qff \in\qff I$\qss
of\qss $\mathcal{S}\dff (\trf I\trf)$\dss with\dss the point\dss 
$v_{\fff i}\qff \in\qff E^{\dff n}$\dnsp.\oss
This\sss identifies\sss the set\sss of\dss vertices of\qss $\mathcal{S}\dff (\trf I\trf)$\dss
with\qss $X\qff \cup\qff V$\nnsp.\oss
For every\sss simplex\dss $\sigma\qff \subset\qff X\qff \cup\qff V$\dss
of\qss $\mathcal{S}\dff (\trf I\trf)$\qss
let\dss us\dss denote by\sss $\Delta\dff(\dff \sigma\dff)$ 
the convex\dss hull\sss of\dss $\sigma$
in\dss $\rrr^{\fff n\dff +\dff 1}$\dnsp.\oss
If\trs we assume\sss that\dss no\dss $n\qff +\qff 1$\dss elements of\qss $X\qff \cup\qff V$\dss
are\sss linearly\sss dependent\qss (which can be achieved\dss by\sss a\sss small\dss
perturbation of\trs $X$\nsp),\oss
then every $\Delta\dff(\dff \sigma\dff)$ is\dss a\sss simplex
and\dss the collection of\trs these simplices seems\sss to be a good candidate\sss
for a\sss triangulation of\trs $E^{\dff n}$\dss behind\dss the proof\halfff.\oss
This\sss idea\dss is\dss supported\dss by\sss a\sss theorem of\qss Scarf\qss to\sss the effect\dss
that\dss the union of\dss simplices $\Delta\dff(\dff \sigma\dff)$
is\dss equal\dss to\dss $E^{\dff n}$\dss
(see\qss \cite{sc3},\oss Theorem\qss 7.2.2).\oss
But\halfff,\pss as an example of\qss Scarf\qss shows,\oss
these simplices do\dss not\sss always\sss form a\sss
triangulation\fff:\pss different\sss $n$\dnsp-simplices $\Delta\dff(\dff \sigma\dff)$
may\dss have common\dss interior\dss points\qss
(see\qss \cite{sc3},\oss the end of\qss Section\qss 7.2).\oss
So,\oss the envelope\sss $\mathcal{S}\dff (\trf I\trf)$\sss 
of\trs $\mathcal{T}(\trf I\trf)$\dss seems\sss to be\sss
the right\sss analogue of\dss a\sss triangulation\dss in\qss Scarf's\qss method.

\mysection{An\qss example\fff:\oss integer\qss points\qss in\qss a\qss simplex}{example}

\myuppar{Cyclic permutations of\dss coordinates.}
Let\sss $n$ be a natural\dss number\halfff,\oss 
$I\off =\off
\{\trf 0\fff,\pff 1\fff,\pff \ldots\fff,\pff n\trf\}$\nnsp,\oss
and\dss let\dss us\dss number\dss the coordinates\sss in\dss $\rrr^{\fff n\dff +\dff 1}$
by\qss $i\qff \in\qff I$\nnsp.\oss
In\dss this section\dss we will\sss often\dss treat\sss elements of\trs $I$\dss
as integers\dss modulo\qss $n\qff +\qff 1$\nnsp.\oss
Actually\halfff,\oss we will\dss need only\dss the resulting\sss cyclic order on\dss $I$\nnsp.\oss
Let\qss
$t\dff \colon\dff
\rrr^{\fff n\dff +\dff 1}
\qff \ttoo\qff
\rrr^{\fff n\dff +\dff 1}$\qss
be\sss the map induced\dss by\dss the cyclic shift\sss of\dss coordinates,\oss
i.e.\qss let\vspace{2.5pt}
\[
\quad
t\qff
(\qff
x_{\trf 0}\fff,\pff x_{\dff 1}\fff,\pff \ldots\fff,\pff x_{\dff n}
\qff)
\off =\off
(\qff
x_{\dff 1}\fff,\pff \ldots\fff,\pff x_{\dff n}\fff,\pff x_{\trf 0} 
\qff)
\]

\vspace{-9.5pt}
and\dss let\qss $t_{\trf i}$\dss be\sss the $i$\dnsp-fold composition\qss
$t\dff \circ\dff t\dff \circ\dff \ldots\dff \circ\dff t$\nnsp.\qff\oss
Then\vspace{2.5pt}
\[
\quad
t_{\trf i}\qff
(\qff
x_{\trf 0}\fff,\pff x_{\dff 1}\fff,\pff \ldots\fff,\pff x_{\dff n}
\qff)
\off =\off
(\qff
x_{\dff i}\fff,\pff x_{\dff i\dff +\dff 1}\fff,\pff \ldots\fff,\pff x_{\dff n\dff +\dff 1}
\qff)
\qff,
\]

\vspace{-9.5pt}
where\sss the addition\dss in\sss subscripts\dss is\dss understood\dss modulo\qss $n\qff +\qff 1$\nnsp.\oss
In\dss particular\halfff,\oss $t_{\trf 0}\off =\off\dff \id$\nnsp.\oss

\myuppar{Lexicographic\sss orders.}
Let\qss $x\fff,\pff y\qff \in\qff \rrr^{\fff n\dff +\dff 1}$\qss and \vspace{2.5pt}
\[
\quad
x
\off =\off
(\qff
x_{\trf 0}\fff,\pff x_{\dff 1}\fff,\pff \ldots\fff,\pff x_{\dff n}
\qff)
\qff,\quad\off
y
\off =\off
(\qff
y_{\dff 0}\fff,\pff y_{\dff 1}\fff,\pff \ldots\fff,\pff y_{\dff n}
\qff)
\qff.
\]

\vspace{-9.5pt}
The inequality\qss
$x\qff <\qff y$\qss 
in\dss the\qss
\emph{lexicographic\dss order}\qss $<$\qss
means\sss that\pss
$x\off \neq\off y$\qss and\qss
$x_{\dff i}\qff <\qff y_{\dff i}$\qss
for\dss the smallest\sss $i$\sss such\dss that\qss
$x_{\dff i}\off \neq\off y_{\dff i}$\nsp.\oss
Together\dss the maps\dss $t_{\trf i}$\dss
the order\qss $<$\qss generates a\qss
\emph{lexicographic\dss order}\qss $<_{\dff i}$\qss
for every\qss $i\qff \in\qff I$\nnsp.\oss
Namely\halfff,\oss the inequality\qss
$x\off <_{\dff i}\off y$\qss
means\sss that\vspace{2.25pt}
\[
\quad
t_{\trf i}\qff(\dff x\trf)\off <\off t_{\trf i}\qff(\dff y\trf)
\]

\vspace{-10pt}
Clearly\halfff,\oss $<_{\dff 0}$\qss is\trs the same order as\qss $<$\nsp.\oss

\myuppar{Integer\dss points in a simplex.}
Let\sss $N$ be a natural\sss number\halfff,\oss 
and\sss let\sss $D$\sss be\sss the set\sss of\dss all\dss points\vspace{2.5pt}
\[
\quad
a
\off =\off
(\qff
a_{\trf 0}\fff,\pff 
a_{\dff 1}\fff,\pff 
\ldots\fff,\pff 
a_{\dff n}
\qff)
\off \in\off\dff
\rrr^{\fff n\dff +\dff 1}
\]

\vspace{-9.5pt}
such\dss that\sss $a_{\dff i}$\sss are integers\qss $\geq\qff 0$\dss
and\dss
$a_{\trf 0}\qff +\qff a_{\dff 1}\qff +\qff \ldots\qff +\qff a_{\dff n}
\off =\off
N$\nnsp.\oss
Equivalently\halfff,\pss
$D$ is\dss the set\sss of\dss integer\dss points\sss
in\dss the $n$\dnsp-simplex\dss
$\Delta$\sss
defined\dss by\dss the equation\qss
$x_{\trf 0}\qff +\qff x_{\dff 1}\qff +\qff \ldots\qff +\qff x_{\dff n}
\off =\off
N$\qss
and\sss inequalities\dss $x_{\dff i}\qff \geq\qff 0$\nnsp.\oss
We will\sss consider $D$ together\dss with\sss
the orders\qss $<_{\dff i}$\qss restricted\dss to $D$\nnsp.\oss

\myuppar{The operators\dss $S_{\fff i}$.}
Let\qss $i\qff \in\qff I$\qss and\dss let\qss
$a
\off =\off
(\qff
a_{\trf 0}\fff,\pff 
a_{\dff 1}\fff,\pff 
\ldots\fff,\pff 
a_{\dff n}
\qff)
\off \in\off\dff
\rrr^{\fff n\dff +\dff 1}$\dnsp.\oss
Then\qss\vspace{2.5pt}
\[
\quad
S_{\fff i}\dff(\dff a\trf)
\off =\off
(\qff
a_{\trf 0}\fff,\pff 
\ldots\fff,\pff
a_{\dff i\dff -\dff 1}\qff +\qff 1\fff,\pff 
a_{\dff i}\qff -\qff 1\fff,\pff 
\ldots\fff,\pff 
a_{\dff n}
\qff)
\qff,
\]

\vspace{-9.5pt}
where\sss the subtraction of\dss $1$\dss in\dss the subscript\dss is\dss understood\dss 
modulo\qss $n\qff +\qff 1$\nnsp.\oss
In other\dss words,\oss
$S_{\fff 0}\dff(\dff a\trf)
\off =\off
(\qff
a_{\trf 0}\qff -\qff 1\fff,\pff 
a_{\dff 1}\fff,\pff 
\ldots\fff,\pff 
a_{\dff n\dff -\dff 1}\fff,\pff
a_{\dff n}\qff +\qff 1
\qff)$\nnsp.\oss

\mypar{Lemma.}{s-operators}
\emph{If\oss $i\qff \geq\qff 1$\nnsp,\oss
then\oss
$a
\off\dff <_{\dff i}\off\dff
S_{\fff 0}\dff(\dff a\trf)$\nnsp.\oss}

\proof
Since\qss
$t_{\trf i}\qff(\dff a\trf)
\off =\off
(\qff
a_{\dff i}\fff,\pff 
a_{\dff i\dff +\dff 1}\fff,\pff
\ldots\fff,\pff 
a_{\dff n}\fff,\pff
\ldots 
\off)$\qss
and\vspace{3pt}
\[
\quad
\hspace*{1em}
t_{\trf i}\qff
\bigl(\qff
S_{\fff 0}\dff(\dff a\trf)
\qff\bigr)
\off =\off
(\qff
a_{\dff i}\fff,\pff 
a_{\dff i\dff +\dff 1}\fff,\pff
\ldots\fff,\pff 
a_{\dff n}\qff +\qff 1\fff,\pff
\ldots 
\off)\qff,
\]

\vspace{-9pt}
this follows directly\dss from\dss the definitions.\oss  \eproof

\myuppar{The $I$\dnsp-cells.}
Our\dss first\dss goal\dss is\dss to prove some necessary\sss conditions for a subset\sss
of\dss $D$\sss to be an $I$\dnsp-cell.\oss
Eventually\dss these necessary\sss conditions will\dss turn out\dss to be also sufficient\halfff.\oss
Suppose\sss that\qss $\sigma\pff \subset\off D$\qss is\dss an $I$\dnsp-cell.\oss
Then\dss the set\sss $\sigma$ consists of\qss $n\qff +\qff 1$\qss elements,\oss
which\dss we will\sss denote by\qss 
$a\trf(\dff 0\dff)\fff,\off a\trf(\dff 1\dff)\fff,\off 
\ldots\fff,\off a\trf(\dff n\dff)$\nnsp.\qff\oss
Let\dss $a_{\dff i}\trf(\dff k\trf)$\dss be\sss the $i${\nnsp}th\dss coordinate
of\trs $a\trf(\dff k\trf)$\nnsp,\oss
i.e.\qss let\vspace{3pt}
\[
\quad
a\trf(\dff k\trf)
\off =\off
\bigl(\qff
a_{\trf 0}\trf(\dff k\trf)\fff,\pff 
a_{\dff 1}\trf(\dff k\trf)\fff,\pff 
\ldots\fff,\pff 
a_{\dff n}\trf(\dff k\trf)
\qff\bigr)
\qff.
\]

\vspace{-9pt}
Without\sss any\dss loss of\dss generality\dss we can assume\sss that\vspace{3pt}
\begin{equation}
\label{ordered}
\quad
a\trf(\dff 0\trf)\off <_{\dff 0}\dff\off
a\trf(\dff 1\trf)\off <_{\dff 0}\dff\off
\ldots\off <_{\dff 0}\dff\off
a\trf(\dff n\trf)
\qff
\end{equation}

\vspace{-9pt}
and\dss hence,\oss in\dss particular\halfff,\oss\vspace{3pt}
\begin{equation}
\label{ordered-0}
\quad
a_{\trf 0}\trf(\dff 0\trf)\off \leq\off
a_{\trf 0}\trf(\dff 1\trf)\off \leq\off
\ldots\off \leq\off
a_{\trf 0}\trf(\dff n\trf)
\qff.
\end{equation}

\vspace{-7pt}
\mypar{Lemma.}{not-equal}
\emph{For each\qss $i\qff \in\qff I$\qss
not\dss all\dss coordinates\dss 
$a_{\trf i}\trf(\dff k\trf)$\nnsp,\oss $0\qff \leq\qff k\qff \leq\qff n$\nnsp,\oss
are equal.\oss}

\proof
Indeed,\oss if\qss all\dss $a_{\trf i}\trf(\dff k\trf)$\dss are equal,\oss
the\qss $<_{\dff i}$\qss and\qss $<_{\dff i\dff +\dff 1}$\qss induce\sss the same orders
on $\sigma$ and\dss hence\qss
$\min\nolimits_{\trf i}\trf \sigma\off =\off \min\nolimits_{\trf i\dff +\dff 1}\trf \sigma$\nnsp.\oss
This contradicts\sss to\sss the fact\dss that\sss $\sigma$ is\dss a\sss $I$\dnsp-cell.\oss  \eproof

\mypar{Lemma.}{differ-by-one}
\emph{For each\dss $i$\dss
the\dss coordinates\dss 
$a_{\trf i}\trf(\dff k\trf)$\nnsp,\pss $0\qff \leq\qff k\qff \leq\qff n$\nnsp,\pss
differ\dss by\dss no more\sss than\dss $1$\nnsp.}

\proof
Suppose\sss that\qss $i\off =\off 0$\qss 
and\dss that\sss
some coordinates\dss $a_{\trf 0}\trf(\dff i\trf)$\dss
differ\dss by\sss at\dss least\sss $2$\nnsp.\oss
Then\qss
(\ref{ordered-0})\qss implies\sss that\qss
$a_{\trf 0}\trf(\dff n\trf)
\off \geq\off 
a_{\trf 0}\trf(\dff 0\trf)\qff +\qff 2$\qss
and\dss hence\qss
$a_{\trf 0}\trf(\dff n\trf)\qff -\qff 1
\off >\off
a_{\trf 0}\trf(\dff 0\trf)$\nnsp.\oss
Let\vspace{3pt}
\[
\quad
a\off\dff =\off\dff S_{\fff 0}\dff\bigl(\qff a\dff(\dff n\trf)\qff\bigr)
\qff.
\]

\vspace{-9pt}
Then\qss
$a_{\trf 0}\trf(\dff n\trf)\qff -\qff 1
\pff >\pff
a_{\trf 0}\trf(\dff 0\trf)$\qss
implies\sss that\dss $a\qff \in\pff D$\dss
and\dss 
$a\trf(\dff 0\trf)
\off <_{\dff 0}\dff\off 
a$\nnsp.\oss
Lemma\qss \ref{s-operators}\qss implies\sss that\qss
$a\trf(\dff n\trf)
\off <_{\dff i}\dff\off 
a$\qss
for\qss $i\qff \geq\qff 1$\nnsp.\oss
Therefore for every\dss $i\qff \in\qff I$\qss
there\dss exists\qss $b\qff \in\qff \sigma$\qss
such\dss that\qss $b\off <_{\dff i}\fff\off a$\qss
and\dss hence\qss
$\min\nolimits_{\trf i}\trf \sigma
\off <_{\dff i}\off\dff
a$\nnsp.\oss
In\dss view of\qss Lemma\qss \ref{minima}\qss
this contradicts\sss to $\sigma$ being an $I$\dnsp-cell.\oss 

This proves\sss the lemma for\qss $i\off =\off 0$\nnsp.\oss
The general\sss case reduces\sss to\sss this one by\sss
a cyclic permutation of\dss coordinates.\oss  \eproof

\mypar{Lemma.}{coordinate-j}
\emph{If\pss
$a_{\trf i\dff -\dff 1}\trf(\trf k\qff)
\off =\off
a_{\trf i\dff -\dff 1}\trf(\trf k\qff -\qff 1\trf)\qff +\qff 1$\nnsp,\oss
then\dss
$a\trf(\trf k\qff)
\off =\off\dff
\min\nolimits_{\trf i}\trf \sigma$\nnsp,\oss
the sequence}\vspace{1.95pt}
\[
\quad
a\trf(\trf k\qff)\fff,\off 
a\trf(\trf k\qff +\qff 1\dff)\fff,\off 
\ldots\fff,\off 
a\trf(\dff n\dff)\fff,\off
a\trf(\dff 0\dff)\fff,\off
a\trf(\dff 1\dff)\fff,\off
\ldots\fff,\off
a\trf(\trf k\qff -\qff 1\dff)
\]

\vspace{-10.05pt}
\emph{is\dss increasing\dss with\dss respect\dss to\sss the order\qss
$<_{\dff i}$\nsp,\oss
and\qss
$a_{\trf i}\trf(\trf k\qff)
\off =\off
a_{\trf i}\trf(\trf k\qff -\qff 1\trf)\qff -\qff 1$\nnsp.\oss}

\proof
Let\dss us\dss consider\dss the case\qss $i\off =\off 1$\qss first\halfff.\oss
Lemma\qss \ref{differ-by-one}\qss together\dss with\dss the inequalities\qss
(\ref{ordered-0})\qss implies\sss that\dss the coordinates\dss
$a_{\trf 0}\trf(\trf l\qff)$\dss with\qss $l\qff \leq\qff k\qff -\qff 1$\qss
are equal,\oss
as also\sss the coordinates\dss
$a_{\trf 0}\trf(\trf l\qff)$\dss with\qss $l\qff \geq\qff k$\nnsp.\oss
In\dss turn,\oss this\sss implies\sss that\dss the orders\qss $<_{\dff 0}$\qss
and\qss $<_{\dff 1}$\qss coincide on\dss the sets\vspace{1.5pt}
\[
\quad 
\bigl\{\qff a\trf(\trf l\trf)
\off \bigl|\off
l\qff \leq\qff k\qff -\qff 1
\off\bigr\}
\hspace*{1.5em}\mbox{and}\hspace*{1.5em}
\bigl\{\qff a\trf(\trf l\trf)
\off \bigl|\off
l\qff \geq\qff k
\off\bigr\}
\qff.
\]

\vspace{-10.5pt}
Together\dss with\qss (\ref{ordered})\qss this\sss implies\sss that\dss 
the sequence\dss $a\trf(\trf l\trf)$\dss
is\dss increasing\dss with\dss respect\dss to\sss the order\qss $<_{\dff 1}$\pss
for\qss $l\qff \leq\pff k\qff -\qff 1$\qss
and\qss for\qss $l\qff \geq\qff k$\nnsp.\oss
This implies,\oss in\dss particular\halfff,\oss
that\dss the minimum\qss $\min\nolimits_{\trf 1}\trf \sigma$\qss
could\dss be equal\sss only\dss to\dss $a\trf(\trf k\qff)$\dss
or\dss $a\trf(\trf 0\trf)$\nnsp.\oss
But\dss
$a\trf(\trf 0\trf)
\off =\off
\min\nolimits_{\trf 0}\trf \sigma$\nnsp,\oss
and\sss the minima with\sss respect\dss to different\sss orders are different\dss
because $\sigma$ is\dff~an\dss $I$\dnsp-cell.\oss
Therefore\qss
$\min\nolimits_{\trf 1}\trf \sigma
\off =\off
a\trf(\trf k\qff)$\nnsp.

Let\dss us\dss prove\sss that\qss
$a\trf(\trf n\trf)\off <_{\dff 1}\dff\off a\trf(\trf 0\trf)$\nnsp.\oss
Lemma\qss \ref{not-equal}\qss together\dss with\dss the inequalities\qss (\ref{ordered-0})\qss
implies\dss that\qss
$a_{\trf 0}\trf(\dff n\trf)
\off =\off
a_{\trf 0}\trf(\dff 0\trf)
\qff +\qff 
1$\nnsp.\oss
Since\sss the sums of\dss coordinates of\dss $a\trf(\dff n\trf)$\dss
and\dss $a\trf(\dff 0\trf)$\dss are both equal\dss to\dss $N$\nnsp,\oss
this\sss implies\sss that\qss
$a_{\fff j}\trf(\dff n\trf)
\off \neq\off
a_{\fff j}\trf(\dff 0\trf)$\qss 
for some\qss $i\qff \geq\qff 1$\nnsp.\oss
Let\dss $j$\dss be\sss the minimal\sss such subscript\halfff.\oss
Suppose\sss that\qss
$a_{\fff j}\trf(\trf 0\trf)\off <\dff\off a_{\fff j}\trf(\trf n\trf)$\qss
and\dss let\vspace{1.5pt}
\[
\quad
a\off\dff =\off\dff S_{\fff 0}\dff\bigl(\qff a\dff(\dff n\trf)\qff\bigr)
\qff.
\]

\vspace{-10.5pt}
We claim\dss that\qss $a\dff(\dff 0\trf)\off <_{\dff 0}\dff\off a$\nnsp.\oss
Indeed,\oss
$a_{\trf 0}
\off =\off 
a_{\trf 0}\trf(\dff n\trf)\qff -\qff 1
\off =\off
a_{\trf 0}\trf(\dff 0\trf)$\qss
and\qss
$a_{\dff m}\off =\off a_{\dff m}\trf(\dff 0\trf)$\qss
for\qss 
$1\leq\qff m\qff \leq\qff j\qff -\qff 1$\qss
by\dss the choice\sss of\dss $j$\nnsp.\oss
Since\qss
$a_{\fff j}\trf(\dff n\trf)
\off \leq\off 
a_{\fff j}$\nsp,\oss
the assumption\qss
$a_{\fff j}\trf(\trf 0\trf)\off <\dff\off a_{\fff j}\trf(\trf n\trf)$\qss
implies\sss that\qss
$a_{\trf i}\trf(\trf 0\trf)\off <\dff\off a_{\trf i}$\nnsp.\oss
It\dss follows\dss that\qss
$a\dff(\dff 0\trf)\off <_{\dff 0}\dff\off a$\nnsp,\oss
proving our claim.\oss
At\dss the same\sss time\qss
Lemma\qss \ref{s-operators}\qss implies\sss that\qss
$a\trf(\dff n\trf)
\off <_{\dff j}\dff\off 
a$\qss
for\qss $j\qff \geq\qff 1$\nnsp.\oss
Therefore for every\dss $j\qff \in\qff I$\qss
there\dss exists\qss $b\qff \in\qff \sigma$\qss
such\dss that\qss $b\off <_{\dff j}\fff\off a$\nnsp.\oss
But\dss this contradicts\sss to $\sigma$ being\sss an $I$\dnsp-cell.\oss
Therefore our assumption\dss is\dss wrong and actually\qss
$a_{\fff j}\trf(\trf n\trf)\off <\dff\off a_{\fff j}\trf(\trf 0\trf)$\nnsp.\oss
By\dss the choice of\dss $j$\sss this implies\sss that\qss
$a\trf(\trf n\trf)\off <_{\dff 1}\dff\off a\trf(\trf 0\trf)$\nnsp.\oss
Together\dss with\dss the observations\sss in\dss the previous paragraph\dss
this implies\sss that\dss for\qss $i\off =\off 1$\qss
the sequence from\dss the\sss lemma\dss is\dss
indeed\sss increasing\halfff.\oss

Since\sss this sequence\dss is\dss increasing\halfff,\oss
Lemma\qss \ref{not-equal}\qss implies\sss that\qss
$a_{\trf 1}\trf(\trf k\qff)
\off <\off
a_{\trf 1}\trf(\trf k\qff -\qff 1\trf)$\nnsp.\oss
Now\qss Lemma\qss \ref{differ-by-one}\qss
implies\sss that\qss
$a_{\trf 1}\trf(\trf k\qff)
\off =\off
a_{\trf 1}\trf(\trf k\qff -\qff 1\trf)
\qff -\qff 1$\nnsp.\oss  
This completes\sss the proof\dss in\dss the case\qss $i\off =\off 1$\nnsp.\oss
The general\sss case reduces\sss to\sss this one by\sss a cyclic permutation
of\dss coordinates.\oss  \eproof

\myuppar{The\sss sets\qss $a\qff +\qff \sigma\dff(\dff \iota\dff)$\nnsp.}
Let\qss $a\qff \in\pff \zzz^{\fff n\dff +\dff 1}$\qss and\dss let\qss
$\iota\dff \colon\dff I\qff \ttoo\qff I$\qss
be a\sss permutation.\oss
Let\dss us\dss define a sequence\qss 
$\alpha\trf(\dff 0\dff)\fff,\off \alpha\trf(\dff 1\dff)\fff,\off 
\ldots\fff,\off \alpha\trf(\dff n\dff)$\qss
by\sss setting\qss $\alpha\trf(\dff 0\dff)\off =\off a$\qss and using\dss 
the recursive rule\vspace{3.75pt}
\begin{equation}
\label{flip}
\quad
\alpha\trf(\dff k\qff)
\off =\off
S_{\fff i}\qff\bigl(\qff \alpha\dff(\trf k\qff -\qff 1\trf)\qff\bigr)
\qff,
\hspace*{1.0em}\mbox{where}\hspace*{1.0em}
i\off =\off\dff \iota\trf(\trf k\qff)
\qff,
\end{equation}

\vspace{-8.25pt}
for\qss $k\qff \geq\qff 1$\nnsp.\oss
Let\dss $\sigma\dff(\trf \iota\trf)$\dss be\sss the set\sss of\trs terms of\dss
such sequence defined\dss by\qss $a\off =\off 0$\qss and\dss the permutation $\iota$\nnsp.\oss
Clearly\halfff,\oss the sequence defined\dss by\qss 
$a\qff \in\pff \zzz^{\fff n\dff +\dff 1}$\qss and $\iota$
differs from\dss the sequence defined\dss by\sss $0$ and $\iota$ 
by\sss adding\sss $a$\sss to each\dss term.\oss
Therefore\sss the set\sss of\trs terms of\trs the sequence\qss
$\alpha\trf(\dff 0\dff)\fff,\off \alpha\trf(\dff 1\dff)\fff,\off 
\ldots\fff,\off \alpha\trf(\dff n\dff)$\qss 
is\dss equal\dss to\qss $a\qff +\qff \sigma\dff(\trf \iota\trf)$\nnsp.\oss\vspace{1.15pt}

\mypar{Theorem.}{top-cells}
\emph{If\pss $\sigma\pff \subset\off D$\qss is\dss an\qss $I$\dnsp-cell,\oss
then\qss $\sigma\off =\off a\qff +\qff \sigma\dff(\trf \iota\trf)$\qss
for some\qss $a\qff \in\qff \sigma$\qss
and\dss some\sss permutation\qss
$\iota\dff \colon\dff I\qff \ttoo\qff I$\nnsp.\oss
Moreover\halfff,\oss one can assume\dss that\qss
$\iota\trf(\dff 0\dff)\off =\off 0$\nnsp.\oss}\vspace{1.15pt}

\proof
As above,\oss we can\sss assume\sss that\qss 
$\sigma
\off =\off
\{
\qff
a\trf(\dff 0\dff)\fff,\off a\trf(\dff 1\dff)\fff,\off 
\ldots\fff,\off a\trf(\dff n\dff)
\qff\}$\qss
and\dss that\dss the inequalities\qss
(\ref{ordered})\qss hold.\oss
Since $\sigma$ is\dss an $I$\dnsp-cell\halfff,\oss
the map\qss
$i\off \longmapsto\off \min\nolimits_{\dff i}\qff \sigma$\qss
is\dss a\sss bijection\qss (see\dss Section\qss \ref{orders}).\oss
Therefore\sss
for every\qss $k\qff \in\qff I$\pss
there\dss is\dss a\sss unique\pss
$i\off =\off \iota\trf(\dff k\qff)$\qss
such\dss that\qss
$a\trf(\dff k\qff)
\off =\off
\min\nolimits_{\trf i}\trf \sigma$\qss
and\dss the map $\iota$ is\dss a\sss bijection\qss
$I\qff \ttoo\qff I$\nnsp.\oss
The inequalities\qss (\ref{ordered})\qss imply\dss that\qss
$a\trf(\dff 0\qff)
\off =\off
\min\nolimits_{\trf 0}\trf \sigma$\qss
and\dss hence\qss
$\iota\trf(\dff 0\dff)\off =\off 0$\nnsp.\oss
Let\qss $a\off =\off a\trf(\dff 0\qff)$\nnsp.\oss
As we will\sss see,\pss
$\sigma\off =\off a\qff +\qff \sigma\dff(\trf \iota\trf)$\qss
for\dss these $a$ and\dss $\iota$\nnsp.\oss

Let\qss $i\qff \in\qff I$\qss and\dss let\dss us\sss consider\dss the $i${\nnsp}th\dss coordinates\dss
$a_{\trf i}\trf(\trf k\qff)$\dss as cyclicly\sss ordered\dss by\qss $k\qff \in\qff I$\nnsp.\oss
By\qss Lemma\qss \ref{coordinate-j}\qss 
the $i${\nnsp}th\dss coordinate decreases\sss
in\dss this cyclic order\sss 
only\sss once and only\dss by\sss $1$\nnsp.\pss
Together\dss with\qss
Lemma\qss \ref{differ-by-one}\qss this\sss implies\sss that\sss
$i${\nnsp}th\dss coordinate\dss increases\sss also only\sss once
and only\dss by\sss $1$\nnsp.\oss
By\qss Lemma\qss \ref{coordinate-j}\pss
decreasing\sss of\trs
$i${\nnsp}th\dss coordinate happens simultaneously\dss
with\dss increasing of\trs
$(\fff i\dff -\dff 1\fff)${\nnsp}th\dss coordinate,\pss
namely\halfff,\oss 
when\sss one\sss passes\sss from\dss $k\qff -\qff 1$\qss to\dss $k$\nnsp,\oss
where\dss $k$\dss is\dss uniquely\sss determined\dss by\dss the condition\qss
$i\off =\off \iota\trf(\dff k\qff)$\nnsp.\oss

Since\sss the map\qss
$k\qff \longmapsto\qff \iota\trf(\dff k\qff)$\qss
is\dss bijective,\oss
the above observations\sss imply\dss that\dss 
$a\trf(\trf k\qff)$\dss differs\sss from\dss
$a\trf(\trf k\qff -\qff 1\qff)$
only\dss in\dss two coordinates,\oss
namely\halfff,\oss in\sss coordinates numbered\dss by\sss $i$\sss
and\qss $i\qff -\qff 1$\nnsp,\oss
where\qss $i\off =\off \iota\trf(\dff k\qff)$\nnsp.\oss
The $i${\nnsp}th\dss coordinate decreases by\sss $1$\sss
and\dss the $(\fff i\dff -\dff 1\fff)${\nnsp}th\dss coordinate increases\sss by\sss $1$\nnsp.\oss
Therefore\qss (\ref{flip})\qss holds for\dss the sequence\dss $a\trf(\dff k\qff)$\dss
in\dss the role of\trs the sequence\dss
$\alpha\trf(\dff k\qff)$\dss and\dss hence\qss
$\sigma\off =\off a\qff +\qff \sigma\dff(\trf \iota\trf)$\nnsp.\oss
This completes\sss the proof\halfff.\oss  \eproof\vspace{1.15pt}

\myuppar{Remarks.}
If\trs 
in\dss the assumption\qss (\ref{ordered})\qss
the order\dss $<_{\dff 0}$\qss is\dss replaced\dss by\dss the order\dss $<_{\dff i}$\nsp,\oss
the above arguments show\dss that\dss
$\sigma\off =\off b\qff +\qff \sigma\dff(\dff \kappa\dff)$\dss
for some\dss $b\qff \in\qff \sigma$\dss
and\sss a\sss permutation\qss
$\kappa\dff \colon\dff I\qff \ttoo\qff I$\qss
such\dss that\qss
$\kappa\trf(\dff i\trf)\off =\off i$\nnsp.\oss
One can also deduce\sss this from\qss Theorem\qss \ref{top-cells}\qss
by\dss using a cyclic permutation of\dss coordinates.\oss
One can check\dss that\sss $\kappa$ is\dss
uniquely\sss determined\dss by\dss the condition\qss
$\kappa\trf(\dff i\trf)\off =\off i$\qss
and\dss that\dss permutations corresponding\dss to different\qss
$i\qff \in\pff I$\qss differ\dss by\sss cyclic permutations of\trs
$I$\nnsp.\oss\vspace{1.15pt}

\myuppar{Toward a converse\sss to\trs Theorem\qss \ref{top-cells}.}
Let\qss
$I_{\trf 0}
\off =\off
\{\qff 1\fff,\pff 2\fff,\pff \ldots\fff,\pff n \qff\}$\nnsp.\oss
We will\dss identify\dss permutations\qss
$\iota\dff \colon\dff I\qff \ttoo\qff I$\qss
such\dss that\qss
$\iota\trf(\dff 0\dff)\off =\off 0$\qss
with\dss the induced\dss by\dss them\dss permutations\qss
$I_{\trf 0}\qff \ttoo\qff I_{\trf 0}$\nnsp.\oss
By\qss Theorem\qss \ref{top-cells}\qss every\sss $I$\dnsp-cell\sss $\sigma$
has\sss the form\dss
$a\qff +\qff \sigma\dff(\trf \iota\trf)$\dss
for some\qss $a\qff \in\qff \sigma$\qss
and some permutation\qss
$\iota\dff \colon\dff
I_{\trf 0}\qff \ttoo\qff I_{\trf 0}$\nnsp.\oss
It\dss turns out\dss that\halfff,\oss conversely\halfff,\oss
every\dss set\sss of\trs the form\qss $a\qff +\qff \sigma\dff(\trf \iota\trf)$\qss
contained\dss in\dss $D$\sss is\dss an $I$\dnsp-cell\halfff.\oss
We\sss will\dss prove\sss this in an\sss indirect\sss way\halfff.\oss
Namely\halfff,\oss it\dss turns out\dss that\sss a simple change of\dss coordinates
turns\dss $D$\sss into\sss the set\sss of\dss vertices of\dss
a canonical\dss triangulation of\dss a\sss standard simplex\halfff.\oss
Moreover\halfff,\oss this change of\dss coordinates\sss turns\sss
the sets of\trs the form\qss $a\qff +\qff \sigma\dff(\trf \iota\trf)$\qss
into sets of\dss vertices of\dss simplices of\dss this\sss triangulation.\oss
Once\sss this\dss is\dss established,\oss a\sss homological\sss
argument\qss (using only\sss chains)\qss implies\sss that\dss the
sets of\trs the form\qss
$a\qff +\qff \sigma\dff(\trf \iota\trf)$\qss
are $I$\dnsp-cells.\oss
The same argument\sss shows\sss that\dss $D$\sss is\dss the set\sss
of\trs the vertices of\dss a\sss triangulation.\oss

\myuppar{Freudenthal\trs triangulations.}
Let\dss us\dss summarize\sss the basic facts about\qss Freudenthal\qss triangulations.\oss
See\qss \cite{i3},\oss Section\qss 9\qss and\dss Appendix\qss 1\qss for\dss the details.\oss
The mentioned above standard simplex\dss is\dss the set\dss
$\Gamma$\dss of\dss all\dss points\qss
$y
\off =\off
(\qff
y_{\dff 1}\fff,\pff 
y_{\trf 2}\fff,\pff
\ldots\fff,\pff 
y_{\dff n}
\qff)
\off \in\off\dff
\rrr^{\fff n}$\qss
such\dss that\vspace{3pt}
\[
\quad
N\off \geq\off
y_{\dff n}\off \geq\off
y_{\dff n\dff -\dff 1}\off \geq\off
\ldots\off \geq\off
y_{\dff 1}\off \geq\off
0
\qff.
\]

\vspace{-9pt}
The set\dss $\Gamma$\dss is\dss an $n$\dnsp-simplex\dss
with\dss the vertices\qss
$N\dff u_{\dff 0}\fff,\off N\dff u_{\dff 1}\trf,\off \ldots\fff,\off N\dff u_{\dff n}
\off \in\off 
\rrr^{\fff n}$\dnsp,\oss
where\vspace{3pt}
\[
\quad
u_{\dff i}
\off =\off
(\dff 1\fff,\pff \ldots\fff,\pff 1\fff,\pff 0\trf,\pff \ldots\fff,\pff 0 \dff)
\off \in\off 
\rrr^{\fff n}
\qff
\]

\vspace{-9pt}
has $1$ as\sss the\dss first\sss $i$ coordinates
and $0$ as\sss the last\sss $n\qff -\qff i$\nnsp.\oss
Let\dss $\Gamma_{\fff 1}$\dss be\sss the $n$\dnsp-simplex\dss 
with\dss the vertices\qss
$u_{\trf 0}\fff,\pff u_{\dff 1}\trf,\pff \ldots\fff,\pff u_{\dff n}$\nsp.\oss
It\dss is\dss equal\dss to $\Gamma$\dss if\qss $N\off =\off 1$\nnsp.\oss
An $(\fff n\dff -\dff 1\fff)$\dnsp-face of\trs $\Gamma$\dss
is\dss defined\dss by\sss one of\trs the\sss following\sss equations\fff:\pss
$N\off =\off y_{\dff n}$\nsp,\oss
$y_{\dff i}\off =\off y_{\trf i\dff -\dff 1}$\qss
with\qss $n\qff \geq\qff i\qff \geq\qff 2$\nnsp,\oss
and\qss
$y_{\dff 1}\off =\off 0$\nnsp.\oss
Let\vspace{3pt}
\[
\quad
e_{\dff i}
\off =\off
(\dff 0\fff,\pff \ldots\fff,\pff 0\fff,\pff 1\fff,\pff 0\trf,\pff \ldots\fff,\pff 0 \dff)
\off \in\off 
\rrr^{\fff n}
\qff
\]

\vspace{-9pt}
be\sss the vector\dss having $1$ as\sss the $i${\nnsp}th\dss coordinate and $0$
as all\sss other coordinates.\oss
Then\vspace{3pt}
\[
\quad
u_{\dff i}
\off =\off
e_{\dff 1}\qff +\qff e_{\dff 2}\qff +\qff \ldots\qff +\qff e_{\dff i}
\qff
\]

\vspace{-9pt}
for every\qss $i\qff \in\pff I$\nnsp,\oss
where\sss the empty\sss sum\dss is\dss interpreted as $0$\nnsp.\oss
For a\sss permutation 
$\omega$ of\trs $I_{\trf 0}$\dss 
let\vspace{3pt}
\[
\quad
L_{\qff \omega}\trf \colon\trf
\rrr^{\fff n}\qff \ttoo\qff \rrr^{\fff n}
\]

\vspace{-9pt}
be\sss the\sss linear\dss map\sss taking\dss $e_{\dff i}$\sss to\sss
$e_{\trf \omega\dff(\dff i\trf)}$\dss
for every\qss $i\qff \in\pff I_{\trf 0}$\nsp,\oss
and\dss let\qss 
$\Gamma\trf(\dff \omega\dff)
\off =\off\dff
L_{\qff \omega}\trf(\qff \Gamma_{\fff 1}\trf)$\nnsp.\oss
Then\dss $\Gamma\trf(\dff \omega\dff)$\dss is\dss an $n$\dnsp-simplex\dss
with\dss the vertices\qss
$u_{\trf 0}\dff(\dff \omega\dff)
\fff,\pff 
u_{\dff 1}\dff(\dff \omega\dff)
\fff,\pff 
\ldots
\fff,\pff 
u_{\dff n}\dff(\dff \omega\dff)$\nsp,\oss
where\vspace{3pt}
\[
\quad
u_{\dff i}\dff(\dff \omega\dff)
\off =\off
e_{\trf \omega\dff(\dff 1\trf)}
\qff +\qff 
e_{\trf \omega\dff(\dff 2\trf)}
\qff +\qff 
\ldots
\qff +\qff 
e_{\trf \omega\dff(\dff i\trf)}
\qff.
\]

\vspace{-9pt}
Clearly\halfff,\oss 
$u\qff +\qff \Gamma\trf(\dff \omega\dff)$\qss
is\dss an $n$\dnsp-simplex\dss
for every\qss $u\qff \in\pff \zzz^{\fff n}$\qss
and every\dss permutation $\omega$ of\trs $I_{\trf 0}$\nsp.\oss
It\dss turns out\dss that\sss $n$\dnsp-simplices\sss
having\dss such\sss form and contained\dss in\dss
$\Gamma$\nsp,\oss
together\dss with\dss their\sss faces,\oss
form a\sss triangulation of\trs
$\Gamma$\nsp.\oss
The vertices of\trs this\sss triangulation are integer\dss points of\trs
$\Gamma$\nnsp.\oss
This\dss triangulation\dss is\trs the\qss \emph{Freudenthal\trs triangulation}\qss 
of\trs $\Gamma$\nnsp.\oss
Freudenthal\pss \cite{f}\qss considered only\dss the case\dss $N\off =\off 2$\nnsp,\oss
but\qss Kuhn\qss and\pss Scarf\pss \cite{sc3}\qss considered\dss the general\sss case also.\oss

\myuppar{The sets of\dss vertices of\sss $n$\dnsp-simplices.} 
Let\qss $u\qff \in\pff \zzz^{\fff n}$\qss and\dss let\sss
$\omega$\sss 
be a\sss permutation of\trs $I_{\trf 0}$\nsp.\oss
Let\qss 
$\upsilon\trf(\dff 0\dff)\fff,\off \upsilon\trf(\dff 1\dff)\fff,\off 
\ldots\fff,\off \upsilon\trf(\dff n\dff)$\qss
be\sss the sequence defined\dss by\qss
$\upsilon\trf(\dff 0\dff)\off =\off u$\qss and\dss  
the recursive rule\vspace{3pt}
\begin{equation*}
\quad
\upsilon\trf(\dff k\qff)
\off =\off\dff
\upsilon\dff(\trf k\qff -\qff 1\trf)\pff +\off e_{\dff i}
\qff,
\hspace*{1.0em}\mbox{where}\hspace*{1.0em}
i\off =\off\dff \omega\trf(\trf k\qff)
\qff
\end{equation*}

\vspace{-9pt}
and\qss $k\qff \geq\qff 1$\nnsp.\oss
Let\dss $\tau\dff(\dff \omega\dff)$\dss be\sss the set\sss of\trs terms of\dss
such sequence defined\dss by\qss $u\off =\off 0$\qss and\dss the permutation $\omega$\nnsp.\oss
Clearly\halfff,\oss the sequence defined\dss by\qss 
$u\qff \in\pff \zzz^{\fff n}$\qss and $\omega$
differs from\dss the sequence defined\dss by\sss $0$ and $\omega$ 
by\sss adding\sss $u$\sss to each\dss term.\oss
Therefore\sss the set\sss of\trs terms of\trs the sequence\qss
$\upsilon\trf(\dff 0\dff)\fff,\off \upsilon\trf(\dff 1\dff)\fff,\off 
\ldots\fff,\off \upsilon\trf(\dff n\dff)$\qss 
is\dss equal\dss to\qss $u\qff +\qff \tau\dff(\dff \omega\dff)$\nnsp.\oss
Clearly\halfff,\pss
$u\qff +\qff \tau\dff(\dff \omega\dff)$\qss
is\dss the set\sss of\dss vertices of\trs the simplex\qss
$u\qff +\qff \Gamma\trf(\dff \omega\dff)$\nnsp.\oss
It\dss follows\dss that\sss a\sss set\sss $\tau$ of\dss
integer\dss points of\trs $\Gamma$\dss
is\sss the set\sss of\dss vertices of\dss an $n$\dnsp-simplex\dss
if\trs and\dss only\trs if\dss 
$\tau$\dss has\sss the form\qss $u\qff +\qff \tau\dff(\dff \omega\dff)$\qss
with $u$ and $\tau$ as above.\oss\vspace{0.25pt}

\myuppar{The change of\dss coordinates.}
Recall\dss that\trs $D$\sss is\dss the set\sss of\dss integer\dss points in\dss
the $n$\dnsp-simplex\dss
$\Delta$\nsp.\oss 
Let\qss
$s\dff \colon\dff
\Delta\qff \ttoo\qff \Gamma$\qss
be\sss the affine map defined\dss by\vspace{3.75pt}
\[
\quad
s\qff
(\qff
x_{\trf 0}\fff,\pff
x_{\trf 1}\fff,\pff
x_{\trf 2}\fff,\pff 
\ldots\fff,\pff 
x_{\dff n}
\qff)
\off =\off
(\qff
y_{\trf 1}\fff,\pff 
y_{\trf 2}\fff,\pff
\ldots\fff,\pff 
y_{\dff n}
\qff) 
\qff,
\]

\vspace{-8.25pt}
where\qss
$y_{\dff i}
\off =\off
x_{\trf 0}\qff +\qff \ldots\qff +\qff x_{\trf i\dff -\dff 1}$\qss
for each\qss $i\qff \in\pff I_{\trf 0}$\nsp.\oss
The map\qss
$t\dff \colon\dff 
\Gamma\qff \ttoo\qff \Delta$\qss
defined\dss by\vspace{3.5pt}
\[
\quad
t\qff
(\qff
y_{\trf 1}\fff,\pff 
y_{\trf 2}\fff,\pff
\ldots\fff,\pff 
y_{\dff n}
\qff) 
\off =\off
(\qff
x_{\trf 0}\fff,\pff
x_{\trf 1}\fff,\pff
x_{\trf 2}\fff,\pff 
\ldots\fff,\pff 
x_{\dff n}
\qff)
\qff,
\]

\vspace{-8.5pt}
where\qss
$x_{\trf 0}\off =\off y_{\trf 1}$\nsp,\oss
$x_{\trf i}\off =\off y_{\trf i\dff +\dff 1}\qff -\pff y_{\trf i}$\pss
for\qss
$1\qff \leq\qff i\qff \leq\qff n\qff -\qff 1$\nnsp,\oss
and\qss
$x_{\dff n}\off =\off N\pff -\pff y_{\dff n}$\nsp,\oss
is\dss the\dss inverse\sss of\dss $s$\nnsp.\oss
Clearly\halfff,\oss $s$\dss and\dss $t$\dss map integer\dss points\sss
to integer\dss points and\dss hence establish a bijection between\dss
the set\dss $D$\sss of\trs integer\dss points in\dss $\Delta$\dss
and\dss the set\sss of\trs integer\dss points in\dss $\Gamma$\dnsp.\oss\vspace{0.25pt}

\mypar{Lemma.}{cells-simplices}
\emph{A subset\qss $\sigma\off \subset\off D$\qss
has\dss the\dss form\qss
$a\qff +\qff \sigma\dff(\trf \iota\trf)$\qss
for\sss a\sss permutation\qss 
$\iota\dff \colon\dff
I_{\trf 0}\qff \ttoo\qff I_{\trf 0}$\oss
if\pss and\dss only\trs if\qss
$s\dff(\dff \sigma\dff)$\sss is\dss the set\sss of\qss vertices
of\qss an $n$\dnsp-simplex.\oss}

\proof
Suppose\sss that\qss
$a\qff \in\pff D$\qss
and\qss
$b
\off =\off
S_{\dff i}\dff(\dff a\trf)$\qss
for some\qss $i\qff \in\qff I_{\trf 0}$\nsp.\oss
Let\vspace{4pt}
\[
\quad
s\qff
(\qff
a_{\trf 0}\fff,\pff
a_{\trf 1}\fff,\pff
\ldots\fff,\pff 
a_{\dff n}
\qff)
\off =\off
(\qff
u_{\trf 1}\fff,\pff 
\ldots\fff,\pff 
u_{\dff n}
\qff) 
\hspace*{1.0em}\mbox{and}\hspace*{1.0em}
s\qff
(\qff
b_{\trf 0}\fff,\pff
b_{\trf 1}\fff,\pff
\ldots\fff,\pff 
b_{\dff n}
\qff)
\off =\off
(\qff
v_{\trf 1}\fff,\pff 
\ldots\fff,\pff 
v_{\dff n}
\qff)
\qff. 
\]

\vspace{-8pt}
Then\qss
$v_{\fff j}\off =\off u_{\trf j}$\qss
for\qss 
$j\qff \leq\qff i\qff -\qff 1$\nnsp,\oss
$v_{\fff i}\off =\off u_{\trf i}\qff +\qff 1$\nnsp,\oss
and\qss
$v_{\fff j}\off =\off u_{\trf j}$\qss
for\qss 
$j\qff \geq\qff i\qff +\qff 1$\nnsp.\oss
Therefore\vspace{3.75pt}
\begin{equation*}
\quad
s\qff(\dff b\trf)
\off =\off
s\qff(\dff a\trf)\qff +\qff e_{\dff i}
\qff.
\end{equation*}

\vspace{-8.25pt}
Therefore applying\dss the operator\dss $S_{\dff i}$\dss to\qss
$a\qff \in\pff D$\qss
corresponds\sss to adding\dss $e_{\dff i}$\dss to\dss $s\qff(\dff a\trf)$\nnsp.\oss
It\dss follows\dss that\dss for every\qss $a\fff,\pff \iota$\qss the map\sss
$s$\sss takes\sss the set\dss
$a\qff +\qff \sigma\dff(\trf \iota\trf)$\dss
to\sss the set\dss
$s\qff(\dff a\trf)\qff +\qff \tau\dff(\trf \iota\trf)$\nnsp.\oss
Con\-verse\-ly\halfff,\oss
for every\qss $u\fff,\pff \omega$\qss the map\sss $t$\dss
takes\sss the set\dss
$u\qff +\qff \tau\dff(\dff \omega\dff)$\dss
to\sss the set\dss
$t\qff(\dff u\trf)\qff +\qff \sigma\dff(\dff \omega\dff)$\nnsp.\oss
Since\sss the sets of\dss integer\dss points\sss in\dss $\Gamma$\sss
having\dss the form\qss
$u\qff +\qff \tau\dff(\dff \omega\dff)$\dss
are exactly\dss the sets of\dss vertices of\dss $n$\dnsp-simplices,\oss
this proves\sss the\sss lemma.\oss  \eproof\vspace{0.25pt}

\mypar{Lemma.}{j-cells}
\emph{Let\qss
$C\pff \subset\pff I$\nnsp.\oss
If\pss $\sigma$ is\dss a\sss $C$\dnsp-cell,\oss
then\qss
$a_{\trf k}\off =\off 0$\qss
for every\qss $a\qff \in\qff \sigma$\dss
and\pss
$k\qff \in\qff I\qff \smallsetminus\pff C$\nnsp.}

\proof
Suppose\sss that\sss $\sigma$ is\dss a\dss $C$\dnsp-cell\fff,\oss
$a\qff \in\qff \sigma$\dnsp,\oss
and\qss
$a_{\trf k}\off \neq\off 0$\pss
for\dss some\qss
$k\qff \in\pff I\qff \smallsetminus\pff C$\nnsp.\oss
Let\sss $j$\trs be\sss the last\sss element\sss of\dss $C$
appearing after\sss $k$\sss in\dss the cyclic order of\qss $I$\nnsp.\oss
Let\qss
$b
\off =\off
(\qff
b_{\trf 0}\fff,\pff 
b_{\dff 1}\fff,\pff 
\ldots\fff,\pff 
b_{\dff n}
\qff)$\nnsp,\oss
where\qss 
$b_{\dff k}\off =\off a_{\trf k}\qff -\qff 1$\nnsp,\pss
$b_{\fff j}\off =\off a_{\fff j}\qff +\qff 1$\nnsp,\oss
and\qss
$b_{\dff i}\off =\off a_{\trf i}$\qss
for\qss $i\off \neq\off k\fff,\pff j$\nnsp.\oss
Then\qss
$b\qff \in\qff D$\qss
and\qss
$a\off <_{\dff i}\dff\off b$\qss
for every\qss $i\qff \in\pff C$\nnsp.\oss
This contradicts\sss to $\sigma$\sss being\sss a\sss $C$\dnsp-cell\halfff.\oss  \eproof

\myuppar{Simplicial\sss complexes and chains.}
For each\qss $i\qff \in\pff I$\qss let\dss 
$\Delta_{\dff i}$\dss 
be\sss the face of\trs
$\Delta$\dss 
defined\dss by\dss the equation\qss
$x_{\dff i}\off =\off 0$\nnsp.\oss
Let\dss 
$\Gamma_{\fff i}$\dss 
be\sss the face of\trs
$\Gamma$\dss 
defined\dss by\dss the equation\qss
$y_{\dff 1}\off =\off 0$\qss if\qss $i\off =\off 0$\nnsp,\oss
the equation\qss
$y_{\dff i\dff +\dff 1}\off =\off y_{\dff i}$\qss
if\qss $1\leq\qff i\qff \leq\qff n\qff -\qff 1$\nnsp,\oss
and\dss the equation\qss
$y_{\dff n}\off =\off N$\qss
if\qss $i\off =\off n$\nnsp.\oss
Then\qss\vspace{3pt}
\[
\quad
\Gamma_{\fff i}\off =\off s\trf(\trf \Delta_{\dff i}\trf)
\qff.
\]

\vspace{-9pt}
Let\qss $C\off \subset\off I$\nnsp.\oss
Recall\dss that\dss
the orders\dss $<_{\dff i}$\dss  
define an abstract\sss simplicial\sss complex\dss
$\mathcal{T} (\trf C\trf)$\dss
having\dss as\sss its\sss vertices 
and simplices
the elements and subsets of\dss $C$\dnsp-cells\sss  respectively\halfff,\oss
and\dss that\dss
the chain\dss $\mathcal{T}\dff\fclass{C}$\dss is\dss defined as\sss
the sum of\dss all\dss
$C$\dnsp-cells.\oss
We are interested\dss mainly\dss in cases\qss $C\off =\off I$\qss
and\qss $C\off =\off I\qff -\qff i$\qss for some\qss $i\qff \in\pff I$\nnsp.\oss 
By\qss Theorem\qss \ref{scarf-pseudo-simplex}\qss\vspace{4.5pt}
\begin{equation}
\label{lemma-boundary-chain}
\quad
\partial\trf \mathcal{T}\dff\fclass{\trf I\trf}
\off =\off
\sum\nolimits_{\qff i\qff \in\qff I}\qff \mathcal{T}\dff\fclass{\trf I\qff -\qff i\trf}
\qff.
\end{equation}

\vspace{-7.5pt}
Let\dss $\mathcal{F}\dff (\trf I\trf)$\dss be\sss the abstract\sss
simplicial\sss complex associated\dss with\qss Freudenthal\qss triangulation
of\trs $\Gamma$\dss
and\dss $\mathcal{F}\fff (\qff I\qff -\qff i\qff)$\dss
be\sss the subcomplex of\dss $\mathcal{F}\dff (\trf I\trf)$\dss
having as\sss its\sss simplices all\sss simplices of\dss $\mathcal{F}\dff (\trf I\trf)$\dss
contained\dss in\trs $\Gamma_{\fff i}$\nsp.\oss
Let\dss us define\sss the chain\dss 
$\mathcal{F}\dff\fclass{\trf I\trf}$\dss as\sss the sum of\dss all\sss 
$n$\dnsp-simplices of\trs $\mathcal{F}\dff (\trf I\trf)$\dss 
and\dss the chain\dss
$\mathcal{F}\dff\fclass{\trf I\qff -\qff i\trf}$\dss
as\sss the sum of\dss all $(\fff n\dff -\dff 1\fff)$\dnsp-simplices
of\sss $\mathcal{F}\dff (\qff I\qff -\qff i\qff)$\nnsp.\oss
Since\sss $\mathcal{F}\dff (\trf I\trf)$ is\dss the abstract\sss simplicial\sss
complex of\dss a\sss 
triangulation of\sss $\Gamma$\dnsp,\vspace{4.5pt}
\begin{equation}
\label{gamma-boundary-chain}
\quad
\partial\trf \mathcal{F}\dff\fclass{\trf I\trf}
\off =\off
\sum\nolimits_{\qff i\qff \in\qff I}\qff \mathcal{F}\dff\fclass{\trf I\qff -\qff i\trf}
\qff.
\end{equation}

\vspace{-7.5pt}
Theorem\qss \ref{top-cells}\qss together\dss with\qss 
Lemma\qss \ref{cells-simplices}\qss 
imply\dss that\dss the map $s$ defines\sss a\sss simplicial\sss map\qss
$\varphi\dff \colon\dff
\mathcal{T}\fff (\trf I\trf)
\qff \ttoo\qff
\mathcal{F}\fff (\trf I\trf)$\nnsp.\oss
Lemma\qss \ref{j-cells}\qss implies\sss that\dss 
for every\qss $i\qff \in\pff I$\qss
the vertices of\trs
$\mathcal{T}\fff (\qff I\qff -\qff i\qff)$\dss
belong\trs to\dss $\Delta_{\dff i}$\dss
and\dss hence\sss $\varphi$
maps\dss $\mathcal{T}\fff (\qff I\qff -\qff i\qff)$\dss to\qss
$\mathcal{F}\fff (\qff I\qff -\qff i\qff)$\nnsp.\oss

\mypar{Theorem.}{isomorphism}
\emph{\dnsp$\varphi\dff \colon\dff
\mathcal{T}\fff (\trf I\trf)
\qff \ttoo\qff
\mathcal{F}\fff (\trf I\trf)$\qss
is\dss an\dss isomorphism of\qss simplicial\dss complexes and}\vspace{4.5pt}
\begin{equation}
\label{chain-sur}
\quad
\varphi_{\dff *}\qff\bigl(\qff \mathcal{T}\dff\fclass{\trf I\trf}\qff\bigr)
\off =\off
\mathcal{F}\dff\fclass{\trf I\trf}
\qff.
\end{equation}

\vspace{-7.5pt}
\proof
Since $\varphi$ is\dss induced\dss by\sss an\sss 
injective affine map,\oss it\dss is\dss an\sss isomorphism\dss if\dss
every $n$\dnsp-simplex
of\dss $\mathcal{F}\fff (\trf I\trf)$\dss is\dss the image of\dss an $n$\dnsp-simplex
of\dss $\mathcal{T}\fff (\trf I\trf)$\nnsp.\oss
Clearly\halfff,\oss this would\dss follow once\qss (\ref{chain-sur})\qss
is\dss proved.\oss
Therefore\sss it\dss is\dss sufficient\dss to prove\qss (\ref{chain-sur}).\oss
We will\dss prove\qss (\ref{chain-sur})\qss using\sss
an\dss induction\dss by\sss $n$\nnsp,\oss
the case\qss $n\off =\off 0$\qss 
(and even\dss the case\qss $n\off =\off 1$\nsp)\qss 
being\dss trivial\halfff.\oss

Suppose\sss that\qss (\ref{chain-sur})\qss holds with\qss
$n\qff -\qff 1$\qss in\dss the role of\sss $n$\nnsp.\oss
Let\dss $i\qff \in\pff I$\nnsp.\oss 
By\sss omitting\dss the coordinate $x_{\dff i}$ we can\dss
identify\dss $\Delta_{\dff i}$\dss with\dss the $(\fff n\dff -\dff 1\fff)$\dnsp-dimensional\sss
version of\dss $\Delta$\nnsp.\oss
Lemma\qss \ref{j-cells}\qss implies\sss that\dss this\sss identification\dss
turns\dss $(\qff I\qff -\qff i\qff)$\dnsp-cells\dss
into\sss the $(\fff n\dff -\dff 1\fff)$\dnsp-dimensional\sss
version of\dss $I$\dnsp-cells.\oss
By\sss omitting\dss the coordinate $y_{\dff i}$ we can\dss
identify\dss $\Gamma_{\fff i}$\dss
with\dss the $(\fff n\dff -\dff 1\fff)$\dnsp-dimensional\sss
version of\dss $\Gamma$\nnsp.\oss
These identifications agree with\sss $s$\sss and\dss hence
our\dss inductive assumption\dss implies\sss that\vspace{4.5pt}
\begin{equation}
\label{surjectivity-1}
\quad
\varphi_{\dff *}\qff\bigl(\qff \mathcal{T}\dff\fclass{\trf I\qff -\qff i\trf}\qff\bigr)
\off =\off
\mathcal{F}\dff\fclass{\trf I\qff -\qff i\trf}
\qff.
\end{equation}

\vspace{-7.5pt}
Since\dss 
$\partial\dff \circ\dff \varphi_{\dff *}
\off =\off
\varphi_{\dff *} \circ\qff \partial$\nnsp,\oss
the equality\qss (\ref{lemma-boundary-chain})\qss
implies\sss that\vspace{4.5pt}
\[
\quad
\partial\qff
\varphi_{\dff *}\qff\bigl(\qff \mathcal{T}\dff\fclass{\trf I\trf}\qff\bigr)
\off\dff =\off\qff
\varphi_{\dff *}\dff\bigl(\qff \partial\trf\mathcal{T}\dff\fclass{\trf I\trf}\qff\bigr)
\off\dff 
=\off\qff
\varphi_{\dff *}\dff
\left(\off 
\sum\nolimits_{\qff i\qff \in\qff I}\pff \mathcal{T}\dff\fclass{\trf I\qff -\qff i\trf}
\off\right)
\qff.
\]

\vspace{-7.5pt}
In\dss view\sss of\qss
(\ref{surjectivity-1})\qss and\qss (\ref{gamma-boundary-chain})\qss
the\sss last\sss expression\dss is\dss equal\dss to\vspace{4.5pt}
\[
\quad
\sum\nolimits_{\qff i\qff \in\qff I}\pff 
\varphi_{\dff *}\qff
\bigl(\qff \mathcal{T}\dff\fclass{\trf I\qff -\qff i\trf}
\off\bigr)
\off\dff =\off\qff
\sum\nolimits_{\qff i\qff \in\qff I}\pff 
\mathcal{F}\dff\fclass{\trf I\qff -\qff i\trf}
\off\dff =\off\qff
\partial\trf \mathcal{F}\dff\fclass{\trf I\trf}
\qff
\]

\vspace{-7.75pt}
and\dss therefore\vspace{1pt}
\begin{equation*}
\quad
\partial\qff
\varphi_{\dff *}\qff\bigl(\qff \mathcal{T}\dff\fclass{\trf I\trf}\qff\bigr)
\off\dff 
=\off\qff
\partial\trf \mathcal{F}\dff\fclass{\trf I\trf}
\qff.
\end{equation*}

\vspace{-9pt}
Since\dss the complex\dss $\mathcal{F}\dff (\trf I\trf)$\dss
is\dss the abstract\sss simplicial\sss complex of\dss a\sss
triangulation of\dss a\sss simplex,\oss it\dss
is\dss non-branching\halfff,\oss strongly\sss connected,\oss
and\dss has non-empty\dss boundary\halfff.\oss
It\dss is\dss well\dss known and easy\dss to see\sss that\dss
these properties imply\dss that\dss
$\mathcal{F}\dff\fclass{\trf I\trf}$\dss
is\dss the only\sss chain\dss with\dss the boundary\dss
$\partial\trf \mathcal{F}\dff\fclass{\trf I\trf}$\nnsp.\oss
Therefore\sss the last\sss equality\dss implies\qss (\ref{chain-sur}).\oss  
This completes\sss the proof\halfff.\oss  \eproof

\mypar{Corollary.}{cells-triangulation}
\emph{The $n$\dnsp-simplices having\qss $I$\dnsp-cells\dss 
as\dss their\dss sets of\qss vertices
are\sss the $n$\dnsp-simplices of\dss a\sss triangulation of\qss $\Delta$\nnsp,\oss
namely\halfff,\oss of\qss the image of\qss Freudenthal\qss triangulation of\pss
$\Gamma$\sss under\dss $t$\nnsp.\oss}

\proof
This\sss follows from\dss the\sss theorem and\dss
the definitions of\dss $\mathcal{T}\fff (\trf I\trf)$\dss and\dss
$\mathcal{F}\fff (\trf I\trf)$\nnsp.\oss  \eproof

\mypar{Corollary.}{description-of-cells}
\emph{A subset\sss of\pss $D$\dss is\dss an\qss $I$\dnsp-cell\qss
if\trs and\dss only\trs if\qss it\dss has\dss
the\dss form\dss
$a\qff +\qff \sigma\dff(\trf \iota\trf)$\nnsp,\oss
where\qss $a\qff \in\pff D$\qss
and\pss
$\iota\dff \colon\dff
I\qff \ttoo\qff I$\qss
is\dss a\sss permutation\dss such\dss that\qss
that\qss
$\iota\trf(\dff 0\dff)\off =\off 0$\nnsp.\oss}

\proof
The\qss ``only\dss if\dff''\qss part\trs 
is\qss Theorem\qss \ref{top-cells}.\oss
Let\dss us\dss prove\sss the\qss ``if\dff''\qss part\halfff.\oss
Let\dss
$\sigma
\off \subset\off D$\dss
be a subset\sss of\trs the\dss form\dss
$a\qff +\qff \sigma\dff(\trf \iota\trf)$\nnsp.\oss
Lemma\qss \ref{cells-simplices}\qss implies\sss that\sss
$s\dff(\dff \sigma\dff)$\sss is\dss
an $n$\dnsp-simplex of\trs $\mathcal{F}\dff (\trf I\trf)$\nnsp,\oss
and\dss hence\qss Theorem\qss \ref{isomorphism}\pss
implies\sss that 
$\sigma$ is\dss an $n$\dnsp-simplex of\dss
$\mathcal{T}\dff (\trf I\trf)$\nnsp,\oss
i.e.\qss is\trs an\sss $I$\dnsp-cell\halfff.\oss  \eproof

\myuppar{Broken\sss and\dss restored\sss symmetry\halfff.}
The set\dss $D$\dss and\dss the family\sss of\dss orders\dss $<_{\dff i}$\dss
are invariant\dss under\dss the cyclic permutation of\dss
coordinates in\dss $\rrr^{\fff n\dff +\dff 1}$\dnsp,\oss
as\dss is\dss the notion of\dss an\sss $I$\dnsp-cell\halfff.\oss
But\dss the description of\trs $I$\dnsp-cells as subsets of\trs $D$\dss 
having\dss the form\qss
$a\qff +\qff \sigma\dff(\trf \iota\trf)$\qss
is\dss not\halfff,\oss
because\sss their definition\dss in\dss terms of\dss sequences\qss
$\alpha\trf(\dff 0\dff)\fff,\off \alpha\trf(\dff 1\dff)\fff,\off 
\ldots\fff,\off \alpha\trf(\dff n\dff)$\qss 
does not\dss require\qss (\ref{flip})\qss to hold\dss for\qss
$k\off =\off 0$\nnsp.\oss
The cyclic symmetry\sss of\trs the order\dss $n\qff +\qff 1$\dss
is\dss also broken\dss in\dss the very\sss definition of\trs $\Gamma$\dnsp.\oss
This symmetry\dss is\dss restored\dss by\dss
Corollary\qss \ref{description-of-cells},\oss
which\dss implies\sss that\dss if\qss 
(\ref{flip})\qss holds for\dss
$k\off \geq\off 1$\nnsp,\oss
then\qss (\ref{flip})\qss automatically\dss holds for\qss $k\off =\off 0$\nnsp.\oss
Since everything\dss is\dss invariant\dss
under\dss translations\dss by\dss integer\sss vectors,\oss
there\dss is\dss even\sss no need\dss to assume\sss that\dss
the\sss terms\dss $\alpha\trf(\dff i\trf)$\dss belong\dss to\sss $D$\nnsp.\oss

\myuppar{Remark.}
Scarf\qss proved\qss Corollaries\qss \ref{cells-triangulation}\qss
and\qss \ref{description-of-cells}\qss 
without\dss explicitly\dss invoking\qss Freudenthal\qss triangulations.\oss
See\qss \cite{sc3},\oss Theorems\qss 7.1.8\qss and\qss 6.2.1\qss respectively\halfff.\oss
But\halfff,\oss for example,\oss 
his formulas on\dss p.\qss 176\qss of\pss
\cite{sc3}\qss amount\dss to using an analogue of\trs the\sss map $s$\nnsp.\oss

\newpage
\mysection{Oriented\qss matroids}{oriented-matroids}

\myuppar{Signed subsets.}
Let\dss us fix a set\dss $M$\nnsp.\oss
A\qss \emph{signed subset}\qss $\sigma$\dss of\dss $M$\dss
is\dss a\sss pair\qss $(\qff \sigma_{\dff +}\dff,\pff  \sigma_{\dff -}\qff)$\qss
of\trs two disjoint\sss subsets\qss
$\sigma_{\dff +}\dff,\pff  \sigma_{\dff -}$\qss of\qss $M$\nnsp.\oss
The\qss \emph{support}\qss of\dss a\sss signed\sss subset\dss $\sigma$\dss
is\dss the union\vspace{4.5pt}
\[
\quad
\underline{\sigma}
\off =\off
\sigma_{\dff +}\qff \cup\qff \sigma_{\dff -}
\off.
\]

\vspace{-7.5pt}
The opposite of\dss a\sss signed\dss subset\qss
$\sigma
\off =\off
(\qff \sigma_{\dff +}\dff,\pff  \sigma_{\dff -}\qff)$\qss
is\dss the signed subset\qss
$(\qff \sigma_{\dff -}\dff,\pff  \sigma_{\dff +}\qff)$\nnsp.\oss
Equiv\-a\-lent\-ly\halfff,\oss the opposite of\dss $\sigma$\dss
is\dss the unique signed subset\dss $\tau$\dss such\dss that\qss\vspace{4.5pt}
\[
\quad
\tau_{\dff +}\off =\off \sigma_{\dff -}
\hspace*{1.2em}\mbox{and}\hspace*{1.2em}
\tau_{\dff -}\off =\off \sigma_{\dff +}
\off.
\]

\vspace{-7.5pt}
The opposite of\dss $\sigma$\dss is\dss usually\sss denoted\dss by\qss
${}-\qff \sigma$\nnsp.\oss
A signed\dss subset\dss $\sigma$\dss can\dss be identified\dss with a map\qss
$\underline{\sigma}
\qff \ttoo\qff
\{\qff +\dff,\pff -\qff\}$\nsp,\oss
which\dss is\dss also denoted\dss by\dss $\sigma$\dss 
and\sss defined\dss by\vspace{4.5pt}
\[
\quad
\sigma\dff(\dff x\trf)
\off =\off\dff 
+
\hspace*{1em}\mbox{if}\hspace*{1em}
x\qff \in\qff \sigma_{\dff +}
\hspace*{1.2em}\mbox{and}\hspace*{1.2em}
\sigma\dff(\dff x\trf)
\off =\off\dff 
-
\hspace*{1em}\mbox{if}\hspace*{1em}
x\qff \in\qff \sigma_{\dff -}
\off.
\]

\vspace{-7.5pt}
This\dss identification\dss matches\dss 
the notation\dss for\dss the opposite signed subset\halfff.\oss

\myuppar{Vector\dss configurations.}
Let\qss $M\off \subset\off \rrr^{\dff n}$\qss for some\dss $n$\nnsp.\oss
Suppose\sss that\qss $0\qff \not\in\qff M$\qss
and\sss that\sss different\sss elements of\trs $M$\dss
are not\dss proportional.\pss  
If\pss $X\off \subset\off M$\dss
is\dss a\sss linearly dependent\sss subset,\oss 
then\vspace{4.5pt}
\[
\quad
\sum\nolimits_{\qff x\qff \in\qff X\vphantom{X^X}}\qff \lambda_{\dff x}\dff x
\off\dff =\off\dff
0
\] 

\vspace{-7.5pt}
for some coefficients\qss $\lambda_{\dff x}\qff \in\qff \rrr$\nnsp.\oss
Such a relation defines a signed subset\sss $\sigma$ of\trs $M$\dss as follows:\oss
the set\sss $\sigma_{\dff +}$\sss is\dss equal\dss to\sss the set\sss
of\dss elements\qss $x\qff \in\qff X$\qss such\dss that\qss
$\lambda_{\dff x}\qff >\qff 0$\nnsp,\oss
and\sss $\sigma_{\dff -}$\sss is\dss equal\dss to\sss the set\sss
of\dss elements\qss $x\qff \in\qff X$\qss such\dss that\qss
$\lambda_{\dff x}\qff <\qff 0$\nnsp.\oss
In\dss general,\pss $\sigma$\dss depends not\sss only on\dss $X$\nnsp,\oss
but\sss also on\dss the relation.\oss
Let\dss us\sss
suppose\sss now\dss that\dss $X$\dss is\dss a minimal\dss linearly\sss dependent\sss subset\halfff,\oss 
i.e.\qss that\dss
no proper subset\sss of\dss $X$\dss
is linearly dependent\halfff.\oss 
In\dss this case\qss $\lambda_{\dff x}\off \neq\off 0$\qss
for all\qss $x\qff \in\qff X$\qss
and\dss the coefficients\dss $\lambda_{\dff x}$\dss
are uniquely\sss determined\dss by\dss $X$\dss up\sss to multiplication\dss by
a common\sss non-zero factor\halfff.\oss
Therefore\sss in\dss this case $\sigma$ is\dss determined\dss by\dss $X$\dss
up\sss to replacing $\sigma$ by\qss
${}-\qff \sigma$\nnsp.\oss
The signed subsets\dss $\sigma$\dss defined\dss by\dss
minimal\dss linearly\sss dependent\sss subsets\qss 
$X\off \subset\off M$\dss
in\dss the above manner\sss are called\dss the\qss \emph{circuits}\qss
of\dss the\qss \emph{vector configuration}\qss $M$\nnsp.\oss

The notion of\dss an\qss \emph{oriented\dss matroid}\pss
is an axiomatization of\dss the basic properties of\dss these circuits.\oss
The first\sss of\trs these basic properties are obvious.\oss
Clearly\halfff,\oss 
$\varnothing$\dss is not a circuit\halfff,\oss and\trs if\dss $\sigma$\dss
is a circuit\sss of\dss $M$\nnsp,\oss
then\dss ${}-\dff \sigma$\dss is also a circuit\halfff.\oss
Since\sss the circuits result\dss from\qss \emph{minimal}\pss
linear dependencies,\oss
if\qss $\sigma\fff,\pff \tau$\qss are circuits,\oss then\dss
the inclusion\qss $\underline{\sigma}\off \subset\off \underline{\tau}$\qss
implies\sss that\sss $\sigma$ is\dss equal\dss to either\qss
$\tau$\nnsp,\oss or\qss
${} -\qff \tau$\nnsp.\oss
The main\dss property\dss of\dss circuits\dss is\dss concerned\dss with an elimination\dss
procedure removing\dss a
element\dss from\dss two given relations.\oss

\myuppar{An elimination\dss procedure.}
Let\qss $X\fff,\pff Y\off \subset\off M$\qss be\sss two\sss linearly dependent\sss subsets.\oss
Then\vspace{4.5pt}
\[
\quad
\sum\nolimits_{\qff x\qff \in\qff X\vphantom{X^X}}\qff \lambda_{\dff x}\dff x
\off\dff =\off\dff
0
\hspace*{1.45em}\mbox{and}\hspace*{1.75em}
\sum\nolimits_{\qff y\qff \in\qff Y\vphantom{X^X}}\qff \mu_{\dff y}\dff y
\off\dff =\off\dff
0
\] 

\vspace{-7.5pt}
for some coefficients\qss $\lambda_{\dff x}\dff,\pff \mu_{\dff y}\qff \in\qff \rrr$\nnsp.\oss
Let\dss $\sigma$\dss and\dss $\tau$\dss be\sss the signed subsets defined\dss
by\dss these relations.\oss
Suppose\sss that\dss $u\qff \in\qff X\qff \cap\qff Y$\dss
and\dss the coefficients\dss $\lambda_{\dff u}$\dss and\dss $\mu_{\dff u}$\dss have opposite signs.\oss
Let\dss us multiply\dss the second\dss relation\dss by\qss
$\alpha
\off =\off
\num{\lambda_{\dff u}\dff/\fff\mu_{\dff u}}$\qss 
and add\dss the resulting\dss
relation\dss to\sss the first\sss one.\oss
This procedure\qss \emph{eliminates}\qss $z$\dss from\dss the given\dss relations
between\dss the elements of\trs $X$\dss and\sss $Y$\dnsp.\oss
Indeed,\oss
the result\trs is\dss a relation\dss between\dss the elements of\trs
the set\qss 
$Z
\off =\off
(\qff X\qff \cup\qff Y\qff)\qff \smallsetminus\qff \{\dff u\qff\}$\nnsp.\oss
Let\dss $\omega$\dss be\sss the signed\sss set\sss defined\dss
by\dss this relation.\oss
The coefficient\sss of\qss $z\qff \in\qff Z$\qss 
is equal\dss to\dss\vspace{3pt}
\[
\quad 
\lambda_{\dff z}
\hspace*{1.2em}\mbox{for}\hspace*{1.2em}
z\qff \in\qff X\qff \smallsetminus\qff Y
\qff,
\]

\vspace{-39pt}
\[
\quad 
\mu_{\dff z}
\hspace*{1.2em}\mbox{for}\hspace*{1.2em}
z\qff \in\qff Y\qff \smallsetminus\qff X
\qff,
\hspace*{1.2em}\mbox{and}\hspace*{1.2em}
\]

\vspace{-39pt}
\[
\quad 
\lambda_{\dff z}\pff +\off \alpha\dff \mu_{\dff z}
\hspace*{1.2em}\mbox{for}\hspace*{1.2em}
z\qff \in\qff
(\qff X\qff \cap\qff Y\qff)\qff \smallsetminus\qff \{\dff u\qff\}
\qff.
\]

\vspace{-9pt}
Clearly\halfff,\oss
the coefficient\qss
$\lambda_{\dff z}\pff +\off \alpha\dff \mu_{\dff z}$\qss
can\sss be positive only\dss if\dss at\dss least\sss one of\trs the coefficients\dss
$\lambda_{\dff z}\dff,\pff \mu_{\dff z}$\dss is\dss positive,\oss
and\dss negative if\dss at\dss least\sss one of\trs them\dss is\dss negative.\oss
It\dss follows\sss that\qss\vspace{3pt}
\begin{equation}
\label{elimination-1}
\quad
\omega_{\dff +}
\off \subset\off\dff 
(\qff \sigma_{\dff +}\qff \cup\qff \tau_{\dff +}\qff)
\qff \smallsetminus\qff 
\{\dff u\qff\}
\hspace*{1.2em}\mbox{and}\hspace*{1.2em}
\omega_{\dff -}
\off \subset\off\dff 
(\qff \sigma_{\dff -}\qff \cup\qff \tau_{\dff -}\qff)
\qff \smallsetminus\qff 
\{\dff u\qff\}
\qff.
\end{equation} 

\vspace{-9pt}
Similarly\halfff,\oss the coefficient\sss of\dss $z$\dss is\dss positive\dss if\pss
$\lambda_{\dff z}\qff >\qff 0$\qss
and\dss either\dss $z\qff \not\in\qff Y$\nnsp,\oss
or\qss $\mu_{\dff z}\qff \geq\qff 0$\nnsp,\oss
and\dss is\dss negative\dss if\pss
$\lambda_{\dff z}\qff <\qff 0$\qss
and\dss either\dss $z\qff \not\in\qff Y$\nnsp,\oss
or\qss $\mu_{\dff z}\qff \leq\qff 0$\nnsp.\oss
It\dss follows\sss that\qss\vspace{3pt}
\begin{equation}
\label{elimination-2}
\quad
\omega_{\dff +}
\off \supset\off 
\sigma_{\dff +}\qff \smallsetminus\qff \tau_{\dff -}
\hspace*{1.2em}\mbox{and}\hspace*{1.2em}
\omega_{\dff -}
\off \supset\off\dff 
\sigma_{\dff -}\qff \smallsetminus\qff \tau_{\dff +}
\off.
\end{equation}

\vspace{-9pt}
\mypar{Theorem.}{vector-elimination}
\emph{Let\qss $\sigma\fff,\pff \tau$\qss be\dss two circuits.\oss
Suppose\sss that\qss}\vspace{3pt}
\[
\quad
u
\off \in\off
(\dff \sigma_{\dff +}\qff \cap\qff \tau_{\dff -}\trf)
\qff \cup\qff
(\dff \sigma_{\dff -}\qff \cap\qff \tau_{\dff +}\trf)
\hspace*{1.2em}\mbox{\emph{and}}\hspace*{1.2em}
\qff
\]
 
\vspace{-39pt}
\[
\quad
v
\off \in\off
(\dff \sigma_{\dff +}\qff \smallsetminus\qff \tau_{\dff -}\trf)
\qff \cup\qff
(\dff \sigma_{\dff -}\qff \smallsetminus\qff \tau_{\dff +}\trf)
\qff.
\]

\vspace{-9pt}
\emph{Then\dss there exists a circuit\dss $\omega$\dss such\dss
that\qss the inclusions\pss \textup{(\ref{elimination-1})}\qss hold\sss
and\qss $v\qff \in\pff \underline{\omega}$\nsp.\oss}

\proof
Let\dss $X\fff,\pff Y$\dss
be minimal\dss linearly\sss dependent\sss subsets defining\dss circuits\dss
$\sigma\fff,\pff \tau$\nnsp.\oss
The elimi\-nation\dss procedure leads\sss to a relation\dss between\dss the elements of\qss
$Z
\off =\off
(\qff X\qff \cup\qff Y\qff)\qff \smallsetminus\qff \{\dff u\qff\}$\dss
such\dss that\dss the corresponding\sss signed subset\dss $\omega$\dss
satisfies\dss (\ref{elimination-1})\dss and\dss (\ref{elimination-2}).\oss
In\dss particular\halfff,\pss
$v\qff \in\pff \underline{\omega}$.\oss
But\dss it\dss may\dss happen\dss that\dss $Z$\dss is\dss not\dss minimal\sss and $\omega$
is\dss not\dss a circuit\halfff.\oss
We will\dss prove\sss that\sss some minimal\sss dependent\sss subset\sss of\trs $Z$\dss
leads\sss to a circuit\dss with\dss required\sss properties.\oss
Let\dss us\dss consider\dss relations
\vspace{3pt}
\begin{equation}
\label{elimination-3}
\quad
\sum\nolimits_{\qff z\qff \in\qff \underline{\omega}\vphantom{X^X}}\qff \alpha_{\dff z}\trf z
\off\dff =\off\dff
0
\end{equation}

\vspace{-9pt}
such\dss that\qss 
$\alpha_{\dff v}\off \neq\off 0$\qss and\dss
the sign of\dss $\alpha_{\dff z}$\dss
is\dss equal\dss to\dss $\omega\dff(\dff z\dff)$\dss
if\qss
$\alpha_{\dff z}\off \neq\off 0$\nnsp.\oss
Obviously\halfff,\oss the relation\dss resulting\dss from\dss the elimination\sss procedure
has\sss these properties.\oss
Suppose\sss that\qss (\ref{elimination-3})\qss
is\dss a\sss relation\dss with\dss the minimal\dss possible number
of\dss non-zero coefficients among\dss the relations with\dss these properties.\oss
Let\qss $V\qff \subset\pff \underline{\omega}$\qss be\sss the set\sss
of\qss $z\qff \in\qff \underline{\omega}$\qss such\dss that\qss
$\alpha_{\dff z}\off \neq\off 0$\nnsp.\oss 

We claim\dss that\trs $V$\dss is\dss a\sss minimal\dss linearly\sss dependent\sss set\halfff.\oss
Suppose\sss that\trs $V$\dss is\dss not\dss minimal.\oss
Then\dss there exists a\sss  
proper\sss subset\qss $W\pff \subset\pff V$\qss
and\sss a relation of\trs the form\vspace{1.75pt}
\begin{equation}
\label{elimination-4}
\quad
\sum\nolimits_{\qff z\qff \in\qff W\vphantom{X^X}}\pff \beta_{\dff z}\trf z
\off\dff =\off\dff
0
\end{equation}

\vspace{-9pt}
such\dss that\dss  
all\sss coefficients\sss $\beta_{\dff z}$ are non-zero.\oss
Moreover\halfff,\oss the set\sss $W$ and\dss
the relation\qss (\ref{elimination-4})\qss
can\dss be chosen\dss in\sss such a way\dss that\qss
$v\qff \not\in\pff W$\nnsp.\oss
Indeed,\oss if\qss $v\qff \in\pff W$\nnsp,\oss
then,\oss after\dss multiplying\qss (\ref{elimination-4})\qss
by\dss ${}-\qff 1$\dss if\trs necessary\halfff,\oss
we may assume\sss that\dss the signs of\dss $\alpha_{\dff v}$ and\dss
$\beta_{\dff v}$ are opposite.\oss
By\sss applying\dss the elimination\dss procedure\sss to\qss
$W\fff,\pff V$\qss and $v$ in\dss the roles of\qss
$X\fff,\pff Y$\qss and $u$ respectively\halfff,\oss
we will\dss get\sss a new\dss relation of\trs the form\qss (\ref{elimination-4}),\oss
this\sss time such\dss that\qss $v\qff \not\in\pff W$\nnsp.\oss

Let\qss 
$\gamma
\off =\off
\num{\alpha_{\dff w}\dff/\dff\beta_{\dff w}}$\qss
be\sss the minimal\dss number among\dss
$\num{\alpha_{\dff z}\dff/\dff\beta_{\dff z}}$\dss
with\qss $z\qff \in\pff W$\nnsp.\oss
Then\vspace{3pt}
\[
\quad 
\num{\gamma\qff \beta_{\dff z}}
\off =\off
\num{\alpha_{\dff w}\qff \beta_{\dff z}\dff/\dff\beta_{\dff w}}
\off \leq\off
\num{\alpha_{\dff z}}
\]

\vspace{-9pt}
for every\qss $z\qff \in\pff W$\nnsp.\oss
It\dss follows\dss that\dss the signs of\qss
$\alpha_{\dff z}\pff +\off \gamma\qff \beta_{\dff z}$\qss
and\dss $\alpha_{\dff z}$\dss are\sss 
the same for every\qss $z\qff \in\pff W$\qss
such\dss that\qss
$\alpha_{\dff z}\pff +\off \gamma\qff \beta_{\dff z}
\off \neq\off
0$\nnsp.\oss
Multiplying\qss (\ref{elimination-4})\qss
by\dss ${}-\qff 1$\dss does not\sss affects $\gamma$\nnsp.\oss
Therefore we may assume\sss that\dss the signs 
of\dss $\alpha_{\dff w}$ and\dss
$\beta_{\dff w}$ are opposite.\oss
Let\dss us multiply\qss (\ref{elimination-4})\qss by\sss
$\alpha$
and add\dss the resulting\dss
relation\dss to\qss (\ref{elimination-3}).\oss
The coefficient\sss of\qss $z\qff \in\pff V$\qss in\dss the
resulting\dss relation\dss is\vspace{2.5pt}
\[
\quad 
\alpha_{\dff z}
\hspace*{1.2em}\mbox{if}\hspace*{1.2em}
z\qff \in\qff V\qff \smallsetminus\qff W
\hspace*{1.2em}\mbox{and}\hspace*{1.2em}
\]

\vspace{-39.5pt}
\[
\quad 
\alpha_{\dff z}\pff +\off \gamma\qff \beta_{\dff z}
\hspace*{1.2em}\mbox{if}\hspace*{1.2em}
z\qff \in\qff W
\qff.
\]

\vspace{-9.5pt}
In\dss particular\halfff,\oss the coefficient\sss of\dss $z$\sss
is\dss either $0$ or\dss has\sss the same sign\dss $\omega\dff(\dff z\dff)$\dss as $\alpha_{\dff z}$
for every\qss $z\qff \in\qff V$\nnsp.\oss
Moreover\halfff,\oss the coefficient\sss of\dss $w$\sss is\dss equal\dss to $0$
and\dss the coefficient\sss of\dss $v$\sss is\dss non-zero.\oss
Hence\dss the new\dss relation has\sss 
the same properties as\qss (\ref{elimination-3}),\oss
but\dss the set\sss of\dss vectors with\dss non-zero coefficient\dss in\dss
the new\dss relation\dss is\dss properly\sss contained\dss in\dss $V$\nnsp.\oss
This contradicts\sss the minimality\sss assumption\sss about\qss (\ref{elimination-3}).\oss
The contradiction shows\sss that\dss $V$\dss is\dss indeed\sss 
a\sss minimal\dss linearly\sss dependent\sss set\halfff.\oss
The corresponding circuit\sss satisfies\sss 
the conditions of\trs the\sss theorem.\oss  \eproof

\mypar{Corollary\halfff.}{weak-vector-elimination}
\emph{Let\qss $\sigma\fff,\pff \tau$\qss be\dss two circuits.\oss 
Suppose\sss that\qss
$\sigma\off \neq\off -\qff \tau$\qss
and\qss
$\sigma\dff(\dff u\trf)\off =\off -\qff \tau\dff(\dff u\trf)$\qss
for some\dss $u$\nnsp.\oss
Then\dss there exists a circuit\dss $\omega$\dss such\dss
that\qss the inclusions\pss \textup{(\ref{elimination-1})}\qss hold.\oss}

\proof
If\trs both differences\qss
$\sigma_{\dff +}\qff \smallsetminus\qff \tau_{\dff -}$\qss
and\qss
$\sigma_{\dff -}\qff \smallsetminus\qff \tau_{\dff +}$\qss
are empty\halfff,\oss
then\qss 
$\underline{\sigma}\off \subset\off \underline{\tau}$\qss
and\dss hence either\qss
$\sigma\off =\off \tau$\qss
or\qss
$\sigma\off =\off -\qff \tau$\nnsp.\oss
The second\sss equality\dss is\dss 
explicitly\sss excluded,\oss
and\dss the first\sss one\dss is\dss
impossible because\qss
$\sigma\dff(\dff u\trf)\off =\off -\qff \tau\dff(\dff u\trf)$\nnsp.\oss
Hence at\dss least\sss one of\trs these differences\dss is\dss non-empty\dss
and\dss hence\qss Theorem\qss \ref{vector-elimination}\qss applies.\oss  \eproof

\myuppar{Oriented\dss matroids.}
An\qss \emph{oriented\dss matroid}\qss is\dss a\sss set\dss $M$\dss
together with a collection of\dss signed subsets of\dss $M$\nnsp,\oss
called\dss its\qss \emph{circuits},\oss
such\dss that\dss the above obvious properties of\dss vector configurations 
together\sss with\dss the property\sss of\qss Corollary\qss \ref{weak-vector-elimination}\qss
hold.\oss
More formally\halfff,\oss
the axioms of\dss oriented\dss matroids are\sss the following\halfff.\oss\vspace{-6pt}
\begin{enumerate}
\item[({\fff}i\fff)]\oss The empty\sss set\dss is\dss not\sss a circuit\halfff.\oss 
\item[({\fff}ii\fff)]\oss If\dss $\sigma$\dss is\dss circuit\halfff,\oss
then\qss ${}-\qff \sigma$\qss is\dss also\sss a circuit\halfff.\oss 
\item[({\fff}iii\fff)]\oss If\qss $\sigma\fff,\pff \tau$\qss are circuits
and\qss $\underline{\sigma}\off \subset\off \underline{\tau}$\nnsp,\oss
then either\qss
$\sigma\off =\off \tau$\qss
or\qss
$\sigma\off =\off -\qff \tau$\nnsp.\oss
\item[({\fff}iv\fff)]\oss If\qss $\sigma\fff,\pff \tau$\qss are circuits
and\qss
$\sigma\dff(\dff u\trf)\off =\off -\qff \tau\dff(\dff u\trf)$\qss
for some\dss $u$\nnsp,\oss
then\dss there exists a circuit\dss $\omega$\dss\\ \hspace*{0.7em}such\dss
that\qss
$\omega_{\dff +}
\off \subset\off\dff 
(\qff \sigma_{\dff +}\qff \cup\qff \tau_{\dff +}\qff)
\qff \smallsetminus\qff 
\{\dff u\qff\}$\qss
and\qss
$\omega_{\dff -}
\off \subset\off\dff 
(\qff \sigma_{\dff -}\qff \cup\qff \tau_{\dff -}\qff)
\qff \smallsetminus\qff 
\{\dff u\qff\}$\nnsp.\oss
\end{enumerate}

\vspace{-3pt}
The property\qss ({\fff}iv\fff)\qss is\dss known as\sss the\qss
\emph{weak elimination property}.\oss
Remarkably\halfff,\oss
in every oriented\dss matroid a stronger\dss form of\dss this property\dss holds.\oss
Namely\halfff,\oss
every oriented\dss matroid\dss has\sss the property\dss proved\dss in\qss
Theorem\qss \ref{vector-elimination}\qss for vector configurations.\oss
This property\dss is\dss known as\sss the\qss
\emph{strong elimination\dss property}.\oss
See\qss \cite{om},\oss Theorem\qss 3.2.5\qss for a proof\halfff,\oss
which\dss is\dss highly\dss non-trivial.\oss
Cf.\pss the situation with\dss the symmetric exchange property of\dss the usual\dss matroids,\oss
which\dss follows\sss from\dss the exchange property\halfff,\oss
but\dss this\dss is\dss far\dss from\dss being obvious.\oss
For\dss the purposes of\trs the present\dss paper\sss one can\dss take\sss the
strong elimination\dss property as another\sss axiom.\oss

Let\dss $M$\dss be an oriented\dss matroid.\oss
Suppose\sss that\qss $X\qff \subset\qff M$\qss and\qss $y\qff \in\qff M$\nnsp.\oss
The element\dss $y$\dss is\dss said\dss to belong\dss to\sss the\dss
\emph{convex hull}\pss of\trs $X$\trs
if\trs either\qss 
$y\qff \in\qff X$\qss
or\dss there exists a circuit\sss $\sigma$\sss
such\dss that\qss
$\sigma_{\dff +}\qff \subset\qff X$\qss
and\qss
$\sigma_{\dff -}\off =\off \{\dff y\qff\}$\nnsp.\oss
Naturally\halfff,\oss the convex\dss hull\sss of\trs $X$\dss is\dss
the set\sss of\dss all\qss $y\qff \in\qff M$\qss belonging\dss to\sss the convex hull of\dss $X$\nnsp.\oss
It\dss is\dss denoted\dss by\dss $\conv X$\nnsp.\oss

A subset\qss $X\qff \subset\qff M$\qss is said\dss to be\qss \emph{independent}\oss
if\dss there is no circuit\sss
$\sigma$\sss
such\dss that\qss $\underline{\sigma}\off \subset\off X$\nnsp.\oss
A maximal\sss independent\sss subset\dss is\dss called\dss a\qss \emph{basis}\qss of\dss $M$\nnsp.\oss
Clearly,\oss if\qss $B$\dss is a basis and\qss
$y\qff \in\qff M\qff \smallsetminus\qff B$\nnsp,\oss
then\dss there exists a circuit\sss
$\sigma$\sss
such\dss that\dss $\sigma\qff \subset\off B\qff \cup\qff \{\dff y\qff\}$\dss
and\dss $y\qff \in\qff \sigma$\nnsp.\oss
One can check\dss that\dss the number of\dss elements 
of\dss a basis depends only\sss on\dss $M$\nnsp.\oss

\myuppar{Todd\halfff's\qss theorem.}
\emph{Let\qss
$\sigma\fff,\pff \tau$\qss
be\sss circuits\sss of\qss an oriented\dss matroid\dss
and\dss let\qss
$w\pff \in\off \underline{\tau}\pff \smallsetminus\off \underline{\sigma}$\nsp.\oss
Suppose\sss that\dss there exists\qss 
$e\pff \in\off \underline{\sigma}\off \cap\off \underline{\tau}$\qss
such\dss that\qss
$\sigma\dff(\trf e\qff)\off =\off -\qff \tau\dff(\trf e\qff)$\nnsp.\oss
Then\dss there\dss is\dss a\sss circuit\sss $\omega$\sss
such\dss that}\vspace{3pt}
\[
\quad
\omega_{\dff +}
\off \subset\off
(\dff \sigma_{\dff +}\qff \cup\off \tau_{\dff +}\trf)
\off \smallsetminus\off
\sigma_{\dff -}\off,
\]

\vspace{-36pt}
\[
\quad
\omega_{\dff -}
\off \subset\off
(\dff \sigma_{\dff -}\qff \cup\off \tau_{\dff -}\trf)
\off \smallsetminus\off
\sigma_{\dff +}
\off,
\]

\vspace{-36pt}
\[
\quad
w\pff \in\off \underline{\omega}\off,
\hspace{1.2em}\mbox{\emph{and}}\hspace{1.5em}
\omega\dff(\trf w\qff) 
\off\qff =\off\qff
\tau\dff(\trf w\qff)
\qff.
\]

\vspace{-9pt}
\emph{Moreover\halfff,\oss
if\qff\oss
$\underline{\tau}
\qff\off \subset\qff\off 
\underline{\sigma}\off \cup\qff \{\qff w\qff\}$\nnsp,\qff\oss
then\dss such\sss a\sss circuit\dss $\omega$\dss is\dss unique\sss and}\qss 
$\underline{\sigma}\off \smallsetminus\qff\off \underline{\tau}
\off\off \subset\pff\off
\underline{\omega}$\nnsp.\oss

\proof
See\qss Appendix\qss for\sss a\sss proof\trs following\pss M.\dss Todd\pss
\cite{t},\oss Theorem\qss 4.2.\oss  \eproof

\mysection{Oriented\qss matroid\qss colorings\fff:\oss the\qss non-degenerate\qss case}{scarf-oriented-matroids}

\myuppar{The oriented\dss matroid\qss framework.}
We need one more notion\dss related\dss to oriented\dss matroids.\oss
An oriented\dss matroid\sss $M$\sss is\dss said\dss to be\qss \emph{acyclic}\pss
if\trs there are no circuits
$\sigma$
such\dss that\qss $\sigma_{\dff -}\off =\off \varnothing$\nnsp.\oss
Let\dss $M$\dss be an acyclic oriented\dss matroid\sss
and\dss let\vspace{1.75pt}
\[
\quad
B\off =\off
\left\{\qff
v_{\fff 0}\fff,\pff v_{\fff 1}\fff,\pff \ldots\fff,\pff v_{\fff n}
\pff\right\}
\]

\vspace{-10.25pt}
be a basis of\dss $M$\nnsp.\oss 
Let\qss $b\qff \in\qff M\qff \smallsetminus\qff B$\dss
and suppose\sss that\sss $b$\sss belongs\sss to\sss the convex hull\sss of\dss $B$\nnsp.\oss
In\dss this section\dss
we\dss will\sss also assume\sss that\dss the pair\dss $M\fff,\pff b$\dss is\qss
\emph{non-degenerate}\pss in\dss the sense\sss that\sss $b$
does\sss not\dss belong\dss to\sss the convex\dss hull\sss of\qss
$<\qff n\qff +\qff 1$\dss elements\sss of\qss
$M\qff -\qff b$\nnsp.\oss
This\dss is\dss a\sss sort\sss of\dss general\dss position
assumption,\oss
and\qss if\qss $M$\dss corresponds\sss to a\sss vector configuration,\oss
it\dss holds,\pss in\dss particular\halfff,\oss 
if\qss $b$\sss is\dss not\sss a\sss linear combination of\trs
$<\qff n\qff +\qff 1$\dss 
elements\sss of\qss
$M\qff -\qff b$\nnsp.\oss

\mypar{Lemma.}{not-basis}
\emph{Let\qss $X\off \subset\off M$\nnsp.\oss
If\oss $\num{X}\off =\off n\qff +\qff 1$\qss
and\dss $X$\dss is\dss not\dss a\sss basis,\oss
then\dss $b$\dss does not\dss belong\dss to\sss the
convex\dss hull\sss of\pss $X$\nnsp.\oss}

\proof
Suppose\sss that\sss $b$\sss belongs\sss to\sss the convex\dss hull\sss of\trs $X$\nnsp,\oss
i.e.\qss there\dss is\dss a circuit\sss $\tau$ such\dss that\qss
$\tau_{\dff +}\qff \subset\pff X$\qss
and\qss
$\tau_{\dff -}\off =\off \{\trf b\qff\}$\nnsp.\oss
Then\dss
$\tau_{\dff +}$\dss has\qss $\geq\qff n\qff +\qff 1$\qss elements
and\dss hence\qss $\tau_{\dff +}\off =\off X$\nnsp.\oss

If\dss $X$\dss is\dss independent\halfff,\oss
then $X$ is\dss contained\dss in\sss some basis
with\qss $>\qff n\qff +\qff 1$\qss elements.\oss
But\dss $B$\dss is\dss a\sss basis consisting of\dss $n\qff +\qff 1$\sss elements,\oss
and\dss the number of\dss elements of\dss a\sss basis depends only\sss on\dss $M$\nnsp.\oss
It\dss follows\dss that\sss $X$ is\dss not\dss independent\halfff,\oss
i.e.\qss there\dss is\dss a\sss circuit\sss $\sigma$\sss such\dss that\pss
$\underline{\sigma}\pff \subset\pff X$\nnsp.\oss
It\dss follows\dss that\qss
$\underline{\sigma}\off \subset\off X\off \subset\off \underline{\tau}$\nnsp.\oss
By\dss the axiom\qss ({\fff}iii\fff)\qss either\qss
$\sigma\off =\off \tau$\qss
or\qss
$\sigma\off =\off -\qff \tau$\nnsp.\oss
But\qss $b\qff \in\pff \underline{\tau}$\qss
and\qss
$b\qff \not\in\pff \underline{\sigma}$\nnsp.\oss
The contradiction shows\sss that\sss $b$\sss does not\dss belong\dss
to\sss the convex\dss hull\sss of\trs $X$\nnsp.\oss  \eproof

\mypar{Lemma.}{exchange}
\emph{If\pss
$w\qff \in\qff M\qff \smallsetminus\qff B$\qss
and\dss $w\off \neq\off b$\nnsp,\oss
then\dss there is a unique\dss
$i\qff \in\qff \{\qff 0\fff,\pff 1\fff,\pff \ldots\fff,\pff n\qff\}$\dss
such\dss that\dss $b$\dss belongs\sss to\sss the convex hull\dss of\oss
$B\qff -\qff v_{\fff i}\qff +\qff w$\nnsp.\oss}

\proof
Since\sss $b$\sss belongs\sss to\sss the convex\dss hull\sss of\trs $B$\dss
and\dss $M\fff,\pff b$\dss is\dss non-degenerate,\oss 
there\dss is\dss a\sss circuit\sss $\sigma$ 
such\dss that\qss
$\sigma_{\dff +}\off =\off B$\qss
and\qss
$\sigma_{\dff -}\off =\off \{\dff b\trf\}$\nnsp.\oss
Since\dss $B$\dss is\dss a\dss basis,\oss 
there\dss is\dss a\sss circuit\sss $\tau$
such\dss that\qss
$w\pff \in\pff \underline{\tau}
\qff\off \subset\dff\off
B\qff +\qff w$\nnsp.\oss
Replacing\dss $\tau$\dss by\dss ${}-\qff \tau$\dss if\trs necessary,\oss
we may assume\sss that\qss $\tau\dff(\dff w\dff)\off =\off +$\nnsp.

We would\dss like\sss to apply\qss Todd's\qss theorem\dss to\qss
$\sigma\fff,\off \tau$\nnsp,\oss and\dss $w$\nnsp.\oss
Let\dss us\sss check\dss its\sss assumptions.\oss
Obviously\halfff,\oss
$w\pff \in\pff 
\underline{\tau}
\off \smallsetminus\off 
\underline{\sigma}$\pss
and\pss
$\underline{\tau}
\dff\off \subset\qff\off 
\underline{\sigma}\off \cup\pff \{\trf w\qff\}$\nnsp.\oss
Also,\oss\vspace{3.625pt}
\[
\quad
\underline{\sigma}\off \cap\off \underline{\tau}
\dff\off =\off\dff 
\underline{\tau}\off -\pff w
\dff\off =\off\dff 
\underline{\tau}\off \cap\pff B
\qff.
\]

\vspace{-8.375pt}
Since\dss $M$\dss is\dss acyclic,\oss $\tau_{\dff -}$\dss is\dss non-empty\halfff.\oss
Since\dss
$w\pff \in\pff \tau_{\dff +}$\nnsp,\oss
the\sss set\dss $\tau_{\dff -}$\dss is\dss contained\dss in\dss $B$\dss
and\dss hence\sss
the intersection\qss
$\tau_{\dff -}\qff \cap\qff \sigma_{\dff +}
\off =\off
\tau_{\dff -}\qff \cap\qff B$\qss is\dss non-empty.\oss
If\qss 
$e\qff \in\qff 
\tau_{\dff -}\qff \cap\qff \sigma_{\dff +}$\nsp,\oss
then\qss
$e\pff \in\pff
\underline{\sigma}\off \cap\off \underline{\tau}$\qss
and\qss
$\sigma\dff(\trf e\qff)\off =\off -\qff \tau\dff(\trf e\qff)$\nnsp.\oss
It\dss follows\dss that\sss all\sss assumptions of\qss Todd's\dss the\-o\-rem hold.\oss

By\qss Todd's\qss theorem\dss
there is a circuit\dss $\omega$\dss
such\dss that\vspace{4.5pt}
\[
\quad
\omega_{\dff +}
\off \subset\off\qff
(\dff \sigma_{\dff +}\qff \cup\off \tau_{\dff +}\trf)
\off \smallsetminus\off
\sigma_{\dff -}
\off =\off\dff 
B\qff +\qff w
\qff,
\]

\vspace{-34.5pt}
\[
\quad
\omega_{\dff -}
\off \subset\off\qff
(\dff \sigma_{\dff -}\qff \cup\off \tau_{\dff -}\trf)
\off \smallsetminus\off
\sigma_{\dff +}
\off =\off\dff 
\{\trf b\qff\}
\qff,
\]

\vspace{-34.5pt}
\[
\quad
w\pff \in\off \underline{\omega}\off,
\hspace{1.2em}\mbox{and}\hspace{1.5em}
\omega\dff(\qff w\qff) 
\off\qff =\off\qff
\tau\dff(\qff w\qff)
\off\qff =\off\qff
+
\qff.
\]

\vspace{-7.5pt}
If\pss $b\pff \not\in\pff \underline{\omega}$\nsp,\oss
then\pss $\underline{\omega}\dff\off =\off \omega_{\dff +}$\nsp,\oss
contrary\qss to\sss the assumption\dss that\dss $M$\dss is\dss acyclic.\dff\oss
If\pss $\underline{\omega}\off \cap\off B$\dss consists of\pss $\leq \qff n$\qss
elements,\oss then\dss $\omega$\dss is a circuit\sss containing\dss $b$\dss
and consisting of\qss $\leq\qff n\qff +\qff 1$\qss elements.\oss
This also contradicts\sss to\sss the assumptions of\trs the\sss
oriented\dss matroid\dss framework.\oss
Therefore\dss $\underline{\omega}$\dss has\sss the form\qss
$B\qff -\qff v_{\fff i}\qff +\qff w\qff -\qff b$\nnsp.\oss
This proves\sss that\sss $i$ with\dss the required\dss properties exists.\oss

In order\dss to prove\dss the uniqueness,\oss
suppose\sss that\qss
$\omega_{\dff +}
\off =\off
B\qff -\qff v_{\fff i}\qff +\qff w$\qss
and\qss
$\omega_{\dff -}
\off =\off
\{\trf b\qff\}$\qss
for some circuit\sss $\omega$\nnsp.\oss
Then\sss $\omega$\sss has all\dss properties from\qss Todd's\qss theorem.\oss
Since\pss
$\underline{\tau}
\dff\off \subset\qff\off 
\underline{\sigma}\off \cup\pff \{\trf w\qff\}$\nnsp,\oss 
Todd's\qss theorem\dss implies\sss the uniqueness of\dss such\dss $\omega$\dss
and\dss hence\sss the uniqueness of\dss $i$\nnsp.\oss  \eproof

\mypar{Lemma.}{two-faces}
\emph{Let\qss
$i\qff \in\qff \{\qff 0\fff,\pff 1\fff,\pff \ldots\fff,\pff n\qff\}$\nnsp.\oss
Suppose\sss that\qss $D\qff \subset\qff M$\qss
and\dss $b$\dss does not\dss belong\sss to\sss the
convex\dss hull\dss of\pss $D$\nnsp.\oss
Then\dss the number of\dss elements\qss 
$d\qff \in\qff D$\qss such\dss that\pss
$D\qff -\qff d\qff +\qff v_{\fff i}$\qss
is\dss a\sss basis containing\dss $b$\dss in\dss its convex hull\dss
is\dss either\dss $2$\dss or\dss $0$\nnsp.\oss}

\proof
Suppose\sss that\qss
$D\qff -\qff d\qff +\qff v_{\fff i}$\qss
is\dss a\sss basis containing\dss $b$\dss in\dss its\sss convex hull.\oss
If\trs we replace\dss $B$\dss by\dss the basis\qss
$D\qff -\qff d\qff +\qff v_{\fff i}$\nsp,\oss
all\sss assumptions of\trs the oriented\dss mathroid\dss framework\sss still\dss hold.\oss
By applying\qss Lemma\qss \ref{exchange}\qss to\sss the basis\qss 
$D\qff -\qff d\qff +\qff v_{\fff i}$\qss
in\dss
the role of\dss $B$\dss and\dss to\qss $w\off =\off d$\nnsp,\oss
we see\sss that\dss there is a unique element\qss
$e\qff \in\qff D\qff -\qff d\qff +\qff v_{\fff i}$\qss
such\dss that\dss $b$\dss belongs\sss to\sss the convex hull of\qss\vspace{3pt}
\[
\quad
D\qff -\qff d\qff +\qff v_{\fff i}\qff -\qff e\qff +\qff d
\off =\off
D\qff -\qff e\qff +\qff v_{\fff i}
\qff.
\]

\vspace{-9pt}
If\qss $e\off =\off v_{\fff i}$\nsp,\oss
then\dss this set\dss is equal\dss to\dss $D$\nnsp.\oss
But\dss $D$\dss does not\sss contain\dss $b$\dss in\dss its convex hull\dss by\dss
the assumption of\dss the\sss lemma.\oss
Hence\qss $e\off \neq\off v_{\fff i}$\qss
and\qss $D\qff -\qff e\qff +\qff v_{\fff i}$\qss
is\dss a\sss basis containing\dss $b$\dss in\dss its\sss convex hull.\oss
It\dss follows\dss that\trs if\qss there\dss is\dss at\dss 
least\dss one element\qss $d\qff \in\qff D$\qss
with\dss the stated\dss property\halfff,\oss
then\dss there are exactly\dss two such elements.\oss
This proves\sss the\sss lemma.\oss  \eproof

\myuppar{Colorings\sss by\sss elements of\dss oriented\dss matroids.}
Let\sss $\mathcal{D}$\sss be\sss a\sss chain-simplex\dss based\sss on\dss $I$\dss 
(see\qss Section\qss \ref{pseudo-simplices})\qss
and\trs $X\off =\off V_{\dff \mathcal{D}}$\dss be\sss the union of\trs the sets
of\dss vertices of\dss complexes $\mathcal{D}\dff(\trf J\dff)$ with\dss
$J\qff \subset\pff I$\nnsp.\oss 
Suppose\sss that\qss $\num{I}\off =\off \num{B}$\nnsp.\oss
Then\dss we can assume\sss that\qss
$I\off =\off \{\qff 0\fff,\pff 1\fff,\pff \ldots\fff,\pff n\qff\}$\qss
and\dss that\dss the map\qss 
$v_{\fff \bullet}\dff \colon\dff
i\qff \longmapsto\qff v_{\fff i}$\qss
is\dss a\sss bijection\qss $I \ttoo B$\nnsp.\oss

Let\dss$A\off =\off M\qff -\qff b$\dss and\dss
let\dss $\Delta\dff(\dff A\dff)$\dss and\dss $\Delta\dff(\dff B\dff)$\dss 
be simplicial\sss complexes
having\halfff,\oss respectively\halfff,\pss $A$ and\sss $B$  
as\sss their\sss sets
of\dss vertices and\sss all\sss subsets of\dss $A$ and\sss $B$
as simplices.\oss 
Let\dss $\partial\dff \Delta\dff(\trf B\trf)$\dss be\sss 
the simplicial\sss complex having\dss $B$\dss as\sss its\sss set\sss
of\dss vertices and\sss all\dss proper
subsets of\trs $B$\dss as\sss simplices.\oss
The pair\qss $\Delta\dff(\dff A\dff)\fff,\pff \partial\dff \Delta\dff(\trf B\trf)$\dss
is\dss an analogue of\trs the pair\qss
$\Delta\dff(\dff I\dff)\fff,\pff \partial\dff \Delta\dff(\trf I\trf)$\qss
in\dss the classical situation.\oss

A\qss \emph{matroid\sss coloring}\pss of\trs $\mathcal{D}$\dss is\dss defined as an arbitrary\dss map\qss
$X\off =\off V_{\dff \mathcal{D}}\qff \ttoo\qff A$\nnsp.\oss
A\qss matroid\sss coloring\dss 
$c\dff \colon\dff
X\dff \ttoo\dff A$\dss 
for every\sss subset\trs $J\qff \subset\pff I$\dss induces
a\sss simplicial\dss map\dss
$\mathcal{D}\dff (\trf J\dff)\dff \ttoo\dff \Delta\dff (\trf A\trf)$\nnsp.\oss
The map $c$
canonically\dss extends\sss to a map\qss
$\varphi\qff \colon\dff
X\qff \cup\pff I
\qff \ttoo\qff A$\qss
equal\dss to\dss $v_{\fff \bullet}$\dss on\dss $I$\nnsp.\oss
Let\sss $\mathcal{E}$\dss be\sss the envelope\sss of\dss $\mathcal{D}$\dss in\dss
the sense of\qss Section\qss \ref{pseudo-simplices}.\oss
Identification of\trs $I$\sss with\dss $B$\dss by\dss the map\dss $v_{\fff \bullet}$\dss
turns $\varphi$ into\sss an extension of\dss $c$\sss in\dss the sense of\qss 
Section\qss \ref{pseudo-simplices}.\oss
In any\sss case,\pss $\varphi$ induces 
a simplicial\dss map\qss
$\mathcal{E}\dff (\trf I\trf)\dff \ttoo\dff \Delta\dff (\trf A\trf)$\nnsp,\oss
which,\oss in\dss turn,\oss
induces an\dss isomorphism\qss
$\partial\dff \Delta\dff(\trf I\trf)
\qff \ttoo\qff 
\partial\dff \Delta\dff(\trf B\trf)$\nnsp.\oss

The action of\dss $\varphi$ on $n$\dnsp-simplices can\dss be described
as follows.\oss 
If\dss $\sigma$ is\dss an $n$\dnsp-simplex of\dss $\mathcal{E}\dff (\trf I\trf)$\dss
and\qss
$\sigma\off =\off \tau\dff *\dff (\pff I\pff \smallsetminus\pff C\qff)$\nnsp,\oss
where\oss $C\off \subset\off I$\oss
and\sss $\tau$ is\dss a\sss simplex of\trs $\mathcal{D}\dff(\trf C\trf)$\nnsp,\oss
then\vspace{2pt}
\begin{equation}
\label{extended-image}
\quad 
\varphi\dff(\dff \sigma\trf)
\off =\off\dff
c\trf(\dff \tau\trf)
\off \cup\off\dff 
v_{\fff \bullet}\trf(\pff I\pff \smallsetminus\pff C\qff)
\qff.
\end{equation}

\vspace{-9pt}
By\qss Lemma\qss \ref{top-simplices}\qss every $n$\dnsp-simplex $\sigma$
of\dss $\mathcal{E}\dff (\trf I\trf)$\dss has\sss such\dss form
and\dss hence\qss (\ref{extended-image})\qss applies\sss to $\sigma$\nnsp.\oss

\myuppar{Cochains\dss $\delta_{\fff i}$\nsp.}
These\sss cochains are analogues of\trs the cochains\dss $\delta_{\fff i}$\dss
from\dss the cochain-based\dss proof\dss of\qss Sperner's\dss lemma\dss
in\qss Section\qss 9\qss of\pss \cite{i2}.\oss
Let\qss $i\qff \in\qff I$\nnsp.\oss
The cochain\dss $\delta_{\fff i}$\dss is\dss an
$(\dff n\dff -\dff 1\dff)$\dnsp-cochain\sss
of\dss $\Delta\dff (\trf A\trf)$\dss defined as\sss follows.\oss
Let\sss $\varepsilon$\sss be\dss an
$(\dff n\dff -\dff 1\dff)$\dnsp-simplex\dss 
of\dss $\Delta\dff(\dff A\dff)$\nnsp.\oss 
Then\vspace{3pt}
\[
\quad
\delta_{\dff i}\trf(\dff \varepsilon\trf)
\off =\off\dff
1
\hspace{1.2em}\mbox{if}\hspace{1.2em}
b\pff \in\pff \conv{\dff \varepsilon\qff +\qff v_{\fff i}\dff}
\hspace{1.2em}\mbox{and}\hspace{1.2em}
\]

\vspace{-39pt}
\[
\quad
\delta_{\dff i}\trf(\dff \varepsilon\trf)
\off =\off\dff
0
\hspace{1.2em}\mbox{otherwise.}\hspace{1.2em}
\]

\vspace{-9pt}
By\qss Lemma\qss \ref{not-basis}\pss if\pss
$b\pff \in\pff \conv{\dff \varepsilon\qff +\qff v_{\fff i}\dff}$\nnsp,\oss
then\qss $\varepsilon\qff +\qff v_{\fff i}$\qss
is\dss a\sss basis.

\mypar{Lemma.}{computing-coboundary}
\emph{Let\dss $\sigma\off \subset\off A$\qss be\sss an\dss
$(\dff n\dff -\dff 1\dff)$\dnsp-simplex\dss
of\qss $\Delta\dff(\dff A\dff)$\nnsp.\oss 
Then}\qss\vspace{3pt}
\[
\quad
\partial^{\fff *}\dff \delta_{\dff i}\trf(\dff \sigma\trf)
\off =\off\dff
1
\]

\vspace{-9pt}
\emph{if\qss and\dss only\qss if\qss $\sigma$\dss is\dss a\sss basis containing\dss $b$\dss in\dss its\sss
convex\dss hull.\oss}

\proof
By\dss the definition of\trs the coboundary\sss operator\halfff,\oss\vspace{3pt}
\[
\quad
\partial^{\fff *}\dff \delta_{\dff i}\trf(\dff \sigma\trf)
\off =\off
\sum_{\dff v\qff \in\qff \sigma\vphantom{X^x}}\qff
\delta_{\dff i}\trf(\trf \sigma\qff -\qff v\qff)
\off.
\]

\vspace{-9pt}
Let\dss $\sigma$\dss be\dss a\sss basis containing\dss $b$\dss in\dss
its convex hull.\oss
Suppose\sss first\dss that\qss $v_{\fff i}\qff \not\in\qff \sigma$\nnsp.\oss
Then\dss Lemma\qss \ref{exchange}\qss
with\dss $\sigma$\dss in\dss the role of\dss $B$\dss
and\dss $v_{\fff i}$\dss in\dss the role of\dss $w$\dss
implies\dss that\dss there is a unique element\qss $v\qff \in\qff s$\qss
such\dss that\qss $\sigma\qff -\qff v\qff +\qff v_{\fff i}$\qss
is\dss a\sss basis containing\dss $b$\dss in\dss its convex hull.\oss
Therefore in\dss this case\qss
$\partial^{\fff *}\dff \delta_{\dff i}\trf(\dff \sigma\trf)
\off =\off
1$\nnsp.\oss
Suppose now\dss that\qss $v_{\fff i}\qff \in\qff \sigma$\nnsp.\oss
Then\qss $\sigma\qff -\qff v\qff +\qff v_{\fff i}$\qss
is\dss defined only\dss if\qss $v\off =\off v_{\fff i}$\qss
and\dss hence\qss
$\delta_{\dff i}\trf(\trf \sigma\qff -\qff v\trf)\off =\off\dff 1$\qss
if\trs and\dss only\trs if\qss $v\off =\off v_{\fff i}$\nsp.\oss
Hence\sss in\dss this case\qss
$\partial^{\fff *}\dff \delta_{\dff i}\trf(\dff \sigma\trf)
\off =\off
1$\qss
also.\oss

By\qss Lemma\qss \ref{not-basis},\oss if\dss $\sigma$\sss is\dss not\sss a\sss basis,\oss
then\dss $b$\dss does not\dss belong\dss to\sss the convex\dss hull\sss of\dss $\sigma$\nnsp.\oss
Therefore\sss it\dss remains\sss to consider\dss the case when\dss $\sigma$\dss
does\sss not\sss contain\dss $b$\dss in\dss its\dss convex\dss hull.\oss
If\dss all\sss summands of\trs the sum\sss representing\qss
$\partial^{\fff *}\dff \delta_{\dff i}\dff(\dff \sigma\trf)$\qss
are\dss $0$\nnsp,\oss
then\qss
$\partial^{\fff *}\dff \delta_{\dff i}\trf(\dff \sigma\trf)
\off =\off
0$\nnsp.\oss
Suppose\sss that\sss one of\trs the summands\dss is\dss non-zero,\oss
i.e.\qss
$\delta_{\dff i}\trf(\dff \sigma\qff -\qff v\qff)
\off =\off
1$\qss
for some\qss $v\qff \in\qff \sigma$\nnsp.\oss
Then\qss
$\sigma\qff -\qff v\qff +\qff v_{\fff i}$\qss
is\dss a\sss basis containing\dss $b$\dss in\dss its convex\dss hull.\oss
Since now\sss $\sigma$\sss does\sss not\sss contain\dss $b$\dss
in\dss its convex hull,\oss
Lemma\qss \ref{two-faces}\qss implies\sss that\sss in\dss this case\sss
there are exactly\sss $2$ elements\qss $v\qff \in\qff \delta$\qss
with\dss this property\halfff,\oss and\dss hence exactly\sss $2$\sss
elements\qss $v\qff \in\qff \sigma$\qss
such\dss that\qss 
$\delta_{\dff i}\trf(\trf \sigma\qff -\qff v\qff)
\off =\off
1$\nnsp.\oss
It\dss follows\dss that\qss 
$\partial^{\fff *}\dff \delta_{\fff i}\dff(\dff \sigma\trf)
\off =\off
0$\pss
if\qss
$\sigma$\dss is\dss not\sss a\sss basis containing\dss $b$\dss in\dss
its convex hull.\oss  \eproof

\mypar{The main\dss theorem\dss for\dss matroid colorings\qss 
({\fff}the non-degenerate case).}{main-non-degenerate} 
\emph{Let\trs 
$c$\dss 
be\dss a coloring.\oss
Then\dss 
there exist\sss a\sss non-empty\sss subset\qss $C\qff \subset\pff I$\qss
and\dss a\sss simplex\sss $\tau$\sss of\pss $\mathcal{D}\dff(\trf C\trf)$\dss
such\dss that\qss
$c\trf(\dff \tau\trf)
\qff \cup\pff 
v_{\fff \bullet}\trf(\qff I\qff \smallsetminus\pff C\qff)$\dss
is\dss a\sss basis containing $b$\ in\dss its convex\dss hull.\oss
The number of\trs such\sss $\tau$\sss is\dss odd.}

\proof
We\sss will\dss use\sss the action of\trs the extended\dss map\dss $\varphi$\dss
on\dss chains and\sss cochains\sss in\sss a\sss manner similar\dss to\sss the\sss cochain-based\dss proof\dss
of\pss Sperner's\qss lemma\dss in\qss \cite{i2}.\oss
Let\dss us\sss fix some\qss $i\qff \in\qff I$\nnsp.\oss
Let\sss $e$\sss be\sss the number
of\dss $n$\dnsp-simplices\dss $\sigma$\dss of\qss
$\mathcal{E}\dff (\trf I\trf)$\qss
such\dss that\dss $\varphi\dff(\dff \sigma\trf)$\sss is\dss a\sss basis
containing\dss $b$\dss in\dss its convex\dss hull.\oss
In\dss view\sss of\qss (\ref{extended-image})\qss
the conclusion of\trs
the\sss theorem\dss means\sss that\sss $e$\sss is\dss odd.\oss 
Recall\trs that\dss
$\mathcal{E}\trf\fclass{\trf I\trf}$\dss
is\dss the sum of\dss all\sss $n$\dnsp-simplices of\dss $\mathcal{E}\dff (\trf I\trf)$\nnsp.\oss
By\qss Lemma\qss \ref{computing-coboundary}\vspace{3pt}
\begin{equation*}
\quad
\partial^{\fff *}\dff \delta_{\fff i}\qff
\bigl(\pff
\varphi_{\fff *}\dff
\bigl(\qff \mathcal{E}\trf\fclass{\trf I\trf}\qff\bigr)
\pff\bigr)
\end{equation*}

\vspace{-9pt}
is\dss equal\dss to\sss $e$\sss modulo\dss $2$\nnsp.\oss 
One\sss the other\dss hand,\oss in view of\trs the identification of\trs $I$\dss with\dss $B$\qss
Lemma\qss \ref{d-of-image}\qss implies\sss that\qss
$\varphi_{\fff *}\dff
\bigl(\qff \mathcal{E}\trf\fclass{\trf I\trf}\qff\bigr)
\off =\off
B$\nnsp,\oss
where\dss $B$\dss is\dss considered as\sss an
$n$\dnsp-simplex of\dss
$\Delta\dff(\dff A\dff)$\nnsp.\oss
It\dss follows\dss that\sss $e$\sss is\dss equal\dss  modulo\dss $2$\dss to\vspace{3pt}
\[
\quad
\partial^{\fff *}\dff \delta_{\fff i}\qff(\trf B\trf)
\off =\off
\delta_{\fff i}\qff(\qff \partial\dff B\trf)
\qff.
\]

\vspace{-9pt}
It\dss remains\sss to prove\sss that\qss
$\delta_{\fff i}\qff 
\bigl(\qff
\partial\trf B
\qff\bigr)
\off =\off\dff 
1$\nnsp.\oss
Every\sss $(\dff n\dff -\dff 1\dff)$\dnsp-face of\trs $B$\dss
has\sss the form\qss $B\qff -\qff v_{\dff k}$\nnsp,\oss
where\qss $k\qff \in\qff B$\nnsp.\oss
Since\qss $B\qff -\qff v_{\dff k}\qff +\qff v_{\dff i}$\qss is\dss not\sss defined\dss
if\pss $k\off \neq\off i$\qss
and\qss
$B\qff -\qff v_{\dff i}\qff +\qff v_{\dff i}\off =\off B$\nnsp,\oss
we see\sss that\pss
$\delta_{\dff i}\trf(\trf B\qff -\qff v_{\dff k}\trf)
\off =\off
1$\pss
if\oss
$k\off =\off i$\pss
and\pss
$\delta_{\dff i}\trf(\trf B\qff -\qff v_{\dff k}\trf)
\off =\off
0$\pss
otherwise.\oss
It\dss follows\dss that\qss
$\delta_{\dff i}\dff(\trf \partial\dff B\trf)
\off =\off\dss
1$\nnsp.\oss
This\sss proves\sss the\dss theorem.\oss  \eproof

\myuppar{A cohomological\dss interpretation\sss of\trs the proof\halfff.}
In\dss the classical\sss situation we are given a\sss triangulation\dss $T$\dss
of\dss an $n$\dnsp-simplex $\Delta$ and a simplicial\dss map\dss
$\varphi\dff \colon\qff T\qff \ttoo\qff \Delta$\dss and\dss
we are interested\sss in\dss the number of\dss simplices of\trs $T$\dss mapped\dss
by\sss $\varphi$\sss onto $\Delta$ considered\dss modulo $2$\nnsp.\oss
A way\dss to count\dss them\dss is\dss to consider\dss the $n$\dnsp-cochain $\delta$
of\dss $\Delta$ equal\dss to $1$ on\dss the only $n$\dnsp-simplex of\dss $\Delta$
and\dss its image\dss $\varphi^{\fff *}(\dff \delta\dff)$\dss with\dss is\dss
equal\dss to $1$ exactly on $n$\dnsp-simplices of\dss $T$\dss mapped\dss by $\varphi$
to $\Delta$\nnsp.\oss
It\dss turns out\dss to be useful\dss to represent\sss $\delta$ as a coboundary\halfff.\oss
Namely\halfff,\pss $\delta\off =\off \partial^{\fff *}\dff \delta_{\dff i}$\nsp,\oss
where $\delta_{\dff i}$ is\dss the cochain equal\dss to $1$ on\dss the $i${\nnsp}th\dss
$(\fff n\dff -\dff 1\fff)$\dnsp-face of\dss $\Delta$\nnsp.\oss
Moreover\halfff,\oss the picture\dss is\dss clarified\dss by\dss
passing\dss from cochains\sss to\sss their cohomology classes.\oss
See\qss \cite{i2},\oss Section\qss 10\qss for\dss the details.\oss
Now\dss we have a simplicial\dss map\qss
$\varphi\dff \colon\dff
\mathcal{E}\dff (\trf I\trf)\dff \ttoo\dff \Delta\dff (\trf A\trf)$\qss
and\dss we are\sss interested\dss in\dss the number of\dss simplices
of\dss $\mathcal{E}\dff (\trf I\trf)$
mapped\dss by\sss $\varphi$\sss onto a\sss basis containing $b$
in\dss its\sss convex\dss hull.\oss 
As in\dss the classical\sss situation,\oss
a way\dss to count\dss them\dss is\dss to consider\dss the $n$\dnsp-cochain $\delta$
of\dss $\Delta\dff (\trf A\trf)$\sss equal\dss to $1$ on each $n$\dnsp-simplex\dss
which\dss is\dss a basis containing $b$ in\dss its\sss convex\dss hull\sss
and\dss equal\dss to $0$ on other $n$\dnsp-simplices.\oss
Then\dss $\varphi^{\fff *}(\dff \delta\dff)$\dss is\dss equal\sss to $1$ on exactly\dss
the $n$\dnsp-simplices we would\dss like\sss to count\halfff.\oss
As\dss before,\oss it\dss is\dss useful\dss to
represent\sss $\delta$ as a coboundary\halfff.\oss
Lemma\qss \ref{computing-coboundary}\qss tell\dss us\sss that\qss
$\delta\off =\off \partial^{\fff *}\dff \delta_{\dff i}$\nsp,\oss
where now $\delta_{\dff i}$ is\dss the cochain defined\dss before\sss this\sss lemma.\oss
This $\delta_{\dff i}$ is\dss an analogue of\trs the classical\sss $\delta_{\dff i}$
in\dss the sense\sss that\dss it\dss works\sss in\dss the same way\halfff.\oss
But\sss a naive analogue 
would\dss be\sss the cochain $\gamma_{\fff i}$ equal\dss to $1$ on\qss 
$B\qff -\qff v_{\fff i}$\qss and\dss to $0$ on
all\sss other $(\fff n\dff -\dff 1\fff)$\dnsp-simplices.\oss
Remarkably\halfff,\oss after\dss passing\dss to cohomology\sss groups\sss these\sss
two analogues\sss of\trs the classical\sss co\-chain\sss $\delta_{\dff i}$\sss 
turn out\dss to be equiv\-a\-lent\halfff.\oss
Let\dss us\sss consider $\varphi$ as a map of\dss pairs\qss\vspace{2.625pt}
\[
\quad
\bigl(\qff \mathcal{E}\dff (\trf I\trf)\fff,\pff 
\partial\dff \Delta\dff(\trf I\trf) \qff\bigr)
\off \ttoo\off
\bigl(\qff \Delta\dff (\trf A\trf)\fff,\pff 
\partial\dff \Delta\dff(\trf B\trf) \qff\bigr)
\]

\vspace{-9.365pt}
and\dss let\dss $\partial\dff \varphi$\dss be\sss the induced\dss map\qss
$\partial\dff \Delta\dff(\trf I\trf)
\qff \ttoo\qff
\partial\dff \Delta\dff(\trf B\trf)$\nnsp.\oss
Consider\dss the following diagram.\oss\vspace{2pt}
\begin{equation*}
\quad
\begin{tikzcd}[column sep=spec, row sep=spe]\dis
H^{\fff n\dff -\dff 1}\dff\bigl(\qff \partial\dff \Delta\dff(\trf B\trf) \qff\bigr)\dff\off
\arrow[r, 
"\dis \partial^{*\nsp *}"]
\arrow[d, 
"\dis (\dff\partial\dff \varphi\dff)^{ *\nsp *}"']
& 
\off H^{\fff n}\dff\bigl(\qff \Delta\dff (\trf A\trf)\fff,\pff 
\partial\dff \Delta\dff(\trf B\trf) \qff\bigr)
\arrow[d,   
"\dis \varphi^{*\nsp *}"']
\\
H^{\fff n\dff -\dff 1}\dff\bigl(\qff \partial\dff \Delta\dff(\trf I\trf) \qff\bigr)\off
\arrow[r,  
"\dis \partial^{*\nsp *}"]
& 
\off H^{\fff n}\dff\bigl(\qff \mathcal{E}\dff (\trf I\trf)\fff,\pff 
\partial\dff \Delta\dff(\trf I\trf) \qff\bigr)
\end{tikzcd}
\end{equation*}

\vspace{-10pt}
Since $\partial\dff \Delta\dff(\trf B\trf)$ and $\partial\dff \Delta\dff(\trf I\trf)$
are\dss the boundaries\sss of\trs
the simplices $\Delta\dff(\trf B\trf)$ and $\Delta\dff(\trf I\trf)$ respectively\halfff,\oss
both cohomology\sss groups in\dss the left\sss column are
isomorphic\sss to\dss $\ftwo$\nsp.\oss
Since $\varphi$ is\dss equal\dss to $v_{\fff \bullet}$ 
on $\partial\dff \Delta\dff(\trf I\trf)$\nnsp,\oss
the map\qss
$\partial\dff \varphi
\dff \colon\dff
\partial\dff \Delta\dff(\trf I\trf)
\qff \ttoo\qff
\partial\dff \Delta\dff(\trf B\trf)$\qss
is\dss an\dss isomorphism of\dss simplicial\sss complexes.\oss
It\dss follows\sss that\dss the\sss left\dss vertical\dss map $(\dff\partial\dff \varphi\dff)^{ *\nsp *}$
is\dss an\dss isomorphism.\oss
While $\Delta\dff (\trf A\trf)$ is\dss a simplex,\oss
it\dss is\dss of\dss uncontrollably\sss big dimension and $\partial\dff \Delta\dff(\trf B\trf)$
is\dss not\dss its boundary\halfff.\oss
Still,\oss the upper\dss map $\partial^{*\nsp *}$ is\dss
an isomorphism.\oss
By an easy\sss standard\dss homological\sss argument\dss
this follows\sss from\dss the fact\dss that\sss all\dss homology\sss groups
of\dss a simplex are $0$\nnsp.\oss
In\dss general,\oss the cohomology\dss group\vspace{3pt}
\begin{equation}
\label{lr-group}
\quad
H^{\fff n}\dff\bigl(\qff \mathcal{E}\dff (\trf I\trf)\fff,\pff 
\partial\dff \Delta\dff(\trf I\trf) \qff\bigr)
\end{equation}

\vspace{-9pt}
is\dss not\dss isomorphic\sss to\dss $\ftwo$\nsp.\oss
If\dss $\mathcal{D}$\sss is\dss
not\sss only\sss a\sss chain-simplex,\oss
but\halfff,\oss moreover\halfff,\oss is\dss a\sss pseudo-simplex\qss
(as\sss it\dss is\dss actually\dss the case in\dss the main examples),\oss
then\qss Theorem\qss \ref{envelope-pseudo}\qss
implies\sss that\dss
the simplicial\sss complex $\mathcal{E}\dff (\trf I\trf)$\dss
is\dss a\sss non-branching and\dss its boundary\dss is\dss equal\dss to\sss
the boundary\sss of\dss $ \Delta\dff(\trf I\trf)$\nnsp.\oss
If\dss $\mathcal{E}\dff (\trf I\trf)$\sss
is\dss also strongly\sss connected,\oss
then\dss the cohomology\dss group\qss (\ref{lr-group})\qss
is\dss isomorphic\sss to\dss $\ftwo$\dss
and\dss the lower\dss map $\partial^{*\nsp *}$ is\dss
an isomorphism\qss
(see\qss \cite{i2},\oss Section\qss 10,\oss Theorem\qss 7).\oss

In\dss general,\pss
$\mathcal{E}\trf\fclass{\trf I\trf}$\dss is\dss a\sss relative cycle of\trs
the pair\dss
$\bigl(\qff \mathcal{E}\dff (\trf I\trf)\fff,\pff 
\partial\dff \Delta\dff(\trf I\trf) \qff\bigr)$\dss
and\dss the evaluation of\dss cohomology\sss classes on\dss this relative cycle defines
a\sss homomorphism\vspace{3pt}
\[
\quad
\varepsilon\dff \colon\dff
H^{\fff n}\dff\bigl(\qff \mathcal{E}\dff (\trf I\trf)\fff,\pff 
\partial\dff \Delta\dff(\trf I\trf) \qff\bigr)
\off \ttoo\off
\ftwo
\qff.
\]

\vspace{-9pt}
Theorem\qss \ref{envelope-chain}\qss implies\sss that\dss $\mathcal{E}$\dss is\dss
a chain-simplex and\dss hence\qss (\ref{chain-simplex-def})\qss implies\dss that\trs
$\partial\qff \mathcal{E}\trf\fclass{\trf I\trf}
\off =\off
\partial\dff I$\nnsp.\oss
In\dss turn,\oss by\sss a standard argument\dss this implies\sss that\dss
$\varepsilon\trf \circ\trf \partial^{*\nsp *}$\dss is\dss an\sss isomorphism.\oss
Also,\oss
Lemma\qss \ref{d-of-image}\qss implies\sss that\qss
$\varphi_{\fff *}\qff
\mathcal{E}\trf\fclass{\trf I\trf}
\off =\off
B$\nnsp,\oss
and\dss by\sss a\sss standard argument\dss this implies\sss that\dss
$\varepsilon\trf \circ\trf \varphi^{*\nsp *}$\dss is\dss an\sss isomorphism.\oss
Let\dss us\dss replace\sss in\dss the above diagram\dss the\sss lower\dss right\sss
cohomology\sss group\sss by\dss $\ftwo$\nsp,\oss
the lower\dss horizontal\dss map\dss $\partial^{*\nsp *}$\dss by\qss
$\varepsilon\trf \circ\trf \partial^{*\nsp *}$\nnsp,\oss
and\dss the right\sss vertical\dss map\dss $\varphi^{*\nsp *}$\trs
by\qss
$\varepsilon\trf \circ\trf \varphi^{*\nsp *}$\nnsp.\oss
One can\dss identify\sss all\sss cohomology\dss groups in\dss the resulting\sss diagram\dss with\dss
$\ftwo$\dss and\dss get\dss the\sss diagram\oss\vspace{1.5pt}
\begin{equation*}
\quad
\begin{tikzcd}[column sep=specc, row sep=spe]\dis
\ftwo 
\arrow[r, 
"\dis \partial^{*\nsp *}"]
\arrow[d, 
"\dis (\dff\partial\dff \varphi\dff)^{ *\nsp *}"']
& 
\ftwo
\arrow[d,   
"\dis \varepsilon\trf \circ\trf \varphi^{*\nsp *}"']
\\
\ftwo
\arrow[r,  
"\dis \varepsilon\trf \circ\trf \partial^{*\nsp *}"]
& 
\ftwo
\end{tikzcd}
\end{equation*}

\vspace{-10pt}
All\dss maps\sss in\dss this diagram are isomorphisms.\oss
Let\dss $\mathbb{1}$\dss be\sss the non-zero element\sss of\dss $\ftwo$\nsp.\oss
Then\vspace{3pt}
\[
\quad
\varepsilon\trf \circ\trf \partial^{*\nsp *}\dff \circ\qff (\dff\partial\dff \varphi\dff)^{ *\nsp *}\dff
(\qff \mathbb{1}\pff)
\off =\off\dff
\mathbb{1}
\]

\vspace{-9pt}
and\dss the commutativity\sss of\trs the diagram\dss implies\sss that\vspace{3pt}
\begin{equation}
\label{comm-corollary}
\quad
\varepsilon\trf \circ\trf \varphi^{*\nsp *}\dff \circ\pff \partial^{*\nsp *}\dff
(\qff \mathbb{1}\pff)
\off =\off\dff
\mathbb{1}
\qff.
\end{equation}

\vspace{-9pt}
Now\dss we have\sss to identify\dss the cohomology\sss class\dss 
$\partial^{*\nsp *}\dff(\qff \mathbb{1}\pff)$\nnsp.\oss
The cohomology\sss class\dss $\hclass{\gamma_{\fff i}\fff}$\dss
of\trs the cochain $\gamma_{\fff i}$ is\dss equal\dss to\dss $\mathbb{1}$\nnsp.\oss
By\dss the definition of\trs the connecting\dss homomorphism\dss
$\partial^{*\nsp *}$\vspace{3pt}
\[
\quad
\partial^{*\nsp *}\dff(\qff \mathbb{1}\pff)
\off =\off\dff
\partial^{*\nsp *}\dff\left(\qff \hclass{\gamma_{\fff i}\fff}\qff\right)
\off =\off\dff
\left[\qff\partial\qff \widetilde{\gamma_{\fff i}}\pff\right]
\qff,
\]

\vspace{-9pt}
where\dss $\widetilde{\gamma_{\fff i}}$\dss
is\dss any\sss extension of\trs the cochain\dss $\gamma_{\fff i}$\dss
to a cochain of\dss $\Delta\dff (\trf A\trf)$\nnsp.\oss
As\sss it\dss is\dss often\dss the case,\oss
the\sss tautological\sss extension\sss by $0$ is\dss not\dss really\dss useful\qss
(cf\halfff.\qss \cite{i2},\oss Section\qss 10).\oss
We claim\dss that\dss the cochain 
$\delta_{\dff i}$ is\dss an extension.\oss
Indeed,\pss
$B\qff -\qff v_{\fff i}\qff +\qff v_{\fff i}
\off =\off
B$\qss is\dss a\sss basis containing $b$ in\dss its convex\dss hull,\oss
and\qss
$B\qff -\qff v_{\fff i}\qff +\qff v_{\fff k}$\qss
is\dss not\sss defined\trs if\qss $k\off \neq\off i$\nnsp.\oss
Hence $\delta_{\dff i}$ is\dss equal\dss to $1$ on\qss
$B\qff -\qff v_{\fff i}$\qss and\dss  to $0$ on
other $(\fff n\dff -\dff 1\fff)$\dnsp-simplices of\dss
$\partial\dff \Delta\dff(\trf B\trf)$\dss
and\dss hence extends\dss $\gamma_{\fff i}$\nsp.\oss
It\dss follows\dss that\vspace{3pt}
\[
\quad
\partial^{*\nsp *}\dff(\qff \mathbb{1}\pff)
\off =\off\dff
\hclass{\dff\partial\dff \delta_{\dff i}\dff}
\off =\off\dff
\hclass{\dff \delta\trf}
\qff,
\]

\vspace{-9pt}
where we used\dss the fact\dss that\qss
$\delta\off =\off \partial^{\fff *}\dff \delta_{\dff i}$\nsp.\oss
Now\qss (\ref{comm-corollary})\qss implies\sss that\qss 
$\varepsilon\trf \circ\trf \varphi^{*\nsp *}\qff\bigl(\qff 
\hclass{\dff \delta\trf}
\qff\bigr)
\off =\off\qff
\mathbb{1}$\nnsp.\oss
But\oss \vspace{3.5pt}
\[
\quad
\varepsilon\trf \circ\trf \varphi^{*\nsp *}\dff\bigl(\qff \hclass{\dff \delta\trf}\qff\bigr)
\off\qff =\off\off
\sum\nolimits\qff \varphi^{\fff *}(\dff \delta\dff)\dff(\dff \sigma\dff)
\qff,
\]

\vspace{-8.5pt}
where\sss the sum\dss is\dss taken over all\sss 
$n$\dnsp-simplices of\trs $\mathcal{E}\dff(\trf I\trf)$\nnsp.\oss
Since\qss $\varphi^{\fff *}(\dff \delta\dff)\dff(\dff \sigma\dff)
\off =\off
1$\qss
if\qss and\dss only\trs if\qss $\varphi\dff(\dff \sigma\dff)$\dss
is\dss a\sss basis containing $b$ in\dss its convex\dss hull,\oss
it\dss follows\dss that\dss the number of\dss simplices $\sigma$ with\dss
this property\dss is\dss odd.\oss
In\dss view\sss of\qss (\ref{extended-image})\qss this\dss is\dss equivalent\dss to\sss
the main\dss theorem.\oss
The key\sss step of\trs the both\dss versions of\trs
the proof\qss is\dss an application of\qss Lemma\qss \ref{computing-coboundary}.\oss
From\dss the cohomological\dss point\sss of\dss view\dss 
this\sss lemma\dss allows\sss to\sss 
identify\dss the cohomology\sss 
class\dss
$\partial^{*\nsp *}\dff\hclass{\gamma_{\fff i}\fff}$\nnsp.\oss

\myuppar{A\sss graph-theoretical\dss interpretation of\trs the proof\halfff.}
Suppose now\dss that\dss $\mathcal{D}$\dss 
is\dss a\sss pseudo-simplex.\oss
In\dss this case\sss the proof\dss of\trs the main\dss theorem admits
a graph-theoretical\dss interpretation\sss leading\dss to a\sss path-following\sss
algorithm in\dss the spirit\sss of\qss Scarf\halfff.\qff\oss
Cf.\oss \cite{i2},\oss Section\qss 9.\oss

A\sss basis of\trs $M$\dss is\dss said\dss to be a\qss
\emph{good\dss basis}\pss if\dss $b$ is\dss contained\dss in\sss its convex\dss hull.\oss
Let\qss $i\qff \in\qff I$\nnsp.\oss
The above\dss  
proof\dss 
can\dss be understood\sss as a
computation of\vspace{4.5pt} 
\[
\quad
\partial^{\fff *}\dff \delta_{\fff i}\qff
\bigl(\pff
\varphi_{\fff *}\dff
\bigl(\qff \mathcal{E}\trf\fclass{\trf I\trf}\qff\bigr)
\pff\bigr)
\off\qff =\off\pff
\varphi^{\dff *}\dff\bigl(\qff \partial^{\fff *}\dff \delta_{\fff i}\qff\bigr)
\dff
\bigl(\qff \mathcal{E}\trf\fclass{\trf I\trf}\qff\bigr)
\pff\bigr)
\qff.
\]

\vspace{-7.5pt}
Since\sss $\mathcal{E}\trf\fclass{\trf I\trf}$ involves
all $n$\dnsp-simplices of\dss $\mathcal{E}\dff (\trf I\trf)$\nnsp,\oss
this\dss is\dss almost\dss the same as computing\sss\vspace{4.5pt}
\[
\quad
\varphi^{\dff *}\dff\bigl(\qff \partial^{\fff *}\dff \delta_{\fff i}\qff\bigr)
\off\qff =\off\pff
\partial^{\fff *}\dff\bigl(\qff \varphi^{\dff *}\dff \delta_{\fff i}\qff\bigr)
\qff.
\]

\vspace{-7.5pt}
Let\dss us\dss treat\sss cochains as sums of\dss simplices.\oss
Then\dss  
$\varphi^{\dff *}\dff \delta_{\fff i}$ 
is\dss the sum
of\dss $(\fff n\dff -\dff 1\fff)$\dnsp-simplices $\tau$ of\dss
$\mathcal{E}\dff (\trf I\trf)$\dss
such\dss that\qss
$\varphi\dff(\dff \tau\trf)\qff +\qff v_{\fff i}$\qss
is\dss a\sss good\dss basis.\oss
The coboundary $\partial^{\fff *}\fff \tau$ is\dss the sum of\dss
all $n$\dnsp-simplices having $\tau$ as face,\pss
and 
$\partial^{\fff *}\dff\left(\qff \varphi^{\dff *}\dff \delta_{\fff i}\qff\right)$
is\dss the sum of\dss coboundaries $\partial^{\fff *}\fff \tau$\nnsp.\pss
Therefore\sss the computation of\dss
$\partial^{\fff *}\dff\bigl(\qff \varphi^{\dff *}\dff \delta_{\fff i}\qff\bigr)$
involves only\dss the $(\fff n\dff -\dff 1\fff)$\dnsp-simplices $\tau$
as\dss above,\pss the $n$\dnsp-simplices having\dss them as\sss faces,\pss 
 and\dss the relation\qss 
``\dnsp$\tau$ is\dss a face of\dss $\sigma$\dnsp''\qss 
between\dss them.\oss

This\sss information can\sss be represented\dss in\dss terms  of\dss
a\sss graph $\mathbb{G}_{\dff i}$\nsp.\oss
Its\sss vertices are $n$\dnsp-simplices $\sigma$ 
of\dss $\mathcal{E}\dff (\trf I\trf)$
such\dss that\qss 
$\varphi\dff(\trf \tau\qff)\qff +\qff v_{\dff i}$\qss 
is\dss a\sss good\dss basis\sss for some
$(\fff n\dff -\dff 1\fff)$\dnsp-face $\tau$ of\trs $\sigma$\dnsp.\oss
The\dss edges of\dss $\mathbb{E}_{\dff i}$ correspond\trs to
$(\fff n\dff -\dff 1\fff)$\dnsp-simplices\dss $\tau$\dss
of\qss $\mathcal{E}\dff (\trf I\trf)$\dss
not\sss contained\dss in $I$
and\sss such\dss that\sss 
$\varphi\dff(\trf \tau\qff)\qff +\qff v_{\dff i}$\dss
is\dss a\sss good\dss basis.\oss
Since\dss $\mathcal{D}$\dss is\dss a\sss pseudo-simplex,\oss
Theorem\qss \ref{envelope-pseudo}\qss implies\sss that\dss every\sss such
$(\fff n\dff -\dff 1\fff)$\dnsp-simplex $\tau$
is\dss a\sss face 
of\trs exactly\dss two $n$\dnsp-simplices of\dss $\mathcal{E}\dff (\trf I\trf)$\nnsp.\oss
By\dss the definition,\oss
both\dss these $n$\dnsp-simplices
are vertices of\qss $\mathbb{G}_{\dff i}$\nsp.\oss
The\sss edge corresponding\dss to $\tau$ connects\sss these\sss two vertices.\oss
The main\sss properties of\trs $\mathbb{G}_{\dff i}$
are proved\dss in\qss Lemma\qss \ref{graph}\qss below\sss
after\dss a\dss preliminary\qss Lemma\qss \ref{starting-vertex}.\oss

\mypar{Lemma.}{starting-vertex}
\emph{Let\dss $\tau$\dss be\dss an\dss $(\fff n\dff -\dff 1\fff)$\dnsp-simplex of\pss
$\mathcal{E}\dff (\trf I\trf)$\dss contained\dss in\dss $I$\nnsp.\oss 
If\qss
$\varphi\dff(\trf \tau\qff)\qff +\qff v_{\dff i}$\dss 
is\trs defined,\oss
then\qss $\tau\off =\off I\qff -\qff i$\nnsp.\oss
Conversely\halfff,\oss
$\varphi\dff(\trf I\qff -\qff i\qff)\qff +\qff v_{\dff i}$\dss
is\dss a\sss good\dss basis.\oss
There\dss is\dss exactly\sss one $n$\dnsp-simplex\dss
of\pss $\mathcal{E}\dff (\trf I\trf)$\dss having\qss
$I\qff -\qff i$\qss as\sss a\sss face.\oss
This\sss simplex\dss $\sigma_{\dff i}$\dss is\dss a\sss vertex of\pss $\mathbb{G}_{\dff i}$\nsp.\oss}

\proof
If\qss $\tau\off \subset\off I$\nnsp,\oss
then\qss $\tau\off =\off I\qff -\qff k$\qss
for some\qss $k\qff \in\qff I$\nnsp.\oss
If\pss $k\off \neq\off i$\nnsp,\oss
then\qss $i\qff \in\qff \tau$\qss
and\dss hence\qss
$v_{\fff i}
\off =\off 
\varphi\dff(\dff i\trf)
\qff \in\pff 
\varphi\dff(\dff \tau\trf)$\nnsp.\oss
In\dss this case\qss
$\varphi\dff(\trf \tau\qff)\qff +\qff v_{\dff i}$\dss
is\dss undefined.\oss
Clearly\halfff,\pss
$\varphi\dff(\trf I\qff -\qff i\qff)\qff +\qff v_{\dff i}
\off =\off
B$\qss
is\dss a good\dss basis.\oss
This proves\sss the first\dss two claims of\trs the lemma.\oss
If\qss $\rho\qff \cup\qff (\trf I\qff -\qff i\trf)$\qss
is\dss an $n$\dnsp-simplex of\dss $\mathcal{E}\dff (\trf I\trf)$\nnsp,\oss
then $\rho$ is\dss a $0$\dnsp-simplex of\trs $\mathcal{D}\dff(\trf \{\trf i\qff\}\trf)$\nnsp.\oss
By\qss Lemma\qss \ref{empty}\qss the complex\dss
$\mathcal{D}\dff(\trf \{\trf i\qff\}\trf)$\dss has only\sss one vertex
and\dss hence only\sss one $0$\dnsp-simplex.\oss
It\dss follows\dss that\sss $\rho$ is\dss equal\dss to\sss this $0$\dnsp-simplex.\oss
In\dss particular\halfff,\pss $\rho$ is\dss unique.\oss
Therefore\sss there\dss is\dss exactly\sss one $n$\dnsp-simplex\dss $\sigma_{\dff i}$\dss 
of\dss $\mathcal{E}\dff (\trf I\trf)$\dss having\qss $I\qff -\qff i$\qss
as\sss a\sss face.\oss  
Since\qss 
$\varphi\dff(\trf I\qff -\qff i\qff)\qff +\qff v_{\dff i}$\qss
is\dss a good\dss basis,\pss
$\sigma_{\dff i}$\dss is\dss a\sss vertex\sss 
of\qss $\mathbb{G}_{\dff i}$\nsp.\oss  \eproof

\mypar{Lemma.}{graph}
\emph{Let\sss $\sigma$ be\sss a\sss vertex of\pss $\mathbb{G}_{\dff i}$\nsp.\oss
If\pss $\varphi\dff(\dff \sigma\trf)$ is\dss a\sss good\dss basis,\oss
then\dss $\sigma$ is\dss an\dss endpoint\sss of\trs exactly\dss $1$ edge\dss
if\pss
$\sigma\off \neq\off \sigma_{\dff i}$\qss
and of\pss $0$ edges\dss if\pss
$\sigma\off =\off \sigma_{\dff i}$\nsp.\oss 
Otherwise\dss $\sigma$ is\dss an\dss endpoint\sss of\trs exactly\dss $2$ edges\dss
if\pss
$\sigma\off \neq\off \sigma_{\dff i}$\qss
and of\pss $1$ edge\dss if\pss
$\sigma\off =\off \sigma_{\dff i}$\nsp.\oss}

\proof
Suppose\sss that\dss $\varphi\dff(\dff \sigma\trf)$\sss is\dss a\sss good\dss basis.\oss
By\sss applying\qss Lemma\qss \ref{exchange}\qss to 
$\varphi\dff(\dff \sigma\trf)$ in\dss the role\sss of\trs $B$\sss
and\sss $v_{\fff i}$\dss in\dss the role of\dss $w$\sss
we see\sss that\dss there\dss is\dss exactly\sss one 
element\qss $u\qff \in\qff \varphi\dff(\dff \sigma\trf)$\qss
such\dss that\qss
$\varphi\dff(\dff \sigma\trf)\qff -\pff u\pff +\pff v_{\fff i}$\qss
is\dss a\sss good\dss basis.\oss
Since\dss $\varphi\dff(\dff \sigma\trf)$\sss is\dss a\sss basis,\oss
$\varphi$\qss is\dss a\dss bijection\sss on $\sigma$ and\sss elements\qss
$u\qff \in\qff \varphi\dff(\dff \sigma\trf)$\dss
correspond\dss to $(\fff n\dff -\dff 1\fff)$\dnsp-faces $\tau$ of\dss $\sigma$
such\dss that\qss
$\varphi\dff(\dff \sigma\trf)\qff -\pff u
\off =\off
\varphi\dff(\dff \tau\trf)$\nnsp.\oss
It\dss follows\dss that\dss there\dss is\dss exactly\sss one
$(\fff n\dff -\dff 1\fff)$\dnsp-face $\tau$ of\dss $\sigma$
such\dss that\qss
$\varphi\dff(\trf \tau\qff)\qff +\qff v_{\dff i}$\dss is\dss a\sss good\dss basis.\oss
If\qss $\sigma\off =\off \sigma_{\dff i}$\nsp,\oss
then\qss $I\qff -\qff i$\qss is\dss such a face and\dss hence\qss
$\tau\off =\off I\qff -\qff i$\nnsp.\oss
In\dss this case $\sigma$ is\dss an endpoint\sss of\dss $0$ edges.\oss
If\qss $\sigma\off \neq\off \sigma_{\dff i}$\nsp,\oss
then $\tau\off \not\subset\off I$\qss  
by\qss
Lemma\qss \ref{starting-vertex}\qss
and $\sigma$ is\dss an endpoint\sss of\dss $1$ edge.\oss
This completes\sss the proof\dss in\dss the case when\dss 
$\varphi\dff(\dff \sigma\trf)$\sss is\dss a\sss good\dss basis.\oss

Suppose now\dss that\dss $\varphi\dff(\dff \sigma\trf)$ is\dss not\sss a\sss good\dss basis.\oss
We claim\dss that\dss there exactly\dss two $(\fff n\dff -\dff 1\fff)$\dnsp-faces 
$\tau$ of\dss $\sigma$ such\dss that\qss
$\varphi\dff(\trf \tau\qff)\qff +\qff v_{\dff i}$\qss 
is\dss a\sss good\dss basis.\oss
Since $\sigma$ is\dss a\sss vertex of\pss $\mathbb{G}_{\dff i}$\nsp,\oss
there\dss is\dss at\dss least\sss one such\dss face.\oss
If\dss $\tau$ is\dss such a face,\oss then\dss the set\qss
$\varphi\dff(\trf \tau\qff)\qff +\qff v_{\dff i}$\qss
is\dss a\sss basis
and\dss hence consists of\qss $n\qff +\qff 1$\qss elements.\oss
It\dss follows\dss that\sss $\varphi\dff(\trf \tau\qff)$\dss consists of\dss $n$ elements.\oss
Since $\tau$ also consists of\dss $n$ elements,\oss
this implies\sss that\sss $\varphi$ is\dss injective on $\tau$\nnsp.\oss
Let\sss $u$ be\sss the vertex of\dss $\sigma$ such\dss that\qss
$\tau\off =\off \sigma\qff -\qff u$\nnsp.\oss
There are\sss two cases\sss to consider\halfff,\oss
depending on\trs if\qss 
$\varphi\dff(\dff \sigma\trf)
\off =\off
\varphi\dff(\trf \tau\qff)$\dss
or\trs not\halfff.\oss\vspace{0.625pt}

Suppose first\dss that\qss
$\varphi\dff(\dff \sigma\trf)
\off =\off
\varphi\dff(\trf \tau\qff)$\nnsp.\oss
Since $\varphi$ is\dss injective on $\tau$\nnsp,\oss
in\dss this case\sss
there\dss is\dss exactly\sss one vertex\qss $v\off \neq\off u$\qss
of\dss $\sigma$ such\dss that
$\varphi\dff(\dff v\trf)
\off =\off
\varphi\dff(\dff u\trf)$
and
$\varphi\dff(\dff \sigma\qff -\qff u\trf)
\off =\off
\varphi\dff(\trf \sigma\qff -\qff v\qff)$\nnsp.\oss
If\dss $w$ is\dss a vertex of\dss $\sigma$ different\dss
from\qss $u\fff,\pff v$\nnsp,\oss
then
$\varphi\dff(\trf \sigma\qff -\qff w\qff)$
has\qss $<\qff n$\qss elements and\dss hence\sss the set\qss
$\varphi\dff(\trf \sigma\qff -\qff w\qff)\qff +\qff v_{\dff i}$\qss
is\dss not\sss a\sss basis.\oss
This implies our claim\dss in\dss the case when\qss
$\varphi\dff(\dff \sigma\trf)
\off =\off
\varphi\dff(\trf \tau\qff)$\nnsp.\oss
Suppose now\dss that\qss
$\varphi\dff(\dff \sigma\trf)
\off \neq\off
\varphi\dff(\trf \tau\qff)$\nnsp.\oss
Then\trs $\varphi\dff(\dff \sigma\trf)$\dss
properly\sss contains\trs
$\varphi\dff(\trf \tau\qff)$\dss
and\dss hence consists of\qss $n\qff +\qff 1$\qss elements.\oss
It\dss follows\dss that\trs $\varphi$\qss is\dss injective on $\sigma$
and\qss
$\varphi\dff(\dff \sigma\trf)\qff -\qff \varphi\dff(\dff u\trf)\qff +\qff v_{\fff i}$\qss
is\dss a\sss good\dss basis.\oss
By\sss applying\qss Lemma\qss \ref{two-faces}\qss
to\trs $\varphi\dff(\dff \sigma\trf)$\dss
in\dss the role of\trs $D$ 
we see\sss that\dss there are $0$ or $2$ elements\qss
$d\qff \in\qff \varphi\dff(\dff \sigma\trf)$\qss
such\dss that\qss
$\varphi\dff(\dff \sigma\trf)\qff -\qff d\qff +\qff v_{\fff i}$\qss
is\dss a\sss good\dss basis.\oss
Of\dss course,\pss $\varphi\dff(\dff u\trf)$ is\dss one of\trs them.\oss
Since $\varphi$ is\dss injective on $\sigma$\nnsp,\oss
the other one\dss is\dss equal\dss to $\varphi\dff(\dff u\trf)$
for uniquely\sss determined\qss $v\qff \in\qff \sigma$\nnsp,\oss
$v\off \neq\off u$\nnsp.\oss
It\dss follows\sss that\qss
$\varphi\dff(\trf \sigma\qff -\qff w\qff)\qff +\qff v_{\dff i}$\qss
is\dss a\sss good\dss basis for\qss $w\off =\off u$\qss or\dss $v$\sss
but\dss not\dss for any\sss other vertex of\dss $\sigma$\nnsp.\oss
This\dss implies our claim\dss in\dss the case when\qss
$\varphi\dff(\dff \sigma\trf)
\off \neq\off
\varphi\dff(\trf \tau\qff)$\nnsp.\oss\vspace{0.625pt}

Therefore\dss if\qss $\varphi\dff(\dff \sigma\trf)$ is\dss not\sss a\sss good\dss basis,\oss
then\dss there exactly\dss two $(\fff n\dff -\dff 1\fff)$\dnsp-faces 
$\tau$ of\dss $\sigma$ such\dss that\qss
$\varphi\dff(\trf \tau\qff)\qff +\qff v_{\dff i}$\qss 
is\dss a\sss good\dss basis.\oss
Suppose\sss that\dss one of\trs them\dss is\dss contained\dss in\dss $I$\nnsp.\oss
In\dss this case\qss
Lemma\qss \ref{starting-vertex}\qss implies\sss that\dss the other does not\sss
and\qss $\sigma\off =\off \sigma_{\dff i}$\nnsp.\oss
Hence in\dss this case $\sigma$ is\dss an endpoint\sss of\dss $1$ edge.\oss
If\trs neither of\trs these\sss two faces\dss is\dss contained\dss in\dss $I$\nnsp,\oss
then $\sigma$ is\dss an endpoint\sss of\dss $2$ edges.\oss
This completes\sss the proof\dss in\dss the case when\dss 
$\varphi\dff(\dff \sigma\trf)$\sss is\trs not\sss a\sss good\dss basis.\oss  \eproof\vspace{2pt}

\prooftitle{A\sss graph-based\dss proof\dss of\qss Theorem\qss \ref{main-non-degenerate}}
Lemma\qss \ref{graph}\qss implies\sss that\sss each
component\dss of\pss the graph\dss $\mathbb{G}_{\dff i}$\dss
either\dss is\dss a\dss cycle,\oss
or\trs is\dss a\sss path connecting\dss two $n$\dnsp-simplices\dss
$\sigma\fff,\pff \sigma'$\dss such\dss that\sss
$\varphi\dff(\dff \sigma\trf)$
and\sss
$\varphi\dff(\dff \sigma'\trf)$\dss
are good\dss bases,\oss
or\trs is\dss a\sss path connecting\sss $\sigma_{\dff i}$\dss
with\dss an $n$\dnsp-simplex\sss $\sigma$
such\dss that\sss $\varphi\dff(\dff \sigma\trf)$
is\trs a\sss good\dss basis.\oss
If\dss $\varphi\dff(\dff \sigma_{\dff i}\trf)$
is\dss a\sss good\dss basis,\oss
then\dss the last\dss path degenerates into a single vertex $\sigma_{\dff i}$\nsp.\oss
If\dss $\varphi\dff(\dff \sigma_{\dff i}\trf)$
is\dss not\sss a\sss good\dss basis,\oss
then\dss $\mathbb{G}_{\dff i}$\dss
contains a unique path connecting\sss $\sigma_{\dff i}$\dss
with\dss an $n$\dnsp-simplex\sss $\sigma$
such\dss that\sss $\varphi\dff(\dff \sigma\trf)$
is\trs a\sss good\dss basis.\oss
Other $n$\dnsp-simplices\sss $\sigma$
such\dss that\sss $\varphi\dff(\dff \sigma\trf)$
is\trs a\sss good\dss basis occur\dss in\dss pairs connected\dss
by\sss a\sss path of\qss $\mathbb{G}_{\dff i}$\nsp.\oss
In\dss particular\halfff,\oss
there\dss is\dss at\sss least\sss one $n$\dnsp-simplex\sss $\sigma$
such\dss that\sss $\varphi\dff(\dff \sigma\trf)$ is\dss a\sss basis
containing $b$ in\dss its convex\dss hull.\oss
Moreover\halfff,\oss the number of\dss such $n$\dnsp-simplices\sss $\sigma$
is\dss odd.\oss
In\dss view\sss of\qss (\ref{extended-image})\qss
this implies\dss Theorem\qss \ref{main-non-degenerate}.\oss  \eproof

\newpage
\mysection{Vector\qss colorings}{vector-colorings}

\myuppar{Introduction.}
The existence part\sss of\qss Theorem\qss \ref{main-non-degenerate}\qss
remains\sss valid\dss without\dss the non-degeneracy\sss assumption\qss
({\fff}the part\sss about\dss the number of\dss cells being odd does not).\oss
In\dss this section\sss we will\dss prove\sss this in\dss 
the classical\sss situation when\dss the oriented\dss matroid\dss $M$\dss
is\dss defined\dss by\sss a vector configuration in\dss $\rrr^{\dff n\dff +\dff 1}$\dnsp\dnsp.\oss
In\dss this case\sss the most\dss natural\sss approach\dss is\dss to perturb\sss this
vector configuration.\oss 
It\dss turns out\dss that\dss it\dss is\dss sufficient\dss to perturb only\dss
the vector\sss $b$\nnsp,\oss as suggested\dss by\qss
Scarf\pss \cite{sc1},\oss \cite{sc3}.\oss
As an application we will\dss prove\sss the classical\qss Scarf\qss theorem\qss
\cite{sc1},\oss \cite{sc3}.\vspace{0.25pt}

\myuppar{The vector\dss framework.}
Let\dss us\sss 
number\dss the\sss standard\sss coordinates of\qss $\rrr^{\fff n\dff +\dff 1}$\dss by\qss
$0\fff,\pff 1\fff,\pff \ldots\fff,\pff n$\nnsp.\oss
Let\trs
$v_{\fff 0}\fff,\pff v_{\fff 1}\fff,\pff \ldots\fff,\pff v_{\fff n}$\qss
be\sss the standard\dss basis of\qss $\rrr^{\fff n\dff +\dff 1}$\dnsp.\oss
Let\qss 
$B\off =\off
\left\{\qff
v_{\fff 0}\fff,\pff v_{\fff 1}\fff,\pff \ldots\fff,\pff v_{\fff n}
\pff\right\}$\qss
and\dss let\qss $M\off \subset\off \rrr^{\fff n\dff +\dff 1}$\qss be a finite set\sss
containing\dss $B$\dss and a\sss vector\dss $b\off \neq\off 0$\dss
such\dss that\sss all\dss its coordinates are non-negative.\oss
The last\sss condition\sss means\sss that\dss there are numbers\qss $b_{\dff i}\qff \geq\qff 0$\qss
such\dss that\vspace{4pt}
\begin{equation}
\label{b-positive}
\quad
b
\off =\off
\sum\nolimits_{\qff i\qff \in\qff I}\qff b_{\dff i}\dff v_{\fff i}
\qff,
\end{equation}

\vspace{-8pt}
The circuits of\trs the vector configuration\dss 
$M$\dss
define an oriented\dss matroid\dss which we also denote by\dss $M$\nnsp.\oss
Let\dss $A\off =\off M\qff -\qff b$\dss and\dss
let\dss us\dss consider\dss the equation\vspace{4pt}
\begin{equation}
\label{boundedness}
\quad
b
\off =\off\qff
\sum\nolimits_{\qff v\qff \in\qff A}\qff
y_{\fff v}\trf v
\end{equation}

\vspace{-8pt}
in\dss real\dss unknowns $y_{\fff v}$\nsp.\oss
A solution\dss 
$\left(\dff y_{\fff v}\dff\right)_{\dff v\qff \in\qff A}$\dss
is\dss  a vector\qss
$y\qff \in\pff \rrr^{\dff A}$\dnsp.\oss 
Following\qss Scarf\qss \cite{sc3},\oss Theorem\qss 4.2.3,\oss
we assume\sss that\dss the set\dss of\dss non-negative solutions 
of\qss (\ref{boundedness})\qss is\dss bounded.\oss\vspace{0.25pt}

\mypar{Lemma.}{vector-scarf-acyclic}
\emph{Under\dss the above assumptions\dss the oriented\dss 
matroid\trs $M$\dss is\dss acyclic.}\vspace{0.25pt}

\proof
If\trs $M$\dss is\dss not\sss acyclic,\oss 
then\dss there\dss is\dss a\sss circuit\sss $\sigma$
such\dss that\qss $\sigma_{\dff -}\off =\off \varnothing$\nnsp.\oss
Equivalently\halfff,\oss $0$\dss is\dss a\sss
non-trivial\dss non-negative\sss linear\dss combination of\dss
vectors\sss in\dss $M$\nnsp.\oss
Using\qss (\ref{b-positive})\qss one can eliminate $b$\sss from\dss this linear
combination\dss while keeping\dss it\dss non-negative.\oss 
Therefore\qss\vspace{4pt}
\[
\quad
0
\off =\off\qff
\sum\nolimits_{\qff v\qff \in\qff A}\qff
z_{\trf v}\trf v
\]

\vspace{-8pt}
for some non-negative numbers\dss $z_{\trf v}$\nsp,\oss not\sss
all\sss of\trs which\dss are\dss $0$\nnsp.\oss
Clearly\halfff,\oss if\dss $y$\dss is\dss a\dss non-negative solution of\qss
(\ref{boundedness}),\oss then\qss
$y\qff +\qff \lambda\dff z$\qss is\dss also a non-negative solution\dss
for every\dss $\lambda\qff >\qff 0$\dss and\dss hence\sss
the set\sss of\dss non-negative solutions\dss is\dss unbounded,\oss
contrary\dss to\sss the assumption.\oss  \eproof\vspace{0.25pt}

\myuppar{The non-degeneracy\sss property\halfff.}
The vector configuration\dss $M$\dss together\dss with\dss $b\qff \in\pff M$\dss
is\dss said\dss to be\qss \emph{non-degenerate}\oss if\qss
the vector\dss 
$b$\dss
cannot\dss be represented\sss as a non-negative\sss
linear combination of\qss $<\qff n\qff +\qff 1$\qss elements of\trs $M\qff -\qff b$\nnsp.\oss
See\qss \cite{sc3},\oss the non-degeneracy\sss assumption\qss 4.3.1.\oss
In\dss this case\sss the corresponding\dss matroid\dss $M$\dss together\dss with\dss $b$\sss
are,\pss obviously\halfff,\oss non-degenerate.\oss

\mypar{The main\dss theorem for\dss vector colorings\qss ({\fff}the general\sss case).}{main-theorem-vector}
\emph{Suppose\sss that\dss $\mathcal{D}$ is\dss a\sss chain-simplex\dss based on\dss
$I\off =\off \{\qff 0\fff,\pff 1\fff,\pff \ldots\fff,\pff n\qff\}$\dss
and\dss let\pss $X\off =\off V_{\dff \mathcal{D}}$\nnsp.\oss
For every\dss coloring\qss
$c\dff \colon\dff X\qff \ttoo\qff A$\qss
there exist\sss a\sss non-empty\dss subset\qss $C\qff \subset\pff I$\qss
and\dss a\dss $d\dff(\trf C\trf)$\dnsp-simplex\trs $\tau$\dss of\oss
$\mathcal{D}\dff(\trf C\trf)$\dss such\dss that}\vspace{2.25pt}
\begin{equation}
\label{scarf-basis-vector}
\quad
c\qff(\trf \tau\qff)\off \cup\off\dff \bigl\{\qff v_{\fff i}\off \bigl|\off
i\qff \in\qff I\qff \smallsetminus\pff C \pff\bigr\}
\end{equation}

\vspace{-9.75pt}
\emph{is\dss a\sss basis of\pss $\rrr^{\fff n\dff +\dff 1}$ 
{\nsp}and\dss $b$\dss is\dss a\sss non-negative\sss linear\dss 
combination of\qss its\dss elements.\oss}\vspace{0.375pt}

\proof
Suppose\sss first\dss that\trs $M\fff,\pff b$\dss is\dss non-degenerate.\oss
Then\dss Theorem\qss \ref{main-non-degenerate}\qss applies.\oss
By\dss this\dss theorem\dss
there exist\sss a\sss non-empty\dss subset\dss $C\qff \subset\pff I$\dss
and\dss a\dss $d\dff(\trf C\trf)$\dnsp-simplex $\tau$ of\dss
$\mathcal{D}\dff(\trf C\trf)$ such\dss that\qss (\ref{scarf-basis-vector})\qss
is\dss a\sss basis of\trs $\rrr^{\fff n\dff +\dff 1}$ and\sss $b$\sss belongs\sss to
its convex\dss hull\dss in\dss the sense of\dss oriented\dss matroids.\oss
By\dss the definition,\oss the latter\dss means\sss that\dss $b$\sss 
is\dss a\sss non-negative\sss linear\dss combination of\dss elements
of\trs the basis\qss (\ref{scarf-basis-vector}).\oss
This proves\sss the\sss theorem\dss in\dss the non-degenerate case.\oss

The degenerate case reduces\sss to\sss the non-degenerate one by\sss perturbing\sss $b$\nnsp.\oss
Obviously\halfff,\oss there exists a sequence\qss
$b_{\dff 1}\dff,\off b_{\dff 2}\dff,\off b_{\dff 3}\dff,\pff \ldots$\qss
of\dss vectors\qss $b_{\dff i}\qff \in\pff \rrr^{\fff n\dff +\dff 1}$\dss
such\dss that\dss $b_{\dff i}$\sss tends\sss to $b$ when\dss $i\qff \ttoo\qff \infty$\dss
and\dss $b_{\dff i}$\dss is\dss not\sss a\sss linear combination of\trs
$<\qff n\qff +\qff 1$\dss 
elements\sss of\qss
$M\qff -\qff b$\dss for every\sss $i$\nnsp.\oss
Let\qss $M_{\dff i}\off =\dss\off M\qff -\qff b\qff +\qff b_{\dff i}$\qss for every\sss $i$\nnsp.\oss
Then\dss the configuration\trs $M_{\dff i}\dff,\off b_{\dff i}$\dss is\dss non-degenerate
for every\sss $i$\dss
and\dss $c\dff \colon\dff X\qff \ttoo\off M\qff -\qff b\off =\off M_{\dff i}\qff -\qff b_{\dff i}$\dss
is\dss a matroid\sss coloring.\oss

By\dss the already\dss proved\dss non-degenerate case of\trs the\sss theorem,\oss
for every\sss $i$\dss there exist\sss a non-empty\sss subset\dss 
$C_{\dff i}\qff \subset\pff I$\dss
and\sss a simplex\sss 
$\tau_{\dff i}$\sss of\trs $\mathcal{D}\dff(\trf C_{\dff i}\dff)$\dss
such\dss that\vspace{2.25pt}
\[
\quad
c\qff(\trf \tau_{\dff i}\qff)\off \cup\off\dff \bigl\{\qff v_{\fff i}\off \bigl|\off
i\qff \in\qff I\qff \smallsetminus\pff C_{\dff i} \off\bigr\}
\]

\vspace{-9.75pt}
is\dss a\sss basis and\dss $b_{\dff i}$\dss 
is\dss a\sss non-negative\sss linear\dss combination of\trs its elements.\oss
Since\sss the set\trs $I$\dss and\dss the complexes\dss $\mathcal{D}\dff(\trf C\trf)$\dss are finite,\oss
after\sss passing\dss to a subsequence we may\sss assume\sss that\dss
$C_{\dff i}\off =\off C$\dss for some subset\dss $C\qff \subset\pff I$\dss
and\dss $\tau_{\dff i}\off =\off \tau$\dss some simplex $\tau$ of\trs $\mathcal{D}\dff(\trf C\trf)$\dss
for all\dss $i$\nnsp.\oss
Then\qss (\ref{scarf-basis-vector})\qss
is\dss a\sss basis\sss and\dss $b_{\dff i}$\dss 
is\dss a\sss non-negative\sss linear\dss combination of\trs its elements
for every\sss $i$\nnsp.\oss
In other\dss words,\oss the coordinates of\trs $b_{\dff i}$\dss in\dss this
basis are non-negative for every\sss $i$\nnsp.\oss
When\dss $i\qff \ttoo\qff \infty$\dss these coordinates\sss tend\dss to\sss
the coordinates of\trs $b$\sss in\dss the same basis.\oss
By\dss passing\dss to\sss the limit\sss we see\sss that\dss the coordinates of\trs $b$\sss
are also non-negative.\oss
The\sss theorem\dss follows.\oss  \eproof\vspace{0.375pt}

\myuppar{Families of\qss linear orders.}
Let\dss us\dss move from\dss the simplex-families\sss of\qss 
Section\qss \ref{pseudo-simplices}\qss
to\sss the families\sss of\qss linear orders of\qss Section\qss \ref{orders}.\oss
So,\oss let\qss $X$\dss be a non-empty\dss finite set\sss
and suppose\sss that\dss for
every\qss $i\qff \in\pff I$\qss a\dss linear\sss order\dss $<_{\dff i}$\dss
on\dss $X$\dss is\dss given.\oss
As before,\oss let\dss us\sss assume\sss that\dss 
$X\qff \cap\pff I\off =\off \varnothing$\nnsp.\oss

In order\dss to\sss be closer\dss to\qss Scarf's\qss form of\trs his\sss theorem,\oss
let\dss us\sss extend\dss the orders\dss $<_{\dff i}$\dss
to\sss linear orders on\dss the set\dss $X\qff \cup\pff I$\nnsp,\oss
which we will\sss still\sss denote by\dss $<_{\dff i}$\nsp.\oss
Following\qss Scarf\halfff,\oss
we will\dss require\sss that\sss every\sss element\dss 
$i\qff \in\qff I$\trs
is\dss minimal\dss with\dss respect\dss to\dss $<_{\dff i}$\dss
and\dss that\trs $x\qff <_{\dff i}\qff k$\dss for every\dss
$x\qff \in\pff X$\qss and every\dss $k\qff \in\pff I$\nnsp,\oss $k\off \neq\off i$\nnsp.\oss
Clearly\halfff,\oss such extensions exist\halfff,\oss
but\sss are not\sss unique in\dss general\halfff:\oss the restriction of\trs the order\dss $<_{\dff i}$\dss to\dss
$I\qff -\qff i$\dss can\dss be completely arbitrary.\oss
Let\dss us\dss fix some choice of\dss such extensions.\oss
Now\dss we can speak about\qss \emph{dominant\sss 
subsets of}\qss $X\qff \cup\pff I$\nnsp.\oss

\mypar{Lemma.}{extended-orders-dominance}
\emph{Let\qss $\tau\qff \subset\pff X$\dss and\qss $C\qff \subset\pff I$\nnsp,\oss
$C\off \neq\off \varnothing$\nnsp.\oss
Then $\tau$ is\dss dominant\dss with respect\dss to\qss $C$\dss
if\qss and\dss only\qss if\pss
$\tau\qff \cup\qff (\qff I\qff \smallsetminus\pff C\trf)$\dss
is\dss dominant\dss with respect\dss to $I$\nnsp.\oss}\vspace{0.875pt}

\proof
If\dss $i\qff \in\qff C$\dss and\qss $k\qff \in\pff I\qff \smallsetminus\pff C$\nnsp,\oss
then\dss $k\off \neq\off i$\qss and\dss hence\dss
$x\qff <_{\dff i}\qff k$\dss for every\qss $x\qff \in\pff X$\nnsp.\oss
It\dss follows\dss that\qss
$\min\nolimits_{\dff i}\qff \tau\qff \cup\qff (\qff I\qff \smallsetminus\pff C\trf)
\off =\off
\min\nolimits_{\dff i}\qff \tau$\qss 
for every\dss $i\qff \in\qff C$\nnsp.\oss
If\dss $i\qff \in\qff I\qff \smallsetminus\pff C$\nnsp,\oss
then\dss $\tau\qff \cup\qff (\qff I\qff \smallsetminus\pff C\trf)$\dss
contains\sss the $<_{\dff i}$\dnsp-minimal\sss element\sss $i$\sss and\dss hence
$\min\nolimits_{\dff i}\qff \tau\qff \cup\qff (\qff I\qff \smallsetminus\pff C\trf)
\off =\off
i$\nnsp.\oss

It\dss follows\dss that\qss
$\tau\qff \cup\qff (\qff I\qff \smallsetminus\pff C\trf)$\dss
is\dss dominant\dss with\dss respect\dss to\dss $I$\dss if\trs there\dss is\dss
no element\trs $y\qff \in\pff X\qff \cup\pff I$\dss such\dss that\qss
$\min\nolimits_{\dff i}\qff \tau\qff <\qff_{\dff i}\off y$\qss 
for every\trs $i\qff \in\qff C$\dss
and\qss
$i\qff <\qff_{\dff i}\off y$\qss 
for every\trs $i\qff \in\pff I\qff \smallsetminus\pff C$\nnsp.\oss
Clearly\halfff,\oss there\dss is\dss no such element\dss $y\qff \in\pff X$\qss
if\trs and\dss only\trs if\trs $\tau$ is\dss dominant\dss with\dss respect\dss to
$C$ and\dss the original\sss orders on\dss $X$\nnsp.\oss 
Therefore,\oss it\dss is\dss sufficient\dss to prove\sss that\dss
there\dss are\dss no such elements\dss
$y\qff \in\pff I$\nnsp.\oss
In\-deed,\oss if\trs $y\qff \in\qff C$\nnsp,\oss
then\dss $y\qff <_{\qff y}\qff \min\nolimits_{\dff y}\qff \tau$\dss
and\dss hence\sss the condition\dss
$\min\nolimits_{\dff i}\qff \tau\qff <\qff_{\dff i}\off y$\dss does not\dss hold,\oss
and\dss if\trs $y\qff \in\pff I\qff \smallsetminus\pff C$\nnsp,\oss
then\dss the condition\qss $y\qff <\qff_{\dff y}\off y$\dss does not\dss hold.\oss
The\sss lemma\sss follows.\oss  \eproof\vspace{1.5pt}

\myuppar{Scarf\qss theorem.}
\emph{Suppose\sss that\sss 
$\varphi\dff \colon\dff X\qff \cup\pff I\qff \ttoo\qff \rrr^{\dff n\dff +\dff 1}$\qss
is\dss a map such\dss that\qss
$\varphi\dff(\dff i\trf)\off =\off v_{\dff i}$\qss for every\trs $i\qff \in\pff I$\nnsp.\oss%
Then\dss there exists a subset\dss $\sigma\qff \subset\pff X\qff \cup\pff I$\dss
dominant\dss with respect\dss to\dss $I$\dss and\dss
such\dss that\trs $\varphi\dff(\dff \sigma\dff)$\sss
is\dss a\sss basis\dss of\oss $\rrr^{\dff n\dff +\dff 1}$\dss
and\dss $b$\dss is\dss a\sss non-negative\sss linear combination of\dss
elements of\qss $\varphi\dff(\dff \sigma\dff)$\nnsp.\oss}\vspace{1.5pt}

\proof
Let\dss $\mathcal{T}$\dss be\sss the simplex-family\sss associated\dss 
in\dss Section\qss \ref{orders}\qss
with\dss $X$\dss
and\dss the orders\dss $<_{\dff i}$\dss on\dss $X$\nnsp.\oss
It\dss is\dss a\sss chain-simplex\dss
by\dss Corollary\qss \ref{scarf-chain-simplex}.\oss
If\trs $C\qff \subset\off I$\nnsp,\oss
then every\dss simplex $\tau$ of\qss $\mathcal{T}(\trf C\trf)$\dss
is\dss a\sss subset\sss of\dss a $C$\dnsp-cell,\oss
i.e.\qss of\trs some\sss set\dss $\sigma\qff \subset\pff X$\dss
dominant\dss with\dss respect\sss to $C$ and such\dss that\dss
$\num{\sigma}\off =\off \num{C}$\nnsp.\oss
If\dss $\tau$ is\dss a $d\dff(\trf C\trf)$\dnsp-simplex,\oss
then\dss $\num{\tau}\off =\off \num{C}$\dss
and\dss hence in\dss this case\dss $\tau\off =\off \sigma$\nnsp.\oss
Therefore every\sss $d\dff(\trf C\trf)$\dnsp-simplex of\trs $\mathcal{T}(\trf C\trf)$\dss
is\dss a\sss subset\sss of\trs $X$\sss dominant\dss 
with\sss respect\dss to $C$\nnsp.\oss

Let\qss $M\off =\off \varphi\dff(\qff X\qff \cup\pff I\qff)\qff +\qff b$\qss
and\dss let\sss $c$\dss be\sss the restriction of\trs $\varphi$\sss to\sss $X$\nnsp.\oss
By\qss Theorem\qss \ref{main-theorem-vector}\qss applied\dss
to\dss $M$\dss and\sss $c$\dss 
there exist\sss a\sss non-empty\dss set\dss $C\qff \subset\pff I$\dss
and\dss a\dss $d\dff(\trf C\trf)$\dnsp-simplex $\tau$ of\oss
$\mathcal{T}(\trf C\trf)$\dss such\dss that\qss (\ref{scarf-basis-vector})\qss
is\dss a\sss basis of\trs $\rrr^{\dff n\dff +\dff 1}$\sss and\sss $b$\sss is\dss 
\dss a\sss non-negative\sss linear\dss 
combination of\qss its\dss elements.\oss
Let\trs $\sigma\off =\off \tau\qff \cup\qff (\qff I\qff \smallsetminus\pff C\trf)$\nnsp.\oss
Clearly\halfff,\qss $\num{\sigma}\off =\off \num{\dff I\dff}$\sss and the image\dss
$\varphi\dff(\dff \sigma\dff)$ is\dss equal\dss to\sss 
the set\qss (\ref{scarf-basis-vector}).\oss 
By\dss the previous paragraph $\tau$
is\dss dominant\dss with\dss respect\dss to $C$ and\dss hence $\sigma$ is\dss
dominant\dss with\sss respect\dss to $I$\dss by\trs
Lemma\qss \ref{extended-orders-dominance}.\oss 
The\dss theorem\dss follows.\oss  \eproof\vspace{0.875pt}

\myuppar{The classical\qss Scarf\qss theorem.}
In\qss 
\cite{sc1}\qss and\qss  
\cite{sc3}\pss
Scarf\qss considered only\sss subsets\sss $X$\sss
of\trs the standard simplex\dss $\Delta^{n}\qff \subset\pff \rrr^{\dff n\dff +\dff 1}$
and\dss the orders\dss $<_{\dff i}$\dss which were almost\sss completely\sss
determined\dss by\dss the standard order\dss $<$\dss on\sss the coordinates of\dss points in\dss
$\rrr^{\dff n\dff +\dff 1}$\dnsp,\oss
i.e.\qss the orders from\qss Section\qss \ref{scarf-brouwer}.\oss
At\dss the same\sss time he was aware\sss that\dss his methods work\dss in\dss more general\dss situations.\oss
See\qss \cite{sc3},\oss Chapter\qss 6\qss and especially\dss 
the discussion at\dss the bottom of\qss p.\dss 146.\oss
In\dss fact\halfff,\oss 
every\dss family\sss of\trs $n\qff +\qff 1$\dss orders on a finite
set\dss $X$\sss can\sss be realized\dss by\sss an  
embedding of\trs $X$\sss
into\dss 
$\Delta^{n}$\dnsp.\oss 
So,\oss the above\sss theorem\dss is\dss more abstract\halfff,\oss
but\dss is\dss not\dss more general\dss than\dss the original\qss Scarf\qss theorem.\oss
While\sss the above proof\trs incorporates a\sss lot\sss of\trs ideas of\pss Scarf\halfff,\oss
it\sss depends on\qss Theorem\qss \ref{main-non-degenerate}\qss
about\sss colorings into colors belonging\dss to an oriented\dss matroid.\oss
The\sss latter\sss notion was still\dss in\dss the future at\dss the\sss time
of\pss \cite{sc1}\qss and\qss \cite{sc3}.\oss
Also,\oss our proof\dss of\pss  Theorem\qss \ref{main-non-degenerate}\qss
is\dss based on\dss topological\dss ideas,\oss
in contrast\dss with related,\oss but\sss different\dss
path-following arguments of\qss Scarf\halfff.\oss

\newpage
\mysection{Oriented\qss matroid\qss colorings\fff:\oss the\qss general\qss case}{scarf-oriented-matroids-general}

\myuppar{Perturbations\sss in\sss oriented\dss matroids.}
Let\trs $M$\dss be\sss an\sss acyclic oriented\dss matroid\sss
and\dss let\vspace{3pt}
\[
\quad
B\off =\off
\left\{\qff
v_{\fff 0}\fff,\pff v_{\fff 1}\fff,\pff \ldots\fff,\pff v_{\fff n}
\pff\right\}
\]

\vspace{-9pt}
be\dss a\sss basis of\dss $M$\nnsp.\oss 
Let\qss $b\qff \in\qff M\qff \smallsetminus\qff B$\dss
and suppose\sss that\sss $b$\sss belongs\sss to\sss the convex hull\sss of\dss $B$\nnsp.\oss
In\dss this\sss section\dss
we\qss \emph{do\dss not\dss assume}\pss that\dss the pair\dss $M\fff,\pff b$\dss
is\dss non-degenerate.\oss

The\sss theory\sss of\dss oriented\dss matroids
provides a combinatorial\sss analogue of\dss the perturbation argument\dss
used\dss in\dss the proof\dss of\qss Theorem\qss \ref{main-theorem-vector}.\oss
Namely\halfff,\oss one can add\dss to\dss $M$\dss
a new element\dss $p$\dss playing\dss the role of\trs a\sss particular\dss
perturbation of\dss $b$\sss in\dss the vector\sss framework.\oss
Let\dss us\dss describe\sss this particular\dss perturbation.\oss
Clearly\halfff,\oss in\dss the vector\sss framework\dss 
$B\qff -\qff v_{\dff i}\qff +\qff b$\qss is\dss a\sss basis for some\dss $i\qff \in\pff I$\qss
(usually\dss for all\sss such $i$\nsp).\oss
Without\sss any\dss loss of\dss generality\sss one can assume\sss that\trs
$B\qff -\qff v_{\dff 0}\qff +\qff b$\qss is\dss a\sss basis.\oss
Then one can\dss take as\sss the perturbation of\dss $b$\dss the vector\qss
$p
\off\dff =\off
b\qff +\qff 
\lambda\dff b_{\dff 1}\qff +\qff 
\lambda^{\dff 2}\dff b_{\dff 2}\qff +\qff 
\ldots\qff +\qff
\lambda^{\dff n}\dff b_{\dff n}$\qss
for any\sss sufficiently\sss small\trs $\lambda\qff >\qff 0$\nnsp.\oss

Adding\sss a new element\dss is\dss a special\sss case of\dss the notion of\dss
an\qss \emph{extension}\qss of\dss an oriented\dss matroid.\oss
Namely\halfff,\oss an oriented\dss matroid\sss $M\fff'$\sss is\dss said\dss to be an\qss
\emph{extension}\qss of\trs $M$\dss if\qss $M\qff \subset\pff M\fff'$\dss and
a\sss signed subset\sss of\trs $M$\trs is\dss a\sss circuit\sss of\trs $M\fff'$\sss
if\trs and\dss only\trs if\trs it\dss is\dss a\sss circuit\sss of\trs $M$\nnsp.\oss
We are interested\dss in\dss the\trs \emph{one point\sss extensions}\qss $M\fff'$\dnsp,\oss
i.e.\qss extensions\sss $M\fff'$ such\dss that\trs
$M\fff'\off =\off M\qff +\qff p$\dss for some\dss $p\qff \not\in\qff M$\nnsp.\oss
Moreover\halfff,\oss we are interested only\dss in\sss one-point\sss extensions\dss
$M\fff'$\dss such\dss that\dss a basis of\trs $M$\dss is\dss
also a basis in\dss $M\fff'$\qss
(if\trs this\dss is\dss true for one basis of\trs $M$\nnsp,\oss
then\dss this\dss is\dss true for every\dss basis of\trs $M$\nsp).\oss

Among\dss such extensions are\sss the\qss
\emph{lexicographic extensions}\qss 
introduced\dss by\dss M.\dss Las\dss Vergnas\qss \cite{lv}.\oss
A\qss \emph{lexicographic extension}\qss $L\off =\off M\qff +\qff p$\dss of\trs $M$\dss
is\dss associated\dss
with every\qss \emph{ordered}\pss basis\qss\vspace{3pt}
\begin{equation}
\label{a-basis}
\quad
b
\off =\off
b_{\trf 0}\dff,\off 
b_{\dff 1}\dff,\off b_{\dff 2}\trf,\off \ldots\dff,\off b_{\dff n}
\end{equation}

\vspace{-9pt}
of\trs $M$\nnsp.\oss
It\dss was realized\dss by\trs M.\dss Todd\qss \cite{t}\qss that\dss 
the new element\sss $p$\sss of\trs $L$\dss can\sss play\dss the role of\trs
the perturbation\qss
$b\qff +\qff
\lambda\dff b_{\dff 1}\qff +\qff 
\lambda^{\dff 2}\dff b_{\dff 2}\qff +\qff 
\ldots\qff +\qff
\lambda^{\dff n}\dff b_{\dff n}$\qss
of\dss $b$\nsp.\oss
In\dss the following\dss two\sss lemmas we state,\oss
following\pss M.\dss Todd\qss \cite{t},\oss
the relevant\dss properties of\trs the lexicographic
extension\dss $L\off =\off M\qff +\qff p$\dss associated\dss with\qss (\ref{a-basis}).\oss
Appendix\qss \ref{extensions}\qss contains some additional\dss details\qss
(not\dss used\dss here).\oss\vspace{1.5pt}

\mypar{Lemma.}{lex-1}
\emph{If\qss $\sigma$\dss is\dss a\sss circuit\sss of\pss $L$\dss such\dss that\dss
$p\qff \in\qff \underline{\sigma}$\nsp,\oss
then\dss $\num{\underline{\sigma}\fff}\qff  
\geq\qff n\qff +\qff 2$\nnsp.\oss
The signed subset\sss 
$\sigma
\off =\off
(\qff \sigma_{\dff +}\dff,\pff  \sigma_{\dff -}\qff)$\nnsp,\oss where\sss
$\sigma_{\dff +}
\off =\off
\{\trf b\dff,\off b_{\dff 1}\trf,\off b_{\dff 2}\trf,\off
\ldots\dff,\off b_{\dff n} \qff\}$\sss
and\qss
$\sigma_{\dff -}\off =\off
\{\trf p\qff\}$\nnsp,\oss is\dss a\sss circuit\halfff.\oss}\vspace{1.5pt}

\mypar{Lemma.}{lex-2}
\emph{Suppose\sss that\trs 
$\sigma$\sss is\dss a\sss circuit\sss of\oss $L$\dss such\dss that\trs
$p\qff \in\qff \sigma_{\dff -}$\dss
and\qss
$b\pff \not\in\off \underline{\sigma}$\nsp.\oss
Then\dss there\sss exists\sss exactly\sss one\sss circuit\qss $\rho$\sss of\oss $M$\sss
such\dss that\pss
$b\pff \in\off \rho_{\dff -}$\nsp,\qff\oss
$\rho_{\dff -}\qff \smallsetminus\qff \{\qff b\pff\}
\pff \subset\dff\off \sigma_{\dff -}$\nsp,\qff\oss 
and\pss 
$\rho_{\dff +}
\pff \subset\dff\off \sigma_{\dff +}$\nsp.\oss}\vspace{1.75pt}

\prooftitle{Proofs}
See\qss \cite{t},\oss Propositions\qss 5.3\qss and\qss 5.4\qss and\dss Corollary\qss 5.6.\oss  \eproof

\mypar{Corollary.}{lex-acyclic}
\emph{If\oss $M$\dss is\dss acyclic,\oss then\qss $L\qff -\qff b$\dss
is\dss also acyclic.\oss}

\proof
Let\sss
$\tau$\sss be\dss a\sss circuit\sss of\trs $L$\sss
such\dss that\dss $\underline{\tau}\pff \subset\pff L\qff -\qff a$\dss
and\trs $\tau_{\dff -}\off =\off \varnothing$\nnsp.\oss
If\pss $p\qff \not\in\off \underline{\tau}$,\oss
then\sss $\tau$\sss is\dss a\sss circuit\sss of\trs $M$\dss
and\dss hence\dss $M$\dss is\dss not\sss acyclic.\oss
Therefore\dss we can assume\sss that\trs
$p\qff \in\off \underline{\tau}\dff\off =\off \tau_{\dff +}$.\oss
Let\trs $\sigma\off =\off -\qff \tau$\nnsp.\oss
Then\dss
$p\qff \in\qff \sigma_{\dff -}$\dss
and\dss
$a\pff \not\in\off \underline{\sigma}$\nsp.\oss
Let\sss $\rho$\sss be\sss the circuit\sss of\trs $M$\sss such as in\dss Lemma\qss \ref{lex-2}.\oss
Then\dss
$\rho_{\dff +}\qff \subset\off \sigma_{\dff +}\pff =\off \tau_{\dff -}\pff =\off \varnothing$\nnsp.\oss
Again,\oss this\dss is\dss impossible for an acyclic\dss $M$\nnsp.\oss  \eproof

\mypar{Corollary.}{lex-3}
\emph{If\pss $p$\dss belongs\sss to\sss the convex\dss hull\sss of\qss a\dss
subset\qss $\varepsilon\qff \subset\pff M\qff -\qff b$\dss in\qss $L$\nnsp,\oss
then\dss $b$\dss belongs\sss to\sss the convex\dss hull\dss of\pss
$\varepsilon$\dss in\qss $M$\nnsp.\oss}\vspace{-0.25pt}

\proof
The first\sss statement\dss immediately\dss follows\dss from\qss
Lemma\qss \ref{lex-1}.\oss
Suppose\sss that\dss $p$\dss belongs\sss to\sss the convex\dss hull\dss 
of\qss $\varepsilon\qff \subset\pff M\qff -\qff b$\nnsp.\oss
Then\dss there\dss is\dss a\sss circuit\sss $\sigma$\sss of\qss $L$\dss
such\dss that\trs
$\sigma_{\dff +}\qff \subset\off \varepsilon$\dss 
and\dss
$\sigma_{\dff -}\off =\off
\{\trf p\qff\}$\nnsp.\oss
Now\trs Lemma\qss \ref{lex-2}\qss implies\sss that\dss there exists
a circuit\sss $\rho$\sss of\qss $M$\dss such\dss that\trs
$b\pff \in\off \rho_{\dff -}$\nsp,\oss 
$\rho_{\dff -}\qff \smallsetminus\qff \{\qff b\pff\}
\pff \subset\dff\off 
\sigma_{\dff -}\qff \cap\qff M
\off =\off
\varnothing$\nnsp,\oss and\pss
$\rho_{\dff +}
\pff \subset\dff\off 
\sigma_{\dff +}
\pff \subset\dff\off
\varepsilon$\nnsp.\oss
Hence\trs
$\rho_{\dff -}
\off =\off
\{\qff b\pff\}$\dss
and\dss 
$\rho_{\dff +}
\pff \subset\dff\off \varepsilon$\nnsp,\oss and\dss 
therefore\sss 
$b$\sss belongs\sss to\sss the convex\dss hull\dss
of\dss $\varepsilon$\nnsp.\oss  \eproof

\mypar{The main\dss theorem\sss for\dss matroid colorings\qss
({\fff}the\sss general\sss case).}{main-theorem-general}
\emph{Let\dss $\mathcal{D}$ is\dss a\sss chain-simplex\dss based on\dss
$I\off =\off \{\qff 0\fff,\pff 1\fff,\pff \ldots\fff,\pff n\qff\}$
and\dss let\qss $X\off =\off V_{\dff \mathcal{D}}$\nnsp.\oss
For every\dss coloring\dss
$c\dff \colon\dff X\qff \ttoo\qff M\qff -\qff b$\qss 
there exist\sss a\sss non-empty\dss subset\qss $C\qff \subset\pff I$\qss
and\dss a\dss $d\dff(\trf C\trf)$\dnsp-simplex\trs $\tau$\dss of\oss
$\mathcal{D}\dff(\trf C\trf)$\dss such\dss that}\vspace{4.5pt}
\begin{equation}
\label{scarf-matroid-basis}
\quad
c\qff(\trf \tau\qff)\off \cup\off\dff \bigl\{\qff v_{\fff i}\off \bigl|\off
i\qff \in\qff I\qff \smallsetminus\pff C \pff\bigr\}
\end{equation}

\vspace{-7.5pt}
\emph{is\dss a\sss basis of\pss $M$\dss 
and\dss $b$\dss is\dss contained\dss in\dss its\sss convex\dss hull.\oss}

\proof
The axioms of\dss oriented\dss matroids immediately\dss imply\dss that\dss
the sets\dss $\underline{\sigma}$\nnsp, where\dss $\sigma$\dss is\dss 
a\sss circuit\sss of\trs $M$\nnsp,\oss
are\sss the circuits of\dss a\sss matroid,\oss which we will\sss denote\sss by\dss
$\underline{M}$.\oss
Moreover\halfff,\oss the oriented\dss matroid\dss $M$\dss and\dss the matroid\dss $\underline{M}$\dss
have\sss the same bases.\oss
By\dss the symmetric exchange property\sss of\trs matroids\sss
the set\trs
$B\qff -\qff v_{\fff i}\qff +\qff b$\dss is\dss a\sss basis of\trs $\underline{M}$,\oss
and\dss hence of\trs $M$\nnsp,\oss for some\dss $i\qff \in\pff I$\nnsp.\oss
Without\sss any\dss loss of\dss generality\dss we can assume\sss that\trs
$B\qff -\qff v_{\trf 0}\qff +\qff b$\dss is\dss a\sss basis of\trs $M$\nnsp.\oss

Let\trs $L\off =\off M\qff +\qff p$\dss be\sss the lexicographic extension of\trs $M$\dss
associated\dss with\dss the basis\qss (\ref{a-basis}),\oss where\qss
$b_{\dff i}\off =\off v_{\fff i}$\qss
for\dss $i\qff \geq\qff 1$\nnsp.\oss
As an unordered\sss set\dss this\dss is\dss the basis\dss
$B\qff -\qff v_{\trf 0}\qff +\qff b$\nnsp.\oss
Let\qss $M\fff'\off =\off L\qff -\qff b$\nnsp.\oss
Then\dss $B$\sss is\dss also a\sss basis of\trs $M\fff'$\dss
and\dss the new element\dss $p$\dss belongs\sss to\dss
$M\fff'\qff \smallsetminus\qff B$\nnsp.\oss
By\qss Corollary\qss \ref{lex-acyclic}\qss the oriented\dss matroid\dss $M\fff'$\dss
is\dss acyclic,\oss and\dss by\qss
Lemma\qss \ref{lex-1}\qss  
the pair\dss $M\fff'\fff,\pff p$\dss
is\dss non-degenerate.\oss
 
Let\dss us\dss check\dss that\dss $p$\sss belongs\sss to\sss the convex\dss hull\sss of\trs
$B$\dss in\trs $M\fff'$\nnsp.\oss
Lemma\qss \ref{lex-1}\qss implies\sss that\dss
the signed subset\sss 
$\sigma
\off =\off
(\qff \sigma_{\dff +}\dff,\pff  \sigma_{\dff -}\qff)$\nnsp,\oss where\sss
$\sigma_{\dff +}
\off =\off
B\qff -\qff v_{\trf 0}\qff +\qff b$\sss
and\qss
$\sigma_{\dff -}\off =\off
\{\qff p\qff\}$\nnsp,\oss is\dss a\sss circuit\halfff.\oss
Since\sss $b$\sss belongs\sss to\sss the convex\dss hull\sss of\trs $B$\nnsp,\oss
there exists a circuit\trs $\tau$\dss such\dss that\dss
$\tau_{\dff -}\off =\off \{\trf b\trf\}$\dss and\dss
$\tau_{\dff +}\qff \subset\qff B$\nnsp.\oss
By\dss the axiom\qss ({\fff}iv{\fff})\qss of\dss oriented\dss matroids\sss 
there\dss exists\dss a\sss circuit\dss $\omega$\dss
such\dss that\qss \vspace{4.5pt}
\[
\quad
\omega_{\dff +}
\off \subset\off\dff 
(\qff \sigma_{\dff +}\qff \cup\qff \tau_{\dff +}\qff)
\qff \smallsetminus\qff 
\{\qff b\qff\}
\hspace{1.2em}\mbox{and}\hspace{1.2em}
\omega_{\dff -}
\off \subset\off\dff 
(\qff \sigma_{\dff -}\qff \cup\qff \tau_{\dff -}\qff)
\qff \smallsetminus\qff 
\{\qff b\qff\}
\qff.
\]

\vspace{-7.5pt}
These inclusions imply\dss that\trs
$\omega_{\dff +}\qff \subset\pff B$\qss and\qss
$\omega_{\dff -}\qff \subset\off  \{\trf p\trf\}$\nnsp.\oss
Since\dss $B$\sss is\dss a\sss basis,\pss $\omega_{\dff -}\off \neq\off \varnothing$\nnsp.\oss
Therefore\qss
$\omega_{\dff -}\off =\off  \{\trf p\trf\}$\dss
and\dss hence\dss $p$\dss belongs\sss to\sss the convex\dss hull of\trs $B$\nnsp.\oss

Clearly\halfff,\oss 
$M\fff'\qff -\qff p
\off =\off
M\qff -\qff b$\dss
and one can consider\sss $c$\sss as a\sss matroid\sss coloring\dss
$X\qff \ttoo\qff M\fff'\qff -\qff p$\nnsp.\oss
By\sss applying\qss
Theorem\qss \ref{main-non-degenerate}\qss to\dss $M\fff'\fff,\pff p$\dss
in\dss the role of\qss $M\fff,\pff b$\nnsp,\oss the same basis\sss $B$\nnsp,\oss
and\dss the coloring\sss $c$\nnsp,\oss we see\sss
that\dss there exist\sss
a non-empty\sss subset\trs $C\qff \subset\pff I$\trs and
a $d\dff(\trf C\trf)$\dnsp-simplex\sss $\tau$ of\dss
$\mathcal{D}\dff(\trf C\trf)$\dss such\dss that\qss (\ref{scarf-matroid-basis})\qss
is\dss a\sss basis of\trs $M\fff'$\dss 
containing\dss $p$\dss in\dss its convex\dss hull.\oss
Let\dss us\dss denote\sss this basis by\sss $\varepsilon$\nnsp.\oss
Since\sss $\varepsilon$\sss 
is\dss contained\dss in\dss 
$M\fff'\qff -\qff p
\off =\off
M\qff -\qff b$\nnsp,\oss
the basis\sss $\varepsilon$\dss is\dss also a\sss basis of\trs $M$\nnsp.\oss
Corollary\qss \ref{lex-3}\qss implies\sss that\sss $b$\sss belongs\sss to\sss
the convex\dss hull\sss of\trs $\varepsilon$\sss in\dss $M$\nnsp.\oss
The\sss theorem\dss follows.\oss  \eproof

\myuppar{Families of\qss linear orders.}
Let\qss $X$\dss be a non-empty\dss finite set\sss
and suppose\sss that\dss for
every\qss $i\qff \in\pff I$\qss a\dss linear\sss order\dss $<_{\dff i}$\dss
on\dss $X$\dss is\dss given.\oss
As before,\oss let\dss us\sss assume\sss that\dss 
$X\qff \cap\pff I\off =\off \varnothing$\nnsp.\oss
Let\dss us\sss extend\dss the orders\dss $<_{\dff i}$\dss
to\sss linear orders on\dss $X\qff \cup\pff I$\dss
subject\dss to\sss the same conditions
as\sss in\dss Section\qss \ref{vector-colorings}.\oss

\myuppar{Generalized\qss Scarf\qss theorem.}
\emph{Let\qss $c\dff \colon\dff X\qff \ttoo\qff M\qff -\qff b$\qss
be a matroid coloring\halfff.\oss 
Let\dss us extend\dss $c$\dss to a map\qss
$\varphi\dff \colon\dff X\qff \cup\pff I\qff \ttoo\qff M\qff -\qff b$\qss
by\dss the rule\qss
$\varphi\dff(\dff i\trf)\off =\off v_{\dff i}$\nsp.\oss
Then\dss there exists a subset\dss $\sigma\qff \subset\pff X\qff \cup\pff I$\dss
consisting of\pss $\num{\dff I\dff}$\dss elements,\pss
dominant\dss with respect\dss to\dss $I$\nnsp,\oss and
such\dss that\trs
$\varphi\dff(\dff \sigma\dff)$\dss
is\dss a\sss basis of\pss $M$\dss containing\dss 
$b$\dss in\dss its convex\dss hull.\oss}

\proof
The proof\dss is\dss completely\sss similar\dss to\sss the proof\dss
of\pss Scarf\qss theorem\dss in\dss Section\qss \ref{vector-colorings}.\oss
The only\sss difference\dss is\trs that\dss one needs\sss to apply\qss
Theorem\qss \ref{main-theorem-general}\qss
instead of\pss Theorem\qss \ref{main-theorem-vector}.\oss  \eproof

\myuppar{Hedgehog colorings.}
One can apply\qss Theorems\qss \ref{main-non-degenerate}\qss
and\qss \ref{main-theorem-general}\qss also\sss to simplex-families
arising\dss from\dss triangulations of\dss geometric simplices.\oss
There\dss is\dss a version of\trs these\sss theorems\dss
which\dss for such simplex-families\dss is\dss 
closer\dss to\sss the lemmas of\qss Alexander\qss and\qss Sperner\halfff.\oss
See\qss Theorem\qss \ref{main-theorem-hedgehog}.\oss
It\dss deals with\dss matroid colorings
called\qss \emph{hedgehog\dss colorings}.\oss
Their definition\dss is\dss based on\sss the notions of\qss
\emph{hyperplanes}\qss and\qss \emph{cocircuits}\pss
discussed\sss in\dss
Appendix\qss \ref{extensions}.\oss See also\qss \cite{om}.\oss
Let\trs $H$\dss be a hyperplane of\trs $M$\nnsp.\oss
If\trs $\eta$\dss is\dss one of\trs two cocircuits corresponding\dss to\dss $H$\nnsp,\oss 
then\vspace{4.5pt}
\[
\quad
\eta\trf(\qff \geq\dff 0\dff)
\off =\off
H
\pff \cup\off 
\bigl\{\qff 
v\qff \in\pff M
\off \bigl|\off
\eta\dff(\dff v\trf)\off =\off\dff +
\off\bigr\} 
\qff
\]

\vspace{-7.5pt}
is\dss called\dss the\qss \emph{closed\dss half-space}\pss
defined\dss by\dss $\eta$\nnsp.\oss
We will\sss denote\sss by\dss $\eta\trf(\qff \leq\dff 0\dff)$\dss
the closed\dss half-space defined\dss by\dss the cocircuit\dss $-\qff \eta$\nnsp.\oss
Closed\dss half-spaces are\qss \emph{convex}\qss in\dss the sense\sss that\dss they\sss
are equal\dss to\sss their convex\dss hulls.\oss
This follows,\oss for example,\oss from\qss \cite{om},\oss Exercise\qss 3.9.\oss
For each\dss $i\qff \in\pff I$\qss let\dss $H_{\dff i}$\dss be\sss
the hyperplane spanned\dss by\dss the set\dss $B\qff -\qff v_{\fff i}$\dss
and\dss let\dss $\eta_{\dff i}$ be\sss one of\trs two cocircuits 
associated\sss with\dss $H_{\dff i}$\nsp,\oss
namely\halfff,\oss the one for\sss which\dss
$\eta_{\dff i}\dff(\dff v_{\fff i}\trf)\off =\off\dff +$\nsp.\oss

A\sss matroid coloring\dss
$c\dff \colon\dff X\qff \ttoo\qff M\qff -\qff b$\dss
is\dss called\sss a\qss \emph{hedgehog\sss coloring}\oss
if\qss 
$C\pff \subset\off I\qff -\qff i$\qss
implies\sss that\trs
$c\dff(\dff v\trf)\qff \in\pff \eta_{\dff i}\trf(\qff \leq\dff 0\dff)$\dss 
for\dss every\dss vertex $v$\dss of\trs $\mathcal{D}\dff(\trf C\trf)$
and every\trs $i\qff \in\pff I$\nnsp.\oss
If\trs $\mathcal{D}$\dss is\dss the simplex-family\sss 
associated\dss with\dss a\sss triangulation
of\dss a\sss geometric simplex,\oss 
then\dss $C\pff \subset\pff D$\dss
implies\sss that\trs
$\mathcal{D}\dff(\trf C\trf)$\dss is\dss
a\sss subcomplex of\qss $\mathcal{D}\dff(\dff D\dff)$\dss
and\dss $c$\dss is\dss a\sss
hedgehog\dss coloring\qss if\trs and\dss only\trs if\qss
$c\dff(\dff v\trf)\qff \in\pff \eta_{\dff i}\trf(\qff \leq\dff 0\dff)$\dss
for every\trs $i\qff \in\pff I$\trs and every\sss 
vertex\sss $v$ of\qss $\mathcal{D}\dff(\trf I\qff -\qff i\trf)$\nnsp.\oss\vspace{0.95pt}

\mypar{Lemma.}{half-spaces}
\emph{If\qss $b$\sss is\dss not\sss contained\dss in\dss the convex\dss hull\sss
of\dss any\dss proper subset\sss of\pss the basis\trs $B$\nnsp,\oss
then\qss 
$b\off \not\in\off \eta_{\dff i}\trf(\qff \leq\dff 0\dff)$\qss
for every\qss $i\qff \in\pff I$\nnsp.\oss}\vspace{0.95pt}

\proof
The point\dss $b$\sss is\dss contained\dss in\dss the convex\dss hull\sss of\trs $B$\dss
and\dss hence belongs\sss to\dss $\eta_{\dff i}\trf(\qff \geq\dff 0\dff)$\nnsp.\oss
Clearly\halfff,\pss
$\eta_{\dff i}\trf(\qff \geq\dff 0\dff)
\qff \cap\qff
\eta_{\dff i}\trf(\qff \leq\dff 0\dff)
\off =\off
H_{\dff i}$\nsp.\oss
Therefore,\oss if\qss $b\off \in\off \eta_{\dff i}\trf(\qff \leq\dff 0\dff)$\nnsp,\oss
then\qss $b\qff \in\qff H_{\dff i}$\nsp.\oss
Hence it\dss is\dss sufficient\dss to show\dss that\dss $b$\dss belongs\sss to\sss
the convex\dss hull\sss of\qss $B\qff -\qff v_{\fff i}$\qss if\qss
$b\qff \in\qff H_{\dff i}$\nsp.\oss
So,\oss suppose\sss that\trs
$b\qff \in\qff H_{\dff i}$\nsp.\oss
Then\dss there exists a circuit\dss $\sigma$\dss such\dss that\dss
$b\qff \in\off \underline{\sigma}$\qss and\qss
$\underline{\sigma}\pff \subset\pff B\qff -\qff v_{\fff i}\qff +\qff b$\nnsp.\oss
Without\sss any\dss loss of\dss generality\dss we can assume\sss that\dss
$b\qff \in\qff \sigma_{\dff +}$\nsp.\oss
On\dss the other\dss hand,\oss there exists a circuit\dss $\tau$\dss such\dss that\dss
$\tau_{\dff -}\off =\off \{\trf b\qff\}$\qss and\qss
$\tau_{\dff +}\qff \subset\pff B$\nnsp.\oss
By\dss the axiom\qss ({\fff}iv\fff)\qss there exists a circuit\sss $\omega$\sss
such\dss that\qss
$\omega_{\dff +}
\off \subset\off\dff 
(\qff \tau_{\dff +}\qff \cup\qff \sigma_{\dff +}\qff)
\qff \smallsetminus\qff 
\{\dff b\qff\}
\off \subset\off
B$\qss
and\qss
$\omega_{\dff -}
\off \subset\off\dff 
(\qff \tau_{\dff -}\qff \cup\qff \sigma_{\dff -}\qff)
\qff \smallsetminus\qff 
\{\dff b\qff\}
\off \subset\off
B$\nnsp.\oss
Hence\trs $\underline{\omega}\pff \subset\pff B$\nnsp.\oss
Since\dss $\underline{\omega}\off \neq\off \varnothing$\dss by\dss the axiom\qss ({\fff}i\fff),\oss
this contradicts\sss to $B$\sss being a basis.\oss
The lemma\sss follows.\oss  \eproof

\mypar{The main\dss theorem\sss for\dss hedgehog colorings.}{main-theorem-hedgehog}
\emph{Suppose\sss that\qss $b$\sss is\dss not\sss contained\dss in\dss the convex\dss hull\sss
of\dss any\dss proper subset\sss of\pss the basis\trs $B$\nnsp.\oss
Then for every\dss hedgehog coloring\dss
$c$\dss there exist\dss a\dss $d\dff(\trf I\trf)$\dnsp-simplex\trs $\tau$\dss of\oss
$\mathcal{D}\dff(\trf I\trf)$\dss such\dss that\dss $c\dff(\dff \tau\dff)$\dss 
is\dss a\sss basis of\pss $M$\dss 
and\dss $b$\dss is\dss contained\dss in\dss its\sss convex\dss hull.\oss
If\qss the pair\dss $M\fff,\pff b$\dss is\dss
non-degenerate,\oss then\dss the number of\qss such $\tau$ is\dss odd.\oss}

\proof
For each\dss $i\qff \in\pff I$\qss let\trs $w_{\dff i}\off =\off v_{\fff i\dff +\dff 1}$\nsp,\oss
where subscripts are\sss treated as integers modulo\dss $n\qff +\qff 1$\nnsp.\oss
We can apply\qss Theorem\qss \ref{main-theorem-general}\qss
to\sss the basis\dss
$w_{\fff 0}\fff,\pff w_{\fff 1}\fff,\pff \ldots\fff,\pff w_{\fff n}$\dss
in\dss the role of\qss
$v_{\fff 0}\fff,\pff v_{\fff 1}\fff,\pff \ldots\fff,\pff v_{\fff n}$\qss
({\fff}it\dss is\dss the same as a set\halfff,\oss but\dss the order\dss matters).\oss
Therefore\sss there exist\sss a non-empty\dss subset\qss $C\qff \subset\pff I$\qss
and\dss a\dss $d\dff(\trf C\trf)$\dnsp-simplex\trs $\tau$\dss of\oss
$\mathcal{D}\dff(\trf C\trf)$\dss such\dss that\vspace{2.5pt}
\[
\quad
Y
\off\off =\off\off
c\qff(\trf \tau\qff)\off \cup\off\dff \bigl\{\qff v_{\fff i\dff +\dff 1}\off \bigl|\off
i\qff \in\qff I\qff \smallsetminus\pff C \pff\bigr\}
\off\off =\off\off
c\qff(\trf \tau\qff)\off \cup\off\dff \bigl\{\qff w_{\fff i}\off \bigl|\off
i\qff \in\qff I\qff \smallsetminus\pff C \pff\bigr\}
\]

\vspace{-9.5pt}
is\dss a\sss basis of\pss $M$\dss 
containing\dss $b$\sss in\dss its\sss convex\dss hull.\oss
Suppose\sss that\qss $C\off \neq\dff\off I$\nnsp.\oss
Then\dss there exists\dss $k\qff \in\pff I$\dss such\dss that\trs
$k\qff \not\in\qff C$\dss and\dss $k\qff -\qff 1\qff \in\qff C$\nnsp.\oss
In\dss particular\halfff,\pss
$C\qff \subset\pff I\qff -\qff k$\nnsp.\oss
Since $c$ is\dss a\sss hedgehog\sss coloring\halfff,\oss
this implies\sss that\trs
$c\dff(\dff v\trf)\qff \in\pff \eta_{\dff k}\trf(\qff \leq\dff 0\dff)$\dss
for every\dss $v\qff \in\qff \tau$\nnsp.\oss
Since\dss $k\qff -\qff 1\qff \in\qff C$\nnsp,\oss\vspace{2.5pt}
\[
\quad
v_{\dff k}
\pff \not\in\off 
\bigl\{\qff v_{\fff i\dff +\dff 1}\off \bigl|\off
i\qff \in\qff I\qff \smallsetminus\pff C \pff\bigr\}
\]

\vspace{-10.8pt}
and\dss hence\qss 
$\{\qff v_{\fff i\dff +\dff 1}\qff \mid\qff
i\qff \in\qff I\qff \smallsetminus\pff C \pff\}
\off\qff \subset\off\qff
H_{\dff k}
\off\qff \subset\off\qff
\eta_{\dff k}\trf(\qff \leq\dff 0\dff)$\nnsp.\oss 
It\dss follows\dss that\trs $Y\qff \subset\pff \eta_{\dff k}\trf(\qff \leq\dff 0\dff)$\nnsp.\oss
Since\sss $b$\sss belongs\sss to\sss the convex\dss hull\sss of\trs $Y$\dss and\dss
$\eta_{\dff k}\trf(\qff \leq\dff 0\dff)$\sss is\dss convex,\oss
this implies\sss that\trs $b\qff \in\pff \eta_{\dff k}\trf(\qff \leq\dff 0\dff)$\nnsp.\oss
The contradiction\dss with\qss Lemma\qss \ref{half-spaces}\qss
shows\sss that\trs $C\off =\off I$\nnsp.\oss
This\sss proves\sss the existence claim.\oss
The claim about\dss the non-degenerate pairs\dss $M\fff,\pff b$\dss similarly\dss follows\dss 
from\trs Theorem\qss \ref{main-non-degenerate}.\oss  \eproof

\myuppar{Classical\dss colorings and\trs matroid\sss colorings.}
Let\dss us\dss choose some element\dss $b\qff \not\in\pff I$\dss
and\dss let\trs $M\off =\off I\qff +\qff b$\nnsp.\oss
Let\dss $\sigma$\dss be\sss the signed\sss subset\sss of\trs $M$\dss
defined\dss by\trs $\sigma_{\dff +}\off =\off I$\dss and\dss
$\sigma_{\dff -}\off =\off \{\trf b\trf\}$\nnsp.\oss
Taking\dss $\sigma$\dss and\dss $-\qff \sigma$\dss
as\sss the only\sss circuits\sss turns\dss $M$\dss into an oriented\dss matroid.\oss
Obviously\halfff,\oss the sets\trs $I$\dss and\trs $I\qff -\qff i\qff +\qff b$\dss with\trs $i\qff \in\pff I$\dss
are bases of\trs $M$\dss and\dss there are no other bases.\oss
Clearly\halfff,\oss the pair\trs $M\fff,\pff b$\dss is\dss non-degenerate.\oss
This oriented\dss matroid structure on\trs $M$\dss
allows\sss us\dss to consider classical\dss colorings\dss
$c\dff \colon\dff X\qff \ttoo\qff I$\trs
as matroid\sss colorings.\oss
For example,\oss an\dss Alexander--Sperner\dss coloring\dss
is\dss a\sss hedgehog\sss col\-or\-ing and\dss hence\qss 
Theorem\qss \ref{as-colorings}\qss
is\dss a\sss special\sss case of\qss Theorem\qss \ref{main-theorem-hedgehog}.\oss

\newpage
\mysection{Scarf's\pss proof\pss of\pss Kakutani's\qss fixed\qss point\pss theorem}{kakutani}

\myuppar{Multi-valued\dss maps.}
Let\qss $X\fff,\pff Y$\qss be\sss two sets.\oss
A\qss \emph{multi-valued\dss map}\qss $F\dff \colon\dff X\qff \ttoo\qff Y$\qss
is\dss simply\sss a\sss usual\dss map\qss 
$X
\qff \ttoo\qff 
\mathcal{P}\dff(\trf Y\trf)
\qff \smallsetminus\qff
\{\trf \varnothing\trf\}$\nnsp,\oss
where\dss $\mathcal{P}\dff(\trf Y\trf)$\dss is\dss the set\sss of\dss all\sss
subsets of\trs $Y$\nnsp.\oss
The\qss \emph{graph}\qss of\dss a multi-valued\dss map\qss
$F\dff \colon\dff X\qff \ttoo\qff Y$\qss
is\dss the subset\sss of\qss $X\dff \times\dff Y$\qss
consisting of\dss all\dss pairs\dss $(\dff x\fff,\qff y\trf)$\dss
such\dss that\qss $y\qff \in\qff F\dff(\dff x\trf)$\nnsp.\oss
A\qss \emph{fixed\dss point}\qss of\qss $F$\dss is\dss a point\qss $x\qff \in\qff X$\qss
such\dss that\qss $x\qff \in\pff F\dff(\dff x\trf)$\nnsp.\oss
Suppose\sss that\qss $X\fff,\pff Y$\qss are\sss topological\sss spaces.\oss
A multi-valued\dss map\qss
$X\qff \ttoo\qff Y$\qss
is\dss said\dss to be\trs \emph{closed}\pss if\trs its\dss graph\dss is\dss
a\sss closed subset\sss of\qss $X\dff \times\dff Y$\nnsp.\oss
Usually\sss closed\dss multi-valued\dss maps are called\qss
\emph{upper semicontinuous},\oss but\dss the\sss term\qss ``closed'',\oss
borrowed\dss from\dss the operator\dss theory\halfff,\oss
seems\sss to\sss invoke more relevant\sss ideas.\oss
If\qss $F\dff \colon\dff X\qff \ttoo\qff Y$\qss is\dss closed,\oss then\qss 
$F\dff(\dff x\trf)$\dss is\dss closed\dss for every\qss $x\qff \in\qff X$\nnsp.\oss

If\trs the\sss topological\sss spaces\qss $X\fff,\pff Y$\qss are sufficiently\dss nice\qss
(for example,\oss are subspaces of\dss $\rrr^{\dff n}$\nsp),\oss
then\dss $F$\dss is\dss closed\trs if\trs and\dss only\trs if\trs
the following condition\dss holds.\oss
Suppose\sss that\qss
$x_{\dff 1}\fff,\pff x_{\trf 2}\fff,\pff  x_{\dff 3}\dff,\pff \ldots$\qss
is\dss a\sss sequence of\dss points of\trs $X$\dss
converging\dss to a\sss point\qss $x\qff \in\qff X$\nnsp.\oss 
Suppose\sss that\dss $y_{\dff i}\qff \in\qff F\dff(\dff x_{\dff i}\trf)$\dss
for every\sss $i$\sss and\dss the sequence\qss
$y_{\dff 1}\fff,\pff y_{\dff 2}\fff,\pff  y_{\dff 3}\dff,\pff \ldots$\qss
converges\sss to a\sss point\qss $y\qff \in\qff Y$\nnsp.\oss
Then\qss $y\qff \in\qff F\dff(\dff x\trf)$\nnsp.\oss

We will\sss consider only\sss closed\dss multi-valued\dss maps\qss
$\Delta^n\qff \ttoo\qff \Delta^n$\dnsp,\oss
where\dss $\Delta^n$\sss is\dss the standard $n$\dnsp-simplex\dss in\dss $\rrr^{\fff n\dff +\dff 1}$\sss
(see\qss Section\qss \ref{scarf-brouwer}).\qff\oss
Let\dss $F$\dss be such a map.\oss

\myuppar{Matroids\sss and colorings related\dss to\dss $F$\nnsp.}
Let\qss
$B\off =\off
\{\trf
v_{\fff 0}\fff,\pff v_{\fff 1}\fff,\pff \ldots\fff,\pff v_{\fff n}\qff\}$\qss
be\sss the standard\dss basis of\qss $\rrr^{\fff n\dff +\dff 1}$\dss
as\sss in\dss the vector\dss framework\dss from\qss Section\qss
\ref{scarf-oriented-matroids-general}.\oss
Let\qss $X\off \subset\off \Delta^{n}$\qss be a finite set\halfff.\oss
It\dss will\dss be used\dss in\dss the same manner as in\qss 
Scarf's\qss proof\dss of\qss Brouwer's\qss fixed\dss
point\dss theorem.\oss
Let\sss $b$\sss be\sss the barycenter
of\dss $\Delta^n$\nnsp,\oss
i.e.\qss the point\sss of\qss $\rrr^{\fff n\dff +\dff 1}$\sss
with\sss all\sss coordinates equal\dss to\sss $1/(\fff n\dff +\dff 1\fff)$\nnsp.\oss

Let\qss
$f\dff \colon\dff
X\qff \ttoo\qff \Delta^n$\qss
be an arbitrary\sss map,\oss initially\sss unrelated\dss to\dss $F$\nnsp.\oss
The\qss \emph{vector coloring associated\trs with}\qss $f$\qss
is\dss the map\qss
$c\dff \colon\dff
X\qff \ttoo\qff \rrr^{\fff n\dff +\dff 1}$\qss
defined\dss by\vspace{3pt}
\[
\quad
c\qff(\dff x\trf)
\off =\off
f\dff(\dff x\trf)\qff -\qff x\qff +\qff b
\qff.
\]

\vspace{-9pt}
for every\qss $x\qff \in\qff X$\nnsp.\oss
The\qss \emph{vector configuration
associated\trs with}\qss $f$\qss is\dss the subset\qss 
$M\off \subset\off \rrr^{\fff n\dff +\dff 1}$\qss
consisting of\qss 
$v_{\fff 0}\fff,\pff v_{\fff 1}\fff,\pff \ldots\fff,\pff v_{\fff n}$\nsp,\oss
the point\sss $b$\nnsp,\oss
and\dss the vectors\dss $c\dff(\dff x\trf)$\dss for all\qss $x\qff \in\qff X$\nnsp.\oss
The circuits of\trs the vector configuration\dss $M$\dss turn\dss $M$\dss
into an oriented\dss matroid.\oss
Let\qss $A\off =\off M\qff -\qff b$\nnsp.\oss

\mypar{Lemma.}{kakutani-acyclic}
\emph{Let\qss
$f\dff \colon\dff
X\qff \ttoo\qff \Delta^n$\qss
be an arbitrary\sss map
and\dss let\dss $M$\dss be\sss the vector configuration associated\trs to\qss $f$\nnsp.\oss
If\oss
$(\qff y_{\fff v}\trf)_{\dff v\qff \in\qff A}$\qss
is\dss a\dss non-negative solution of\pss the equation\pss
\textup{(\ref{boundedness})},\oss
then\qss $y_{\fff v}\qff \leq\pff 1$\qss for every\qss $v\qff \in\qff A$\nnsp.\oss
In\dss particular\halfff,\oss the set\sss of\dss non-negative solutions\dss 
is\qss bounded\dss and\dss hence\dss the vector configuration\sss $M$\dss
is\trs acyclic as a matroid.\oss}

\proof
In\dss the present\sss situation\dss the equation\qss (\ref{boundedness})\qss
takes\sss the form\vspace{3pt}
\[
\quad
b
\off =\off\qff
\sum\nolimits_{\qff i\qff \in\qff I}\qff y_{\dff i}\trf v_{\fff i}
\off\off +\off\off
\sum\nolimits_{\qff x\qff \in\qff X}\qff
y_{\dff x}\qff c\qff(\dff x\trf)
\qff.
\]

\vspace{-9pt}
Let\dss us\dss take\sss the sum of\dss all\sss coordinates of\dss each\dss
term of\trs this equation.\oss
Obviously\halfff,\oss 
the sum of\dss coordinates of\dss $b$\dss
and of\dss each\sss 
$v_{\fff i}$\sss is\dss $1$\nnsp.\oss
For every\dss $x\qff \in\qff X$\qss the sums of\trs the coordinates of\trs
$x$\dss and\dss $f\dff(\dff x\trf)$\dss are\sss both equal\dss to 
$1$\sss
and\dss hence\sss the sum of\dss coordinates of\dss
$c\dff(\dff x\trf)$\dss is\dss equal\dss to\sss that\sss of\dss $b$\nnsp,\oss
i.e.\qss to $1$\nnsp.\oss
Now\sss the above equation\dss implies\sss that\vspace{3pt}
\[
\quad
1
\off =\off\qff
\sum\nolimits_{\qff i\qff \in\qff I}\qff y_{\dff i}\trf 
\off\off +\off\off
\sum\nolimits_{\qff x\qff \in\qff X}\qff
y_{\fff x}
\off.
\]

\vspace{-9pt}
Since all\dss numbers\dss $y_{\dff i}$\dss and\dss $y_{\dff x}$\dss
are non-negative,\oss
this\sss implies\sss that\sss all\sss of\trs them are\qss $\leq\qff 1$\nnsp.\oss  \eproof

\mypar{Lemma.}{subsequences}
\emph{Suppose\sss that\dss for\sss every\dss natural\dss number $k$
a\dss finite subset\qss 
$X_{\trf k}\off \subset\off \Delta^n$\qss
and\sss a\sss map\qss
$f_{\dff k}\qff \colon\trf
X_{\dff k}\qff \ttoo\qff \Delta^n$\qss
are given.\oss
Let\dss $c_{\dff k}$\dss be\sss the vector coloring associated\dss
with\dss $f_{\dff k}$\nnsp.\oss
Then\dss for some subset\pss $C\off \subset\off I$\qss
the following\dss property\dss holds.\oss
For an infinite set\sss of\trs natural\dss numbers\sss $k$\sss
there exist\dss a\dss $C$\dnsp-cell\pss 
$\sigma_{\dff k}\off \subset\off X_{\trf k}$\qss
such\dss that}\vspace{3pt}
\begin{equation}
\label{kakutani-basis}
\quad
c_{\dff k}\qff(\trf \sigma_{\dff k}\qff)
\off \cup\off\dff 
\bigl\{\qff 
v_{\fff i}\off \bigl|\off
i\qff \in\qff I\qff \smallsetminus\pff C 
\off\bigr\}
\end{equation}

\vspace{-9pt}
\emph{is\dss a\sss basis of\pss $\rrr^{\fff n\dff +\dff 1}$\sss 
containing\sss $b$\dss in\dss its convex\dss hull.\oss}

\proof
Let\dss $M_{\trf k}$\dss be\sss the vector configuration associated\dss with\dss $f_{\dff k}$\nsp.\oss
Lemma\qss \ref{kakutani-acyclic}\qss 
implies\sss that\dss $M_{\trf k}$\dss together\dss with\dss $b$\sss 
satisfies\dss the assumptions of\trs
the vector\qss Scarf\qss framework.\oss
By\dss the classical\qss Scarf\qss theorem\qss there are\sss subsets\qss
$\sigma_{\dff k}\off \subset\off X_{\dff k}$\qss
and\pss
$C_{\dff k}\off \subset\off\dff I$\qss
such\dss that\dss $\sigma_{\dff k}$\dss is\trs 
a\qss $C_{\dff k}$\nnsp-cell\trs
and\dss the set\qss (\ref{kakutani-basis})\qss with\qss
$C\off =\off\dff C_{\dff k}$\qss
is\dss a\sss basis of\qss $\rrr^{\fff n\dff +\dff 1}$\sss 
containing\sss $b$\dss in\dss its convex\dss hull.\oss
Since\sss there\dss is\dss only\sss a\sss finite number of\dss subsets of\trs $I$\nnsp,\oss
the same subset\sss of\trs $I$\dss occurs as\dss $C_{\dff k}$\dss
an\sss infinite number of\trs times.\oss 
Clearly\halfff,\oss one can\dss take\sss such a subset\dss as\dss $C$\nnsp.\oss  \eproof

\myuppar{Kakutani's\qss theorem.}
\emph{Let\pss
$F\dff \colon\dff \Delta^n\qff \ttoo\qff \Delta^n$\qss
be a closed\dss map.\oss
If\qss the set\pss $F\dff(\dff x\trf)$\sss
is\dss convex\dss for\dss every\qss
$x\qff \in\qff \Delta^n$\nnsp,\oss
then\trs $F$\dss has a\sss fixed\dss point\halfff,\oss
i.e.\qss there exists\qss
$z\qff \in\qff \Delta^n$\qss
such\dss that\qss
$z\qff \in\pff F\dff(\trf z\trf)$\nnsp.\oss}

\proof
The proof\dss follows\sss the outline of\qss Scarf's\qss proof\dss
of\qss Brouwer's\qss fixed\dss point\dss theorem\dss from\qss Section\qss \ref{scarf-brouwer},\oss
and\dss we will\dss use\dss the notions introduced\dss
in\dss that\sss section.\oss
Let\qss
$X_{\dff 1}\fff,\pff X_{\dff 2}\fff,\pff  X_{\dff 3}\dff,\pff \ldots$\qss
be a\sss sequence of\dss finite subsets of\dss $\Delta^{n}$\dnsp.\oss
Suppose\sss that\dss the sets\dss $X_{\dff k}$\dss are
chosen\sss in such a\sss way\dss that\dss $X_{\dff k}$\dss
is\dss $\varepsilon_{\dff k}$\hnsp\dnsp-dense\sss in\dss $X_{\dff k}$\dss
for some sequence\qss
$\varepsilon_{\dff k}\qff \ttoo\qff 0$\nnsp.\oss
Let\dss us\dss choose\sss for every\dss natural\dss number $k$ a\sss map\qss
$f_{\dff k}\qff \colon\trf
X_{\dff k}\qff \ttoo\qff \Delta^n$\qss
such\dss that\qss\vspace{3pt}\vspace{-0.375pt}
\[
\quad 
f_{\dff k}\dff(\dff x\trf)\qff \in\pff F\trf(\trf x\trf)
\]

\vspace{-9pt}\vspace{-0.375pt}
for every\dss $x\qff \in\qff X_{\dff k}$\dss
and\dss let\dss $c_{\dff k}$\dss be\sss the vector coloring associated\dss
with\dss $f_{\dff k}$\nnsp.\oss

Lemma\qss \ref{subsequences}\qss implies\sss that\sss after\dss passing\dss
to a subsequence we can assume\sss that\dss there exist\sss a subset\qss
$C\off \subset\off I$\qss independent\sss of\dss $k$\dss
and\dss $C$\dnsp-cells\qss 
$\sigma_{\dff k}\off \subset\off\dff X_{\trf k}$\qss
such\dss that\qss (\ref{kakutani-basis})\qss
is\dss a\sss basis of\pss $\rrr^{\fff n\dff +\dff 1}$\sss 
containing\sss $b$\dss in\dss its convex\dss hull.\oss
Then\sss
every\dss 
$\sigma_{\dff k}$\dss consists of\trs $\num{\fff C\halfff}$\dss elements.\oss
Moreover\halfff,\oss Lemma\qss \ref{minima}\qss implies\sss that\dss maps\qss
$i
\off \ttoo\off
\min\nolimits_{\trf i}\trf \sigma_{\dff k}$\qss 
are\dss bijections\pss
$C\qff \ttoo\qff \sigma_{\dff k}$\nsp.\oss

Let\qss 
$\dis
c_{\trf k}\trf(\dff i\trf)
\off =\off
c_{\trf k}\qff\bigl(\qff \min\nolimits_{\trf i}\trf \sigma_{\dff k} \qff\bigr)$\nnsp.\oss
Now\dss the set\qss (\ref{kakutani-basis})\qss takes\sss the form\vspace{4.5pt}
\[
\quad
\bigl\{\qff 
c_{\trf k}\trf(\dff i\trf)\off \bigl|\off
i\qff \in\pff C 
\off\bigr\}
\off \cup\off
\bigl\{\qff 
v_{\fff i}\off 
\bigl|\off
i\qff \in\qff I\qff \smallsetminus\pff C 
\off\bigr\}
\qff.
\]

\vspace{-7.5pt}
Since\sss this set\dss is\sss a basis containing $b$ 
in\sss its convex\dss hull,\oss\vspace{4.5pt}
\begin{equation}
\label{b-at-a-step}
\quad
b
\off\dff =\off\qff
\sum\nolimits_{\qff i\qff \in\qff I\qff \smallsetminus\pff C}\pff 
y_{\trf k}\trf(\dff i\trf)\qff v_{\fff i}
\off\off +\off\off
\sum\nolimits_{\qff i\qff \in\pff C}\pff
y_{\trf k}\trf(\dff i\trf)\qff c_{\trf k}\trf(\dff i\trf)
\qff,
\end{equation}

\vspace{-7.5pt}
for some non-negative coefficients\dss $y_{\trf k}\trf(\dff i\trf)$\nnsp.\oss

Let\qss 
$\Delta_{\dff k}
\off =\off\dff
\Delta\dff(\dff \sigma_{\dff k}\dff,\off C\trf)$\nnsp.\oss
Since $X_{\dff k}$\dss
is\dss $\varepsilon_{\dff k}$\hnsp\dnsp-dense\sss in\dss $X_{\dff k}$\dss
and\qss
$\varepsilon_{\dff k}\qff \ttoo\qff 0$\nnsp,\oss
Lemma\qss \ref{dominant-dense}\qss implies\sss that\dss the diameters
of\trs the simplices\dss $\Delta_{\dff k}$\dss tend\dss to $0$\nnsp.\oss
After\dss passing\dss to a subsequence\sss one can assume\sss
that\dss the simplices\sss $\Delta_{\dff k}$\sss converge\sss to a\sss
point\qss $z\qff \in\qff \Delta^{n}$\dnsp,\oss
i.e.\qss that\sss every\sss sequence of\dss points\qss
$x\dff(\dff k\trf)\qff \in\qff \Delta_{\dff k}$\qss
converges\sss to\dss $z$\nnsp.\oss
In\dss particular\halfff,\oss
$\min\nolimits_{\trf i}\trf \sigma_{\dff k}$\qss
converges\sss to $z$\sss for\sss every\sss $i$\nnsp.\oss
After\dss passing\dss to a further subsequence one can
assume\sss that\dss for every\qss $i\qff \in\qff I$\qss the sequence\qss\vspace{3pt}
\[
\quad
f_{\trf k}\qff\bigl(\qff \min\nolimits_{\trf i}\trf \sigma_{\dff k} \qff\bigr)
\qff
\]

\vspace{-9pt}
converges\sss to a\trs limit\qss 
$f\dff(\dff i\trf)\qff \in\qff \Delta^n$\dnsp,\oss
and\dss hence\dss $c_{\trf k}\trf(\dff i\trf)$\dss
converges\sss to\qss
$f\dff(\dff i\trf)\qff -\qff z\qff +\qff b$\nnsp.\oss
For each\sss $k$\sss the numbers\qss
$y_{\dff i}
\off =\off 
y_{\trf k}\trf(\dff i\trf)$\qss
form a solution of\qss (\ref{boundedness})\qss
and\dss hence\qss Lemma\qss \ref{kakutani-acyclic}\qss implies\sss that\qss
$y_{\trf k}\trf(\dff i\trf)\qff \leq\qff 1$\qss
for every\qss $k\fff,\pff i$\nnsp.\oss
Therefore,\oss
after\dss passing\dss to a subsequence once more one can
assume\sss that\dss for every\qss $i\qff \in\qff I$\qss the sequence\dss
$y_{\trf k}\trf(\dff i\trf)$\dss
converges\sss
to a\trs limit\qss $y\dff(\dff i\trf)\qff \geq\qff 0$\nnsp.\oss
Now one can\dss pass\dss to\sss the limit\dss in\dss the equation\qss (\ref{b-at-a-step})\qss
and conclude\sss that\vspace{4.5pt}
\begin{equation}
\label{k-limit}
\quad
b
\off\dff =\off\qff
\sum\nolimits_{\qff i\qff \in\qff I\qff \smallsetminus\pff C}\pff 
y\dff(\dff i\trf)\qff v_{\fff i}
\off\off +\off\off
\sum\nolimits_{\qff i\qff \in\pff C}\pff
y\dff(\dff i\trf)\qff 
\left(\qff
f\dff(\dff i\trf)\qff -\qff z\qff +\qff b
\qff\right)
\qff.
\end{equation}

\vspace{-2.625pt}
\textsc{Claim.\hspace*{0.4em}}
\emph{$y\dff(\dff i\trf)\off =\off 0$\qss
for every\qss
$i\qff \in\qff I\qff \smallsetminus\pff C$\qss and}\vspace{4.5pt}
\begin{equation}
\label{subsum-coeff}
\quad
\sum\nolimits_{\qff i\qff \in\pff C}\pff
y\dff(\dff i\trf)
\off\dff =\off\qff
1
\qff.
\end{equation}

\vspace{-2.625pt}
\emph{Proof\qss of\qss the\dss claim.\hspace*{0.5em}}
As we saw\dss in\qss Section\qss \ref{scarf-brouwer},\oss
the intersection of\trs the simplex\qss
$\Delta\dff(\dff \sigma_{\dff k}\dff,\off C\trf)$\sss
with\dss the face of\dss $\Delta^{n}$\sss
defined\dss by\dss the equations\qss
$x_{\dff i}\off =\off 0$\qss
with\qss $i\qff \in\pff I\qff \smallsetminus\pff C$\qss
is\dss non-empty\halfff.\oss
Since\sss the diameters
of\trs the simplices\dss 
$\Delta_{\dff k}$\dss tend\dss to $0$\nnsp,\oss
this\sss implies\sss that\sss $z$\dss is\dss contained\dss in\dss this\sss face,\oss
i.e.\qss that\dss the coordinates\sss $z_{\dff i}$\sss of\dss $z$\sss
with\qss $i\qff \in\pff I\qff \smallsetminus\pff C$\qss
are equal\dss to $0$\nnsp.\oss
Since\sss
the sum of\dss coordinates of\dss every\dss point\sss of\dss $\Delta^n$\dss
is\dss $1$\nnsp,\oss
taking\dss the sum of\dss coordinates\sss in\qss (\ref{k-limit})\qss
shows\sss that\vspace{4.5pt}
\begin{equation}
\label{k-limit-coeff}
\quad
1
\off\dff =\off\qff
\sum\nolimits_{\qff i\qff \in\pff I}\pff
y\dff(\dff i\trf)
\qff.
\end{equation}

\vspace{-7.5pt}
Next\halfff,\oss let\dss us\sss consider\dss the sums of\dss not\sss
all\sss coordinates,\oss but\sss only\sss of\dss coordinates numbered\dss 
by\sss elements of\dss $C$\nnsp.\oss
For $b$\sss this sum\dss is\dss equal\dss to\dss $m\dff/(\dff n\qff +\qff 1\dff)$\nnsp,\oss
where\qss $m\off =\off \num{\dff C}$\nnsp.\oss
For $v_{\fff i}$ with\dss $i\qff \in\qff I\qff \smallsetminus\pff C$\dss
this\sss sum\dss is\dss equal\dss to\sss $0$\nnsp.\oss
Since $z_{\dff i}\off =\off 0$\sss if\dss
$i\qff \in\qff I\qff \smallsetminus\pff C$\nnsp,\oss
for $z$ this\sss sum\dss is\dss equal\dss to $1$\nnsp.\oss
Finally\halfff,\oss 
since\qss $f\dff(\dff i\trf)\qff \in\qff \Delta^n$\dnsp,\oss
for\dss $f\dff(\dff i\trf)$\dss this\sss sum\dss is\qss $\leq\qff 1$\nnsp.\oss
Therefore\qss (\ref{k-limit})\qss implies\sss that\vspace{4.5pt}
\[
\quad
\frac{m}{n\qff +\qff 1}
\off\dff \leq\off\qff
\sum\nolimits_{\qff i\qff \in\pff C}\pff
y\dff(\dff i\trf)\qff
\frac{m}{n\qff +\qff 1}
\qff.
\]

\vspace{-7.5pt}
In\dss turn,\oss this implies\sss that\vspace{4.5pt}
\begin{equation}
\label{estimate}
\quad
1
\off\dff \leq\off\qff
\sum\nolimits_{\qff i\qff \in\pff C}\pff
y\dff(\dff i\trf)
\qff.
\end{equation}

\vspace{-7.5pt}
Since\dss $y\dff(\dff i\trf)$\dss are non-negative,\oss 
the\sss inequalities\qss 
(\ref{estimate})\qss and\qss
(\ref{k-limit-coeff})\qss
together imply\sss that\qss
$y\dff(\dff i\trf)\off =\off 0$\qss
for every\dss
$i\qff \in\pff I\pff \smallsetminus\pff C$\dss
and\dss hence\qss (\ref{subsum-coeff})\qss holds.\oss
This\sss completes\sss the proof\dss of\trs the claim.\oss  \esubproof

\vspace{6pt}
In\dss view of\trs this claim\dss the first\sss sum\dss in\qss (\ref{k-limit})\qss 
is\dss equal\dss to $0$ and\dss hence\qss (\ref{k-limit})\qss
implies\sss that\vspace{4.5pt}
\[
\quad
b
\off\dff =\off\qff
\sum\nolimits_{\qff i\qff \in\pff C}\pff
y\dff(\dff i\trf)\qff 
\left(\qff
f\dff(\dff i\trf)\qff -\qff z\qff +\qff b
\qff\right)
\qff.
\]

\vspace{-7.5pt}
Together\dss with\qss (\ref{subsum-coeff})\qss this implies\sss that\vspace{4.5pt}
\begin{equation*}
\quad
z
\off\dff =\off\qff
\sum\nolimits_{\qff i\qff \in\pff C}\pff
y\dff(\dff i\trf)\qff 
f\dff(\dff i\trf)
\end{equation*}

\vspace{-7.5pt}
and\sss $z$\sss is\dss a\sss convex\sss combination of\trs
the vectors\dss $f\dff(\dff i\trf)$\dss with\qss $i\qff \in\pff C$\nnsp.\oss
Since\dss $F$\dss is\dss a\sss closed\dss map,\oss
$f\dff(\dff i\trf)\qff \in\pff F\dff(\trf z\trf)$\qss
for every\qss $i\qff \in\pff I$\nnsp.\oss
Since\dss $F\dff(\trf z\trf)$\dss is\dss convex,\oss 
this implies\sss that\qss
$z\qff \in\pff F\dff(\trf z\trf)$\nnsp.\oss  \eproof\vspace{4.5pt}

\myuppar{A geometric interpretation of\pss Scarf's\qss proof\halfff.}
Scarf\qss wrote\sss that\qss \emph{``there\dss is\dss an\dss illuminating\dss
geometrical\dss interpretation''}\qss of\dss his\sss theorem\dss in\dss the case
when\dss $X$\dss is\dss a subset\sss of\trs the standard simplex\qss 
$S\off =\off\Delta^n$\dss
as\sss in\dss this section and\qss Section\qss \ref{scarf-brouwer},\oss
and\dss the vector\sss $b$\dss and\dss the colors\dss $c\trf(\dff x\trf)$\nnsp,\oss
where\dss $x\qff \in\qff X$\nnsp,\oss
also belong\dss to\sss the standard simplex.\oss
Scarf\qss suggests\sss to consider\dss the latter standard simplex\dss $S'$\sss 
as\sss different\dss from\sss $S$\sss and\dss to think about\dss the coloring\sss $c$\sss
as a sort\sss of\dss a map\dss $S\qff \ttoo\qff S'$\nnsp.\oss 
He wrote\sss that\qss 
\emph{``the suggestion\dss that\dss the\sss theorem\dss is\dss concerned\dss with\qss
``inverting''\qss an arbitrary\dss mapping\dss is\dss both accurate and\dss
illuminating}\halfff.''\oss
See\qss \cite{sc3},\oss pp.\qss 76\dff--77.\oss

But\dss in\qss Scarf's\qss proof\dss of\qss Kakutani\qss theorem\dss his\sss theorem\dss
is\dss used\dss in\sss a\sss somewhat\sss different\dss manner\halfff.\oss
Namely\halfff,\oss
it\dss is\dss only\dss natural\dss to\sss think\sss about\dss the vectors\dss 
$f\dff(\dff x\trf)\qff -\pff x$\dss
as\sss tangent\dss vectors\sss to $\Delta^n$\dnsp.\oss
Then\dss the colors\dss
$c\qff(\dff x\trf)
\off =\off
f\dff(\dff x\trf)\qff -\qff x\qff +\qff b$\dss
are simply\dss these\sss tangent\dss vectors with\dss their origins 
moved\dss from $0$\sss to\sss $b$\nnsp,\oss
and\dss both\dss the coloring\dss $c$\dss and\dss the map\dss
$x\off \longrightarrow\off f\dff(\dff x\trf)\qff -\pff x$\dss
are interpreted as\sss a\sss vector\dss field\sss on\sss $\Delta^n$\dnsp.\oss
Since\sss $f$\sss maps\dss $\Delta^n$\dss into\dss $\Delta^n$\dnsp,\oss
this\sss vector\dss field\trs is\dss directed\dss inside of\trs the simplex $\Delta^n$
on\dss its\sss boundary\halfff,\oss
and\dss the proof\dss amounts\sss to showing\dss that\qss
(under a continuity\sss assumption)\qss such a vector\dss field always has a zero.\oss

In other\dss words,\oss in\qss Scarf's\qss proof\dss of\qss Kakutani\qss theorem\qss
Scarf\qss theorem\dss is\dss used\sss not\dss to invert\sss a\sss mapping\halfff,\oss
but\dss to find a zero of\dss a\sss vector\dss field.\oss
This agrees with\qss Scarf's\qss intuition\dss that\dss the colors $c\dff(\dff x\trf)$
should\dss be\sss thought\sss as belonging not\dss to $S$\nnsp,\oss
but\dss to another simplex $S'$\dnsp.\oss

\newpage
\mysection{Chains\qss in\pss $\rrr^{\dff n}$\qss and\pss their\qss applications}{scarf-theorem-geometric}

\myuppar{Chains\sss in\dss $\rrr^{\dff n}$\nsp\dnsp.}
One can consider $m$\dnsp-chains in $\rrr^{\dff n}$ without\dss introducing a simplicial\sss
complex in advance or even\dss later\halfff.\oss 
Let\dss us\dss define an \emph{$m$\dnsp-simplex\dss in}\dss $\rrr^{\dff n}$
as\sss a\sss subset\sss of\sss $\rrr^{\dff n}$ consisting of\dss
$m\qff +\qff 1$ points,\oss
and an \emph{$m$\dnsp-chain\dss in}\dss $\rrr^{\dff n}$\dss
as\sss a\sss finite formal\sss sum of\dss $m$\dnsp-simplices in\sss $\rrr^{\dff n}$\dnsp.\oss
The boundary\sss operator $\partial$ acting on $m$\dnsp-chains in\sss $\rrr^{\dff n}$ 
is\dss defined\sss as before,\oss and\dss the identity\dss
$\partial\dff \circ\dff \partial\off =\off 0$\dss still\dss holds,\oss
with\dss the same proof\halfff.\oss
A simplex $\sigma$ is\dss said\dss to\sss be a\qss
\emph{simplex of\trs the chain}\sss $c$\dss if\dss $\sigma$ occurs in\dss the formal\sss
sum $c$ with\dss the coefficient\sss $1$\nnsp.\oss
If\dss $\sigma$ is\dss an $m$\dnsp-simplex $\sigma$ in\sss $\rrr^{\dff n}$\dnsp,\oss 
then\sss
$\conv{\sigma}$\sss denotes\sss the convex\dss hull\sss of\trs the vertices\qss
({\fff}i.e.\qss elements)\qss of\dss $\sigma$\nnsp.\oss
An $m$\dnsp-simplex $\sigma$ is\dss said\dss to be\qss \emph{generic}\qss
if\trs its\dss the vertices are affinely\dss in\-de\-pen\-dent\halfff.\oss
In\dss this case $\conv{\sigma}$ is\dss a\sss geometric $m$\dnsp-simplex\sss in\sss $\rrr^{\dff n}$\dnsp.\oss
In general\halfff,\qss $\conv{\sigma}$ is\dss a\sss polyhedron of\dss dimension\dss $\leq\qff m$\nnsp.\oss
An $m$\dnsp-chain $c$ in\sss $\rrr^{\dff n}$ is\dss said\dss to be\qss
\emph{generic}\pss if\trs every\sss simplex of\dss $c$ is\dss generic.\oss

One can\dss introduce an abstract\sss simplicial\sss complex\dss having\sss
$\rrr^{\dff n}$ as its set\sss of\dss vertices and all\dss finite subsets of\sss $\rrr^{\dff n}$
as its simplices.\oss 
The $m$\dnsp-simplices and $m$\dnsp-chains in\sss $\rrr^{\dff n}$ 
are $m$\dnsp-simplices and $m$\dnsp-chains of\dss this simplicial\sss complex.\oss
This point\sss of\dss view\sss includes\sss these notions in\dss the framework\sss of\trs
the usual\dss theory\halfff,\oss but\trs is\dss not\dss particularly\sss useful\sss otherwise.\oss

\myuppar{Maps\sss to\sss $\rrr^{\dff n}$\dnsp.}
Suppose\sss that $S$ is\dss an abstract\sss simplicial\sss complex\dss and $V$ is\dss its
set\sss of\dss vertices.\oss
Let\dss 
$\varphi\dff \colon\dff
V\qff \ttoo\qff \rrr^{\dff n}$\dss
be a arbitrary\dss map.\oss
For an $m$\dnsp-simplex $\sigma$ of\dss $S$\sss let\qss\vspace{0pt}
\[
\quad
\varphi_{\dff *}\dff(\dff \sigma\dff)
\off =\off
\varphi\dff(\dff \sigma\dff)\qss
\]

\vspace{-12pt}
if\dss $\varphi\dff(\dff \sigma\dff)$ is\dss an $m$\dnsp-simplex\sss in\sss $\rrr^{\dff n}$\dnsp,\oss
i.e.\qss if\dss $\varphi$ is\dss injective on $\sigma$\dnsp,\pss
and\dss let\dss
$\varphi_{\dff *}\dff(\dff \sigma\dff)
\off =\off
0$\dss
otherwise.\oss
Let\dss us\sss extend\sss $\varphi_{\dff *}$\sss to a map from $m$\dnsp-chains of\dss $S$\sss
to $m$\dnsp-chains in\sss $\rrr^{\dff n}$\sss by\dss linearity\halfff.\oss
The resulting\dss map,\pss still\sss denoted\dss by\sss $\varphi_{\dff *}$\nsp,\pss
commutes with\dss the boundary\sss operators in\dss the sense\sss that\dss
$\partial\qff \circ\qff \varphi_{\dff *}
\off =\off\dff
\varphi_{\dff *}\fff \circ\qff \partial$\nnsp.\oss
The proof\dss is\dss the same as for\dss the simplicial\dss maps between
simplicial\sss complexes.\oss
See\qss \cite{i2},\oss Section\qss 1,\oss Theorem\qss 1,\oss for example.\oss

\myuppar{General\dss position and\dss the intersection\dss numbers.}
We will\dss need\sss only\dss the simplest\dss instances of\trs these notions.\oss
We will\dss identify $0$\dnsp-simplices in\sss $\rrr^{\dff n}$ with\dss
the corresponding\dss points in\sss $\rrr^{\dff n}$\dnsp.\oss

Let\sss $\sigma$\dss be an $n$\dnsp-simplex\sss in\sss $\rrr^{\dff n}$
and\dss let\sss $z\qff \in\pff \rrr^{\dff n}$\dnsp.\oss
The point\sss $z$ is\dss said\dss to be\qss 
\emph{in\dss general\dss position\dss with\dss respect\dss to}\dss $\sigma$
if\dss $z$ does not\dss belong\dss to\dss $\conv{\tau}$\dss for
every $(\fff n\dff -\dff 1\fff)$\dnsp-face $\tau$ 
of\trs $\sigma$\nnsp.\oss
In\dss this case\sss
the\qss \emph{intersection\dss number}\dss $\sigma\dff \cdot\trf z$\dss
is\dss defined\dss as\dss $1\qff \in\pff \ftwo$\qss if\qss $z\qff \in\qff \conv{\sigma}$\qss
and\dss $0\qff \in\pff \ftwo$\qss otherwise.\oss
If\dss $c$\sss is\dss an $n$\dnsp-chain\dss
in\sss $\rrr^{\dff n}$\sss and $d$\dss is\dss a $0$\dnsp-chain\sss
in\sss $\rrr^{\dff n}$\dnsp,\oss
then\dss $c\fff,\pff d$\dss are said\dss to be\qss
\emph{in\dss general\dss position}\pss if\dss every\sss simplex of\dss $d$\dss is\dss
in\dss general\dss position\dss with\dss respect\dss to every\sss simplex of\dss $c$\nnsp.\oss
In\dss this case\sss
the\qss \emph{intersection\dss number}\dss $c\dff \cdot\dff d$\dss
is\dss defined as\sss the sum of\trs the intersection\dss numbers\dss 
$\sigma\dff \cdot\trf z$\dss
over all\dss pairs\dss $\sigma\fff,\pff z$\dss such\dss that\sss $\sigma$ is\dss
a simplex of\dss $c$ and $z$\sss is\dss a simplex of\dss $d$\nnsp.\oss

Next\halfff,\oss let\sss $\tau$\dss be an 
$(\fff n\dff -\dff 1\fff)$\dnsp-simplex\dss
in\sss $\rrr^{\dff n}$\dnsp,\oss
and\dss let\sss $\omega$\dss be a $1$\dnsp-simplex.\oss
Since $\omega$ is\dss a\sss two-points subset\sss of\trs $\rrr^{\dff n}$\dnsp,\oss
the convex\dss hull\sss $\conv{\omega}$ is\dss a\sss segment\halfff.\oss
The $1$\dnsp-simplex $\omega$ is\dss said\dss to be\qss 
\emph{in\dss general\dss position\dss with\dss respect\dss to}\dss $\tau$
if\trs both vertices of\dss $\omega$ do not\dss belong\dss to $\conv{\tau}$
and $\conv{\omega}$ is\dss disjoint\dss from $\conv{\upsilon}$
for every\sss $(\fff n\dff -\dff 2\fff)$\dnsp-face $\upsilon$ 
of\trs $\tau$\nnsp.\oss 
If\trs this\dss is\dss the case,\oss then
the\qss \emph{intersection\dss number}\dss $\tau\dff \cdot\dff \omega$\dss
is\dss defined\dss as\dss $1\qff \in\pff \ftwo$\qss if\qss 
$\conv{\omega}\pff \cap\qff \conv{\tau}\off \neq\off \varnothing$\qss
and\dss $0\qff \in\pff \ftwo$\qss otherwise.\oss
If\dss $c$\sss is\dss an $(\fff n\dff -\dff 1\fff)$\dnsp-chain\dss
in\sss $\rrr^{\dff n}$\sss and $d$\dss is\dss a $1$\dnsp-chain\sss
in\sss $\rrr^{\dff n}$\dnsp,\oss
then\dss $c\fff,\pff d$\dss are said\dss to be\qss
\emph{in\dss general\dss position}\pss if\dss every\sss simplex of\dss $d$\dss is\dss
in\dss general\dss position\dss with\dss respect\dss to every\sss simplex of\dss $c$\nnsp.\oss
In\dss this case\sss
the\qss \emph{intersection\dss number}\dss $c\dff \cdot\dff d$\dss
is\dss defined as\sss the sum of\trs the intersection\dss numbers\dss 
$\tau\dff \cdot\qff \omega$\dss
over all\dss pairs\dss $\tau\fff,\pff \omega$\dss such\dss that\sss $\tau$ is\dss
a simplex of\dss $c$ and $\omega$ is\dss a simplex of\dss $d$\nnsp.\oss\vspace{-0.6pt}

Finally\halfff,\oss let $\sigma$\dss be an\sss $n$\dnsp-simplex\sss in\sss $\rrr^{\dff n}$
and\sss $\omega$\dss be a $1$\dnsp-simplex.\oss
Then $\omega$ is\dss said\dss to be\qss 
\emph{in\dss general\dss position\dss with\dss respect\dss to}\qss $\sigma$\dss
if\dss every\sss vertex of\dss $\omega$\dss is\dss in\dss general\dss position\dss
with\dss respect\dss to $\sigma$ and\sss $\omega$\dss is\dss in\dss general\dss position\dss
with\dss respect\dss to every $(\fff n\dff -\dff 1\fff)$\dnsp-face of\dss $\sigma$\nnsp.\oss
If\dss $c$\sss is\dss an $n$\dnsp-chain\dss
in\sss $\rrr^{\dff n}$\sss and $d$\dss is\dss a $1$\dnsp-chain\sss
in\sss $\rrr^{\dff n}$\dnsp,\oss
then\dss $c\fff,\pff d$\dss are said\dss to be\qss
\emph{in\dss general\dss position}\pss if\dss every\sss simplex of\dss $d$\dss is\dss
in\dss general\dss position\dss with\dss respect\dss 
to every\sss simplex of\dss $c$\nnsp.\oss\vspace{-0.6pt}

\mypar{Lemma.}{intersections-generic}
\emph{Suppose\sss that $\sigma$ and $\omega$ are\dss a\sss generic
$n$\dnsp-simplex\dss and\dss a\dss $1$\dnsp-simplex\dss in\qss $\rrr^{\dff n}$ re\-spec\-tive\-ly\halfff.\oss
If\qss $\sigma\fff,\pff \omega$\qss are\sss in\dss general\trs position,\oss then\oss
$\partial\dff \sigma\qff \cdot\qff \omega
\off =\off 
\sigma\qff \cdot\pff \partial\dff \omega$\nsp.\oss}\vspace{-0.6pt}

\proof
Being a $1$\dnsp-simplex,\pss $\omega$ is\dss always generic\sss because\sss
two different\dss points are automatically\sss affinely\dss independent\halfff.\oss
Since $\sigma$ is\dss generic,\qss
$\conv{\sigma}$ is\dss a\sss geometric $n$\dnsp-simplex.\vspace{-0.6pt}

If\qss both vertices of\dss $\omega$ belong\dss to $\conv{\sigma}$\nnsp,\oss
then\dss 
$\sigma\qff \cdot\qff \partial\dff \omega
\off =\off
2\off =\off 0$\nnsp.\oss
Also,\oss in\dss this case\dss $\conv{\omega}$\dss is\dss contained\dss in\dss  
$\conv{\sigma}$\dss
and\dss hence $\conv{\omega}$ is\dss disjoint\dss from\dss the boundary\sss of\dss
$\conv{\sigma}$\nnsp.\oss
It\dss follows\dss that\dss in\dss this case\dss 
$\partial\dff \sigma\qff \cdot\qff \omega
\off =\off
0$\dss
and\dss hence\dss
$\partial\dff \sigma\qff \cdot\qff \omega
\off =\off 
\sigma\qff \cdot\qff \partial\dff \omega$\nnsp.\oss\vspace{-0.6pt}

If\dss only\sss one vertex of\dss $\omega$ is\dss contained\dss in $\conv{\sigma}$\nnsp,\oss
then\dss 
$\sigma\qff \cdot\qff \partial\dff \omega
\off =\off
1$\nnsp.\oss
Also,\oss in\dss this case\sss the segment\dss $\conv{\omega}$\dss intersects
exactly\sss one $(\fff n\dff -\dff 1\fff)$\dnsp-face of\trs the
geometric $n$\dnsp-simplex $\conv{\sigma}$\nnsp.\oss
It\dss follows\dss that\dss in\dss this case\dss
$\partial\dff \sigma\qff \cdot\qff \omega
\off =\off
1$\dss
and\dss hence\dss
$\partial\dff \sigma\qff \cdot\qff \omega
\off =\off 
\sigma\qff \cdot\qff \partial\dff \omega$\nnsp.\oss\vspace{-0.6pt}

If\qss neither of\trs the vertices of\dss $\omega$ is\dss contained\dss in $\conv{\sigma}$\nnsp,\oss
then\dss 
$\sigma\qff \cdot\qff \partial\dff \omega
\off =\off
0$\nnsp.\oss
The intersection\dss
$\conv{\omega}\qff \cap\qff \conv{\sigma}$\dss
is\dss convex and\dss contained\dss in\dss the segment\dss $\conv{\omega}$\nnsp.\oss
There\-fore\dss
$\conv{\omega}\qff \cap\qff \conv{\sigma}$\dss 
is\dss either empty\halfff,\oss
or\dss is\dss a\sss segment\halfff,\oss
or\dss consists of\dss one point\halfff.\oss
In\dss the first\sss case\dss
$\partial\dff \sigma\qff \cdot\qff \omega
\off =\off
0
\off =\off
\sigma\qff \cdot\qff \partial\dff \omega$\nnsp,\oss
in\dss the second case\dss $\conv{\omega}$ intersects\sss two $(\fff n\dff -\dff 1\fff)$\dnsp-faces
of\dss $\conv{\sigma}$ and\dss hence also\dss
$\partial\dff \sigma\qff \cdot\qff \omega
\off =\off
0$\nnsp.\oss
In\dss the\sss last\sss case 
$\conv{\omega}$ intersects an $(\fff n\dff -\dff 2\fff)$\dnsp-face
of\dss $\conv{\sigma}$\nnsp,\oss 
which\dss is\dss impossible because\dss $\sigma\fff,\off \omega$\dss
are assumed\dss to be in\sss general\dss position.\oss
The lemma follows.\oss  \eproof\vspace{-0.6pt}

\mypar{Lemma.}{intersection-point-cycle}
\emph{Suppose\sss that\sss $c$\sss is\sss an
{\nsp}$n$\dnsp-chain\dss in\qss $\rrr^{\dff n}$\dss
such\dss that\qss
$c\qff \cdot\pff \partial\dff \kappa
\off =\off 0$\qss
for\dss every\dss $1$\dnsp-simplex\dss $\kappa$\dss in\pss $\rrr^{\dff n}$\dss in\dss general\dss position\dss
with\dss respect\dss to $c$\nnsp.\oss
Then\qss
$c\qff \cdot\qff z\off =\off 0$\qss
for\dss every\dss point\qss
$z\qff \in\pff \rrr^{\dff n}$\trs
in\dss general\trs 
position\dss with\dss respect\dss to\sss $c$\nnsp.\oss}\vspace{-0.6pt}

\proof
Let\sss $C$\sss be\sss the union of\dss the geometric simplices\dss $\conv{\sigma}$\dss
over all\sss $n$\dnsp-simplices $\sigma$ of\dss $c$\nnsp.\oss
Since $c$ is\dss a\sss finite sum,\oss the set\sss $C$\sss is\dss bounded.\oss
Clearly\halfff,\oss there exists a\sss point\dss $u\qff \in\pff \rrr^{\dff n}$
such\dss that\trs $u\qff \not\in\qff C$\dss and\dss the segment\sss connecting\sss
$z$\sss and\dss $u$\sss is\dss disjoint\dss from\dss $\conv{\upsilon}$\dss
for every\dss $(\fff n\dff -\dff 2\fff)$\dnsp-face of\dss every\sss $n$\dnsp-simplex of\dss $c$\nnsp.\oss
Then\dss the $1$\dnsp-simplex\dss
$\kappa\off =\off \{\qff z\fff,\qff u\qff\}$\dss is\dss in\dss general\dss position\dss
with\sss respect\dss to\sss $c$\nnsp.\oss
By\dss the assumptions of\trs the lemma,\pss
$c\qff \cdot\pff \partial\dff \kappa
\off =\off 
0$\nsp.\oss
But\trs $\partial\dff \kappa\off =\off z\qff +\qff u$\dss
and\dss hence\dss
$c\qff \cdot\pff \partial\dff \kappa
\off =\off
c\qff \cdot\qff z\pff +\pff c\qff \cdot\qff u$\nnsp.\oss
Clearly\halfff,\pss $u\qff \not\in\qff C$\dss implies\sss that\trs
$c\qff \cdot\qff u\off =\off 0$\nnsp.\oss
Together\dss with\dss
$c\qff \cdot\pff \partial\dff \kappa
\off =\off
0$\dss
this\sss implies\sss that\trs $c\qff \cdot\qff z\off =\off 0$\nnsp.\oss  \eproof\vspace{-0.75pt}

\mypar{Lemma\qss (Carath\'{e}odory\fff).}{caratheodory}
\emph{Let\trs $\sigma$ be\sss an\sss $l$\dnsp-simplex\dss in\dss $\rrr^{\dff n}$\sss
for some $l$\dss and\dss let\dss $m$\dss be\sss the dimension of\qss $\conv{\sigma}$\nnsp.\oss
Then\sss every\dss point\sss of\qss $\conv{\sigma}$ belongs\dss to\dss
$\conv{\tau}$\dss for some $m$\dnsp-face\dss $\tau$ of\qss $\sigma$\nnsp.\oss}\vspace{-0.75pt}

\proof
Since\sss $m$ is\dss the dimension of\dss $\conv{\sigma}$\nnsp,\oss
the convex\dss hull\sss $\conv{\sigma}$\sss is\dss contained\dss in\dss
an affine hyperplane $A$ of\dss dimension $m$\nnsp.\oss
Let\trs $x\qff \in\qff \conv{\sigma}$\nnsp.\oss
By\dss a\sss classical\dss theorem of\qss
Carath\'{e}odory\qss applied\dss to\sss the subset\sss $\sigma$ of\dss $A$\nnsp,\oss
there exists a subset\trs
$\tau\qff \subset\qff \sigma$\dss
such\dss that\trs $\num{\tau}\off =\off m\qff +\qff 1$\dss and\dss 
$x\qff \in\qff \conv{\tau}$\nnsp.\oss
Clearly\halfff,\qss $\tau$ is\dss an\sss $m$\dnsp-face of\dss $\sigma$\nnsp.\oss  \eproof\vspace{-0.75pt}

\mypar{Lemma.}{intersections-non-generic}
\emph{Suppose\sss that $\sigma$ and\dss $\omega$ are\dss an\sss 
$n$\dnsp-simplex\dss and\dss a\dss $1$\dnsp-simplex\dss in\qss $\rrr^{\dff n}$ re\-spec\-tive\-ly\halfff.\oss
If\pss $\sigma\fff,\pff \omega$\dss are\sss in\dss general\trs position
and $\sigma$ is\dss not\dss generic,\oss then\oss
$\partial\dff \sigma\qff \cdot\qff \omega
\off =\off 
\sigma\qff \cdot\pff \partial\dff \omega
\off =\off
0$\nnsp.}\vspace{-0.75pt}

\proof
Since $\sigma$ is\dss not\dss generic,\oss
the vertices of\dss $\sigma$ are affinely\sss dependent\sss
and $\conv{\sigma}$ is\dss a\sss polyhedron of\dss dimension\qss $\leq\qff n\qff -\qff 1$\nnsp.\oss
Suppose first\dss that\dss the dimension of\dss $\conv{\sigma}$ is\qss $\leq\qff n\qff -\qff 2$\nnsp.\oss
Let\trs $x\qff \in\qff \conv{\sigma}$\nnsp.\oss
Lemma\qss \ref{caratheodory}\qss implies\sss that\trs
$x\qff \in\qff \conv{\upsilon\fff}$\dss for some\dss $m\qff \leq\qff n\qff -\qff 2$\dss
and some $m$\dnsp-face $\upsilon$ of\dss $\sigma$\nnsp.\oss
In\dss turn,\oss this implies\sss that\sss $\upsilon$ is\dss an $m$\dnsp-face
of\dss some $(\fff n\dff -\dff 1\fff)$\dnsp-face $\tau$ of\dss $\sigma$\nnsp.\oss
Since\dss $\sigma\fff,\pff \omega$\dss are in\dss general\dss position,\oss
this implies\sss that\trs 
$\conv{\upsilon\fff}\qff \cap\qff \conv{\omega}
\off =\off
\varnothing$\dss
and\dss hence\dss $x\qff \not\in\qff \conv{\omega}$\nnsp.\oss
Since $x$ was an arbitrary\dss point\sss of\dss $\conv{\sigma}$\nnsp,\oss
we see\sss that\trs
$\conv{\sigma}\qff \cap\qff \conv{\omega}
\off =\off
\varnothing$\dss
and\dss hence\qss
$\partial\dff \sigma\qff \cdot\qff \omega
\off =\off 
\sigma\qff \cdot\qff \partial\dff \omega
\off =\off
0$\nnsp.\oss
This proves\sss the\sss theorem\dss when\dss
the dimension of\dss $\conv{\sigma}$ is\qss $\leq\qff n\qff -\qff 2$\nnsp.\oss

Suppose\sss that\dss the dimension of\dss $\conv{\sigma}$ is\qss
$=\off n\qff -\qff 1$\nnsp.\oss
Then $\conv{\sigma}$ is\dss contained\sss in an affine plane\dss
$H\pff \subset\pff \rrr^{\dff n}$ of\dss dimension\dss $n\qff -\qff 1$\nnsp.\oss
The segment\sss $\omega$ is\dss either\sss disjoint\dss from\sss $H$\nnsp,\oss
or\dss intersects $H$ in one of\trs the\sss two vertices of\dss $\omega$\nnsp,\oss
or\dss intersects $H$ in an interior\sss point\halfff,\oss
or\dss is\dss contained\sss in\sss $H$\nnsp.\oss

If\dss $\conv{\omega}$ is\dss disjoint\dss from $H$\nnsp,\oss
then\dss
$\conv{\omega}\qff \cap\qff\conv{\sigma}
\off =\off
\varnothing$\dss
and\dss hence\qss
$\partial\dff \sigma\qff \cdot\qff \omega
\off =\off 
\sigma\qff \cdot\qff \partial\dff \omega
\off =\off
0$\nnsp.\oss

Suppose\sss that\sss $\conv{\omega}$ intersects $H$ in a vertex\sss $z$ of\dss $\omega$\nnsp.\oss
Lemma\qss \ref{caratheodory}\qss implies\sss that\sss
$\conv{\sigma}$ is\dss equal\dss to\sss the union of\dss the sets
$\conv{\tau}$ over all\sss $(\fff n\dff -\dff 1\fff)$\dnsp-faces $\tau$ of\dss $\sigma$\nnsp.\oss
Since\dss $\sigma\fff,\pff \omega$\dss 
are in\dss general\dss position,\oss
this implies\sss that\sss 
$z\qff \not\in\qff \conv{\sigma}$\nnsp.\oss
In\dss turn,\oss this implies\sss that\qss
$\partial\dff \sigma\qff \cdot\qff \omega
\off =\off 
\sigma\qff \cdot\qff \partial\dff \omega
\off =\off
0$\nnsp.\oss

Suppose now\sss that\sss $\conv{\omega}$ intersects $H$ only\dss in a point\sss $u$\sss
which\dss is\dss not\sss
a vertex of\dss $\omega$\nnsp.\oss
In\dss this case\qss
$\sigma\qff \cdot\qff \partial\dff \omega
\off =\off
0$\nnsp.\oss
Also,\oss in\dss this case $u$ is\dss in\sss general\dss position\dss
with\sss respect\dss to $\partial\dff \sigma$ in $H$ and\qss
$\partial\dff \sigma\qff \cdot\qff \omega
\off =\off
\partial\dff \sigma\dff \cdot\qff u$\nnsp,\oss
where\dss the second\dss intersection\sss number\dss is\dss understood\dss in\sss $H$\nnsp.\oss
Let\dss us\dss identify\sss $H$\sss with\sss $\rrr^{\dff n\dff -\dff 1}$\dnsp.\oss
Let\trs $c\off =\off \partial\dff \sigma$\dss
and\dss let\sss $\kappa$\sss be a $1$\dnsp-simplex\sss in\sss $\rrr^{\dff n\dff -\dff 1}$\nnsp.\oss
Arguing\dss by\dss induction,
we can assume\sss that\dss the\sss theorem\dss is\dss true for\dss the intersections
in\sss $\rrr^{\dff n\dff -\dff 1}$\nnsp.\oss
Then\dss\vspace{1.5pt}
\[
\quad
c\qff \cdot\pff \partial\dff \kappa
\off =\off
\partial\dff c\qff \cdot\pff \kappa
\off =\off
\partial\dff \circ\dff \partial\qff (\dff \sigma\dff)\qff  \cdot\pff \kappa
\off =\off
0\qff \cdot\qff \kappa
\off =\off
0
\qff
\]

\vspace{-10.5pt}
and\dss hence\trs Lemma\qss \ref{intersection-point-cycle}\qss
with\dss $n\qff -\qff 1$\dss in\dss the role of\dss $n$
implies\sss that\trs
$c\qff \cdot\qff u
\off =\off
0$\nnsp.\oss
It\dss follows\dss that\trs
$\partial\dff \sigma\qff \cdot\qff \omega
\off =\off
\partial\dff \sigma\qff \cdot\qff u
\off =\off
c\qff \cdot\qff u
\off =\off
0$\nnsp.\oss
Therefore\qss
$\partial\dff \sigma\qff \cdot\qff \omega
\off =\off 
\sigma\qff \cdot\qff \partial\dff \omega
\off =\off
0$\qss
in\dss this case also.\oss

It\dss remains\sss to consider\dss the case when\trs $\conv{\omega}\qff \subset\pff H$\nnsp.\oss
Since\dss $\sigma\fff,\pff \omega$\dss are in\dss general\dss position,\pss
$\omega$ is\dss in\dss general\dss position with
respect\dss to every\sss $(\fff n\dff -\dff 1\fff)$\dnsp-face $\tau$ of\dss $\sigma$\nnsp.\oss
This means\sss that\dss the vertices of\dss $\omega$ are not\sss contained\sss in\sss
$\conv{\tau}$ and $\conv{\omega}$ is\dss disjoint\dss from $(\fff n\dff -\dff 2\fff)$\dnsp-faces\sss of\dss
$\conv{\tau}$ and\dss hence from\dss the boundary\sss of\dss $\conv{\tau}$\nnsp.\oss
It\dss follows\dss that\sss $\conv{\omega}$ is\dss disjoint\dss from $\conv{\tau}$\nnsp.\oss
By\qss Lemma\qss \ref{caratheodory}\qss this implies\sss that\sss
$\conv{\omega}$ is\dss disjoint\dss from $\conv{\sigma}$
and\dss hence\qss
$\partial\dff \sigma\qff \cdot\qff \omega
\off =\off 
\sigma\qff \cdot\qff \partial\dff \omega
\off =\off
0$\nnsp.\oss  \eproof

\mypar{Theorem.}{intersections-general}
\emph{Suppose\sss that $c$ and $d$ are\dss an
$n$\dnsp-chain\dss and\dss a\dss $1$\dnsp-chain\dss in\qss $\rrr^{\dff n}$ re\-spec\-tive\-ly\halfff.\oss
If\qss the chains\dss $c\fff,\pff d$\qss are\sss in\dss general\trs position,\oss then\oss
$\partial\dff c\qff \cdot\qff d\off =\off c\qff \cdot\pff \partial\dff d$\nsp.\oss}

\proof
In view of\trs the definition of\trs the intersection\dss numbers of\dss chains,\oss
it\dss is\dss sufficient\dss to prove\sss that\qss
$\partial\dff \sigma\qff \cdot\qff \omega\off =\off \sigma\qff \cdot\qff \partial\dff \omega$\qss
if\dss $\sigma$ is\dss an $n$\dnsp-simplex,\pss
$\omega$ is\dss a\sss $1$\dnsp-simplex,\oss
and\dss $\sigma\fff,\pff \omega$\dss are\sss in\dss general\dss position.\oss
If\dss $\sigma$ is\dss a\sss generic $n$\dnsp-simplex,\oss
then\dss this equality\dss follows\dss from\qss Lemma\qss \ref{intersections-generic},\oss
if\trs it\dss is\dss not\halfff,\oss
then\dss from\qss Lemma\qss \ref{intersections-non-generic}.\oss  \eproof

\mypar{Corollary.}{cycle-zero}
\emph{Suppose\sss that\sss $c$\sss is\sss an
{\nsp}$n$\dnsp-chain\dss in\qss $\rrr^{\dff n}$\dss
and\qss
$z\qff \in\pff \rrr^{\dff n}$\trs
is\dss a\sss  point\dss
in\dss general\trs 
position\dss with\dss respect\dss to\sss $c$\nnsp.\oss
If\qss $c$ is\dss a\sss cycle,\oss
i.e.\qss if\oss $\partial\dff c\off =\off 0$\nnsp,\oss
then\qss $c\qff \cdot\qff z\off =\off 0$\nnsp.\oss}

\proof
Let\sss $\kappa$\sss be a $1$\dnsp-simplex\dss in\pss 
$\rrr^{\dff n}$\dss in\dss general\dss position\dss
with\dss respect\dss to $c$\nnsp.\oss
Then\qss Theorem\qss \ref{intersections-general}\qss implies\sss that\trs 
$c\qff \cdot\pff \partial\dff \kappa
\off =\off
\partial\dff c\qff \cdot\qff \kappa
\off =\off
0\qff \cdot\qff \kappa
\off =\off
0$\nnsp.\oss
It\dss remains\sss to apply\qss Lemma\qss \ref{intersection-point-cycle}.\oss  \eproof

\myuppar{The affine framework.}
Let\sss $\mathcal{D}$ be\sss a\sss chain-simplex\dss based on  
$I\pff =\pff \{\qff 0\fff,\pff 1\fff,\pff \ldots\fff,\pff n\qff\}$\nnsp,\pss
$X\off =\off V_{\dff \mathcal{D}}$\nnsp,\oss
and\dss let $\mathcal{E}$ 
be\sss the envelope of\dss $\mathcal{D}$\dnsp.\oss
An\qss \emph{affine coloring}\pss of\trs $\mathcal{D}$\dss is\dss simply\sss
a map\dss
$c\dff \colon\dff
X\qff \ttoo\qff \rrr^{\dff n}$\nsp\dnsp.\oss

Let\trs
$w_{\fff 0}\fff,\pff w_{\fff 1}\fff,\pff \ldots\fff,\pff w_{\fff n}$\dss
be affinely\dss independent\dss points in\dss $\rrr^{\dff n}$\sss
and\dss $\Gamma$\dss be\sss the geometric\sss the $n$\dnsp-simplex\dss
in\dss $\rrr^{\dff n}$ having\trs these points
as\dss its\dss vertices.\oss

\mypar{Lemma.}{chain-cover}
\emph{Let\qss $\varphi\dff \colon\dff X\qff \cup\pff I\qff \ttoo\qff \rrr^{\dff n}$\qss
be a map such\dss that\qss
$\varphi\dff(\dff i\trf)\off =\off w_{\fff i}$\qss for every\trs $i\qff \in\pff I$\nnsp.\oss
If\trs $z\qff \in\pff \Gamma$\dss and\dss
$z$\sss 
is\dss in\sss general\dss position\dss with respect\dss to\sss
the chain\dss
$\varphi_{\fff *}\qff
\mathcal{E}\trf\fclass{\trf I\qff}$\dss
and\dss the $n$\dnsp-simplex\dss $\varphi\dff(\trf I\trf)$\nnsp,\oss
then\dss $z\qff \in\qff
\conv{\dff\varphi\dff(\trf \sigma\trf)\fff}$\trs
for some $n$\dnsp-simplex\sss $\sigma$ of\pss $\mathcal{E}\dff(\trf I\trf)$\nnsp.\oss}

\proof
Theorem\qss \ref{envelope-chain}\qss implies\sss that\trs
$\partial\qff \mathcal{E}\trf\fclass{\trf I\qff}
\off =\off
\partial\trf I$\dss
and\dss hence\vspace{4.125pt}
\[
\quad
\partial\qff \left(\qff
\varphi_{\fff *}\qff
\mathcal{E}\trf\fclass{\trf I\qff}
\qff\right)
\off =\off
\varphi_{\fff *}\qff
\bigl(\trf
\partial\trf I
\qff\bigr)
\off =\off
\partial\qff
\bigl(\trf
\varphi_{\fff *}\dff(\trf I\trf)
\qff\bigr)
\qff.
\]

\vspace{-7.875pt}
In\dss turn,\oss this implies\sss that\qss
$\varphi_{\fff *}\qff
\mathcal{E}\trf\fclass{\trf I\qff}
\off +\off
\varphi_{\fff *}\dff(\trf I\trf)$\qss
is\dss a\sss cycle.\oss
Now\trs Corollary\qss \ref{cycle-zero}\qss implies\sss that\vspace{4pt}
\[
\quad
\varphi_{\fff *}\qff
\mathcal{E}\trf\fclass{\trf I\qff}\dff \cdot\qff z
\off =\off
\varphi_{\fff *}\dff(\trf I\trf)\dff \cdot\qff z
\qff.
\]

\vspace{-8pt}
The simplex\trs $\Gamma$\dss is\dss the convex\dss hull\trs $\Gamma$\sss of\trs
$\varphi\dff(\trf I\trf)
\off =\off
\varphi_{\fff *}\dff(\trf I\trf)$\nnsp,\oss
and\dss hence\dss
$z\qff \in\pff \Gamma$\dss implies\sss that\trs
$\varphi_{\fff *}\dff(\trf I\trf)\dff \cdot\qff z
\off =\off
1$\nnsp.\oss
This,\oss in\dss turn,\oss implies\sss that\trs
$\varphi_{\fff *}\qff
\mathcal{E}\trf\fclass{\trf I\qff}\dff \cdot\qff z
\off =\off
1$\dss
and\dss hence\dss
$z\qff \in\qff
\conv{\dff\varphi\dff(\trf \sigma\trf)\fff}$\dss
for some $n$\dnsp-simplex\sss $\sigma$ of\pss 
$\mathcal{E}\dff(\trf I\trf)$\nnsp.\oss  \eproof

\mypar{The main\dss theorem\sss for\dss affine colorings.}{main-theorem-vector-new}
\emph{For\dss every\trs map\dss 
$c\dff \colon\dff X\dff \ttoo\qff  \rrr^{\dff n}$\sss
and\dss every\dss $z\qff \in\pff \Gamma$\dss 
there exist\sss a\sss non-empty\dss subset\qss $C\qff \subset\pff I$\qss
and\trs a\sss $d\dff(\trf C\trf)$\dnsp-simplex\dss $\sigma$\dss of\oss
$\mathcal{D}\dff(\trf C\trf)$\dss such\dss that\dss the set}\vspace{4.5pt}
\begin{equation}
\label{affine-colorings-simplex}
\quad
c\qff(\trf \sigma\trf)\off \cup\off\dff \bigl\{\qff w_{\fff i}\off \bigl|\off
i\qff \in\qff I\qff \smallsetminus\pff C \pff\bigr\}
\end{equation}

\vspace{-7.5pt}
\emph{consists of\qss $n\qff +\qff 1$\dss elements,\oss
is\dss an\sss affinely\dss independent\halfff,\oss
and contains\sss $z$\sss in\dss its convex\dss hull.\oss}

\proof
Let\dss us\dss extend $c$\sss to a\sss map\dss
$\varphi\dff \colon\dff
X\qff \cup\pff I\qff \ttoo\qff \rrr^{\dff n}$\dss
by\dss the rule\dss
$\varphi\dff(\dff i\trf)\off =\off w_{\fff i}$\nnsp.\oss
By\qss Lemma\qss \ref{top-simplices}\pss it\trs is\dss sufficient\dss to
prove\sss that\dss there exists an $n$\dnsp-simplex $\sigma$
of\dss $\mathcal{E}\dff(\trf I\trf)$ such\dss that\dss the set\sss $\varphi\dff(\dff \sigma\dff)$
is\dss affinely\sss independent\sss and\sss contains $z$ in\dss its convex\dss hull.\oss
Cf.\qss the proof\dss of\qss Theorem\qss \ref{main-non-degenerate}.\oss
Since\dss $z\qff \in\pff \Gamma$\nnsp,\oss
there exists a sequence\dss
$z_{\dff 1}\dff,\off z_{\dff 2}\dff,\off z_{\dff 3}\dff,\pff \ldots$\dss
of\dss points\dss $z_{\dff i}\qff \in\pff \Gamma$\sss
such\dss that\dss $z_{\dff i}$\sss tends\sss to $z$ when\dss $i\qff \ttoo\qff \infty$\dss
and\dss for every\sss $i$\sss the point\sss $z_{\dff i}$\dss is\dss not\sss 
an\sss affine\sss linear combination of\trs
$<\qff n\qff +\qff 1$\dss 
elements\sss of\trs the union of\qss 
$\varphi\dff(\trf X\trf)\qff \cup\pff
\{\trf w_{\fff 0}\fff,\pff w_{\fff 1}\fff,\pff \ldots\fff,\pff w_{\fff n} \qff\}$\nnsp.\oss
This implies\sss that\sss $z_{\dff i}$\sss is\dss in\sss the general\dss position\dss
with respect\dss to\sss
$\varphi_{\fff *}\qff
\mathcal{E}\trf\fclass{\trf I\qff}$\dss
and\dss $\varphi\dff(\trf I\trf)$\sss
for every\sss $i$\nnsp.\oss
Since\trs\vspace{4pt}
\[
\quad
\varphi\dff(\trf I\trf)
\off =\off 
\bigl\{\qff
w_{\fff 0}\fff,\pff w_{\fff 1}\fff,\pff \ldots\fff,\pff w_{\fff n}
\qff\fff\bigr\}
\qff,
\]

\vspace{-8pt}
Lemma\qss \ref{chain-cover}\qss implies\sss that\dss for every\sss $i$\sss
there exists an $n$\dnsp-simplex $\sigma$
of\dss $\mathcal{E}\dff(\trf I\trf)$ such\dss that\dss the set\sss $\varphi\dff(\dff \sigma\dff)$
contains $z_{\dff i}$ in\dss its convex\dss hull.\oss
Since\dss $\mathcal{E}\dff(\trf I\trf)$\dss has only\sss a\sss finite number
of\dss $n$\dnsp-simplices,\oss we can assume\sss that\dss this
simplex $\sigma$ is\dss the same for\sss all\sss $i$\nnsp.\oss
By\dss the choice of\dss $z_{\dff i}$\sss the image
$\varphi\dff(\dff \sigma\dff)$ consists of\trs $n\qff +\qff 1$\dss points
and\dss is\dss affinely\dss independent\halfff.\oss
Passing\dss to\sss the\sss limit\sss shows\sss that\sss $z$\sss
is\dss contained\sss in\dss the convex\dss hull\sss of\trs
$\varphi\dff(\dff \sigma\dff)$\nnsp.\oss
The\dss theorem\dss follows.\oss  \eproof

\myuppar{Families of\qss linear orders.}
So,\oss let\qss $X$\dss be a non-empty\dss finite set\sss
and suppose\sss that\dss for
every\qss $i\qff \in\pff I$\qss a\dss linear\sss order\dss $<_{\dff i}$\dss
on\dss $X$\dss is\dss given.\oss
As in\qss Section\qss \ref{vector-colorings},\oss
we assume\sss that\dss $X\qff \cap\pff I\off =\off \varnothing$\dss
and extend\dss the orders\dss $<_{\dff i}$\dss to\dss
$X\qff \cup\pff I$\nnsp,\oss subject\dss to\sss the same conditions.\oss

\myuppar{The affine\sss form\sss of\qss Scarf\qss theorem.}
\emph{Let\qss $\varphi\dff \colon\dff X\qff \cup\pff I\qff \ttoo\qff \rrr^{\dff n}$\qss
be a map such\dss that\qss
$\varphi\dff(\dff i\trf)\off =\off w_{\fff i}$\qss for every\trs $i\qff \in\pff I$\trs
and\dss let\trs $z\qff \in\pff \Gamma$\dnsp.\oss
Then\dss there exists a subset\dss $\sigma\qff \subset\pff X\qff \cup\pff I$\dss
dominant\dss with respect\dss to\dss $I$\dss and
such\dss that\qss 
$z\qff \in\qff
\conv{\dff\varphi\dff(\trf \sigma\trf)\fff}$\dss
and\qss 
$\varphi\dff(\dff \sigma\dff)$\sss
is\dss an\sss $n$\dnsp-simplex\dss in\qss $\rrr^{\dff n}$\dnsp.\oss}

\proof
This\sss can\dss be deduced\dss from\qss Theorem\qss \ref{main-theorem-vector-new}\qss
in\dss exactly\dss the same way\sss as\trs Generalized\qss Scarf\qss Theorem\qss
was deduced\dss from\qss Theorem\qss \ref{main-theorem-general}.\oss  \eproof

\prooftitle{Another\dss proof\trs of\pss Scarf\qss theorem} 
See\dss Section\qss \ref{vector-colorings}.\oss
Let\qss $M\off =\off \varphi\dff(\qff X\qff \cup\pff I\qff)\qff +\qff b$\qss
and\dss let\sss $\mathbf{M}$\sss be\sss the convex\dss hull\sss of\trs $M$\dss 
in\dss the usual\dss geometric sense,\oss i.e.\qss 
the set\sss of\dss all\dss linear\dss combinations\oss\vspace{3pt} 
\[
\quad
\sum\nolimits_{\qff v\qff \in\qff M}\qff
x_{\fff v}\trf v
\]

\vspace{-9pt}
such\dss that\dss $x_{\dff v}\qff \geq\qff 0$\dss for all\trs $v\qff \in\pff M$\trs
and\oss 
$\sum\nolimits_{\qff v\qff \in\qff M}\qff x_{\dff v}\off =\off 1 $\nnsp.\oss
Lemma\qss \ref{vector-scarf-acyclic}\qss implies\sss that\sss $0$\sss
is\dss not\sss a\sss non-negative linear combination of\dss elements of\trs $M$\nnsp.\oss
In\dss particular\halfff,\pss $0$\dss is\dss not\dss contained\dss in\dss
the convex\dss hull\dss $\mathbf{M}$\nnsp.\oss
By\sss a classical\sss separation\dss theorem\dss there exists a\sss linear\sss functional\trs
$L\qff \colon\dff
\rrr^{\dff n\dff +\dff 1}\qff \ttoo\qff \rrr$\qss
such\dss that\trs
$L\dff(\dff v\trf)\qff >\pff 0$\qss 
for every\trs $v\qff \in\pff M$\nnsp.\oss
Let\qss
$r\dff(\dff v\trf)
\off =\off
v\dff/\qff L\dff(\dff v\trf)$\qss
for every\dss $v\qff \in\pff M$\nnsp.\oss
Then each\dss vector\sss $r\dff(\dff v\trf)$\sss belongs\sss to\sss the hyperplane\dss
$L^{\dff -\dff 1}\dff(\dff 1\dff)$\dss
and\qss (\ref{b-positive})\qss implies\sss that\dss
$r\dff(\dff b\trf)$\dss is\dss a\sss convex\dss combination of\trs
the vectors\dss $r\dff(\dff v_{\fff i}\trf)$\nnsp.\oss
Let\qss 
$f\qff \colon\qff
L^{\dff -\dff 1}\dff(\dff 1\dff)
\qff \ttoo\pff
\rrr^{\dff n}$\qss
be an affine isomorphism,\oss and\dss let\trs
$w_{\fff i}\off =\off f\dff \circ\dff r\trf(\dff v_{\fff i}\trf)$\nnsp,\pss
$z\off =\off  f\dff \circ\dff r\trf(\dff b\trf)$\nnsp.\oss
Then\dss the points $v_{\fff i}$ are affinely\sss in\-de\-pen\-dent\qss
and\sss $z$ is\dss a\sss convex\dss combination of\trs these points,\oss
i.e.\qss $z$\dss belongs\sss to\sss the $n$\dnsp-simplex\dss $\Gamma$\dss
with\dss the vertices\dss
$w_{\fff 0}\fff,\pff w_{\fff 1}\fff,\pff \ldots\fff,\pff w_{\fff n}$\nsp.\oss
Therefore\sss the affine form of\qss Scarf\qss theorem applies\sss to\sss 
the map\dss $f\dff \circ\trf \varphi$\dss in\dss the role of\trs $\varphi$
and\dss the point\sss $z$\nnsp.\oss
Clearly\halfff,\oss every\sss simplex $\sigma$ satisfying\dss
the conclusion of\trs this\sss affine form satisfies also\sss the conclusion
of\pss Scarf\qss theorem.\oss  \eproof

\myuppar{The simplex\dss $\Gamma$ and\dss related coordinates.}
The vertices\dss
$w_{\fff 0}\fff,\pff w_{\fff 1}\fff,\pff \ldots\fff,\pff w_{\fff n}$\dss of\qss $\Gamma$\dss
are affinely\dss independent\sss and
every\dss point\trs $x\qff \in\pff \rrr^{\dff n}$\dss
has a unique presentation\dss in\dss the form\vspace{3pt}
\[
\quad
x
\off =\off
a_{\trf 0}\dff(\dff x\trf)\dff w_{\dff 0}\qff +\qff
a_{\trf 1}\dff(\dff x\trf)\dff w_{\dff 1}\qff +\qff
\ldots\qff +\qff
a_{\dff n}\dff(\dff x\trf)\dff w_{\fff n}
\qff
\]

\vspace{-9pt}
where coefficients $a_{\dff i}\dff(\dff x\trf)$\nnsp,\oss
the\qss \emph{barycentric coordinates}\pss of\sss $x$\nnsp,\oss
are such\dss that\trs\vspace{3pt}
\[
\quad
a_{\trf 0}\dff(\dff x\trf)\qff +\qff
a_{\trf 1}\dff(\dff x\trf)\qff +\qff
\ldots\qff +\qff
a_{\dff n}\dff(\dff x\trf)
\off =\off
1
\qff.
\]

\vspace{-9pt}
Clearly\halfff,\pss
$a_{\dff i}\dff(\dff w_{\fff k}\trf)\off =\off 1$\qss
if\qss $i\off =\off k$\dss and\dss
$a_{\dff i}\dff(\dff w_{\fff k}\trf)\off =\off 0$\qss
if\qss $i\off \neq\off k$\nnsp.\oss
Also,\oss
the $(\fff n\dff -\dff 1\fff)$\dnsp-face\dss $\Gamma_{\dff i}$\dss of\qss $\Gamma$\sss
opposite\sss to $w_{\fff i}$ is\dss contained\sss in\dss the
hyperplane\dss
$H_{\dff i}$\dss defined\dss by\dss the equation\dss
$a_{\dff i}\dff(\dff x\trf)\off =\off 0$\nnsp.\oss
The half-space\dss
$\eta_{\dff i}\trf(\qff \geq\dff 0\dff)$\dss 
bounded\dss by\dss $H_{\dff i}$\dss and containing\sss
$w_{\fff i}$\sss is\dss defined\dss by\dss the inequality\trs
$a_{\dff i}\dff(\dff x\trf)\qff \geq\qff 0$\dss
and\dss the other\dss half-space\dss
$\eta_{\dff i}\trf(\qff \leq\dff 0\dff)$\dss 
is\dss defined\dss by\dss
$a_{\dff i}\dff(\dff x\trf)\qff \leq\qff 0$\nnsp.\oss
Let\trs\vspace{3pt}
\[
\quad
l_{\dff i}\dff(\dff x\trf)
\off =\off 
a_{\dff i}\dff(\dff x\trf)
\qff -\pff
a_{\dff i}\dff(\dff 0\dff)
\]

\vspace{-9pt}
for every\sss $i\qff \in\pff I$\nnsp.\oss
Then\dss
$l_{\dff i}\dff \colon\dff
\rrr^{\dff n}\qff \ttoo\qff \rrr$\dss
is\dss a\sss linear\dss map\sss 
and\dss
$l_{\dff i}\dff(\dff y\trf)\pff \leq\pff l_{\dff i}\dff(\dff x\trf)$\dss
if\qss and\dss only\qss if\qss
$a_{\dff i}\dff(\dff y\trf)\pff \leq\pff a_{\dff i}\dff(\dff x\trf)$\nnsp.\oss

For\dss the rest\sss of\trs this section we will\sss assume 
that\sss $0$\sss is\dss the barycenter of\trs $\Gamma$\dnsp,\oss
or\halfff,\oss equivalently\halfff,\oss
that\trs
$w_{\fff 0}\qff +\qff w_{\fff 1}\qff +\qff \ldots\qff +\qff w_{\fff n}
\off =\off 0$\nnsp.\oss
Then\qss
$a_{\dff i}\dff(\dff 0\dff)
\off =\off
1\dff/(\dff n\qff +\qff 1\dff)$\dss and\dss hence\vspace{3pt}
\begin{equation}
\label{l-values}
\quad
l_{\dff i}\dff(\dff w_{\fff i}\trf)
\off =\off 
\frac{n}{n\qff +\qff 1}
\hspace*{1.5em}\mbox{and}\hspace*{1.5em}
l_{\dff i}\dff(\dff w_{\fff k}\trf)
\off =\off 
-\qff 
\frac{1}{n\qff +\qff 1}
\end{equation}

\vspace{-9pt}
for every\sss $i$\sss and every\dss $k\off \neq\off i$\nnsp.\oss\vspace{-1.75pt}

\myuppar{Triangulations and\dss vector\sss colorings.}
The\sss
simplex-families arising\sss from\sss 
triangulations of\trs $\Gamma$\dss 
are\sss pseudo-simplices and\dss hence are chain-simplices,\pss
as we saw\sss in\qss Section\qss \ref{pseudo-simplices}.\oss
We\sss leave\sss the\sss task of\dss stating\dss the corresponding\sss special\sss case of\qss
Theorem\qss \ref{main-theorem-vector-new}\qss to\sss the reader\halfff.\oss
We will\sss discuss in\sss details\sss two special\sss
cases dealing\dss with\qss \emph{vector\dss hedgehog colorings}\qss 
and\dss with\qss \emph{inward\dss tangent\sss colorings},\oss
to be defined\dss in\sss a\sss moment\halfff.\oss
The colors of\qss the inward\dss tangent\sss colorings\sss should\dss 
be interpreted as\sss tangent\dss
vectors\sss to a simplex.\oss
Cf.\qss remarks at\dss the end of\qss Section\qss \ref{kakutani}.\oss

Suppose\dss that\trs $T$\dss is\dss a\sss triangulation of\qss $\Gamma$ 
and\dss let\trs $X$\dss be\sss the set\sss of\dss vertices of\trs $T$\dnsp.\oss
Let\dss $\mathcal{D}_{\qff T}$\dss be\sss the simplex-family\sss associated\dss to\dss $T$\dnsp.\oss
Then\dss $\mathcal{D}_{\qff T}$\dss is\dss a\sss pseudo-simplex and\dss hence\dss 
is\dss a\sss chain-simplex.\oss
A map\dss
$c\dff \colon\dff
X\qff \ttoo\pff \rrr^{\dff n}$\dss is\dss called\sss a\qss
\emph{vector\sss hedgehog coloring}\oss if\vspace{3pt}
\[
\quad
c\dff(\dff x\trf)
\qff \in\pff 
\eta_{\dff i}\trf(\qff \leq\dff 0\dff)
\]

\vspace{-9pt}
for every\trs $i\qff \in\pff I$\dss and\trs
$x\qff \in\pff X\pff \cap\pff \Gamma_{\dff i}$\nsp,\oss
and an\qss \emph{inward\dss tangent\sss coloring}\oss
if\qss\vspace{2pt}
\begin{equation*}
\quad
x\qff +\qff c\dff(\dff x\trf)
\qff \in\pff 
\eta_{\dff i}\trf(\qff \geq\dff 0\dff)
\end{equation*}

\vspace{-9pt}
for every\trs $i\qff \in\pff I$\dss and\trs
$x\qff \in\pff X\pff \cap\pff \Gamma_{\dff i}$\nsp.\oss
Since\dss $a_{\dff i}\dff(\dff x\trf)\off =\off 0$\dss for\sss such $x$\nnsp,\oss
the last\sss condition\dss 
holds\dss if\qss and\dss only\qss if\qss
$a_{\dff i}\dff(\dff x\qff +\qff c\dff(\dff x\trf)\trf)
\pff \geq\pff
a_{\dff i}\dff(\dff x\trf)$\nnsp,\oss
or\halfff,\oss equivalently\halfff,\oss
$l_{\dff i}\dff(\dff x\qff +\qff c\dff(\dff x\trf)\trf)
\pff \geq\pff
l_{\dff i}\dff(\dff x\trf)$\dss
for every\trs $i\qff \in\pff I$\dss and\trs
$x\qff \in\pff X\pff \cap\pff \delta_{\dff i}$\nsp.\oss
Since $l$\dss is\dss linear\halfff,\oss
this condition\dss is\dss equivalent\dss to\dss the condition\vspace{3pt}
\[
\quad
l_{\dff i}\trf(\trf c\dff(\dff x\trf)\trf)\off \geq\dff\off 0
\]

\vspace{-9pt}
for every\trs $i\qff \in\pff I$\trs 
and\trs 
$x\qff \in\pff X\pff \cap\pff \Gamma_{\dff i}$\nsp.\oss

\mypar{The main\dss theorem\sss for\dss vector\dss hedgehog\sss colorings.}{main-theorem-vector-hedgehog}
\emph{For\dss every vector\dss hedgehog\sss coloring\qss
$c\dff \colon\dff
X\qff \ttoo\pff \rrr^{\dff n}$\qss
and\dss every\dss $z\qff \in\pff \Gamma$\dss
there exists an $n$\dnsp-simplex\dss of\pss $T$\dss
such\dss that\trs if\qss $\sigma$\dss is\dss its set\dss of\qss vertices,\oss
then\dss $z$\dss is\dss contained\dss in\dss the convex\dss hull\sss of\qss
$c\dff(\dff \sigma\dff)$\nnsp.\oss}

\proof
Suppose first\dss that\sss $z$ is\dss contained\sss in\dss the interior of\trs $\Gamma$\dnsp.\oss
Let\dss us\sss apply\qss Theorem\qss \ref{main-theorem-vector-new}\qss
to\dss $\mathcal{D}_{\qff T}$\nsp,\oss
the simplex\trs $\Gamma$\nnsp,\oss
the point\sss $z\qff \in\pff \Gamma$\nnsp,\oss
and\dss the points\dss
$w_{\fff 1}\fff,\pff v_{\fff 2}\fff,\pff \ldots\fff,\pff v_{\fff n}\fff,\pff v_{\fff 0}$\dss
in\dss the role of\trs
$v_{\fff 0}\fff,\pff v_{\fff 1}\fff,\pff \ldots\fff,\pff v_{\fff n\dff -\dff 1}\fff,\pff v_{\fff n}$\dss
respectively\halfff.\oss
By\dss this\dss theorem\dss there exists a non-empty\sss subset\trs
$C\qff \subset\pff I$\trs
and\dss a\sss $d\dff(\trf C\trf)$\dnsp-simplex\dss $\sigma$\dss of\pss
$\mathcal{D}_{\qff T}\dff(\trf C\trf)$\dss such\dss that\dss 
the set\qss\vspace{3pt}\vspace{0.125pt}
\begin{equation*}
\quad
Y
\off\off =\off\off
c\qff(\trf \sigma\trf)\off \cup\off\dff \bigl\{\qff v_{\fff i\dff +\dff 1}\off \bigl|\off
i\qff \in\qff I\qff \smallsetminus\pff C \pff\bigr\}
\end{equation*}

\vspace{-9pt}\vspace{0.125pt}
consists of\qss $n\qff +\qff 1$\dss elements,\oss
is\dss an\sss affinely\dss independent\halfff,\oss
and contains\sss $z$\sss in\dss its convex\dss hull.\oss 
If\qss $C\off \neq\dff\off I$\nnsp,\oss
then\dss there exists\dss $k\qff \in\pff I$\dss such\dss that\trs
$k\qff \not\in\qff C$\dss and\dss $k\qff -\qff 1\qff \in\qff C$\nnsp.\oss
By\sss arguing exactly\sss as in\dss the proof\dss of\qss 
Theorem\qss \ref{main-theorem-hedgehog}\qss we conclude\sss that\trs
$Y\qff \subset\pff \eta_{\dff k}\trf(\qff \leq\dff 0\dff)$\dss
and\dss hence\dss
$z\qff \in\qff \eta_{\dff k}\trf(\qff \leq\dff 0\dff)$\nnsp.\oss
This contradicts\sss to\sss the assumption\dss that\sss $z$ is\dss
contained\sss in\dss the interior of\trs $\Gamma$\dnsp.\oss
It\dss follows\dss that\trs $C\off =\off I$\nnsp.\oss
In\dss turn,\oss this implies\sss that\sss 
$\sigma$ is\dss the set\sss of\dss
vertices of\dss an $n$\dnsp-simplex of\trs $T$\dss
and\dss $c\qff(\trf \sigma\trf)$\dss contains $z$ in\dss its convex\dss hull.\oss
This proves\sss the\sss theorem\sss in\dss the case when\dss
$z$ is\dss contained\sss in\dss the interior of\trs $\Gamma$\dnsp.\oss
But\dss the union of\trs the convex\dss hulls of\trs the sets\dss
$c\dff(\dff \sigma\dff)$\dss
with $\sigma$\sss being\dss the set\sss of\trs vertices of\dss an $n$\dnsp-simplex\sss of\trs
$T$\dss is\dss obviously\sss closed.\oss
Therefore\sss this special\sss case implies\sss the\sss theorem\dss
for arbitrary\dss $z\qff \in\pff \Gamma$\dnsp.\oss  \eproof

\mypar{The main\dss theorem\sss for\dss inward\dss tangent\sss colorings.}{main-theorem-inward-tangent}
\emph{For\dss every\dss inward\dss tangent\sss coloring\qss
$c\dff \colon\dff
X\qff \ttoo\pff \rrr^{\dff n}$\qss
there exists an $n$\dnsp-simplex\dss of\pss $T$\dss
such\dss that\trs if\qss $\sigma$\dss is\dss its set\dss of\qss vertices,\oss
then\dss $0$\dss is\dss contained\dss in\dss the convex\dss hull\sss of\qss
$c\dff(\dff \sigma\dff)$\nnsp.\oss}

\proof
Let\dss us\sss apply\qss Theorem\qss \ref{main-theorem-vector-new}\qss
to\dss $\mathcal{D}_{\qff T}$\nsp,\oss
the simplex\trs $\Gamma$\nnsp,\oss
the points\dss
$w_{\fff 0}\fff,\pff w_{\fff 1}\fff,\pff \ldots\fff,\pff w_{\fff n}$\nsp,\oss
and\dss $z\off =\off 0$\nnsp.\oss
By\dss this\dss theorem\dss there exists a non-empty\sss subset\trs
$C\qff \subset\pff I$\trs
and\dss a\sss $d\dff(\trf C\trf)$\dnsp-simplex\dss $\tau$\dss of\pss
$\mathcal{D}_{\qff T}\dff(\trf C\trf)$\dss such\dss that\dss 
the set\qss (\ref{affine-colorings-simplex})\qss
contains $0$ in\dss its convex\dss hull.\oss
Then\vspace{3pt}
\begin{equation}
\label{convex-hull}
\quad
\sum\nolimits_{\qff x\qff \in\qff \tau}\qff a_{\dff x}\qff c\dff(\dff x\trf)
\off +\off
\sum\nolimits_{\qff k\qff \in\qff I\pff \smallsetminus\qff C}\qff a_{\trf k}\qff w_{\fff k}
\off =\off
0
\end{equation}

\vspace{-9pt}
for some non-negative coefficients\dss $a_{\dff x}\dff,\pff a_{\trf i}$\dss
such\dss that\oss\vspace{3.5pt} 
\[
\quad
\sum\nolimits_{\qff x\qff \in\qff \tau}\qff a_{\dff x}
\off +\off
\sum\nolimits_{\qff k\qff \in\qff I\pff \smallsetminus\qff C}\qff a_{\trf k}
\off =\off
1
\qff.
\]

\vspace{-8.5pt}
Suppose\sss that\trs 
$a_{\trf i}\off =\off 0$\dss
for every\trs $i\qff \in\pff I\pff \smallsetminus\qff C$\qss
({\fff}this happens,\oss in\dss particular\halfff,\oss
when\qss $C\off =\off I$\nsp).\oss
In\dss this case $0$\sss belongs\sss to\sss the convex\dss hull\sss of\trs
$c\dff(\dff \tau\dff)$\nnsp.\oss
The simplex $\tau$ is\dss a\sss simplex of\trs the abstract\sss
simplicial\sss complex associated\dss with\dss $T$\dss and\dss hence\dss
is\dss a\sss face of\dss an $n$\dnsp-simplex $\sigma$ of\trs that\sss complex.\oss
It\dss follows\sss that\trs $c\dff(\dff \tau\dff)\qff \subset\qff c\dff(\dff \sigma\dff)$\dss
and\dss hence $0$\sss belongs\sss to\sss the convex\dss hull\sss
of\dss
$c\dff(\dff \sigma\dff)$\nnsp.\oss
This proves\sss the\sss theorem\sss in\dss the case when\dss 
$a_{\trf i}\off =\off 0$\dss
for every\trs $i\qff \in\pff I\pff \smallsetminus\qff C$\nnsp.\oss

Suppose now\dss that\dss
$a_{\trf i}\off \neq\off 0$\dss
for some\trs $i\qff \in\pff I\pff \smallsetminus\qff C$\nnsp.\oss
Let\dss us\dss choose some\trs $i\qff \in\pff I\pff \smallsetminus\qff C$\dss such\dss that\trs
$a_{\dff i}\qff \geq\qff a_{\dff k}$\dss for all\trs
$k\qff \in\pff I\pff \smallsetminus\qff C$\nnsp.\oss
Then\dss $a_{\dff i}\qff >\qff 0$\nnsp.\oss
By applying\dss $l_{\dff i}$\sss to\qss  
(\ref{convex-hull})\qss we conclude\sss that\vspace{3.5pt}
\begin{equation}
\label{convex-hull-l}
\quad
\sum\nolimits_{\qff x\qff \in\qff \tau}\qff 
a_{\dff x}\qff l_{\dff i}\trf(\trf c\dff(\dff x\trf)\trf)
\off +\off
\sum\nolimits_{\qff k\qff \in\qff I\pff \smallsetminus\qff C}\qff 
a_{\trf k}\qff l_{\dff i}\trf(\trf w_{\fff k}\trf)
\off =\off
0
\qff.
\end{equation}

\vspace{-8.5pt}
Since $c$ is\dss an\sss inward\dss tangent\sss coloring\halfff,\oss
the first\sss sum at\dss the left\dss hand side\dss is\dss non-negative.\oss

Let\trs $K\off =\off (\qff I\pff \smallsetminus\qff C \trf)\qff -\qff i$\nnsp.\oss
Then\dss the second sum\dss is\dss equal\dss to\vspace{3pt}
\[
\quad
a_{\trf i}\qff l_{\dff i}\trf(\trf w_{\fff i}\trf)
\off +\off
\sum\nolimits_{\qff k\qff \in\qff K}\qff 
a_{\trf k}\qff l_{\dff i}\trf(\trf w_{\fff k}\trf)
\qff.
\]

\vspace{-7.5pt}
In view of\pss (\ref{l-values})\qss the last\sss expression\dss is\dss equal\dss to\vspace{3pt}
\[
\quad
a_{\trf i}\qff \frac{n}{n\qff +\qff 1}
\off -\off
\sum\nolimits_{\qff k\qff \in\qff K}\qff 
a_{\trf k}\qff \frac{1}{n\qff +\qff 1}
\off \geq\off
a_{\trf i}\qff \frac{n}{n\qff +\qff 1}
\off -\off
\sum\nolimits_{\qff k\qff \in\qff K}\qff 
a_{\trf i}\qff \frac{1}{n\qff +\qff 1}
\]

\vspace{-26.25pt}
\[
\quad
\phantom{a_{\trf i}\qff \frac{n}{n\qff +\qff 1}
\off -\off
\sum\nolimits_{\qff k\qff \in\qff K}\qff 
a_{\trf k}\qff \frac{1}{n\qff +\qff 1}
\off }
=\off
a_{\trf i}\qff \frac{n}{n\qff +\qff 1}
\off -\off
a_{\trf i}\qff \frac{\num{K}}{n\qff +\qff 1}
\off =\off
a_{\trf i}\qff \frac{n\qff -\qff \num{K}}{n\qff +\qff 1}
\qff.
\]

\vspace{-7.5pt}
Since\sss $C$\sss is\dss non-empty\halfff,\oss
$\num{\dff I\pff \smallsetminus\qff C\halfff}\qff \leq\qff n$\nnsp.\oss
Since\dss
$i\qff \in\pff I\pff \smallsetminus\qff C$\nnsp,\oss
this implies\sss that\trs
$\num{K}\qff \leq\qff n\qff -\qff 1$\dss
and\dss hence\dss
$n\qff -\qff \num{K}\qff >\qff 0$\nnsp.\oss
It\dss follows\sss that\dss the second sum\dss in\qss (\ref{convex-hull-l})\qss
is\dss $>\qff 0$\nnsp.\oss
Since\sss the first\sss sum\dss is\dss $\qff \geq\qff 0$\dnsp,\pss
the left\dss hand side of\qss (\ref{convex-hull-l})\qss
is\dss $>\qff 0$\nnsp,\oss
contrary\dss to\qss (\ref{convex-hull-l}).\oss
The contradiction shows\sss that\trs 
$a_{\trf i}\off =\off 0$\dss
for every\trs $i\qff \in\pff I\pff \smallsetminus\qff C$\nnsp.\oss
But\dss we already\sss saw\sss that\dss the\sss theorem\dss holds
in\dss this case.\oss   \eproof

\prooftitle{Another\dss proof\trs of\pss Kakutani\qss theorem}
Let\trs $\Gamma$\dss be as above,\oss
and\dss let\qss
$F\dff \colon\dff 
\Gamma\qff \ttoo\qff \Gamma$\qss
be a closed\dss multi-valued\dss map.\oss
Let\dss us\dss choose a\sss map\qss
$f\dff \colon\dff 
\Gamma\qff \ttoo\qff \Gamma$\qss
(not\dss assumed\dss to be continuous)\qss
such\dss that\trs $f\dff(\dff x\trf)\qff \in\pff F\dff(\dff x\trf)$\dss
for every\trs $x\qff \in\pff \Gamma$\dnsp.\oss
Let\dss $c\trf(\dff x\trf)\off =\off f\dff(\dff x\trf)\qff -\pff x$\nnsp.\oss

Let\qss
$\varepsilon_{\dff 1}\dff,\off
\varepsilon_{\dff 2}\dff,\off
\varepsilon_{\dff 3}\dff,\off
\ldots$\qss
be\dss positive numbers such\dss that\trs
$\varepsilon_{\dff k}\qff \ttoo\qff 0$\dss when $k\ttoo \infty$\nnsp.\oss
By\sss considering\dss triangulations of\trs $\Gamma$\dss into simplices of\dss
diameter\qss $<\qff \varepsilon_{\dff k}$\qss and applying\qss
Theorem\qss \ref{main-theorem-inward-tangent},\oss
we see\sss that\dss there exist\sss subsets\dss
$\sigma_{k}\qff \subset\pff \Gamma$\dss such\dss that\trs
$\num{\sigma_k}\off =\off n\qff +\qff 1$\nnsp,\oss
the diameter of\dss $\sigma_{k}$\dss is\dss $<\qff \varepsilon_{\dff k}$\dss
and\dss the image\dss $c\trf(\dff \sigma_{k}\trf)$\dss contains $0$ in\dss its\sss
convex\dss hull.\oss
Let\trs 
$w_{\dff 0}\dff(\dff k\trf)\fff,\off
w_{\dff 1}\dff(\dff k\trf)\fff,\off
\ldots\fff,\off
w_{\fff n}\dff(\dff k\trf)$\dss
be\sss the elements\sss of\dss $\sigma_k$\nsp.\oss
After\sss passing\dss to a subsequence we can assume\sss
that\dss every\sss sequence\dss
$w_{\fff i}\dff(\dff 1\trf)\dff,\off
w_{\fff i}\dff(\dff 2\trf)\dff,\off
w_{\fff i}\dff(\dff 3\trf)\dff,\off
\ldots$\dss
converges\sss to a\dss limit\sss 
when\dss $k \ttoo \infty$\nnsp.\oss 
Since\sss the diameters of\dss simplices $\sigma_k$\sss tend\dss to $0$\nnsp,\oss
these\sss limits are,\oss in\dss fact\halfff,\oss equal.\oss
Let\sss $u$\sss be\sss the common\dss values of\trs these\sss limits.\oss 
After\dss passing\dss to a\sss further subsequence we can assume\dss
that\sss each sequence\vspace{3pt}
\[
\quad
f\trf (\dff w_{\fff i}\dff(\dff 1\trf)\trf)\dff,\off\off
f\trf (\dff w_{\fff i}\dff(\dff 2\trf)\trf)\dff,\off\off
f\trf (\dff w_{\fff i}\dff(\dff 3\trf)\trf)\dff,\off\off
\ldots
\]

\vspace{-9pt}
converges\sss to a\sss limit\dss which we will\sss denote\sss by\dss $w_{\dff i}$\nsp.\oss
It\dss follows\dss that\sss each sequence\vspace{3pt}
\[
\quad
c\trf (\dff w_{\fff i}\dff(\dff 1\trf)\trf)\dff,\off\off
c\trf (\dff w_{\fff i}\dff(\dff 2\trf)\trf)\dff,\off\off
c\trf (\dff w_{\fff i}\dff(\dff 3\trf)\trf)\dff,\off\off
\ldots
\]

\vspace{-9pt}
converges\sss to\sss the limit\trs 
$c_{\dff i}
\off =\off
w_{\dff i}\qff -\qff u$\nnsp.\oss
Since every\sss set\dss $c\trf(\dff \sigma_{k}\trf)$\dss contains $0$ in\dss its\sss
convex\dss hull,\oss
the set\trs
$\{\qff c_{\dff 0}\dff,\pff  c_{\dff 1}\dff,\pff \ldots\dff,\pff  c_{\dff n}\pff\}$\trs
also contains $0$ in\dss its convex\dss hull.\oss
Cleary\halfff,\oss every\sss convex combination of\trs the points $c_{\dff i}$\dss
has\sss the form\dss $w\qff -\qff u$\nnsp,\oss
where $w$ is\dss a\sss convex combination of\trs the points $w_{\dff i}$\nsp.\oss
It\dss follows\dss that\sss $u$ is\dss a convex combination of\trs the points $w_{\dff i}$\nsp.\oss
Since\sss $F$\sss is\dss a\sss closed\dss multi-valued\dss map,\oss
the choice of\trs the map $f$ implies\sss that\trs
$w_{\dff i}\qff \in\pff F\dff(\dff u\trf)$\dss for every\trs $i\qff \in\pff I$\nnsp.\oss
Since\dss $F\dff(\dff u\trf)$ is\dss convex,\oss this implies\sss that\trs
$u\qff \in\pff F\dff(\dff u\trf)$\nnsp,\oss
i.e.\qss that\sss $u$ is\dss a\sss fixed\dss point\sss of\trs $F$\nnsp.\oss \eproof

\myuppar{Another version of\qss Scarf\qss theorem.}
The\sss theorem called\qss
\emph{Scarf\qss theorem}\qss above\trs is\qss Theorem\qss 4.2.3\qss from\qss \cite{sc3}\qss
written\sss in\sss a\sss somewhat\sss different\dss language.\oss
In\qss Theorem\qss 7.1\qss of\pss  
\cite{sc4}\pss Scarf\qss
replaced\dss the assumption\dss that\dss
the set\dss of\dss non-negative solutions 
of\pss (\ref{boundedness})\qss is\dss bounded\qss
(which\dss is\dss a\sss part\dss of\dss our definition of\dss
a vector framework)\qss by\dss the following\halfff.\oss
Suppose\sss that\vspace{4.5pt}
\[
\quad
x
\off =\off
\sum\nolimits_{\qff v\qff \in\pff M\qff -\qff b}\qff
y_{\fff v}\trf v
\qff,
\]

\vspace{-7.5pt}
where\dss
$y_{\fff v}\qff \geq\qff 0$\dss for all\qss $v\qff \in\pff M\qff -\qff b$\nnsp.\oss
Then,\oss if\qss all\sss coordinates of\trs the vector $x$ are non-positive,\oss
all\sss coefficients $y_{\fff v}$ should\dss be equal\dss to $0$\nnsp.\oss
This assumption\dss immediately\sss implies\sss that\dss $0$\dss is\dss 
not\sss a\sss non-negative linear combination of\dss elements\sss of\dss $M$\nnsp.\oss
Therefore\sss the above proof\dss works under\dss this assumption also,\oss
and,\pss moreover\halfff,\pss there\dss is\dss no need\dss to use\qss
Lemma\qss \ref{vector-scarf-acyclic}.\oss

\myuppar{Kannai's\qss Generalized\trs Sperner\qss lemma.}
Theorem\dss \ref{main-theorem-vector-hedgehog}\qss generalizes\qss
\emph{Generalized\trs Sperner\qss lem\-ma}\qss of\pss Kannai\qss \cite{kan},\oss 
who assumed\dss that\dss $c\dff(\dff x\trf)$\dss belongs\sss
to\sss the affine hyperplane spanned\dss by\dss the face\dss $\Gamma_{\dff i}$\dss
for every\trs $i\qff \in\pff I$\dss and\trs
$x\qff \in\pff X\pff \cap\pff \Gamma_{\dff i}$\nsp,\oss
in contrast\dss with our\dss weaker assumption\dss
$c\dff(\dff x\trf)
\qff \in\pff 
\eta_{\dff i}\trf(\qff \leq\dff 0\dff)$\nnsp.\oss
The\sss assumption\dss used\dss by\trs Kannai\qss is\dss forced\dss
by\dss his\sss method of\dss proof\dss 
based on an induction\sss by\sss $n$ similar\dss to\trs
Sperner's\trs arguments.\oss
In contrast\dss with\qss Kannai\qss \cite{kan},\oss
our\sss proof\dss does not\dss require any\dss form of\trs
the simplicial\sss approximation\dss theorem.\oss
Using\dss the\sss latter\dss is\dss often considered\dss as a\sss topological,\pss
and\dss hence non-elementary\sss
and\dss non-combinatorial,\pss argument\halfff.\oss
See\trs Kannai\qss \cite{kan},\pss Section\qss 4\qss and\qss Ziegler\qss \cite{z},\oss
Introduction,\oss for example.\oss

\newpage
\myappend{Todd's\qss theorem}{todd}

\myuppar{Todd\halfff's\qss theorem.}
\emph{Let\qss
$\sigma\fff,\pff \tau$\qss
be\sss circuits\sss of\qss an oriented\dss matroid\dss
and\dss let\qss
$w\pff \in\off \underline{\tau}\pff \smallsetminus\off \underline{\sigma}$\nsp.\oss
Suppose\sss that\dss there exists\qss 
$e\pff \in\off \underline{\sigma}\off \cap\off \underline{\tau}$\qss
such\dss that\qss
$\sigma\dff(\trf e\qff)\off =\off -\qff \tau\dff(\trf e\qff)$\nnsp.\oss
Then\dss there\dss is\dss a\sss circuit\sss $\omega$\sss
such\dss that}\vspace{3pt}
\begin{equation}
\label{todd-inclusion-1}
\quad
\omega_{\dff +}
\off \subset\off
(\dff \sigma_{\dff +}\qff \cup\off \tau_{\dff +}\trf)
\off \smallsetminus\off
\sigma_{\dff -}\off,
\end{equation}

\vspace{-36pt}
\begin{equation}
\label{todd-inclusion-2}
\quad
\omega_{\dff -}
\off \subset\off
(\dff \sigma_{\dff -}\qff \cup\off \tau_{\dff -}\trf)
\off \smallsetminus\off
\sigma_{\dff +}
\off,
\end{equation}

\vspace{-36pt}
\begin{equation}
\label{todd-inclusion-3}
\quad
w\pff \in\off \underline{\omega}\off,
\hspace{1.2em}\mbox{\emph{and}}\hspace{1.5em}
\omega\dff(\trf w\qff) 
\off\qff =\off\qff
\tau\dff(\trf w\qff)
\qff.
\end{equation}

\vspace{-9pt}
\emph{Moreover\halfff,\oss
if\qff\oss
$\underline{\tau}
\dff\off \subset\qff\off 
\underline{\sigma}\off \cup\qff \{\trf w\qff\}$\nnsp,\qff\oss
then\dss such\sss a\sss circuit\dss $\omega$\dss is\dss unique\sss and}\qss 
$\underline{\sigma}\off \smallsetminus\qff\off \underline{\tau}
\off\off \subset\pff\off
\underline{\omega}$\nnsp.\oss

\proof
Let\qss
$d\trf(\trf \sigma\fff,\pff \tau\trf)
\off =\off
\num{(\trf \sigma_{\dff +}\qff \cap\qff \tau_{\dff -}\trf)
\pff \cup\pff
(\trf \sigma_{\dff -}\qff \cap\qff \tau_{\dff +}\trf)}$\nnsp.\oss
Let\dss us\dss prove\sss the first\sss statement\dss 
using\sss an\dss induction\dss by\dss
$d\trf(\trf \sigma\fff,\pff \tau\trf)$\nnsp.\oss
By\dss the assumptions of\trs the\sss theorem,\pss
$d\trf(\trf \sigma\fff,\pff \tau\trf)\qff \geq\qff 1$\nnsp.\oss
Without\sss any\dss loss of\trs generality\dss we can assume\sss
that\qss
$e\qff \in\qff \sigma_{\dff +}\qff \cap\qff \tau_{\dff -}$\qss
and\dss hence\qss
$\num{\sigma_{\dff +}\qff \cap\qff \tau_{\dff -}}\off \geq\off 1$\nnsp.\oss

Suppose\sss that\qss
$d\trf(\trf \sigma\fff,\pff \tau\trf)\off =\off 1$\nnsp.\oss
In\dss this case\qss
$\num{\sigma_{\dff +}\qff \cap\qff \tau_{\dff -}}\off =\off 1$\qss
and\qss
$\sigma_{\dff -}\qff \cap\qff \tau_{\dff +}\off =\off \varnothing$\nnsp.\oss
It\dss follows\dss that\qss 
$e\qff \in\qff 
\sigma_{\dff +}\qff \cap\qff \tau_{\dff -}$\nsp.\oss
Let\dss us\sss apply\dss the strong elimination\sss property\sss
as stated\dss in\qss Theorem\qss \ref{vector-elimination}\qss
to\qss $u\off =\off e$\qss 
and\qss $v\off =\off w$\nnsp.\oss
We see\sss that\dss
there exists a circuit\sss $\omega$ such\dss that\qss
$w\pff \in\off \underline{\omega}$\nsp,\vspace{3pt}
\[
\quad
\omega_{\dff +}
\off \subset\off\dff 
(\qff \sigma_{\dff +}\qff \cup\qff \tau_{\dff +}\qff)
\qff \smallsetminus\qff 
\{\trf e\qff\}
\qff,\hspace{1.2em}\mbox{and}\hspace{1.2em}
\omega_{\dff -}
\off \subset\off\dff 
(\qff \sigma_{\dff -}\qff \cup\qff \tau_{\dff -}\qff)
\qff \smallsetminus\qff 
\{\trf e\qff\}
\qff.
\]

\vspace{-9pt}
Since\sss
$\sigma_{\dff -}$\sss is\dss disjoint\dss from\sss $\tau_{\dff +}$\sss
and\sss $\sigma_{\dff +}$\nsp,\vspace{3pt}
\[
\quad
\omega_{\dff +}
\off \subset\off\dff 
\sigma_{\dff +}\qff \cup\qff \tau_{\dff +}
\off \subset\off
(\dff \sigma_{\dff +}\qff \cup\off \tau_{\dff +}\trf)
\off \smallsetminus\off
\sigma_{\dff -}
\off.
\]

\vspace{-9pt}
On\dss the other\dss hand,\pss $\sigma_{\dff +}$\sss intersects\sss $\tau_{\dff -}$\sss
only\dss by\sss $e$\sss and\dss is\dss disjoint\dss from\sss $\sigma_{\dff -}$\nsp.\oss
Therefore\vspace{3pt}
\[
\quad
\omega_{\dff -}
\off \subset\off\dff 
(\qff 
\sigma_{\dff -}\qff \cup\qff \tau_{\dff -}
\qff)
\qff \smallsetminus\qff 
\{\trf e\qff\}
\off \subset\off
(\dff \sigma_{\dff -}\qff \cup\off \tau_{\dff -}\trf)
\off \smallsetminus\off
\sigma_{\dff +}
\off.
\]

\vspace{-9pt}
It\dss follows\dss that\qss the properties\qss 
(\ref{todd-inclusion-1})\qss and\qss (\ref{todd-inclusion-2})\qss hold.\oss
Since$w\pff \in\off \underline{\omega}$\nsp,\oss
these properties imply\dss that\qss
$\omega\dff(\trf w\qff) 
\off\qff =\off\qff
\tau\dff(\trf w\qff)$\nnsp.\oss
This completes\sss the proof\dss of\trs the first\sss statement\dss for\qss
$d\trf(\trf \sigma\fff,\pff \tau\trf)\off =\off 1$\nnsp.\oss

Suppose\sss now\dss that\dss
$k\qff >\qff 1$\qss
and\dss that\dss the first\dss statement\dss is\dss proved\dss for\qss
$d\trf(\trf \sigma\fff,\pff \tau\trf)\off <\off k$\nnsp.\oss
Suppose\sss that\qss
$d\trf(\trf \sigma\fff,\pff \tau\trf)\off =\off k$\qss
and\sss apply\dss the strong elimination\dss property\sss as in\dss the
previous paragraph.\oss
If\trs the properties\qss (\ref{todd-inclusion-1})\qss and\qss (\ref{todd-inclusion-2})\qss hold,\oss
then\dss the property\qss (\ref{todd-inclusion-3})\qss also holds and\dss we are done.\oss
Otherwise at\dss least\sss one of\trs the sets\qss
$\sigma_{\dff +}\qff \cap\pff \omega_{\dff -}$\qss
and\qss
$\sigma_{\dff -}\qff \cap\pff \omega_{\dff +}$\qss
is\dss non-empty\halfff.\oss
Let\dss $f$\dss be an element\sss either of\trs these sets.\oss
The inclusions\qss (\ref{elimination-1})\qss imply\dss that\vspace{3pt}
\[
\quad
\sigma_{\dff +}\qff \cup\off \omega_{\dff +}
\off \subset\off
\sigma_{\dff +}\qff \cup\off \tau_{\dff +}
\hspace{1.2em}\mbox{and}\hspace{1.2em}
\sigma_{\dff -}\qff \cup\off \omega_{\dff -}
\off \subset\off
\sigma_{\dff -}\qff \cup\off \tau_{\dff -}
\off.
\]

\vspace{-9pt}
Therefore,\oss if\qss the first\sss statement\dss holds for\qss
$\sigma\fff,\off \omega$\nnsp,\pss and\dss $f$\dss
in\dss the role of\qss
$\sigma\fff,\off \tau$\nnsp,\pss and\dss $e$\dss re\-spec\-tive\-ly\halfff,\oss
then\dss it\dss holds for\qss
$\sigma\fff,\off \tau$\nnsp,\pss and\dss $e$\sss also.\oss
But\vspace{3pt}\vspace{-0.25pt}
\[
\quad
\sigma_{\dff +}\qff \cap\qff \omega_{\dff -}
\off \subset\off
(\qff
\sigma_{\dff +}\qff \cap\qff \tau_{\dff -}
\qff)
\qff \smallsetminus\qff 
\{\trf e\qff\}
\]

\vspace{-9pt}\vspace{-0.25pt}
and\dss hence\qss
$\num{\sigma_{\dff +}\qff \cap\qff \omega_{\dff -}}
\off \leq\off
\num{\sigma_{\dff +}\qff \cap\qff \tau_{\dff -}}
\qff -\qff
1$\nnsp.\oss
On\dss the other hand,\pss\vspace{3pt}\vspace{-0.25pt}
\[
\quad
\sigma_{\dff -}\qff \cap\qff \omega_{\dff +}
\off \subset\off
\sigma_{\dff -}\qff \cap\qff \tau_{\dff +}
\]

\vspace{-9pt}\vspace{-0.25pt}
and\dss hence\qss
$\num{\sigma_{\dff -}\qff \cap\qff \omega_{\dff +}}
\off \leq\off
\num{\sigma_{\dff -}\qff \cap\qff \tau_{\dff +}}$\nsp.\oss
It\dss follows\sss that\qss
$d\trf(\trf \sigma\fff,\pff \omega\trf)
\off \leq\off
d\trf(\trf \sigma\fff,\pff \tau\trf)
\qff -\qff
1$\qss
and\dss hence\sss the first\sss statement\dss holds for\qss
$\sigma\fff,\off \omega$\nnsp,\pss and\dss $f$\dss
by\dss the inductive assumption.\oss
This completes\sss the step of\trs the induction and\dss hence\sss
the proof\dss of\trs the first\sss statement\halfff.\oss

Let\dss us\dss prove\sss the second\sss statement\halfff.\oss
If\oss
(\ref{todd-inclusion-1})\qss --\qss (\ref{todd-inclusion-3})\pss hold\sss and\pss
$\underline{\tau}
\dff\off \subset\qff\off 
\underline{\sigma}\off \cup\pff \{\trf w\qff\}$\nnsp,\qff\oss
then\vspace{3pt}\vspace{-0.25pt}
\[
\quad
\underline{\omega}
\off\qff \subset\off\qff
\underline{\sigma}\off \cup\off \underline{\tau}
\off\qff \subset\off\qff
\underline{\sigma}\off \cup\off 
\left(\pff
\underline{\sigma}\off \cup\pff \{\trf w\qff\}
\qff\right)
\off =\off\qff
\underline{\sigma}\off \cup\pff \{\trf w\qff\}
\qff.
\]

\vspace{-9pt}\vspace{-0.25pt}
By\sss applying\dss the elimination\dss property\dss to\qss 
$\tau\fff,\pff -\qff \omega$\nnsp,\pss and\sss $w$\dss
we\sss get\sss a\sss circuit\sss $\pi$\sss such\dss that\vspace{3pt}\vspace{-0.25pt}
\[
\quad
\underline{\pi}
\off\qff \subset\off\qff
\left(\pff
\underline{\tau}\off \cup\off \underline{\omega} 
\pff\right)
\off \smallsetminus\off 
\{\trf w\qff\}
\off\qff \subset\off\qff
\left(\pff
\underline{\sigma}\off \cup\off \{\trf w\qff\} 
\qff\right)
\off \smallsetminus\off 
\{\trf w\qff\}
\off\dff =\off\pff
\underline{\sigma}
\off.
\]

\vspace{-9pt}\vspace{-0.25pt}
The axiom\qss ({\fff}iii\fff)\qss implies\sss that\qss
$\pi\off =\off \sigma$\qss or\qss ${}-\qff \sigma$\qss
and\dss hence\qss
$\underline{\pi}\off\qff =\off\qff \underline{\sigma}$\nsp.\oss
It\dss follows\dss that\vspace{3pt}\vspace{-0.25pt}
\[
\quad
\underline{\sigma}
\off\qff =\off\qff
\left(\pff
\underline{\tau}\off \cup\off \underline{\omega} 
\pff\right)
\off \smallsetminus\off 
\{\trf w\qff\}
\]

\vspace{-9pt}\vspace{-0.25pt}
and\dss hence\qss
$\underline{\sigma}
\off\qff \subset\off\qff
\underline{\tau}\off \cup\off \underline{\omega}$\nsp.\oss
In\dss turn,\oss this implies\sss that\qss
$\underline{\sigma}
\off \smallsetminus\off
\underline{\tau}
\off\qff \subset\off\qff
\underline{\omega}$\nsp.\oss 

It\dss remains\sss to prove\sss the uniqueness of\dss $\omega$\dss
with\dss the properties\qss
(\ref{todd-inclusion-1})\qss --\qss (\ref{todd-inclusion-3})\qss
under\dss the assumption\qss
$\underline{\tau}
\qff\off \subset\qff\off 
\underline{\sigma}\off \cup\qff \{\trf w\qff\}$\nnsp.\qff\oss
Suppose\sss that\qss 
$\alpha
\off \neq\off
\omega$\qss 
also satisfies\sss these properties.\oss
Since\qss
$\alpha\trf(\trf w\qff) 
\off\qff =\off\qff
\omega\dff(\trf w\qff)$\qss
and\dss hence\qss
$\alpha
\off\dff \neq\off
-\qff \omega$\nnsp,\oss
the axiom\qss ({\fff}iii\fff)\qss implies\sss that\trs\vspace{3pt}\vspace{-0.25pt}
\[
\quad
\underline{\omega}
\off \smallsetminus\dff\off
\underline{\alpha}
\off\pff \neq\off\qff
\varnothing
\hspace{1.5em}\mbox{and}\hspace{1.7em}
\underline{\alpha}
\off \smallsetminus\dff\off
\underline{\omega}
\off\pff \neq\off\qff
\varnothing
\qff.
\]

\vspace{-9pt}\vspace{-0.25pt}
By\sss applying\dss the elimination\sss property\dss to\qss
$\omega\fff,\off {}-\qff \alpha$\nnsp,\oss and\dss $w$\dss
we get\sss a\sss circuit\dss $\pi$\dss such\dss that\vspace{3pt}\vspace{-0.25pt}
\[
\quad
\underline{\pi}
\off\qff \subset\off\qff
\left(\pff
\underline{\omega}\off \cup\off \underline{\alpha} 
\pff\right)
\off \smallsetminus\off 
\{\trf w\qff\}
\off\qff \subset\off\qff
\left(\pff
\underline{\sigma}\off \cup\off \{\trf w\qff\} 
\qff\right)
\off \smallsetminus\off 
\{\trf w\qff\}
\off\dff =\off\pff
\underline{\sigma}
\off.
\]

\vspace{-9pt}\vspace{-0.25pt}
It\dss follows\dss that\qss
$\underline{\pi}\off\qff =\off\qff \underline{\sigma}$\qss
and\qss
$\underline{\sigma}
\off\qff =\off\qff
\left(\pff
\underline{\omega}\off \cup\off \underline{\alpha} 
\pff\right)
\off \smallsetminus\off 
\{\trf w\qff\}$\nnsp.\qff\oss
If\oss
$x\qff \in\pff
\underline{\omega}
\off \smallsetminus\off
\underline{\alpha}$\nsp,\oss
then\vspace{3pt}\vspace{-0.25pt}
\[
\quad
\pi\dff(\dff x\trf)
\off =\off
\omega\dff(\dff x\trf)
\off =\off
\sigma\dff(\dff x\trf)
\]

\vspace{-9pt}\vspace{-0.25pt}
and\dss hence\qss $\pi\off =\off \sigma$\nnsp.\oss
A similar argument\sss using\sss some\dss
$y\qff \in\pff
\underline{\alpha}
\off \smallsetminus\off
\underline{\omega}$\qss
leads\sss to\sss the conclusion\dss that\qss
$\pi\off =\off {}-\qff \sigma$\nnsp.\oss
The contradiction shows\sss that\sss $\omega$\sss is\dss unique.\oss  \eproof

\newpage
\myappend{Cocircuits\qss and\pss lexicographic\qss extensions}{extensions}

\myuppar{Cocircuits.}
Let\dss $M$\dss be an oriented\dss matroid.\oss
It\trs is\dss known\dss that\dss
the number of\dss elements in a base of\trs $M$\dss depends only\sss on $M$\nnsp.\oss
It\dss is\dss called\dss the\qss \emph{rank}\qss of\trs $M$\nnsp.\oss
We will\dss denote it\dss by\sss $n$\nnsp.\oss
The\dss \emph{span}\qss of\dss a subset\dss $X\qff \subset\qff M$\dss
is\dss defined as\sss the set\sss of\dss all\sss elements\dss $a\qff \in\qff M$\dss
such\dss that\dss there exists a circuit\sss $\sigma$ of\trs $M$\sss such\dss that\dss
$a\qff \in\qff \underline{\sigma}$\dss and\dss
$\underline{\sigma}\qff -\qff a\qff \subset\qff X$\nnsp.\oss
A subset\sss of\trs $M$\sss is\dss called
a\dss \emph{hyperplane}\pss if\trs it\dss is\dss equal\dss to\sss the span of\dss 
some independent\sss subset\dss $X\qff \subset\qff M$\dss such\dss that\dss 
$\num{X}\off =\off n\qff -\qff 1$\nnsp.\oss

Let\dss $X\qff \subset\qff M$\dss be an independent\sss subset\sss
such\dss that\dss 
$\num{X}\off =\off n\qff -\qff 1$\nnsp,\oss
and\dss let\sss $H$\sss be\sss the hyperplane spanned\dss by\dss $X$\nnsp.\oss
Clearly\halfff,\oss there exists an element\sss $e\qff \in\qff M$\dss such\dss that\dss
$X\qff +\qff e$\dss is\dss a\sss basis.\oss
For every\dss $u\qff \in\qff M\pff \smallsetminus\qff (\trf X\qff +\qff e\dff)$\dss
there\dss is\sss a\sss unique circuit\dss $\sigma$\sss such\dss that\dss
$u\qff \in\qff \underline{\sigma}$\dss
and\dss
$\underline{\sigma}\qff -\qff u\off \in\off X\qff +\qff e$\nnsp.\oss
By\dss the definition of\trs $H$\nnsp,\oss
if\trs $u\qff \in\qff M\pff \smallsetminus\qff H$\nnsp,\oss
then\dss $e\qff \in\qff \underline{\sigma}$\nsp.\oss
Let\trs $\tau\dff(\dff u\trf)\off =\off\dff +$\dss if\trs
$e\qff \in\qff \sigma_{\dff +}$\dss
and\dss $\tau\dff(\dff u\trf)\off =\off\dff -$\dss if\trs
$e\qff \in\qff \sigma_{\dff -}$\nsp.\oss
Then $\tau$ is\dss a\sss signed subset\sss of\trs $M$\sss such\dss that\trs
$\underline{\tau}\off =\off M\pff \smallsetminus\qff H$\nnsp.\oss
One can\dss prove\sss that\halfff,\oss up\sss to replacing $\tau$ by\sss $-\qff \tau$\nnsp,\oss
the signed subset\sss $\tau$ depends only\sss on\sss $H$\nnsp.\oss
It\dss is\dss called a\qss \emph{cocircuit}\pss of\trs $M$\sss corresponding\dss to\sss $H$\nnsp.\oss
A signed subset\sss of\trs $M$\sss is\dss called\sss a\qss \emph{cocircuit}\pss
if\trs it\sss is\dss a\sss cocircuit\sss corresponding\dss to a hyperplane.\oss

Two signed subsets\dss $\sigma\fff,\pff \tau$\dss are said\dss to be\qss
\emph{ortogonal}\pss if\trs either\dss
$\underline{\sigma}\qff \cap\qff \underline{\tau}\off =\off \varnothing$\nnsp,\oss
or\dss there exist\dss two elements\dss
$u\fff,\pff v\qff \in\qff \underline{\sigma}\qff \cap\qff \underline{\tau}$\dss
such\dss that\dss
$\sigma\dff(\dff u\trf)\off =\off \tau\dff(\dff u\trf)$\dss
and\dss
$\sigma\dff(\dff v\trf)\off =\off -\qff \tau\dff(\dff v\trf)$\nnsp.\oss
If\trs this\dss is\dss the case,\oss then\dss we write\dss
$\sigma\trf \perp\qff \tau$\nnsp.\oss
It\dss turns out\dss that\sss a signed subset\sss $\tau$ of\trs $M$\sss 
is\dss a\sss cocircuit\trs if\trs and\dss only\trs if\qss
$\sigma\trf \perp\qff \tau$\dss for every\sss circuit\sss $\sigma$\sss
of\trs $M$\nnsp.\oss

\myuppar{Las\dss Vergnas\qss ({\fff}lexicographic)\qss extensions.}
An oriented\dss matroid\sss $M\fff'$\sss is\dss said\dss to be an\qss
\emph{extension}\qss of\trs $M$\dss if\qss $M\qff \subset\pff M\fff'$\dss and
a\sss signed subset\sss of\trs $M$\trs is\dss a\sss circuit\sss of\trs $M\fff'$\sss
if\trs and\dss only\trs if\trs it\dss is\dss a\sss circuit\sss of\trs $M$\nnsp.\oss
Clearly\halfff,\oss the rank of\trs $M\fff'$\dss is\dss $\geq\qff n$\nnsp.\oss
We are interested\dss in\dss the\trs \emph{one point\sss extensions}\qss $M\fff'$\dnsp,\oss
i.e.\qss extensions\sss $M\fff'$ such\dss that\trs
$M\fff'\off =\off M\qff +\qff p$\dss for some\dss $p\qff \not\in\qff M$\nnsp.\oss
Up\sss to isomorphism,\oss there\dss is\dss only\sss one such extension\dss
with\dss the rank\dss $n\qff +\qff 1$\nnsp,\oss
and\sss we are interested\sss in\dss the ones of\dss rank $n$\nnsp.\oss
Among\dss such extensions are\sss the\qss
\emph{lexicographic extensions}\qss 
introduced\dss by\dss M.\dss Las\dss Vergnas\qss \cite{lv}.\oss
Let\dss us\dss describe\sss their\dss basic properties  
following\pss M.\dss Todd\qss \cite{t},\oss Theorem\qss 5.1.

Recall\dss that\sss a signed set\sss $\tau$\sss can\dss be considered as a map\dss
$\underline{\tau}
\qff \ttoo\qff
\{\qff +\dff,\pff -\qff\}$\nsp.\oss
This allows\sss to speak about\dss the restriction of\dss $\tau$\sss to a 
subset\sss of\dss $\underline{\tau}$\nsp.\oss
Let\trs $\{\qff a_{\dff 1}\dff,\off a_{\dff 2}\trf,\off \ldots\dff,\off a_{\dff k} \qff\}
\qff \subset\pff M$\trs be an independent\sss set\sss
and\dss let\trs $p\qff \not\in\qff M$\nnsp.\oss
Then\dss there\dss is\dss a\sss unique structure of\dss an oriented\dss matroid\sss on\dss
$M\fff'\off =\off M\qff +\qff p$\dss such\dss that\dss $M\fff'$\sss is\dss an extension of\trs $M$\dss of\trs
the rank\sss $n$\sss and\dss the following\dss two properties hold.\oss\vspace{2.3pt}

({\fff}a{\fff})\oss If\dss $\tau$\sss is\dss a\sss cocircuit\sss of\qss $M$\dss and\dss
$\underline{\tau}
\qff \cap\qff
\{\trf a_{\dff 1}\dff,\off a_{\dff 2}\trf,\off \ldots\dff,\off a_{\dff k} \qff\}
\off =\off
\varnothing$\nnsp,\oss
then\sss $\tau$\sss is\dss a\sss cocircuit\sss of\qss $M\fff'$\nnsp.\vspace{2.3pt}

({\fff}b{\fff})\oss Let\dss $\tau$\sss is\dss a\sss signed subset\sss of\qss $M\fff'$\dss
such\dss that\trs $p\qff \in\qff \underline{\tau}$\dss and\dss
$\underline{\tau}
\qff \cap\qff
\{\qff a_{\dff 1}\dff,\off a_{\dff 2}\trf,\off \ldots\dff,\off a_{\dff k} \qff\}
\off \neq\off
\varnothing$\nnsp,\oss
and\dss 
\hspace*{2.2em}let\dss $i$\dss be\sss the minimal\dss number such\dss that\dss
$e_{\dff i}\qff \in\qff \underline{\tau}$\nsp.\oss
Then $\tau$ is\dss a\sss cocircuit\sss of\trs $M\fff'$\dss if\trs and\dss 
\hspace*{2.2em}only\trs if\trs the restriction of\dss $\tau$\sss to\dss
$\underline{\tau}\qff \cap\qff M$\dss is\dss a\sss cocircuit\sss of\trs $M$\dss
and\dss $\tau\dff(\dff p\dff)\off =\off \tau\dff(\dff e_{\dff i}\dff)$\nnsp.\oss\vspace{2.3pt}

Moreover\halfff,\oss every\sss cocircuit\sss $\tau$ of\trs $M\fff'$\sss
such\dss that\trs $p\qff \in\qff \underline{\tau}$\dss has\sss the form
described\dss in\qss ({\fff}b{\fff}).\oss

If\trs $M$\sss is\dss a subset\sss of\trs $\rrr^{\dff n}$\sss considered as oriented\dss matroid,\oss
then\dss $M\fff'\off =\off M\qff +\qff p$\nnsp,\oss
were\vspace{3pt}
\[
\quad
p
\off\dff =\off
a_{\dff 1}\qff +\qff 
\lambda\dff a_{\dff 2}\qff +\qff 
\lambda^{\dff 2}\dff a_{\dff 3}\qff +\qff 
\ldots\qff +\qff
\lambda^{\dff k\dff -\dff 1}\dff a_{\dff k}
\qff
\]

\vspace{-9pt}
for a sufficiently\sss small\dss $\lambda\qff >\qff 0$\nnsp,\oss
considers as an oriented\dss matroid,\oss
has\sss the above properties and so\dss is\dss a\sss 
lexicographic extension.\oss
In\dss this case\sss the vector $p$\sss is\dss a\sss perturbation of\trs the vector $a_{\dff 1}$\nsp.\oss
It\dss turns out\dss that\dss the added element\sss $p$\sss of\dss a\sss
lexicographic extension can\dss play\dss the role of\dss a\sss perturbation
of\trs the vector $a_{\dff 1}$ also\sss in\dss general\sss case.\oss
This was realized\dss by\qss M.\dss Todd\qss \cite{t}.\oss

\vspace*{12pt}
\begin{flushright}

October\qss 27,\oss 2019.\off\oss
Preface:\oss July\qss 21,\oss 2022
 
https\halfff:/\!/\hspace*{-0.06em}nikolaivivanov.com

E-mail\halfff:\oss nikolai.v.ivanov{\fff}@{\dff}icloud.com,\oss ivanov{\fff}@{\dff}msu.edu

Department\sss of\qss Mathematics,\oss Michigan\sss State\sss University

\end{flushright}

\end{document}